\documentclass[a4paper,12pt]{article} 
\usepackage{amssymb}
\usepackage{amscd}
\usepackage{amsmath} 
\usepackage{graphicx}
\usepackage{latexsym} 
\usepackage{enumerate}

\usepackage{theorem}
\theoremstyle{plain}
\theorembodyfont{\itshape}
\newtheorem{theorem}{Theorem}[section]

\newtheorem{lemma}[theorem]{Lemma}
\newtheorem{corollary}[theorem]{Corollary}
 
\theorembodyfont{\rmfamily}

\newtheorem{definition}[theorem]{Definition}
\newtheorem{remark}[theorem]{Remark}
\newtheorem{problem}[theorem]{Problem}
\newtheorem{example}[theorem]{Example}

\def\Kl{\mathrm{Kl}}

\def\PVI{\mathrm{P}_{\mathrm{VI}}}

\def\RH{\mathrm{RH}} 
\def\rh{\mathrm{rh}} 

\def\Tr{\mathrm{Tr}}

\def\E{\mathcal{E}}

\def\K{\mathcal{K}}
\def\M{\mathcal{M}}

\def\R{\mathcal{R}}
\def\Y{\mathcal{Y}}
\def\Sol{\mathcal{S}}
\def\C{\mathbb{C}}

\def\P{\mathbb{P}}
\def\bR{\mathbb{R}}
\def\Z{\mathbb{Z}}
\def\Wall{\mathbf{Wall}} 
\def\ga{\gamma}
\def\k{\kappa}
\def\lam{\lambda}
\def\si{\sigma}
\def\th{\theta}

\def\Th{\Theta}
\def\ve{\varepsilon}
\def\carl{\circlearrowleft}
\def\car{\curvearrowright}

\def\ol{\overline}

\def\dfrac#1#2{{\displaystyle\frac{#1}{#2}}}

\def\la{\langle}
\def\ra{\rangle}

\def\wt{\widetilde} 
\def\-{\phantom{-}}
\title{\bf Finite Branch Solutions to Painlev\'e VI \\ 
Around a Fixed Singular Point} 
\author{Katsunori Iwasaki \\ \\
Faculty of Mathematics, Kyushu University \\
6-10-1 Hakozaki, Higashi-ku, Fukuoka 812-8581 
Japan\thanks{E-mail address: {\tt iwasaki@math.kyushu-u.ac.jp}}} 
\begin{document}
\date{Dedicated to Professor Kazuo Okamoto on his sixtieth birthday}
\maketitle
\begin{abstract} 
Every finite branch local solution to the sixth Painlev\'{e} 
equation around a fixed singular point is an algebraic branch 
solution. 
In particular a global solution is an algebraic solution if 
and only if it is finitely many-valued globally. 
The proof of this result relies on algebraic geometry of 
Painlev\'{e} VI, Riemann-Hilbert correspondence, geometry 
and dynamics on cubic surfaces, resolutions of Kleinian 
singularities, and power geometry of algebraic differential 
equations. 
In the course of the proof we are also able to classify all 
finite branch solutions up to B\"{a}cklund transformations. 
\end{abstract} 
%%%%%%%%%%%%%%%%%%%%%%%%%% sec:intro %%%%%%%%%%%%%%%%%%%%%%%%%%%%%%
\section{Introduction} \label{sec:intro}
%%%%%%%%%%%%%%%%%%%%%%%%%%%%%%%%%%%%%%%%%%%%%%%%%%%%%%%%%%%%%%%%%%%
We are interested in a {\sl finite branch} local solution to 
the sixth Painlev\'e equation around a fixed singular point. 
We show that every such solution is in fact an {\sl algebraic 
branch} solution (see Definition \ref{def:finite} for the 
terminology). 
In particular a global solution is an {\sl algebraic solution} 
if and only if it is {\sl finitely many-valued} globally. 
Although the problem under study is local in nature, our 
solution to it relies on an effective combination of 
some global technologies and some local tools. 
The former includes the algebraic geometry of the sixth Painlev\'e 
equation, Riemann-Hilbert correspondence, geometry and dynamics on 
cubic surfaces, Kleinian singularities and their minimal resolutions 
\cite{IIS0,IIS1,IIS2,IIS3,IU}, while the latter includes 
the power geometry of algebraic differential equation 
\cite{Bruno1,Bruno2,BG}, which is a method of constructing formal 
solutions by means of Newton polygons, and the theory of nonlinear 
differential equations of ``regular singular type" 
\cite{Gerard,GS}, which discusses the convergence of formal 
solutions. 
\par 
Let us describe our main results in more detail. 
First we recall that the sixth Painlev\'e equation $\PVI(\k)$ 
is a Hamiltonian system of nonlinear differential equations 
%%%%%%%%%%%%%%%%%%%%%%%%%%% eqn:PVI %%%%%%%%%%%%%%%%%%%%%%%%%%%
\begin{equation} \label{eqn:PVI} 
\dfrac{d q}{d z} = \dfrac{\partial H(\k)}{\partial p}, 
\qquad 
\dfrac{d p}{d z} = -\dfrac{\partial H(\k)}{\partial q}, 
\end{equation}
%%%%%%%%%%%%%%%%%%%%%%%%%%%%%%%%%%%%%%%%%%%%%%%%%%%%%%%%%%%%%%%
with time variable $z \in Z := \P^1-\{0,1,\infty\}$ and unknown 
functions $q = q(z)$ and $p = p(z)$, depending on complex 
parameters $\k = (\k_0,\k_1,\k_2,\k_3,\k_4)$ in the 
$4$-dimensional affine space
%%%%%%%%%%%%%%%%%%%%%%%%%%% eqn:K %%%%%%%%%%%%%%%%%%%%%%%%%%%%%
\begin{equation} \label{eqn:K}
\K := \{\, \k = (\k_0,\k_1,\k_2,\k_3,\k_4) \in \C^5 \,:\, 
2 \k_0 + \k_1 + \k_2 + \k_3 +\k_4 = 1 \,\}, 
\end{equation}
%%%%%%%%%%%%%%%%%%%%%%%%%%%%%%%%%%%%%%%%%%%%%%%%%%%%%%%%%%%%%%%
where the Hamiltonian $H(\k) = H(q,p,z;\k)$ is given by 
%%%%%%%%
\[
\begin{array}{rcl}
z(z-1) H(\k) &=&  
(q_0q_1q_z) p^2 - \{\k_1q_1q_z 
+ (\k_2-1)q_0q_1 + \k_3q_0q_z \} p 
+ \k_0(\k_0+\k_4) q_z, 
\end{array}
\]
%%%%%%%%
with $q_{\nu} := q - \nu$ for $\nu \in \{0, 1, z \}$. 
Each of the points $0$, $1$, $\infty$ is called a 
{\sl fixed singular point}. 
%%%%%%%%%
\par 
It is well known that equation (\ref{eqn:PVI}) has the analytic 
Painlev\'{e} property, that is, any meromorphic solution germ 
at a base point $z \in Z$ can be continued {\sl meromorphically} 
along any path in $Z$ emanating from $z$. 
Thus a solution can branch only around a fixed singular point. 
We are interested in finite branch solutions around it, 
by which we mean the following. 
%%%%%%%%%%%%%%%%%%%%%%% def:finite %%%%%%%%%%%%%%%%%%%%%%%%%%%%%
\begin{definition} \label{def:finite} 
A {\sl finite branch} solution to equation (\ref{eqn:PVI}), say, 
around $z = 0$ is a local solution $(q(z),p(z))$ on a punctured 
disk $D^{\times} = D-\{0\}$ centered at $z = 0$ such that its 
lift $(\tilde{q}(\tilde{z}),\tilde{p}(\tilde{z}))$ along some 
finite branched covering 
$\varphi : (\tilde{D}, \tilde{0}) \to (D, 0)$, $\tilde{z} 
\mapsto z = \tilde{z}^n$ around $z = 0$ is a single-valued 
meromorphic function on $\tilde{D}^{\times} = \tilde{D} - 
\{\tilde{0}\}$. 
Such a solution is said to be an {\sl algebraic branch} 
solution if it can be represented by a convergent 
Puiseux-Laurent expansion 
%%%%%%%%%%%%%%%%%%%%%%%%%% eqn:PL %%%%%%%%%%%%%%%%%%%%%%%%%%%%%%
\begin{equation} \label{eqn:PL} 
q(z) = \sum_{i \gg -\infty} a_i \, z^{i/n}, 
\qquad 
p(z) = \sum_{i \gg -\infty} b_i \, z^{i/n}, 
\end{equation} 
%%%%%%%%%%%%%%%%%%%%%%%%%%%%%%%%%%%%%%%%%%%%%%%%%%%%%%%%%%%%%%%%%
with $a_i = b_i = 0$ for all sufficiently small $i \ll 0$, 
namely, if the lift $(\tilde{q}(\tilde{z}),\tilde{p}(\tilde{z}))$ 
is a single-valued meromorphic function on $\tilde{D}$ with at 
most pole at the origin $\tilde{z} = \tilde{0}$. 
\end{definition} 
%%%%%%%%%%%%%%%%%%%%%%%%%% prob:main %%%%%%%%%%%%%%%%%%%%%%%%%
\begin{problem} \label{prob:main} 
Is any finite branch solution to $\PVI(\k)$ an algebraic 
branch solution ? 
\end{problem}
%%%%%%%%%%%%%%%%%%%%%%%%%%%%%%%%%%%%%%%%%%%%%%%%%%%%%%%%%%%%%%
\par 
In this article we settle this problem in the affirmative 
as is stated in the following. 
%%%%%%%%%%%%%%%%%%%%%%%%%% thm:main %%%%%%%%%%%%%%%%%%%%%%%%%%%%
\begin{theorem} \label{thm:main}
Any finite branch solution to Painlev\'e VI around a fixed 
singular point is an algebraic branch solution. 
In particular a global solution is an algebraic solution if 
and only if it is finitely many-valued globally. 
These results are valid for all parameters $\k \in \K$. 
\end{theorem} 
%%%%%%%%%%%%%%%%%%%%%%%%%%%%%%%%%%%%%%%%%%%%%%%%%%%%%%%%%%%%%%%%
\par 
It is an interesting problem to consider algebraic solutions to 
Painlev\'{e} VI. 
Many algebraic solutions have been constructed in 
\cite{AK,Boalch1,Boalch2,DM,Hitchin1,Hitchin2,Kitaev,Mazzocco}, 
but a complete classification seems to be outstanding. 
We hope that Theorem \ref{thm:main} will play an important part 
in discussing this issue. 
The following remark explains what Theorem \ref{thm:main} 
signifies and why it is remarkable. 
%%%%%%%%%%%%%%%%%%%%%%%%%%% rem:main %%%%%%%%%%%%%%%%%%%%%%%%%%%
\begin{remark} \label{rem:main} 
Logically, according to Definition \ref{def:finite}, 
a finite branch solution $(q(z),p(z))$ around 
$z = 0$ may have a very transcendental singularity at $z = 0$, 
to the effect that its lift $(\tilde{q}(\tilde{z}),
\tilde{p}(\tilde{z}))$ may have infinitely many poles in 
$\tilde{D}^{\times}$ accumulating to the origin 
$\tilde{z} = \tilde{0}$, or even if such an accumulation 
phenomenon does not occur, it may have an essential 
singularity at $\tilde{z} = \tilde{0}$. 
Rather surprisingly, however, Theorem \ref{thm:main} excludes 
the possibility for a finite branch solution to admit such 
transcendental phenomena. 
This result becomes more intriguing if we recall that 
wild behaviors of a generic solution to Painlev\'e VI have 
been observed in \cite{Garnier,Guzzetti,Jimbo,Shimomura,Takano} 
and examples of solutions with infintely many poles 
accumulating to $\tilde{z} = 0$ are given in 
\cite{Guzzetti,Shimomura}; such a distribution of poles may 
be expected for a generic solution, though it is not 
rigorously verified yet to the author's knowledge. 
Thus we can think that a finite branch solution is quite 
distinguished from generic solutions, necessarily being an 
algebraic branch solution. 
\end{remark} 
%%%%%%%%%%%%%%%%%%%%%%%%%%%%%%%%%%%%%%%%%%%%%%%%%%%%%%%%%%%%%%%%
\par
\begin{figure}[t]
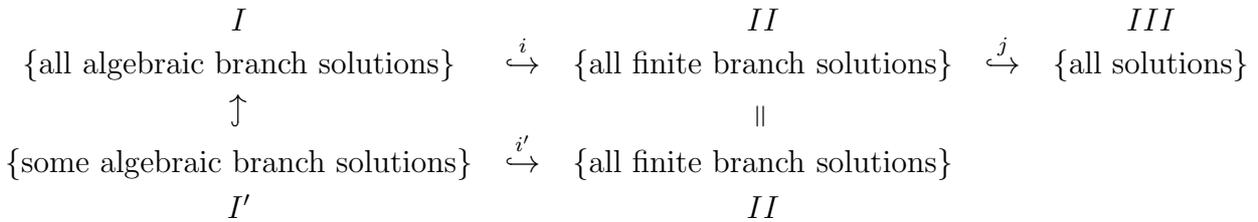

\begin{center}
\begin{tabular}{ccccc} 
  $I$ & & $II$ & & $III$ \\[-1mm]
$\{\mbox{all algebraic branch solutions}\}$  
& $\overset{i}{\hookrightarrow}$ &
$\{\mbox{all finite branch solutions}\}$ &
$\overset{j}{\hookrightarrow}$ & $\{\mbox{all solutions}\}$ \\[1mm] 
$\mbox{\rotatebox{90}{$\hookrightarrow$}}$ & & 
$\mbox{\rotatebox{90}{$=$}}$ & & \\[-1mm]
$\{\mbox{some algebraic branch solutions}\}$ 
& $\overset{i'}{\hookrightarrow}$ &
$\{\mbox{all finite branch solutions}\}$ & & \\[1mm] 
 $I'$ & & $II$ & & \\
\end{tabular}
\end{center}
\caption{Main idea for the proof of Theorem \ref{thm:main}}
\label{fig:idea}
\end{figure}

 
The main idea for the proof of Theorem \ref{thm:main} is 
presented in Figure \ref{fig:idea}. 
We have natural inclusions $i$ and $j$ in the top line of 
Figure \ref{fig:idea} and we wish to show that the injection 
$i$ is in fact a surjection. 
Our strategy consists of the ``upper bound part" and the 
``lower bound part". 
%%%%%%%%%%%%%%%%%%%%%%%%%%%%%% fig:idea %%%%%%%%%%%%%%%%%%%%%%%%%%
%%%%%%%%%%%%%%%%%%%%%%%%%%%%%%%%%%%%%%%%%%%%%%%%%%%%%%%%%%%%%%%
\begin{enumerate} 
\item {\sl Upper bound part}: 
In this part we investigate the inclusion $j : II 
\hookrightarrow III$ in Figure \ref{fig:idea}, 
considering how the locus of finite branch solutions is 
included in the moduli space of all solutions. 
In other words, we make a confinement of the locus $II$ 
in the entire space $III$. 
What we shall really do is not an upper bound 
estimation of this locus but rather a pinpoint 
identification of it. 
This is the main part of the article and we use the 
algebraic geometry of Painlev\'e VI, Riemann-Hilbert 
correspondence, geometry and dynamics on cubic surfaces, 
and minimal resolutions of Kleinian singularities 
\cite{IIS0,IIS1,IIS2,IIS3,Iwasaki,IU}. 
\item {\sl Lower bound part}: In this part we fill in 
the diagram of Figure \ref{fig:idea} by adding the 
bottom line to the top one. 
We try to construct as many algebraic branch 
solutions as possible in order to make the set $I'$ 
as large as possible. 
The construction is based on the power geometry 
technique developed in \cite{Bruno1,Bruno2,BG} and 
the convergence arguments in \cite{Gerard,GS}. 
We are done if the set $I'$ is large enough to 
show that the injection $i' : I' \hookrightarrow II$ 
is in fact a surjection. 
This does not mean that we verify the equality 
$I' = II$ directly. 
(If such a direct approach were feasible, then our 
problem would not be difficult from the beginning!) 
Instead, we prove it very indirectly based on the 
following idea. 
\item {\sl Key trick}: Suppose that a component 
$A$ of $I'$ injects into a component $B$ of $II$. 
If the cardinalities of $A$ and $B$ are finite and 
the same, then the injection $i' : A \hookrightarrow 
B$ is in fact a surjection. 
If $A$ and $B$ are biholomorphic to $\C$ and 
the injection $i' : A \hookrightarrow B$ is 
holomorphic, then it must be a surjection 
because any holomorphic injection 
$\C \hookrightarrow \C$ is a surjection (use 
Casorati-Weierstrass or Picard's little theorem). 
The same argument holds true if $\C$ is replaced 
by $\C^{\times}$, since any holomorphic injection 
$\C^{\times} \hookrightarrow \C^{\times}$ is a 
surjection (lift it to the universal covering 
$\C \hookrightarrow \C$). 
These tricks enable us to identify the component 
$A \subset I'$ with the component $B \subset II$. 
We show that each component involved is either of 
the three types mentioned above. 
Then we make this kind of argument componentwise to 
get an identification $I' = II$, which leads to 
the desired coincidence $I' = I = II$. 
\end{enumerate} 
%%%%%%%%%%%%%%
\par 
In view of the way in which Theorem \ref{thm:main} is 
established, the power geometry technique provides us with 
an efficient method of identifying all finite branch 
solutions (up to B\"{a}cklund transformations), which have 
now turned out to be algebraic branch solutions, by 
determining the leading terms of their Puiseux-Laurent 
expansions. 
\par
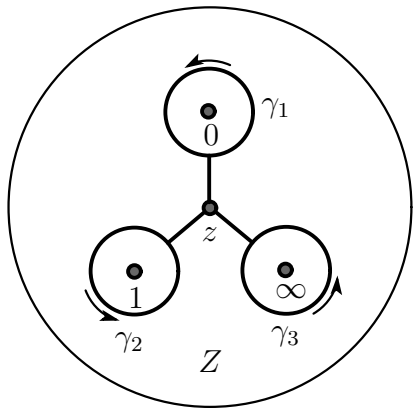
\begin{figure}[t] 
\begin{center}
%WinTpicVersion2.15
\unitlength 0.1in
\begin{picture}(20.86,20.86)(3.60,-23.16)
% CIRCLE 1 0 3 0
% 4 1403 1673 2323 2164 2323 2164 2323 2164
% 
\special{pn 13}%
\special{ar 1403 1273 1043 1043  0.0000000 6.2831853}%
% CIRCLE 0 0 3 0
% 4 1399 1177 1403 1404 1403 1404 1403 1404
% 
\special{pn 20}%
\special{ar 1399 777 227 227  0.0000000 6.2831853}%
% CIRCLE 0 0 0 0
% 4 1403 1677 1424 1702 1424 1702 1424 1702
% 
\special{pn 20}%
\special{sh 0.600}%
\special{ar 1403 1277 33 33  0.0000000 6.2831853}%
% LINE 0 0 3 0
% 2 1399 1404 1399 1648
% 
\special{pn 20}%
\special{pa 1399 1004}%
\special{pa 1399 1248}%
\special{fp}%
% CIRCLE 0 0 3 0
% 4 1012 2007 1182 1857 1182 1857 1182 1857
% 
\special{pn 20}%
\special{ar 1012 1607 227 227  0.0000000 6.2831853}%
% LINE 0 0 3 0
% 2 1185 1860 1371 1702
% 
\special{pn 20}%
\special{pa 1185 1460}%
\special{pa 1371 1302}%
\special{fp}%
% CIRCLE 0 0 3 0
% 4 1798 2003 1621 1861 1621 1861 1621 1861
% 
\special{pn 20}%
\special{ar 1798 1603 227 227  0.0000000 6.2831853}%
% LINE 0 0 3 0
% 2 1623 1858 1437 1702
% 
\special{pn 20}%
\special{pa 1623 1458}%
\special{pa 1437 1302}%
\special{fp}%
% CIRCLE 0 0 0 0
% 4 1395 1173 1416 1198 1416 1198 1416 1198
% 
\special{pn 20}%
\special{sh 0.600}%
\special{ar 1395 773 33 33  0.0000000 6.2831853}%
% CIRCLE 0 0 0 0
% 4 1012 2009 1033 2034 1033 2034 1033 2034
% 
\special{pn 20}%
\special{sh 0.600}%
\special{ar 1012 1609 33 33  0.0000000 6.2831853}%
% CIRCLE 0 0 0 0
% 4 1794 2001 1815 2026 1815 2026 1815 2026
% 
\special{pn 20}%
\special{sh 0.600}%
\special{ar 1794 1601 33 33  0.0000000 6.2831853}%
% STR 2 0 3 0
% 3 1370 1308 1370 1350 2 0
% $0$
\put(13.7000,-9.5000){\makebox(0,0)[lb]{$0$}}%
% STR 2 0 3 0
% 3 980 2158 980 2200 2 0
% $1$
\put(9.8000,-18.0000){\makebox(0,0)[lb]{$1$}}%
% STR 2 0 3 0
% 3 1740 2108 1740 2150 2 0
% $\infty$
\put(17.4000,-17.5000){\makebox(0,0)[lb]{$\infty$}}%
% STR 2 0 3 0
% 3 1360 1808 1360 1850 2 0
% $z$
\put(13.6000,-14.5000){\makebox(0,0)[lb]{$z$}}%
% STR 2 0 3 0
% 3 1343 2491 1343 2533 2 0
% $Z$
\put(13.4300,-21.3300){\makebox(0,0)[lb]{$Z$}}%
% CIRCLE 1 0 3 0
% 4 1397 1176 1507 932 1507 932 1309 928
% 
\special{pn 13}%
\special{ar 1397 776 268 268  4.3714100 5.1359243}%
% VECTOR 1 0 3 0
% 2 1318 920 1288 937
% 
\special{pn 13}%
\special{pa 1318 520}%
\special{pa 1288 537}%
\special{fp}%
\special{sh 1}%
\special{pa 1288 537}%
\special{pa 1356 522}%
\special{pa 1334 511}%
\special{pa 1336 487}%
\special{pa 1288 537}%
\special{fp}%
% CIRCLE 1 0 3 0
% 4 1008 2014 759 2108 759 2108 895 2251
% 
\special{pn 13}%
\special{ar 1008 1614 266 266  2.0157068 2.7806232}%
% VECTOR 1 0 3 0
% 2 884 2251 916 2260
% 
\special{pn 13}%
\special{pa 884 1851}%
\special{pa 916 1860}%
\special{fp}%
\special{sh 1}%
\special{pa 916 1860}%
\special{pa 857 1823}%
\special{pa 865 1846}%
\special{pa 846 1861}%
\special{pa 916 1860}%
\special{fp}%
% CIRCLE 1 0 3 0
% 4 1804 2000 1939 2229 1939 2229 2057 2072
% 
\special{pn 13}%
\special{ar 1804 1600 266 266  0.2772553 1.0381186}%
% VECTOR 1 0 3 0
% 2 2060 2082 2062 2050
% 
\special{pn 13}%
\special{pa 2060 1682}%
\special{pa 2062 1650}%
\special{fp}%
\special{sh 1}%
\special{pa 2062 1650}%
\special{pa 2038 1715}%
\special{pa 2059 1703}%
\special{pa 2078 1718}%
\special{pa 2062 1650}%
\special{fp}%
% STR 2 0 3 0
% 3 910 2388 910 2430 2 0
% $\gamma_2$
\put(9.1000,-20.3000){\makebox(0,0)[lb]{$\gamma_2$}}%
% STR 2 0 3 0
% 3 1720 2358 1720 2400 2 0
% $\gamma_3$
\put(17.2000,-20.0000){\makebox(0,0)[lb]{$\gamma_3$}}%
% STR 2 0 3 0
% 3 1670 1168 1670 1210 2 0
% $\gamma_1$
\put(16.7000,-8.1000){\makebox(0,0)[lb]{$\gamma_1$}}%
\end{picture}%
\end{center}
\caption{Three basic loops $\ga_1$, $\ga_2$, $\ga_3$ in 
$Z = \P^1-\{0,1,\infty\}$} 
\label{fig:selem} 
\end{figure}

 
In some sense this article is a counterpart of the 
previous paper \cite{IU} where an ergodic study of 
Painlev\'e VI is developed (see also the survey \cite{IU2}). 
Put $z_1 = 0$, $z_2 = 1$, $z_3 = \infty$. 
For each $\{i,j,k\} = \{1,2,3\}$, let $\ga_i$ be a loop in 
$Z$ surrounding $z_i$ once anti-clockwise and leaving 
$z_j$ and $z_k$ outside as in Figure \ref{fig:selem}. 
Then the fundamental group $\pi_1(Z,z)$ is represented as 
%%%%%%%%%%%%%%%%%%%%%%% eqn:fundgr %%%%%%%%%%%%%%%%%%%%%%%%%%
\begin{equation} \label{eqn:fundgr} 
\pi_1(Z,z) = 
\la \, \ga_1, \ga_2, \ga_3 \,|\, \ga_1\ga_2\ga_3 = 1 \, \ra. 
\end{equation} 
%%%%%%%%%%%%%%%%%%%%%%%%%%%%%%%%%%%%%%%%%%%%%%%%%%%%%%%%%%%%%%
A loop $\ga \in \pi_1(Z,z)$ is said to be {\sl elementary} 
if it is conjugate to $\ga_i^m$ for some $i \in \{1,2,3\}$ 
and $m \in \Z$; otherwise, it is said to be 
{\sl non-elementary}. 
The main theme of \cite{IU} is the dynamics of the nonlinear 
monodromy of $\PVI(\k)$ along a given loop $\ga$. 
It is shown there that, along every non-elementary loop, the 
nonlinear monodromy is chaotic and the number of its periodic 
points grows exponentially as the period tends to infinity. 
On the other hand, it is Liouville integrable along an 
elementary loop, in the sense that it preserves a 
Lagrangian fibration. 
Now we notice that from the dynamical point of view the main 
problem of this article is nothing other than discussing 
the periodic points of the nonlinear monodromy along the 
basic loop $\ga_i$, which is of course an elementary loop. 
In view of its integrable character, one may doubt if there 
is something very deep with this issue. 
As Theorem \ref{thm:main} and Remark \ref{rem:main} show, 
however, this issue is actually quite interesting from the 
function-theoretical point of view. 
%%%%%%%%%%%%%%%%%%%%%% fig:selem %%%%%%%%%%%%%%%%%%%%%%%%%%%%%%
%%%%%%%%%%%%%%%%%%%%%%%%%%%%%%%%%%%%%%%%%%%%%%%%%%%%%%%%%%%%%%
\par 
The plan of this article is as follows. 
In \S\ref{sec:PhaseSpace} the phase space of Painlev\'e VI 
is introduced as a moduli space of stable parabolic 
connections. 
In \S\ref{sec:RHC} the Riemann-Hilbert correspondence 
from the moduli space to an affine cubic surface is 
formulated and its character as an analytic minimal 
resolution of Kleininan singularities is stated. 
In \S\ref{sec:cubic} the dynamical system on the cubic 
surface representing the nonlinear monodromy of 
Painlev\'e VI is formulated and some preliminary 
properties of it are given. 
In \S\ref{sec:backlund} we briefly review B\"acklund 
transformations and their relation to the 
Riemann-Hilbert correspondence. 
In \S\ref{sec:fixed} fixed points and periodic 
points of the dynamical system are discussed. 
A stratification of the parameter space $\K$ 
is also introduced in order to describe the 
singularities of the cubic surfaces. 
In \S\ref{sec:case} a case-by-case study of 
fixed points and periodic points is made according 
to the stratification, thereby a pinpoint 
identification of finite branch solutions is made 
on each stratum. 
In \S\ref{sec:power} power geometry of algebraic 
differential equations is applied to 
Painlev\'e VI in order to construct as many 
algebraic branch solutions as possible. 
In \S\ref{sec:surjection} we consider the inclusion 
of those solutions constructed in \S\ref{sec:power} 
into the moduli space of all finite branch solutions. 
After some preliminaries on Riccati solutions, 
we show that this inclusion is in fact a surjection, 
thereby complete the proof of Theorem \ref{thm:main}. 
%%%%%%%%%%%%%%%%%%%%% sec:PhaseSpace %%%%%%%%%%%%%%%%%%%%%%%%
\section{Phase Space} \label{sec:PhaseSpace}
%%%%%%%%%%%%%%%%%%%%%%%%%%%%%%%%%%%%%%%%%%%%%%%%%%%%%%%%%%%%%%
Equation (\ref{eqn:PVI}) is only a fragmentary appearance of a 
more intrinsic object constructed algebro-geometrically 
\cite{IIS1,IIS2,IIS3}. 
We review this construction following the expositions of 
\cite{IU,IU2}. 
The sixth Painlev\'e dynamical system $\PVI(\k)$ is formulated 
as a holomorphic, uniform, transversal foliation on a fibration 
of certain smooth quasi-projective rational surfaces 
%%%%%%%%%% 
\begin{equation*}
\pi_{\k} : \M(\k) \to Z := \P^1-\{0,1,\infty\}, 
\end{equation*}
%%%%%%%%%%
whose fiber $\M_z(\k) := \pi_{\k}^{-1}(z)$ over $z \in Z$, 
called the {\sl space of initial conditions} at time $z$, is 
realized as a moduli space of stable parabolic connections. 
The total space $\M(\k)$ is called the {\sl phase space} of 
$\PVI(\k)$. 
In this formulation, the uniformity of the Painlev\'e foliation, 
in other words, the geometric Painlev\'e property of it is a 
natural consequence of a solution to the Riemann-Hilbert problem 
(see Theorem \ref{thm:SolRHP}), especially of the properness of 
the Riemann-Hilbert correspondence \cite{IIS1}. 
Then equation (\ref{eqn:PVI}) is just a coordinate expression of 
the foliation on an affine open subset of $\M(\k)$ and the analytic 
Painlev\'e property for equation (\ref{eqn:PVI}) is an immediate 
consequence of the geometric Painlev\'e property for the 
foliation and the algebraicity of the phase space $\M(\k)$. 
Moreover there exists a natural compactification 
$\M_z(\k) \hookrightarrow \ol{\M}_z(\k)$ 
of the moduli space $\M_z(\k)$ into a moduli space 
$\ol{\M}_z(\k)$ of stable parabolic phi-connections. 
\par
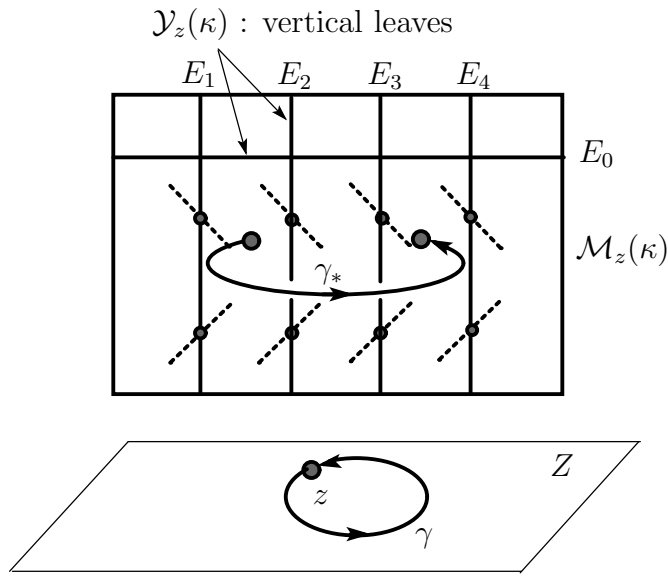
\begin{figure}[t]
\begin{center}
%WinTpicVersion2.15
\unitlength 0.1in
\begin{picture}(32.66,29.53)(1.30,-30.79)
% LINE 0 0 3 0
% 2 2517 990 2517 2546
% 
\special{pn 20}%
\special{pa 2517 590}%
\special{pa 2517 2146}%
\special{fp}%
% LINE 0 0 3 0
% 2 1118 983 1118 2538
% 
\special{pn 20}%
\special{pa 1118 583}%
\special{pa 1118 2138}%
\special{fp}%
% BOX 0 0 3 0
% 2 667 983 2992 2546
% 
\special{pn 20}%
\special{pa 667 583}%
\special{pa 2992 583}%
\special{pa 2992 2146}%
\special{pa 667 2146}%
\special{pa 667 583}%
\special{fp}%
% LINE 0 0 3 0
% 2 667 1310 3000 1310
% 
\special{pn 20}%
\special{pa 667 910}%
\special{pa 3000 910}%
\special{fp}%
% LINE 0 2 3 0
% 2 1436 1458 1740 1776
% 
\special{pn 20}%
\special{pa 1436 1058}%
\special{pa 1740 1376}%
\special{dt 0.054}%
\special{pa 1740 1376}%
\special{pa 1740 1376}%
\special{dt 0.054}%
% LINE 0 2 3 0
% 2 1896 1450 2199 1768
% 
\special{pn 20}%
\special{pa 1896 1050}%
\special{pa 2199 1368}%
\special{dt 0.054}%
\special{pa 2199 1368}%
\special{pa 2199 1368}%
\special{dt 0.054}%
% LINE 0 2 3 0
% 2 2370 1465 2673 1784
% 
\special{pn 20}%
\special{pa 2370 1065}%
\special{pa 2673 1384}%
\special{dt 0.054}%
\special{pa 2673 1384}%
\special{pa 2673 1384}%
\special{dt 0.054}%
% LINE 0 2 3 0
% 2 962 1458 1266 1776
% 
\special{pn 20}%
\special{pa 962 1058}%
\special{pa 1266 1376}%
\special{dt 0.054}%
\special{pa 1266 1376}%
\special{pa 1266 1376}%
\special{dt 0.054}%
% LINE 0 2 3 0
% 2 1273 2079 962 2382
% 
\special{pn 20}%
\special{pa 1273 1679}%
\special{pa 962 1982}%
\special{dt 0.054}%
\special{pa 962 1982}%
\special{pa 962 1982}%
\special{dt 0.054}%
% LINE 0 2 3 0
% 2 1732 2079 1421 2382
% 
\special{pn 20}%
\special{pa 1732 1679}%
\special{pa 1421 1982}%
\special{dt 0.054}%
\special{pa 1421 1982}%
\special{pa 1421 1982}%
\special{dt 0.054}%
% LINE 0 2 3 0
% 2 2207 2072 1896 2375
% 
\special{pn 20}%
\special{pa 2207 1672}%
\special{pa 1896 1975}%
\special{dt 0.054}%
\special{pa 1896 1975}%
\special{pa 1896 1975}%
\special{dt 0.054}%
% LINE 0 2 3 0
% 2 2665 2064 2354 2367
% 
\special{pn 20}%
\special{pa 2665 1664}%
\special{pa 2354 1967}%
\special{dt 0.054}%
\special{pa 2354 1967}%
\special{pa 2354 1967}%
\special{dt 0.054}%
% STR 2 0 3 0
% 3 3062 1799 3062 1877 2 0
% $\mathcal{M}_z(\kappa)$
\put(30.6200,-14.7700){\makebox(0,0)[lb]{$\mathcal{M}_z(\kappa)$}}%
% CIRCLE 0 0 0 0
% 4 2517 1621 2517 1648 2517 1648 2517 1648
% 
\special{pn 20}%
\special{sh 0.600}%
\special{ar 2517 1221 27 27  0.0000000 6.2831853}%
% CIRCLE 0 0 0 0
% 4 2525 2212 2525 2238 2525 2238 2525 2238
% 
\special{pn 20}%
\special{sh 0.600}%
\special{ar 2525 1812 26 26  0.0000000 6.2831853}%
% CIRCLE 0 0 0 0
% 4 2051 2227 2051 2255 2051 2255 2051 2255
% 
\special{pn 20}%
\special{sh 0.600}%
\special{ar 2051 1827 28 28  0.0000000 6.2831853}%
% CIRCLE 0 0 0 0
% 4 2059 1628 2059 1656 2059 1656 2059 1656
% 
\special{pn 20}%
\special{sh 0.600}%
\special{ar 2059 1228 28 28  0.0000000 6.2831853}%
% CIRCLE 0 0 0 0
% 4 1592 1636 1592 1664 1592 1664 1592 1664
% 
\special{pn 20}%
\special{sh 0.600}%
\special{ar 1592 1236 28 28  0.0000000 6.2831853}%
% CIRCLE 0 0 0 0
% 4 1592 2227 1592 2255 1592 2255 1592 2255
% 
\special{pn 20}%
\special{sh 0.600}%
\special{ar 1592 1827 28 28  0.0000000 6.2831853}%
% CIRCLE 0 0 0 0
% 4 1118 1628 1118 1656 1118 1656 1118 1656
% 
\special{pn 20}%
\special{sh 0.600}%
\special{ar 1118 1228 28 28  0.0000000 6.2831853}%
% CIRCLE 0 0 0 0
% 4 1118 2227 1118 2255 1118 2255 1118 2255
% 
\special{pn 20}%
\special{sh 0.600}%
\special{ar 1118 1827 28 28  0.0000000 6.2831853}%
% LINE 2 0 3 0
% 2 145 3471 2781 3471
% 
\special{pn 8}%
\special{pa 145 3071}%
\special{pa 2781 3071}%
\special{fp}%
% LINE 2 0 3 0
% 2 744 2794 3396 2794
% 
\special{pn 8}%
\special{pa 744 2394}%
\special{pa 3396 2394}%
\special{fp}%
% STR 2 0 3 0
% 3 869 618 869 696 2 0
% $\mathcal{Y}_z(\kappa)$ : vertical leaves
\put(8.6900,-2.9600){\makebox(0,0)[lb]{$\mathcal{Y}_z(\kappa)$ : vertical leaves}}%
% VECTOR 2 0 3 0
% 2 1210 744 1559 1109
% 
\special{pn 8}%
\special{pa 1210 344}%
\special{pa 1559 709}%
\special{fp}%
\special{sh 1}%
\special{pa 1559 709}%
\special{pa 1527 647}%
\special{pa 1522 670}%
\special{pa 1498 675}%
\special{pa 1559 709}%
\special{fp}%
% STR 2 0 3 0
% 3 2929 2895 2929 2973 2 0
% $Z$
\put(29.2900,-25.7300){\makebox(0,0)[lb]{$Z$}}%
% CIRCLE 0 0 0 0
% 4 1693 2942 1708 2981 1708 2981 1708 2981
% 
\special{pn 20}%
\special{sh 0.600}%
\special{ar 1693 2542 42 42  0.0000000 6.2831853}%
% STR 2 0 3 0
% 3 1700 3042 1700 3120 2 0
% $z$
\put(17.0000,-27.2000){\makebox(0,0)[lb]{$z$}}%
% LINE 2 0 3 0
% 2 744 2794 130 3471
% 
\special{pn 8}%
\special{pa 744 2394}%
\special{pa 130 3071}%
\special{fp}%
% LINE 2 0 3 0
% 2 3388 2794 2774 3479
% 
\special{pn 8}%
\special{pa 3388 2394}%
\special{pa 2774 3079}%
\special{fp}%
% ELLIPSE 0 0 3 0
% 4 1926 3082 1560 2880 1677 2942 1740 2834
% 
\special{pn 20}%
\special{ar 1926 2682 366 202  4.3196555 6.2831853}%
\special{ar 1926 2682 366 202  0.0000000 3.9369308}%
% VECTOR 0 0 3 0
% 2 1802 2903 1756 2912
% 
\special{pn 20}%
\special{pa 1802 2503}%
\special{pa 1756 2512}%
\special{fp}%
\special{sh 1}%
\special{pa 1756 2512}%
\special{pa 1825 2519}%
\special{pa 1808 2502}%
\special{pa 1818 2480}%
\special{pa 1756 2512}%
\special{fp}%
% VECTOR 0 0 3 0
% 2 1896 3277 1965 3277
% 
\special{pn 20}%
\special{pa 1896 2877}%
\special{pa 1965 2877}%
\special{fp}%
\special{sh 1}%
\special{pa 1965 2877}%
\special{pa 1898 2857}%
\special{pa 1912 2877}%
\special{pa 1898 2897}%
\special{pa 1965 2877}%
\special{fp}%
% STR 2 0 3 0
% 3 2230 3273 2230 3350 2 0
% $\gamma$
\put(22.3000,-29.5000){\makebox(0,0)[lb]{$\gamma$}}%
% STR 2 0 3 0
% 3 3077 1278 3077 1357 2 0
% $E_0$
\put(30.7700,-9.5700){\makebox(0,0)[lb]{$E_0$}}%
% STR 2 0 3 0
% 3 1016 866 1016 944 2 0
% $E_1$
\put(10.1600,-5.4400){\makebox(0,0)[lb]{$E_1$}}%
% STR 2 0 3 0
% 3 1514 874 1514 951 2 0
% $E_2$
\put(15.1400,-5.5100){\makebox(0,0)[lb]{$E_2$}}%
% STR 2 0 3 0
% 3 1981 874 1981 951 2 0
% $E_3$
\put(19.8100,-5.5100){\makebox(0,0)[lb]{$E_3$}}%
% STR 2 0 3 0
% 3 2440 874 2440 951 2 0
% $E_4$
\put(24.4000,-5.5100){\makebox(0,0)[lb]{$E_4$}}%
% ELLIPSE 0 0 3 0
% 4 1818 1862 1156 1698 1397 1746 2400 1761
% 
\special{pn 20}%
\special{ar 1818 1462 662 164  5.6717678 6.2831853}%
\special{ar 1818 1462 662 164  0.0000000 3.9798100}%
% CIRCLE 0 0 0 0
% 4 1382 1746 1397 1784 1397 1784 1397 1784
% 
\special{pn 20}%
\special{sh 0.600}%
\special{ar 1382 1346 41 41  0.0000000 6.2831853}%
% CIRCLE 0 0 0 0
% 4 2260 1737 2276 1776 2276 1776 2276 1776
% 
\special{pn 20}%
\special{sh 0.600}%
\special{ar 2260 1337 42 42  0.0000000 6.2831853}%
% VECTOR 0 0 3 0
% 2 1752 2025 1877 2025
% 
\special{pn 20}%
\special{pa 1752 1625}%
\special{pa 1877 1625}%
\special{fp}%
\special{sh 1}%
\special{pa 1877 1625}%
\special{pa 1810 1605}%
\special{pa 1824 1625}%
\special{pa 1810 1645}%
\special{pa 1877 1625}%
\special{fp}%
% VECTOR 0 0 3 0
% 4 2400 1784 2330 1753 2330 1753 2330 1753
% 
\special{pn 20}%
\special{pa 2400 1384}%
\special{pa 2330 1353}%
\special{fp}%
\special{sh 1}%
\special{pa 2330 1353}%
\special{pa 2383 1398}%
\special{pa 2379 1375}%
\special{pa 2399 1362}%
\special{pa 2330 1353}%
\special{fp}%
\special{pa 2330 1353}%
\special{pa 2330 1353}%
\special{fp}%
% STR 2 0 3 0
% 3 1708 1892 1708 1970 2 0
% $\gamma_*$
\put(17.0800,-15.7000){\makebox(0,0)[lb]{$\gamma_*$}}%
% LINE 0 0 3 0
% 2 1589 998 1589 1948
% 
\special{pn 20}%
\special{pa 1589 598}%
\special{pa 1589 1548}%
\special{fp}%
% LINE 0 0 3 0
% 2 1589 2054 1589 2543
% 
\special{pn 20}%
\special{pa 1589 1654}%
\special{pa 1589 2143}%
\special{fp}%
% LINE 0 0 3 0
% 2 2050 988 2050 1958
% 
\special{pn 20}%
\special{pa 2050 588}%
\special{pa 2050 1558}%
\special{fp}%
% LINE 0 0 3 0
% 2 2050 2054 2050 2543
% 
\special{pn 20}%
\special{pa 2050 1654}%
\special{pa 2050 2143}%
\special{fp}%
% VECTOR 2 0 3 0
% 2 1210 752 1354 1288
% 
\special{pn 8}%
\special{pa 1210 352}%
\special{pa 1354 888}%
\special{fp}%
\special{sh 1}%
\special{pa 1354 888}%
\special{pa 1356 818}%
\special{pa 1340 836}%
\special{pa 1317 829}%
\special{pa 1354 888}%
\special{fp}%
\end{picture}%
\end{center}
\caption{Nonlinear monodromy $\ga_* : \M_z(\k) \carl$ 
along a loop $\ga \in \pi_1(Z,z)$} 
\label{fig:sps} 
\end{figure}

%%%%%%%%%%%%%%%%%%%%%%%%%%%% tab:riemann %%%%%%%%%%%%%%%%%%%%%%%%
\begin{table}[t]
\begin{center} 
\begin{tabular}{|c||c|c|c|c|}
\hline
\vspace{-4mm} &         &         &         &              \\
singularities & $t_1=0$ & $t_2=z$ & $t_3=1$ & $t_4=\infty$ \\[1mm]
\hline
\vspace{-4mm} &         &         &         &              \\
first exponent & $-\lam_1$ & $-\lam_2$ & $-\lam_3$ & $-\lam_4$ \\[1mm]
\hline
\vspace{-4mm} &         &         &         &              \\
second exponent & $\lam_1$ & $\lam_2$  & $\lam_3$ & $\lam_4-1$ \\[1mm]
\hline
\vspace{-4mm} &         &         &         &              \\
difference  & $\k_1$ & $\k_2$ & $\k_3$ & $\k_4$ \\[1mm]  
\hline
\end{tabular}
\end{center}
\caption{Riemann scheme: $\k_i$ is the difference of the second 
exponent from the first.} 
\label{tab:riemann}
\end{table}
%%%%%%%%%%%%%%%%%%%%%%%%%%%%%%%%%%%%%%%%%%%%%%%%%%%%%%%%%%%%%%%%%%%
Here we include a very sketchy explanation of the terminology 
used in the last paragraph. 
A {\sl stable parabolic connection} is a Fuchsian connection 
equipped with a parabolic structure on a (rank $2$) vector bundle 
over $\P^1$ having a Riemann scheme as in Table \ref{tab:riemann}, 
where the parabolic structure corresponds to the first exponents,  
which satisfies a sort of stability condition in geometric 
invariant theory. 
Here the parameter $\k_i$ stands for the difference of the second 
exponent from the first one at the regular singular point $t_i$. 
On the other hand, a {\sl stable parabolic phi-connection} is 
a variant of stable parabolic connection allowing a 
``matrix-valued Planck constant" called a phi-operator $\phi$ 
such that the generalized Leibniz rule 
%%%%
\[
\nabla(f s) = df \otimes \phi(s) + f \nabla(s) 
\]
%%%%
is satisfied, where the key point here is that the field 
$\phi$ may be degenerate or simi-classical. 
Then the moduli space $\M_z(\k)$ can be compactified by 
adding some semi-classical objects, that is, some stable 
parabolic phi-connections with degenerate phi-operator $\phi$. 
There is the following characterization of our moduli spaces 
(see Figure \ref{fig:sps}).
%%%%%%%%%%%%%%%%%%%%%% thm:moduli %%%%%%%%%%%%%%%%%%%%%%%%%%%%%
\begin{theorem}[\cite{IIS1,IIS2,IIS3}] \label{thm:moduli} 
$\phantom{a}$ 
\begin{enumerate}
\item The compactified moduli space $\ol{\M}_z(\k)$ is isomorphic 
to an $8$-point blow-up of the Hirzebruch surface 
$\varSigma_2 \to \P^1$ of degree $2$. 
\item $\ol{\M}_z(\k)$ has a unique effective anti-canonical 
divisor $\Y_z(\k)$, which is given by 
%%%%%%%%%%%%%%%%%%%%%%%% eqn:anti-canonical %%%%%%%%%%%%%%%%%%%%%
\begin{equation} \label{eqn:anti-canonical} 
\Y_z(\k) = 2 E_0 + E_1 + E_2 + E_3 + E_4, 
\end{equation}
%%%%%%%%%%%%%%%%%%%%%%%%%%%%%%%%%%%%%%%%%%%%%%%%%%%%%%%%%%%%%%%%%%
where $E_0$ is the strict transform of the section at infinity and 
$E_i$ $(i = 1,2,3,4)$ is the strict transform of the fiber 
over the point $t_i \in \P^1$ of the Hirzebruch 
surface ${\mit\Sigma}_2 \to \P^1$. 
\item The support of the divisor $\Y_z(\k)$ is exactly the locus 
where the phi-operator $\phi$ is degenerate, with the coefficients 
of formula $(\ref{eqn:anti-canonical})$ being the ranks of 
degeneracy of $\phi$. 
In particular, 
%%%%%
\[
\M_z(\k) = \ol{\M}_z(\k) - \Y_z(\k). 
\]
\end{enumerate}
\end{theorem}
%%%%%
This theorem implies that $\M_z(\k)$ is a moduli-theoretical 
realization of the space of initial conditions for $\PVI(\k)$ 
constructed ``by hands" in \cite{Okamoto}, 
$\ol{\M}_z(\k)$ is a generalized Halphen surface of type 
$D_4^{(1)}$ in \cite{Sakai} and $(\ol{\M}_z(\k), \Y_z(\k))$ 
is an Okamoto-Painlev\'e pair of type $\wt{D}_4$ in \cite{STT}. 
%%%%%%%%%%%%%%%%%%%%%%%%%% fig:sps %%%%%%%%%%%%%%%%%%%%%%%%%
%%%%%%%%%%%%%%%%%%%%%%%%%%%%%%%%%%%%%%%%%%%%%%%%%%%%%%%%%%%%%%%%
\par  
Since the Painlev\'e foliation has the geometric 
Painlev\'e property \cite{IIS1}, each loop $\ga \in \pi_1(Z,z)$ 
admits global horizontal lifts along the foliation and 
induces an automorphism 
%%%%%%%%%%%%%%%%%%%%%%%%% eqn:nm %%%%%%%%%%%%%%%%%%%%%%%%%%%%%%%%%
\begin{equation} \label{eqn:nm} 
\ga_* : \M_z(\k) \to \M_z(\k), \quad Q \mapsto Q', 
\end{equation}
%%%%%%%%%%%%%%%%%%%%%%%%%%%%%%%%%%%%%%%%%%%%%%%%%%%%%%%%%%%%%%
called the {\sl nonlinear monodromy} along the loop $\ga$ 
(see Figure \ref{fig:sps}). 
Note that a fixed point or a periodic point of the map 
$\ga_* : \M_z(\k) \carl$ can be identified with a solution 
germ at $z$ which is single-valued or finitely many-valued 
along the loop $\ga$, respectively. 
%%%%%%%%%%%%%%%%%%%%%%%%% sec:RHC %%%%%%%%%%%%%%%%%%%%%%%%%%%%
\section{Riemann-Hilbert Correspondence} \label{sec:RHC} 
\begin{figure}[t] 
\begin{center}
%WinTpicVersion2.15
\unitlength 0.1in
\begin{picture}(27.93,16.78)(2.10,-18.78)
% CIRCLE 0 0 3 0
% 4 483 880 483 607 476 600 203 887
% 
\special{pn 20}%
\special{ar 483 480 273 273  3.1165979 4.6873942}%
% LINE 0 0 3 0
% 2 210 880 210 2000
% 
\special{pn 20}%
\special{pa 210 480}%
\special{pa 210 1600}%
\special{fp}%
% CIRCLE 0 0 3 0
% 4 483 2000 210 2000 210 2000 490 2280
% 
\special{pn 20}%
\special{ar 483 1600 273 273  1.5458015 3.1415927}%
% CIRCLE 0 0 3 0
% 4 1603 1440 1603 1727 1603 1727 1603 1727
% 
\special{pn 20}%
\special{ar 1603 1040 287 287  0.0000000 6.2831853}%
% CIRCLE 0 0 3 0
% 4 2450 1447 2450 1734 2450 1734 2450 1734
% 
\special{pn 20}%
\special{ar 2450 1047 287 287  0.0000000 6.2831853}%
% LINE 0 0 3 0
% 2 490 600 2730 600
% 
\special{pn 20}%
\special{pa 490 200}%
\special{pa 2730 200}%
\special{fp}%
% LINE 0 0 3 0
% 2 483 2273 2730 2273
% 
\special{pn 20}%
\special{pa 483 1873}%
\special{pa 2730 1873}%
\special{fp}%
% CIRCLE 0 0 3 0
% 4 2730 2000 2730 2273 2723 2273 3003 2000
% 
\special{pn 20}%
\special{ar 2730 1600 273 273  6.2831853 6.2831853}%
\special{ar 2730 1600 273 273  0.0000000 1.5964317}%
% LINE 0 0 3 0
% 2 763 600 763 1160
% 
\special{pn 20}%
\special{pa 763 200}%
\special{pa 763 760}%
\special{fp}%
% LINE 0 0 3 0
% 2 1610 600 1610 1153
% 
\special{pn 20}%
\special{pa 1610 200}%
\special{pa 1610 753}%
\special{fp}%
% LINE 0 0 3 0
% 2 2443 600 2443 1153
% 
\special{pn 20}%
\special{pa 2443 200}%
\special{pa 2443 753}%
\special{fp}%
% LINE 0 0 3 0
% 2 3003 887 3003 2007
% 
\special{pn 20}%
\special{pa 3003 487}%
\special{pa 3003 1607}%
\special{fp}%
% CIRCLE 0 0 3 0
% 4 2723 880 2744 607 2996 887 2723 600
% 
\special{pn 20}%
\special{ar 2723 480 274 274  4.7123890 6.2831853}%
\special{ar 2723 480 274 274  0.0000000 0.0256354}%
% CIRCLE 0 0 3 0
% 4 763 1454 763 1741 763 1741 763 1741
% 
\special{pn 20}%
\special{ar 763 1054 287 287  0.0000000 6.2831853}%
% CIRCLE 1 0 0 0
% 4 2450 1447 2471 1482 2471 1482 2471 1482
% 
\special{pn 13}%
\special{sh 0.600}%
\special{ar 2450 1047 41 41  0.0000000 6.2831853}%
% CIRCLE 0 0 0 0
% 4 1610 1440 1631 1475 1631 1475 1631 1475
% 
\special{pn 20}%
\special{sh 0.600}%
\special{ar 1610 1040 41 41  0.0000000 6.2831853}%
% CIRCLE 0 0 0 0
% 4 763 1447 784 1482 784 1482 784 1482
% 
\special{pn 20}%
\special{sh 0.600}%
\special{ar 763 1047 41 41  0.0000000 6.2831853}%
% VECTOR 1 0 3 0
% 2 707 740 707 1020
% 
\special{pn 13}%
\special{pa 707 340}%
\special{pa 707 620}%
\special{fp}%
\special{sh 1}%
\special{pa 707 620}%
\special{pa 727 553}%
\special{pa 707 567}%
\special{pa 687 553}%
\special{pa 707 620}%
\special{fp}%
% VECTOR 1 0 3 0
% 2 1554 733 1554 1013
% 
\special{pn 13}%
\special{pa 1554 333}%
\special{pa 1554 613}%
\special{fp}%
\special{sh 1}%
\special{pa 1554 613}%
\special{pa 1574 546}%
\special{pa 1554 560}%
\special{pa 1534 546}%
\special{pa 1554 613}%
\special{fp}%
% VECTOR 1 0 3 0
% 2 2380 733 2380 1013
% 
\special{pn 13}%
\special{pa 2380 333}%
\special{pa 2380 613}%
\special{fp}%
\special{sh 1}%
\special{pa 2380 613}%
\special{pa 2400 546}%
\special{pa 2380 560}%
\special{pa 2360 546}%
\special{pa 2380 613}%
\special{fp}%
% VECTOR 1 0 3 0
% 2 826 1013 826 740
% 
\special{pn 13}%
\special{pa 826 613}%
\special{pa 826 340}%
\special{fp}%
\special{sh 1}%
\special{pa 826 340}%
\special{pa 806 407}%
\special{pa 826 393}%
\special{pa 846 407}%
\special{pa 826 340}%
\special{fp}%
% VECTOR 1 0 3 0
% 2 1673 1006 1673 733
% 
\special{pn 13}%
\special{pa 1673 606}%
\special{pa 1673 333}%
\special{fp}%
\special{sh 1}%
\special{pa 1673 333}%
\special{pa 1653 400}%
\special{pa 1673 386}%
\special{pa 1693 400}%
\special{pa 1673 333}%
\special{fp}%
% VECTOR 1 0 3 0
% 2 2506 1006 2506 733
% 
\special{pn 13}%
\special{pa 2506 606}%
\special{pa 2506 333}%
\special{fp}%
\special{sh 1}%
\special{pa 2506 333}%
\special{pa 2486 400}%
\special{pa 2506 386}%
\special{pa 2526 400}%
\special{pa 2506 333}%
\special{fp}%
% VECTOR 1 0 3 0
% 2 728 1741 791 1741
% 
\special{pn 13}%
\special{pa 728 1341}%
\special{pa 791 1341}%
\special{fp}%
\special{sh 1}%
\special{pa 791 1341}%
\special{pa 724 1321}%
\special{pa 738 1341}%
\special{pa 724 1361}%
\special{pa 791 1341}%
\special{fp}%
% VECTOR 1 0 3 0
% 2 1568 1727 1638 1720
% 
\special{pn 13}%
\special{pa 1568 1327}%
\special{pa 1638 1320}%
\special{fp}%
\special{sh 1}%
\special{pa 1638 1320}%
\special{pa 1570 1307}%
\special{pa 1585 1325}%
\special{pa 1574 1347}%
\special{pa 1638 1320}%
\special{fp}%
% VECTOR 1 0 3 0
% 2 2422 1734 2478 1741
% 
\special{pn 13}%
\special{pa 2422 1334}%
\special{pa 2478 1341}%
\special{fp}%
\special{sh 1}%
\special{pa 2478 1341}%
\special{pa 2414 1313}%
\special{pa 2425 1334}%
\special{pa 2409 1353}%
\special{pa 2478 1341}%
\special{fp}%
% VECTOR 1 0 3 0
% 2 1673 2273 1610 2273
% 
\special{pn 13}%
\special{pa 1673 1873}%
\special{pa 1610 1873}%
\special{fp}%
\special{sh 1}%
\special{pa 1610 1873}%
\special{pa 1677 1893}%
\special{pa 1663 1873}%
\special{pa 1677 1853}%
\special{pa 1610 1873}%
\special{fp}%
% STR 2 0 3 0
% 3 756 1797 756 1867 5 0
% $C_1$
\put(7.5600,-14.6700){\makebox(0,0){$C_1$}}%
% STR 2 0 3 0
% 3 1638 1811 1638 1881 5 0
% $C_2$
\put(16.3800,-14.8100){\makebox(0,0){$C_2$}}%
% STR 2 0 3 0
% 3 2471 1818 2471 1888 5 0
% $C_3$
\put(24.7100,-14.8800){\makebox(0,0){$C_3$}}%
% STR 2 0 3 0
% 3 1645 2070 1645 2140 5 0
% $C_4$
\put(16.4500,-17.4000){\makebox(0,0){$C_4$}}%
% STR 2 0 3 0
% 3 763 1510 763 1580 5 0
% $0$
\put(7.6300,-11.8000){\makebox(0,0){$0$}}%
% STR 2 0 3 0
% 3 1617 1517 1617 1587 5 0
% $z$
\put(16.1700,-11.8700){\makebox(0,0){$z$}}%
% STR 2 0 3 0
% 3 2471 1510 2471 1580 5 0
% $1$
\put(24.7100,-11.8000){\makebox(0,0){$1$}}%
\end{picture}%
\end{center}
\caption{Four loops in $\P^1 - \{0,z,1,\infty\}$}
\label{fig:sloop}
\end{figure}
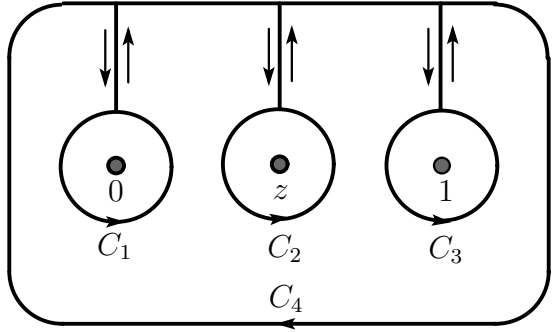

%%%%%%%%%%%%%%%%%%%%%%%%%%%%%%%%%%%%%%%%%%%%%%%%%%%%%%%%%%%%%%
Generally speaking, a Riemann-Hilbert correspondence is the 
map from a moduli space of flat connections to a moduli 
space of monodromy representations, sending a connection 
to its monodromy. 
In our situation an appropriate Riemann-Hilbert correspondence 
%%%%%%%%%%%%%%%%%%%%%%%%%% eqn:RHka %%%%%%%%%%%%%%%%%%%%%%%%%%%
\begin{equation} \label{eqn:RHka}
\RH_{z,\k} : \M_z(\k) \to \R_z(a), \qquad 
Q \mapsto \rho,
\end{equation}
%%%%%%%%%%%%%%%%%%%%%%%%%%%%%%%%%%%%%%%%%%%%%%%%%%%%%%%%%%%%%%
is formulated in \cite{IIS1,IIS2,IIS3}. 
For each $a = (a_1,a_2,a_3,a_4) \in A := \C^4$,  
let $\R_z(a)$ denote the moduli space of Jordan equivalence 
classes of linear monodromy representations 
%%%%%%%%
\[
\rho : \pi_1(\P^1-\{0,z,1,\infty\},*) \to SL_2(\C), 
\]
%%%%%%%
with the prescribed local monodromy data 
$\Tr \, \rho(C_i) = a_i$ $(i = 1,2,3,4)$, where $C_i$ is 
a loop as in Figure \ref{fig:sloop}. 
Any stable parabolic connection $Q \in \M_z(\k)$, restricted 
to $\P^1-\{0,z,1,\infty\}$, induces a flat connection and 
determines the Jordan equivalence class $\rho \in \R_z(a)$ 
of its monodromy representations, where the correspondence 
of parameters $\k \mapsto a$ is described as follows. 
If 
%%%%%%%%%%%%%%%%%%%%%%% eqn:bi %%%%%%%%%%%%%%%%%%%%%%%%%%%%%%
\begin{equation} \label{eqn:bi}
b_i = \left\{ \begin{array}{rl}
\exp(\sqrt{-1} \pi \k_i) \qquad & (i = 0,1,2,3), \\[2mm]
-\exp({\sqrt{-1} \pi \k_4}) \qquad & (i = 4), 
\end{array} \right.
\end{equation}
%%%%%%%%%%%%%%%%%%%%%%%%%%%%%%%%%%%%%%%%%%%%%%%%%%%%%%%%%%%%
then $b = (b_0,b_1,b_2,b_3,b_4)$ belongs to the 
multiplicative space 
%%%%%
\[
B := \{\,b = (b_0,b_1,b_2,b_3,b_4) \in (\C^{\times})^5 
\,:\, b_0^2 b_1b_2b_3b_4 = 1 \, \}. 
\]
%%%%%
The Riemann scheme in Table \ref{tab:riemann} then implies 
that the monodromy matrix $\rho(C_i)$ has an eigenvalue 
$b_i$ for each $i = 1,2,3,4$. 
Since $\rho(C_i) \in SL_2(\C)$, its 
trace $a_i = \Tr \, \rho(C_i)$ is given by 
%%%%%%%%%%%%%%%%%%%%%%%% eqn:ai %%%%%%%%%%%%%%%%%%%%%%%%%%%
\begin{equation}  \label{eqn:ai}
a_i = b_i + b_i^{-1} \qquad (i = 1,2,3,4). 
\end{equation}
%%%%%%%%%%%%%%%%%%%%%%%%%%%%%%%%%%%%%%%%%%%%%%%%%%%%%%%%%%%
\par 
%%%%%%%%%%%%%%%%%%%%%%%%%%% fig:sloop %%%%%%%%%%%%%%%%%%%%%%
%%%%%%%%%%%%%%%%%%%%%%%%%%%%%%%%%%%%%%%%%%%%%%%%%%%%%%%%%%%%
Given any $\th = (\th_1,\th_2,\th_3,\th_4) \in \Th := 
\C^4_{\th}$, consider the affine cubic surface 
%%%%
\[
\Sol(\th) = \{x \in \C^3_x \,:\, 
f(x,\th) := x_1x_2x_3 + x_1^2 + x_2^2 + x_3^2 
- \th_1 x_1 - \th_2 x_2 - \th_3 x_3 + \th_4 = 0 \}.  
\]
%%%%%%%
Then there exists an isomorphism of affine algebraic surfaces 
%%%%%%%
\[
\R_z(a) \to \Sol(\th), \quad 
\rho \mapsto x = (x_1,x_2,x_3), \quad 
\mbox{with} \quad x_i = \Tr\, \rho(C_jC_k)
\]
%%%%%%%
for $\{i,j,k\} = \{1,2,3\}$, where the correspondence of 
parameters $a \mapsto \th$ is given by 
%%%%%%%%%%%%%%%%%%%%%%%%%%% eqn:thi %%%%%%%%%%%%%%%%%%%%%%%
\begin{equation} \label{eqn:thi} 
\theta_i = \left\{
\begin{array}{ll}
a_i a_4 + a_j a_k \qquad & (\{i,j,k\}=\{1,2,3\}), \\[2mm]
a_1 a_2 a_3 a_4 + a_1^2 + a_2^2 + a_3^2 + a_4^2 - 4 
\qquad & (i = 4). 
\end{array}
\right. 
\end{equation} 
%%%%%%%%%%%%%%%%%%%%%%%%%%%%%%%%%%%%%%%%%%%%%%%%%%%%%%%%%%%%%%
The composition of the sequence 
$\k \mapsto b \mapsto a \mapsto \th$ of the three maps 
(\ref{eqn:bi}), (\ref{eqn:ai}) and (\ref{eqn:thi}) 
is referred to as the {\sl Riemann-Hilbert correspondence 
in the parameter level} \cite{IIS1} and is denoted by 
%%%%%%%%%%%%%%%%%%%%%%%%%% eqn:rh %%%%%%%%%%%%%%%%%%%%%%%%%%%%%
\begin{equation} \label{eqn:rh}
\rh : \K \to \Th. 
\end{equation}
%%%%%%%%%%%%%%%%%%%%%%%%%%%%%%%%%%%%%%%%%%%%%%%%%%%%%%%%%%%%%%%%
Then the Riemann-Hilbert correspondence 
(\ref{eqn:RHka}) is reformulated as a holomorphic map 
%%%%%%%%%%%%%%%%%%%%%%%%% eqn:RHkth %%%%%%%%%%%%%%%%%%%%%%%%%%%%
\begin{equation} \label{eqn:RHkth}
\RH_{z,\k} : \M_z(\k) \to \Sol(\th) 
\qquad \mbox{with} \quad \th = \rh(\k). 
\end{equation} 
%%%%%%%%%%%%%%%%%%%%%%%%%%%%%%%%%%%%%%%%%%%%%%%%%%%%%%%%%%%%%%%
\par
\begin{figure}[t]
\begin{center}
%WinTpicVersion2.15
\unitlength 0.1in
\begin{picture}(18.90,10.40)(1.40,-14.40)
% CIRCLE 0 0 0 0
% 4 403 1026 403 1072 403 1072 403 1072
% 
\special{pn 20}%
\special{sh 0.600}%
\special{ar 403 626 46 46  0.0000000 6.2831853}%
% CIRCLE 0 0 0 0
% 4 1171 1018 1171 1065 1171 1065 1171 1065
% 
\special{pn 20}%
\special{sh 0.600}%
\special{ar 1171 618 47 47  0.0000000 6.2831853}%
% CIRCLE 0 0 0 0
% 4 410 1794 410 1840 410 1840 410 1840
% 
\special{pn 20}%
\special{sh 0.600}%
\special{ar 410 1394 46 46  0.0000000 6.2831853}%
% CIRCLE 0 0 0 0
% 4 1178 1794 1178 1840 1178 1840 1178 1840
% 
\special{pn 20}%
\special{sh 0.600}%
\special{ar 1178 1394 46 46  0.0000000 6.2831853}%
% CIRCLE 0 0 0 0
% 4 787 1410 787 1456 787 1456 787 1456
% 
\special{pn 20}%
\special{sh 0.600}%
\special{ar 787 1010 46 46  0.0000000 6.2831853}%
% LINE 0 0 3 0
% 4 433 1072 756 1388 764 1388 764 1388
% 
\special{pn 20}%
\special{pa 433 672}%
\special{pa 756 988}%
\special{fp}%
\special{pa 764 988}%
\special{pa 764 988}%
\special{fp}%
% LINE 0 0 3 0
% 4 817 1449 1140 1763 1140 1763 1140 1763
% 
\special{pn 20}%
\special{pa 817 1049}%
\special{pa 1140 1363}%
\special{fp}%
\special{pa 1140 1363}%
\special{pa 1140 1363}%
\special{fp}%
% LINE 0 0 3 0
% 2 1140 1056 817 1372
% 
\special{pn 20}%
\special{pa 1140 656}%
\special{pa 817 972}%
\special{fp}%
% LINE 0 0 3 0
% 2 748 1449 457 1748
% 
\special{pn 20}%
\special{pa 748 1049}%
\special{pa 457 1348}%
\special{fp}%
% STR 2 0 3 0
% 3 140 893 140 970 2 0
% $w_1$
\put(1.4000,-5.7000){\makebox(0,0)[lb]{$w_1$}}%
% STR 2 0 3 0
% 3 1232 896 1232 972 2 0
% $w_2$
\put(12.3200,-5.7200){\makebox(0,0)[lb]{$w_2$}}%
% STR 2 0 3 0
% 3 176 1889 176 1966 2 0
% $w_3$
\put(1.7600,-15.6600){\makebox(0,0)[lb]{$w_3$}}%
% STR 2 0 3 0
% 3 1208 1901 1208 1978 2 0
% $w_4$
\put(12.0800,-15.7800){\makebox(0,0)[lb]{$w_4$}}%
% STR 2 0 3 0
% 3 692 1229 692 1306 2 0
% $w_0$
\put(6.9200,-9.0600){\makebox(0,0)[lb]{$w_0$}}%
% STR 2 0 3 0
% 3 2030 1786 2030 1886 2 0
% $C = \left(\begin{array}{rrrrr} -2&1&1&1&1\\ 1&-2&0&0&0\\ 1&0&-2&0&0\\ 1&0&0&-2&0\\ 1&0&0&0&-2 \end{array}\right)$
\put(20.3000,-14.8600){\makebox(0,0)[lb]{$C = \left(\begin{array}{rrrrr} -2&1&1&1&1\\ 1&-2&0&0&0\\ 1&0&-2&0&0\\ 1&0&0&-2&0\\ 1&0&0&0&-2 \end{array}\right)$}}%
\end{picture}%
\hspace{4.5cm} $\phantom{a}$
\end{center}
\caption{Dynkin diagram and Cartan matrix of type $D_4^{(1)}$} 
\label{fig:sdynkin} 
\end{figure}
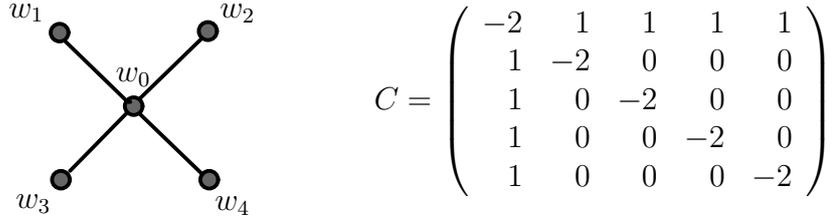

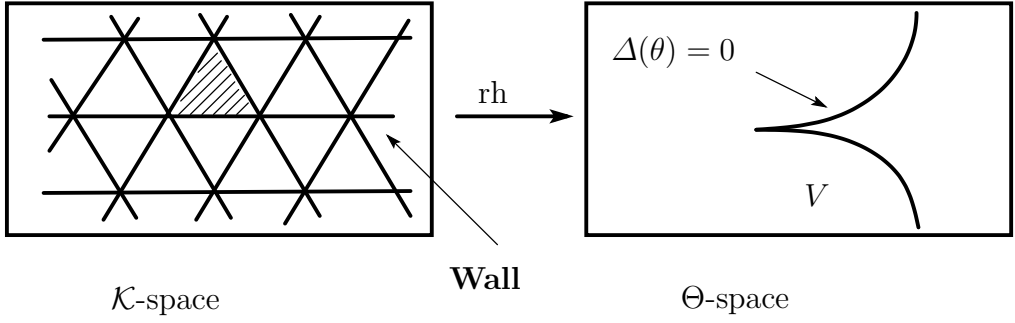
\begin{figure}[t] 
\begin{center}
%WinTpicVersion2.15
\unitlength 0.1in
\begin{picture}(52.00,14.60)(2.00,-16.70)
% BOX 0 0 3 0
% 2 200 610 2400 1820
% 
\special{pn 20}%
\special{pa 200 210}%
\special{pa 2400 210}%
\special{pa 2400 1420}%
\special{pa 200 1420}%
\special{pa 200 210}%
\special{fp}%
% BOX 0 0 3 0
% 2 3200 610 5400 1820
% 
\special{pn 20}%
\special{pa 3200 210}%
\special{pa 5400 210}%
\special{pa 5400 1420}%
\special{pa 3200 1420}%
\special{pa 3200 210}%
\special{fp}%
% LINE 0 0 3 0
% 2 420 1200 2200 1200
% 
\special{pn 20}%
\special{pa 420 800}%
\special{pa 2200 800}%
\special{fp}%
% LINE 0 0 3 0
% 2 420 1380 880 690
% 
\special{pn 20}%
\special{pa 420 980}%
\special{pa 880 290}%
\special{fp}%
% LINE 0 0 3 0
% 2 740 680 1360 1730
% 
\special{pn 20}%
\special{pa 740 280}%
\special{pa 1360 1330}%
\special{fp}%
% LINE 0 0 3 0
% 2 1340 680 700 1740
% 
\special{pn 20}%
\special{pa 1340 280}%
\special{pa 700 1340}%
\special{fp}%
% LINE 0 0 3 0
% 2 1200 670 1820 1720
% 
\special{pn 20}%
\special{pa 1200 270}%
\special{pa 1820 1320}%
\special{fp}%
% LINE 0 0 3 0
% 2 1820 680 1180 1740
% 
\special{pn 20}%
\special{pa 1820 280}%
\special{pa 1180 1340}%
\special{fp}%
% LINE 0 0 3 0
% 2 1670 670 2290 1720
% 
\special{pn 20}%
\special{pa 1670 270}%
\special{pa 2290 1320}%
\special{fp}%
% LINE 0 0 3 0
% 2 2280 700 1640 1760
% 
\special{pn 20}%
\special{pa 2280 300}%
\special{pa 1640 1360}%
\special{fp}%
% LINE 0 0 3 0
% 2 870 1730 430 1010
% 
\special{pn 20}%
\special{pa 870 1330}%
\special{pa 430 610}%
\special{fp}%
% LINE 0 0 3 0
% 2 390 800 2290 790
% 
\special{pn 20}%
\special{pa 390 400}%
\special{pa 2290 390}%
\special{fp}%
% LINE 0 0 3 0
% 2 390 1600 2290 1590
% 
\special{pn 20}%
\special{pa 390 1200}%
\special{pa 2290 1190}%
\special{fp}%
% LINE 2 0 3 0
% 16 1380 1020 1210 1190 1360 980 1150 1190 1340 940 1090 1190 1310 910 1100 1120 1290 870 1190 970 1400 1060 1270 1190 1430 1090 1330 1190 1450 1130 1390 1190
% 
\special{pn 8}%
\special{pa 1380 620}%
\special{pa 1210 790}%
\special{fp}%
\special{pa 1360 580}%
\special{pa 1150 790}%
\special{fp}%
\special{pa 1340 540}%
\special{pa 1090 790}%
\special{fp}%
\special{pa 1310 510}%
\special{pa 1100 720}%
\special{fp}%
\special{pa 1290 470}%
\special{pa 1190 570}%
\special{fp}%
\special{pa 1400 660}%
\special{pa 1270 790}%
\special{fp}%
\special{pa 1430 690}%
\special{pa 1330 790}%
\special{fp}%
\special{pa 1450 730}%
\special{pa 1390 790}%
\special{fp}%
% STR 2 0 3 0
% 3 2640 1030 2640 1130 2 0
% $\rh$
\put(26.4000,-7.3000){\makebox(0,0)[lb]{$\rh$}}%
% STR 2 0 3 0
% 3 3330 850 3330 950 2 0
% $\varDelta(\theta) = 0$
\put(33.3000,-5.5000){\makebox(0,0)[lb]{$\varDelta(\theta) = 0$}}%
% STR 2 0 3 0
% 3 740 2140 740 2240 2 0
% $\mathcal{K}$-space
\put(7.4000,-18.4000){\makebox(0,0)[lb]{$\mathcal{K}$-space}}%
% STR 2 0 3 0
% 3 3690 2130 3690 2230 2 0
% $\Theta$-space
\put(36.9000,-18.3000){\makebox(0,0)[lb]{$\Theta$-space}}%
% STR 2 0 3 0
% 3 2490 2000 2490 2100 2 0
% $\mathbf{Wall}$
\put(24.9000,-17.0000){\makebox(0,0)[lb]{$\mathbf{Wall}$}}%
% VECTOR 2 0 3 0
% 2 4080 970 4440 1150
% 
\special{pn 8}%
\special{pa 4080 570}%
\special{pa 4440 750}%
\special{fp}%
\special{sh 1}%
\special{pa 4440 750}%
\special{pa 4389 702}%
\special{pa 4392 726}%
\special{pa 4371 738}%
\special{pa 4440 750}%
\special{fp}%
% SPLINE 0 0 3 0
% 6 4080 1270 4500 1200 4800 980 4910 670 4910 660 4910 660
% 
\special{pn 20}%
\special{pa 4080 870}%
\special{pa 4112 868}%
\special{pa 4145 865}%
\special{pa 4177 863}%
\special{pa 4209 860}%
\special{pa 4241 857}%
\special{pa 4273 853}%
\special{pa 4305 849}%
\special{pa 4337 843}%
\special{pa 4368 837}%
\special{pa 4399 831}%
\special{pa 4430 823}%
\special{pa 4461 814}%
\special{pa 4491 803}%
\special{pa 4521 792}%
\special{pa 4550 779}%
\special{pa 4579 765}%
\special{pa 4607 749}%
\special{pa 4635 732}%
\special{pa 4662 714}%
\special{pa 4688 694}%
\special{pa 4713 674}%
\special{pa 4737 652}%
\special{pa 4760 628}%
\special{pa 4781 604}%
\special{pa 4801 578}%
\special{pa 4820 551}%
\special{pa 4837 523}%
\special{pa 4853 495}%
\special{pa 4867 465}%
\special{pa 4879 435}%
\special{pa 4889 404}%
\special{pa 4898 372}%
\special{pa 4904 341}%
\special{pa 4908 309}%
\special{pa 4910 277}%
\special{pa 4910 260}%
\special{sp}%
% SPLINE 0 0 3 0
% 5 4080 1270 4590 1340 4840 1540 4920 1780 4920 1780
% 
\special{pn 20}%
\special{pa 4080 870}%
\special{pa 4113 870}%
\special{pa 4145 871}%
\special{pa 4178 871}%
\special{pa 4210 872}%
\special{pa 4242 873}%
\special{pa 4275 875}%
\special{pa 4307 877}%
\special{pa 4339 880}%
\special{pa 4370 883}%
\special{pa 4402 888}%
\special{pa 4434 894}%
\special{pa 4465 900}%
\special{pa 4496 908}%
\special{pa 4527 917}%
\special{pa 4557 927}%
\special{pa 4587 939}%
\special{pa 4617 952}%
\special{pa 4647 967}%
\special{pa 4675 983}%
\special{pa 4703 1001}%
\special{pa 4730 1020}%
\special{pa 4755 1041}%
\special{pa 4779 1063}%
\special{pa 4801 1087}%
\special{pa 4821 1111}%
\special{pa 4838 1138}%
\special{pa 4854 1165}%
\special{pa 4867 1194}%
\special{pa 4879 1223}%
\special{pa 4889 1254}%
\special{pa 4898 1285}%
\special{pa 4906 1316}%
\special{pa 4913 1348}%
\special{pa 4920 1380}%
\special{sp}%
% VECTOR 2 0 3 0
% 2 2740 1870 2170 1300
% 
\special{pn 8}%
\special{pa 2740 1470}%
\special{pa 2170 900}%
\special{fp}%
\special{sh 1}%
\special{pa 2170 900}%
\special{pa 2203 961}%
\special{pa 2208 938}%
\special{pa 2231 933}%
\special{pa 2170 900}%
\special{fp}%
% VECTOR 0 0 3 0
% 2 2530 1200 3100 1200
% 
\special{pn 20}%
\special{pa 2530 800}%
\special{pa 3100 800}%
\special{fp}%
\special{sh 1}%
\special{pa 3100 800}%
\special{pa 3033 780}%
\special{pa 3047 800}%
\special{pa 3033 820}%
\special{pa 3100 800}%
\special{fp}%
% STR 2 0 3 0
% 3 4330 1570 4330 1670 2 0
% $V$
\put(43.3000,-12.7000){\makebox(0,0)[lb]{$V$}}%
\end{picture}%
\end{center}
\caption{Riemann-Hilbert correspondence in the parameter level} 
\label{fig:sparam} 
\end{figure}

 
The map (\ref{eqn:rh}) admits a remarkable affine Weyl group 
structure \cite{IIS1,Iwasaki}, from which the 
B\"acklund transformations of Painlev\'e VI emerge \cite{IIS0}. 
In view of formula (\ref{eqn:K}) the affine space $\mathcal{K}$ 
can be identified with the linear space $\mathbb{C}^4$ by the 
forgetful isomorphism $\mathcal{K} \to \mathbb{C}^4$, $\kappa = 
(\kappa_0,\kappa_1,\kappa_2,\kappa_3,\kappa_4) 
\mapsto (\kappa_1,\kappa_2,\kappa_3,\kappa_4)$, 
where the latter space $\mathbb{C}^4$ is equipped with the 
standard (complex) Euclidean inner product. 
For each $i \in \{0,1,2,3,4\}$, let $w_i : \mathcal{K} \to 
\mathcal{K}$, $\k \mapsto \k'$, be the orthogonal reflection 
in the hyperplane 
$\{\, \kappa \in \mathcal{K}\,:\, \kappa_i =0\}$, which is 
explicitly represented as 
%%%%%%%%%%%%%%%%%%%%%%%%%%%% eqn:cartan %%%%%%%%%%%%%%%%%%%%%
\begin{equation} \label{eqn:cartan}
\k'_j = \k_j + \k_i c_{ij} \qquad (i,j \in \{0,1,2,3,4\}), 
\end{equation}
%%%%%%%%%%%%%%%%%%%%%%%%%%%%%%%%%%%%%%%%%%%%%%%%%%%%%%%%%%%%%
where $C = (c_{ij})$ is the Cartan matrix of type $D_4^{(1)}$ 
given in Figure \ref{fig:sdynkin}. 
Then the group generated by $w_0$, $w_1$, $w_2$, $w_3$, 
$w_4$ is an affine Weyl group of type $D_4^{(1)}$, 
%%%%%%
\[
W(D_4^{(1)}) = \langle w_0, w_1, w_2, w_3, w_4 \rangle 
\curvearrowright \mathcal{K}.   
\]
%%%%%%
corresponding to the Dynkin diagram in 
Figure \ref{fig:sdynkin}. 
%%%%%%%%%%%%%%%%%%%%%%%%% fig:sdynkin %%%%%%%%%%%%%%%%%%%%%%%
%%%%%%%%%%%%%%%%%%%%%%%%%%%%%%%%%%%%%%%%%%%%%%%%%%%%%%%%%%%%%%%
The reflecting hyperplanes of all reflections in the 
group $W(D_4^{(1)})$ are given by affine linear relations 
%%%%%%%%
\[
\kappa_i = m, \qquad \kappa_1 \pm \kappa_2 \pm \kappa_3 
\pm \kappa_4 = 2m+1 
\qquad (i \in \{1,2,3,4\}, \, m \in \mathbb{Z}), 
\]
%%%%%%%%
where the signs $\pm$ may be chosen arbitrarily. 
Let $\mathbf{Wall}$ be the union of all these hyperplanes. 
Then the affine Weyl group structure on (\ref{eqn:rh}) is 
stated as follows \cite{IIS1} (see Figure \ref{fig:sparam}).  
%%%%%%%%%%%%%%%%%%%%%%% lem:discriminant %%%%%%%%%%%%%%%%%%%%%%
\begin{lemma} \label{lem:discriminant} 
In terms of $b \in B$, the discriminant $\varDelta(\theta)$ 
of the cubic surfaces $\mathcal{S}(\theta)$ factors as 
%%%%%%%%%%%%%%%%%%%%%%%%%% eqn:vD %%%%%%%%%%%%%%%%%%%%%%%%%%%%%%%
\begin{equation} \label{eqn:vD}
\varDelta(\theta) = \prod_{l=1}^4(b_l-b_l^{-1})^2 
\prod_{\varepsilon \in \{\pm1\}^4}(b^{\varepsilon}-1), 
\end{equation}
%%%%%%%%%%%%%%%%%%%%%%%%%%%%%%%%%%%%%%%%%%%%%%%%%%%%%%%%%%%%%%%%%
where we put 
$b^{\varepsilon} = b_1^{\varepsilon_1}b_2^{\varepsilon_2}
b_3^{\varepsilon_3}b_4^{\varepsilon_4}$ 
for each quadruple sign 
$\varepsilon = 
(\varepsilon_1,\varepsilon_2,\varepsilon_3,\varepsilon_4) 
\in \{\pm1\}^4$. 
The Riemann-Hilbert correspondence in the parameter level 
$(\ref{eqn:rh})$ is a branched $W(D_4^{(1)})$-covering 
ramifying along $\mathbf{Wall}$ and mapping $\mathbf{Wall}$ 
onto the discriminant locus $\varDelta(\theta) = 0$ 
in $\Theta$. 
\end{lemma} 
%%%%%%%%%%%%%%%%%%%%%%%%%%%%%%%%%%%%%%%%%%%%%%%%%%%%%%%%%%%%%%%
%%%%%%%%%%%%%%%%%%%%%%%%% fig:sparam %%%%%%%%%%%%%%%%%%%%%%%%%%
%%%%%%%%%%%%%%%%%%%%%%%%%%%%%%%%%%%%%%%%%%%%%%%%%%%%%%%%%%%%%%%
\par
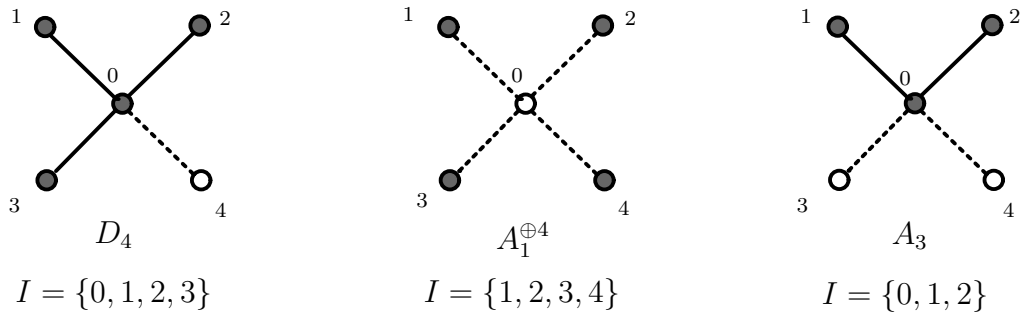
\begin{figure}[b]
\begin{center}
%WinTpicVersion2.15
\unitlength 0.1in
\begin{picture}(51.80,15.20)(2.60,-17.70)
% CIRCLE 0 0 0 0
% 4 448 852 448 900 448 900 448 900
% 
\special{pn 20}%
\special{sh 0.600}%
\special{ar 448 452 48 48  0.0000000 6.2831853}%
% CIRCLE 0 0 0 0
% 4 1248 844 1248 892 1248 892 1248 892
% 
\special{pn 20}%
\special{sh 0.600}%
\special{ar 1248 444 48 48  0.0000000 6.2831853}%
% CIRCLE 0 0 0 0
% 4 456 1652 456 1700 456 1700 456 1700
% 
\special{pn 20}%
\special{sh 0.600}%
\special{ar 456 1252 48 48  0.0000000 6.2831853}%
% CIRCLE 0 0 3 0
% 4 1256 1652 1256 1700 1256 1700 1256 1700
% 
\special{pn 20}%
\special{ar 1256 1252 48 48  0.0000000 6.2831853}%
% CIRCLE 0 0 0 0
% 4 848 1252 848 1300 848 1300 848 1300
% 
\special{pn 20}%
\special{sh 0.600}%
\special{ar 848 852 48 48  0.0000000 6.2831853}%
% LINE 0 0 3 0
% 4 480 900 816 1228 824 1228 824 1228
% 
\special{pn 20}%
\special{pa 480 500}%
\special{pa 816 828}%
\special{fp}%
\special{pa 824 828}%
\special{pa 824 828}%
\special{fp}%
% LINE 0 2 3 0
% 4 880 1292 1216 1620 1216 1620 1216 1620
% 
\special{pn 20}%
\special{pa 880 892}%
\special{pa 1216 1220}%
\special{dt 0.054}%
\special{pa 1216 1220}%
\special{pa 1216 1220}%
\special{dt 0.054}%
\special{pa 1216 1220}%
\special{pa 1216 1220}%
\special{dt 0.054}%
% LINE 0 0 3 0
% 2 1216 884 880 1212
% 
\special{pn 20}%
\special{pa 1216 484}%
\special{pa 880 812}%
\special{fp}%
% LINE 0 0 3 0
% 2 808 1292 504 1604
% 
\special{pn 20}%
\special{pa 808 892}%
\special{pa 504 1204}%
\special{fp}%
% STR 2 0 3 0
% 3 270 750 270 830 2 0
% $\scriptstyle{1}$
\put(2.7000,-4.3000){\makebox(0,0)[lb]{$\scriptstyle{1}$}}%
% STR 2 0 3 0
% 3 1350 760 1350 840 2 0
% $\scriptstyle{2}$
\put(13.5000,-4.4000){\makebox(0,0)[lb]{$\scriptstyle{2}$}}%
% STR 2 0 3 0
% 3 260 1740 260 1820 2 0
% $\scriptstyle{3}$
\put(2.6000,-14.2000){\makebox(0,0)[lb]{$\scriptstyle{3}$}}%
% STR 2 0 3 0
% 3 1328 1764 1328 1844 2 0
% $\scriptstyle{4}$
\put(13.2800,-14.4400){\makebox(0,0)[lb]{$\scriptstyle{4}$}}%
% STR 2 0 3 0
% 3 768 1060 768 1140 2 0
% $\scriptstyle{0}$
\put(7.6800,-7.4000){\makebox(0,0)[lb]{$\scriptstyle{0}$}}%
% CIRCLE 0 0 0 0
% 4 2536 852 2536 900 2536 900 2536 900
% 
\special{pn 20}%
\special{sh 0.600}%
\special{ar 2536 452 48 48  0.0000000 6.2831853}%
% CIRCLE 0 0 0 0
% 4 3336 844 3336 892 3336 892 3336 892
% 
\special{pn 20}%
\special{sh 0.600}%
\special{ar 3336 444 48 48  0.0000000 6.2831853}%
% CIRCLE 0 0 0 0
% 4 2544 1652 2544 1700 2544 1700 2544 1700
% 
\special{pn 20}%
\special{sh 0.600}%
\special{ar 2544 1252 48 48  0.0000000 6.2831853}%
% CIRCLE 0 0 0 0
% 4 3344 1652 3344 1700 3344 1700 3344 1700
% 
\special{pn 20}%
\special{sh 0.600}%
\special{ar 3344 1252 48 48  0.0000000 6.2831853}%
% CIRCLE 0 0 3 0
% 4 2936 1252 2936 1300 2936 1300 2936 1300
% 
\special{pn 20}%
\special{ar 2936 852 48 48  0.0000000 6.2831853}%
% LINE 0 2 3 0
% 4 2568 900 2904 1228 2912 1228 2912 1228
% 
\special{pn 20}%
\special{pa 2568 500}%
\special{pa 2904 828}%
\special{dt 0.054}%
\special{pa 2904 828}%
\special{pa 2904 828}%
\special{dt 0.054}%
\special{pa 2912 828}%
\special{pa 2912 828}%
\special{dt 0.054}%
% LINE 0 2 3 0
% 4 2968 1292 3304 1620 3304 1620 3304 1620
% 
\special{pn 20}%
\special{pa 2968 892}%
\special{pa 3304 1220}%
\special{dt 0.054}%
\special{pa 3304 1220}%
\special{pa 3304 1220}%
\special{dt 0.054}%
\special{pa 3304 1220}%
\special{pa 3304 1220}%
\special{dt 0.054}%
% LINE 0 2 3 0
% 2 3304 884 2968 1212
% 
\special{pn 20}%
\special{pa 3304 484}%
\special{pa 2968 812}%
\special{dt 0.054}%
\special{pa 2968 812}%
\special{pa 2968 812}%
\special{dt 0.054}%
% LINE 0 2 3 0
% 2 2896 1292 2592 1604
% 
\special{pn 20}%
\special{pa 2896 892}%
\special{pa 2592 1204}%
\special{dt 0.054}%
\special{pa 2592 1204}%
\special{pa 2592 1204}%
\special{dt 0.054}%
% STR 2 0 3 0
% 3 2300 740 2300 820 2 0
% $\scriptstyle{1}$
\put(23.0000,-4.2000){\makebox(0,0)[lb]{$\scriptstyle{1}$}}%
% STR 2 0 3 0
% 3 3450 760 3450 840 2 0
% $\scriptstyle{2}$
\put(34.5000,-4.4000){\makebox(0,0)[lb]{$\scriptstyle{2}$}}%
% STR 2 0 3 0
% 3 2370 1710 2370 1790 2 0
% $\scriptstyle{3}$
\put(23.7000,-13.9000){\makebox(0,0)[lb]{$\scriptstyle{3}$}}%
% STR 2 0 3 0
% 3 3416 1756 3416 1836 2 0
% $\scriptstyle{4}$
\put(34.1600,-14.3600){\makebox(0,0)[lb]{$\scriptstyle{4}$}}%
% STR 2 0 3 0
% 3 2864 1060 2864 1140 2 0
% $\scriptstyle{0}$
\put(28.6400,-7.4000){\makebox(0,0)[lb]{$\scriptstyle{0}$}}%
% CIRCLE 0 0 0 0
% 4 4552 850 4552 898 4552 898 4552 898
% 
\special{pn 20}%
\special{sh 0.600}%
\special{ar 4552 450 48 48  0.0000000 6.2831853}%
% CIRCLE 0 0 0 0
% 4 5352 842 5352 890 5352 890 5352 890
% 
\special{pn 20}%
\special{sh 0.600}%
\special{ar 5352 442 48 48  0.0000000 6.2831853}%
% CIRCLE 0 0 3 0
% 4 4560 1650 4560 1698 4560 1698 4560 1698
% 
\special{pn 20}%
\special{ar 4560 1250 48 48  0.0000000 6.2831853}%
% CIRCLE 0 0 3 0
% 4 5360 1650 5360 1698 5360 1698 5360 1698
% 
\special{pn 20}%
\special{ar 5360 1250 48 48  0.0000000 6.2831853}%
% CIRCLE 0 0 0 0
% 4 4952 1250 4952 1298 4952 1298 4952 1298
% 
\special{pn 20}%
\special{sh 0.600}%
\special{ar 4952 850 48 48  0.0000000 6.2831853}%
% LINE 0 0 3 0
% 4 4584 898 4920 1226 4928 1226 4928 1226
% 
\special{pn 20}%
\special{pa 4584 498}%
\special{pa 4920 826}%
\special{fp}%
\special{pa 4928 826}%
\special{pa 4928 826}%
\special{fp}%
% LINE 0 2 3 0
% 4 4984 1290 5320 1618 5320 1618 5320 1618
% 
\special{pn 20}%
\special{pa 4984 890}%
\special{pa 5320 1218}%
\special{dt 0.054}%
\special{pa 5320 1218}%
\special{pa 5320 1218}%
\special{dt 0.054}%
\special{pa 5320 1218}%
\special{pa 5320 1218}%
\special{dt 0.054}%
% LINE 0 0 3 0
% 2 5320 882 4984 1210
% 
\special{pn 20}%
\special{pa 5320 482}%
\special{pa 4984 810}%
\special{fp}%
% LINE 0 2 3 0
% 2 4912 1290 4608 1602
% 
\special{pn 20}%
\special{pa 4912 890}%
\special{pa 4608 1202}%
\special{dt 0.054}%
\special{pa 4608 1202}%
\special{pa 4608 1202}%
\special{dt 0.054}%
% STR 2 0 3 0
% 3 4360 750 4360 830 2 0
% $\scriptstyle{1}$
\put(43.6000,-4.3000){\makebox(0,0)[lb]{$\scriptstyle{1}$}}%
% STR 2 0 3 0
% 3 5440 750 5440 830 2 0
% $\scriptstyle{2}$
\put(54.4000,-4.3000){\makebox(0,0)[lb]{$\scriptstyle{2}$}}%
% STR 2 0 3 0
% 3 4340 1740 4340 1820 2 0
% $\scriptstyle{3}$
\put(43.4000,-14.2000){\makebox(0,0)[lb]{$\scriptstyle{3}$}}%
% STR 2 0 3 0
% 3 5408 1762 5408 1842 2 0
% $\scriptstyle{4}$
\put(54.0800,-14.4200){\makebox(0,0)[lb]{$\scriptstyle{4}$}}%
% STR 2 0 3 0
% 3 4872 1074 4872 1154 2 0
% $\scriptstyle{0}$
\put(48.7200,-7.5400){\makebox(0,0)[lb]{$\scriptstyle{0}$}}%
% STR 2 0 3 0
% 3 704 1924 704 2004 2 0
% $D_4$
\put(7.0400,-16.0400){\makebox(0,0)[lb]{$D_4$}}%
% STR 2 0 3 0
% 3 2784 1972 2784 2052 2 0
% $A_1^{\oplus 4}$
\put(27.8400,-16.5200){\makebox(0,0)[lb]{$A_1^{\oplus 4}$}}%
% STR 2 0 3 0
% 3 4832 1930 4832 2010 2 0
% $A_3$
\put(48.3200,-16.1000){\makebox(0,0)[lb]{$A_3$}}%
% STR 2 0 3 0
% 3 296 2252 296 2332 2 0
% $I = \{0,1,2,3\}$
\put(2.9600,-19.3200){\makebox(0,0)[lb]{$I = \{0,1,2,3\}$}}%
% STR 2 0 3 0
% 3 2408 2252 2408 2332 2 0
% $I = \{1,2,3,4\}$
\put(24.0800,-19.3200){\makebox(0,0)[lb]{$I = \{1,2,3,4\}$}}%
% STR 2 0 3 0
% 3 4470 2260 4470 2340 2 0
% $I = \{0,1,2\}$
\put(44.7000,-19.4000){\makebox(0,0)[lb]{$I = \{0,1,2\}$}}%
\end{picture}%
\end{center}
\caption{Examples of $W(D_4^{(1)})$-strata} 
\label{fig:sstrata}
\vspace{0.5cm}
\end{figure}

 
The singularity structure of the cubic surfaces $\Sol(\th)$ 
can be described in terms of the stratification of 
$\K$ by proper Dynkin subdiagrams, which we now define. 
%%%%%%%%%%%%%%%%%%% def:stratification %%%%%%%%%%%%%%%%%%%%%%%%
\begin{definition} \label{def:stratification} 
Let $\mathcal{I}$ be the set of all {\sl proper} subsets of 
$\{0,1,2,3,4\}$ including the empty set $\emptyset$. 
For each element $I \in \mathcal{I}$, we put 
%%%%%%%%
\[
\begin{array}{rcl}
\ol{\K}_I &=& \mbox{the $W(D_4^{(1)})$-translates of the set} \,\, 
\{\, \k \in \K \,:\, \k_i = 0 \,\, (i \in I)\,\}, \\[2mm] 
D_I &=& \mbox{the Dynkin subdiagram of $D_4^{(1)}$ that has nodes 
$\bullet$ exactly in $I$}. 
\end{array}
\]
%%%%%%%%
Let $\K_I$ be the set obtained from $\ol{\K}_I$ by removing 
the sets $\ol{\K}_J$ with $\# J = \# I + 1$. 
Then it turns out that we have either $\K_I = \K_{I'}$ or 
$\K_I \cap \K_{I'} = \emptyset$ for any distinct subsets $I$, 
$I' \in \mathcal{I}$ (see Remark \ref{rem:stratification}). 
So we can think of the stratification of $\K$ by the subsets 
$\K_I$ $(I \in \mathcal{I})$, called the 
$W(D_4^{(1)})$-{\sl stratification}, where each $\K_I$ 
is referred to as a $W(D_4^{(1)})$-{\sl stratum}. 
For example, if $I = \emptyset$, one has the {\sl big open} 
$\K_{\emptyset} = \K - \Wall$. 
Other examples of $W(D_4^{()1})$-strata are given in 
Figure \ref{fig:sstrata}. 
The diagram $D_I$ encodes not only its underlying 
{\sl abstract} Dynkin type but also the inclusion pattern 
$D_I \hookrightarrow D_4^{(1)}$, a kind of marking. 
The abstact Dynkin type of $D_I$ is denoted by 
$\mathrm{Dynk}(I)$. 
All the feasible abstract Dynkin types are tabulated in 
Table \ref{tab:type}. 
\end{definition}
%%%%%%%%%%%%%%%%%%%%%%%%%%%%%%%%%%%%%%%%%%%%%%%%%%%%%%%%%%%%%%
\par 
There is a mistake in the definition of $\K_I$ in 
\cite[Definition 9.3]{IIS1} and \cite{IU2}, which is now 
corrected in Definition \ref{def:stratification}. 
(As for \cite{IIS1}, correction may be possible 
before it is published.) 
%%%%%%%%%%%%%%%%%%%%%%%%% fig:sstrata %%%%%%%%%%%%%%%%%%%%%%%%
%%%%%%%%%%%%%%%%%%%%%%%%%%%%%%%%%%%%%%%%%%%%%%%%%%%%%%%%%%%%%%
%%%%%%%%%%%%%%%%%%%%%%%%% rem:stratification %%%%%%%%%%%%%%%%%%
\begin{remark} \label{rem:stratification} 
Let $I$ and $I'$ be distinct elements of $\mathcal{I}$. 
If $\mathrm{Dynk}(I) \neq \mathrm{Dynk}(I')$, then 
$\K_I \cap \K_{I'} = \emptyset$. 
On the other hand, if $\K_I = \K_{I'}$ then 
$\mathrm{Dynk}(I) = \mathrm{Dynk}(I')$ must be of abstract 
type $A_1$ or $A_2$. 
\begin{enumerate} 
\item There is a unique $W(D_4^{(1)})$-stratum of abstract 
type $\emptyset$, or $A_1$, or $A_2$, or $A_1^{\oplus 4}$. 
\item There are six $W(D_4^{(1)})$-strata of abstract 
type $A_1^{\oplus 2}$ or $A_3$. 
\item There are four $W(D_4^{(1)})$-strata of abstract 
type $A_1^{\oplus 3}$ or $D_4$. 
\end{enumerate}
\end{remark} 
%%%%%%%%%%%%%%%%%%%%%%%%%%%%%%%%%%%%%%%%%%%%%%%%%%%%%%%%%%%%%%%
%%%%%%%%%%%%%%%%%%%%%%%%%%% tab:type %%%%%%%%%%%%%%%%%%%%%%%%%
\begin{table}[t]
\begin{center}
\begin{tabular}{|c||c|c|c|c|c|}
\hline
\vspace{-4mm} &    &    &    &    &  \\
number of nodes & 4 & 3 & 2 & 1 & 0 \\[1mm]
\hline
\vspace{-4mm} &    &    &    &    &  \\
abstract & $D_4$ & $A_3$ & $A_2$ & $A_1$ & 
$\emptyset$\\[1mm] 
\cline{2-6}
\vspace{-3mm} &    &    &    &    &  \\ 
Dynkin type & $A_1^{\oplus 4}$ & $A_1^{\oplus 3}$ & 
$A_1^{\oplus 2}$ & $-$ & $-$\\[1mm] 
\hline 
\end{tabular}
\end{center}
\vspace{5mm}
\caption{Feasible abstract Dynkin types} 
\label{tab:type}
\end{table}
%%%%%%%%%%%%%%%%%%%%%%%%%%%%%%%%%%%%%%%%%%%%%%%%%%%%%%%%%%%%%%%
%%%%%%%%%%%%%%%%%%%%%%%%%% ex:sing %%%%%%%%%%%%%%%%%%%%%%%%%%
\begin{example} \label{ex:sing} 
We consider the $W(D_4^{(1)})$-strata of abstract types 
$A_1^{\oplus 4}$ and $D_4$. 
%%%%%%%
\begin{enumerate}
\item The unique $W(D_4^{(1)})$-stratum of abstract type 
$A_1^{\oplus 4}$ exactly corresponds to the value 
$\th = (0,0,0,-4)$. 
A parameter $\k \in \K$ lies in this stratum if and only 
if either 
%%%%%%
\begin{enumerate}
\item $\k_1$, $\k_2$, $\k_3$, $\k_4 \in \Z$,  
$\k_1+\k_2+\k_3+\k_4 \in 2\Z$; or 
\item $\k_1$, $\k_2$, $\k_3$, $\k_4 \in \Z+1/2$. 
\end{enumerate}
%%%%%%
\item The four $W(D_4^{(1)})$-strata of abstract type $D_4$ 
exactly correspond to the values 
$\th = (8\ve_1,8\ve_2,8\ve_3,28)$, where 
$\ve = (\ve_1,\ve_2,\ve_3) \in \{\pm\}^3$ ranges over all 
triple signs such that $\ve_1\ve_2\ve_3 = 1$. 
A parameter $\k \in \K$ lies in the union of these 
$W(D_4^{(1)})$-strata if and only if 
%%%%%
\[
\k_1, \k_2, \k_3, \k_4 \in \Z, \quad 
\k_1+\k_2+\k_3+\k_4 \in 2\Z+1. 
\]
%%%%%%
\end{enumerate}
\end{example}
%%%%%%%%%%%%%%%%%%%%%%%%%%%%%%%%%%%%%%%%%%%%%%%%%%%%%%%%%%%%%
\par 
With this stratification, we have a very neat solution 
to the Riemann-Hilbert problem. 
%%%%%%%%%%%%%%%%%%%%%%%%% thm:SolRHP %%%%%%%%%%%%%%%%%%%%%%%%%%
\begin{theorem}[\cite{IIS1,IIS2,IIS3}]  \label{thm:SolRHP} 
Given any $\k \in \K$, put $\th = \rh(\k) \in \Th$. 
Then, 
%%%%%%%
\begin{enumerate} 
\item if $\k \in \K_I$ then $\Sol(\th)$ has Kleinian 
singularities of Dynkin type $D_I$, 
\item the Riemann-Hilbert correspondence $(\ref{eqn:RHkth})$ 
is a proper surjective map that is an analytic minimal 
resolution of Kleinian singularities. 
\end{enumerate}
\end{theorem}
%%%%%%%%%%%%%%%%%%%%%%%%%%%%%%%%%%%%%%%%%%%%%%%%%%%%%%%%%%%%%%
\par
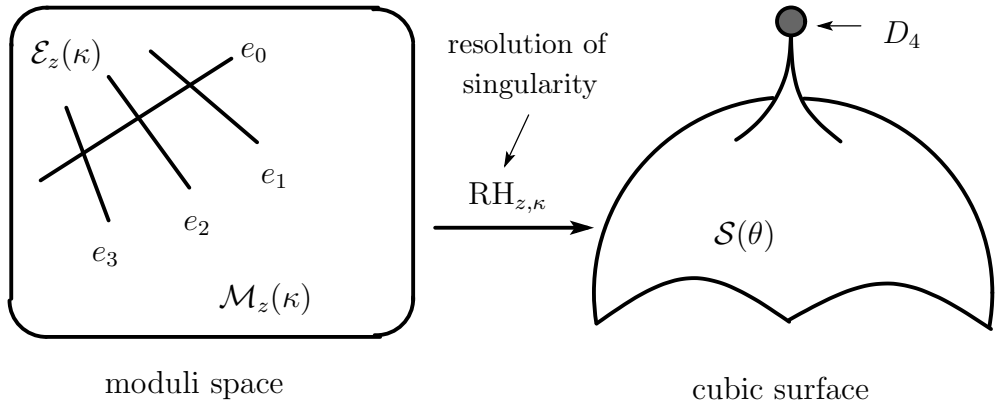
\begin{figure}[t] 
\begin{center} 
%WinTpicVersion2.15
\unitlength 0.1in
\begin{picture}(50.89,18.73)(2.00,-20.72)
% LINE 0 0 3 0
% 2 390 600 2109 600
% 
\special{pn 20}%
\special{pa 390 200}%
\special{pa 2109 200}%
\special{fp}%
% LINE 0 0 3 0
% 2 390 2319 2119 2319
% 
\special{pn 20}%
\special{pa 390 1919}%
\special{pa 2119 1919}%
\special{fp}%
% LINE 0 0 3 0
% 2 209 790 209 2120
% 
\special{pn 20}%
\special{pa 209 390}%
\special{pa 209 1720}%
\special{fp}%
% LINE 0 0 3 0
% 2 2290 790 2290 2129
% 
\special{pn 20}%
\special{pa 2290 390}%
\special{pa 2290 1729}%
\special{fp}%
% CIRCLE 2 0 3 0
% 4 390 790 390 790 209 780 209 790
% 
\special{pn 8}%
\special{ar 390 390 0 0  3.1415927 3.1967852}%
% CIRCLE 0 0 3 0
% 4 390 780 399 600 399 600 200 799
% 
\special{pn 20}%
\special{ar 390 380 180 180  3.0419240 4.7623474}%
% CIRCLE 0 0 3 0
% 4 390 2129 200 2129 200 2129 390 2319
% 
\special{pn 20}%
\special{ar 390 1729 190 190  1.5707963 3.1415927}%
% CIRCLE 0 0 3 0
% 4 2100 2120 2100 2310 2090 2329 2290 2110
% 
\special{pn 20}%
\special{ar 2100 1720 190 190  6.2306022 6.2831853}%
\special{ar 2100 1720 190 190  0.0000000 1.6186068}%
% CIRCLE 0 0 3 0
% 4 2100 799 2290 780 2290 780 2081 591
% 
\special{pn 20}%
\special{ar 2100 399 191 191  4.6212956 6.1835167}%
% LINE 0 0 3 0
% 2 1349 866 361 1502
% 
\special{pn 20}%
\special{pa 1349 466}%
\special{pa 361 1102}%
\special{fp}%
% STR 2 0 3 0
% 3 1397 790 1397 885 2 0
% $e_0$
\put(13.9700,-4.8500){\makebox(0,0)[lb]{$e_0$}}%
% STR 2 0 3 0
% 3 1501 1445 1501 1540 2 0
% $e_1$
\put(15.0100,-11.4000){\makebox(0,0)[lb]{$e_1$}}%
% STR 2 0 3 0
% 3 1102 1692 1102 1787 2 0
% $e_2$
\put(11.0200,-13.8700){\makebox(0,0)[lb]{$e_2$}}%
% STR 2 0 3 0
% 3 627 1854 627 1949 2 0
% $e_3$
\put(6.2700,-15.4900){\makebox(0,0)[lb]{$e_3$}}%
% STR 2 0 3 0
% 3 694 2547 694 2642 2 0
% moduli space
\put(6.9400,-22.4200){\makebox(0,0)[lb]{moduli space}}%
% SPLINE 0 0 3 0
% 4 4247 672 4218 957 3990 1271 3962 1290
% 
\special{pn 20}%
\special{pa 4247 272}%
\special{pa 4247 304}%
\special{pa 4246 337}%
\special{pa 4245 369}%
\special{pa 4243 401}%
\special{pa 4241 433}%
\special{pa 4237 465}%
\special{pa 4232 496}%
\special{pa 4226 528}%
\special{pa 4218 558}%
\special{pa 4207 589}%
\special{pa 4195 619}%
\special{pa 4181 648}%
\special{pa 4166 677}%
\special{pa 4149 705}%
\special{pa 4130 732}%
\special{pa 4110 758}%
\special{pa 4089 783}%
\special{pa 4066 807}%
\special{pa 4042 829}%
\special{pa 4018 850}%
\special{pa 3992 870}%
\special{pa 3966 888}%
\special{pa 3962 890}%
\special{sp}%
% SPLINE 0 0 3 0
% 4 4247 672 4323 1119 4503 1299 4503 1299
% 
\special{pn 20}%
\special{pa 4247 272}%
\special{pa 4247 306}%
\special{pa 4247 341}%
\special{pa 4247 375}%
\special{pa 4248 409}%
\special{pa 4250 442}%
\special{pa 4252 475}%
\special{pa 4255 508}%
\special{pa 4260 540}%
\special{pa 4266 571}%
\special{pa 4273 602}%
\special{pa 4282 632}%
\special{pa 4293 660}%
\special{pa 4305 688}%
\special{pa 4320 715}%
\special{pa 4337 740}%
\special{pa 4357 764}%
\special{pa 4378 788}%
\special{pa 4400 811}%
\special{pa 4424 833}%
\special{pa 4449 854}%
\special{pa 4474 875}%
\special{pa 4500 896}%
\special{pa 4503 899}%
\special{sp}%
% CIRCLE 0 0 3 0
% 4 4247 2088 4152 1071 4152 1071 3240 2249
% 
\special{pn 20}%
\special{ar 4247 1688 1021 1021  2.9830536 4.6192473}%
% CIRCLE 0 0 3 0
% 4 4275 2088 5273 2268 5273 2268 4313 1090
% 
\special{pn 20}%
\special{ar 4275 1688 1014 1014  4.7504467 6.2831853}%
\special{ar 4275 1688 1014 1014  0.0000000 0.1784423}%
% SPLINE 0 0 3 0
% 4 3240 2249 3743 2012 4237 2240 4237 2240
% 
\special{pn 20}%
\special{pa 3240 1849}%
\special{pa 3269 1829}%
\special{pa 3298 1808}%
\special{pa 3327 1788}%
\special{pa 3356 1769}%
\special{pa 3385 1749}%
\special{pa 3414 1731}%
\special{pa 3442 1713}%
\special{pa 3471 1697}%
\special{pa 3500 1681}%
\special{pa 3529 1667}%
\special{pa 3558 1654}%
\special{pa 3587 1642}%
\special{pa 3616 1633}%
\special{pa 3645 1624}%
\special{pa 3674 1618}%
\special{pa 3703 1614}%
\special{pa 3732 1612}%
\special{pa 3761 1612}%
\special{pa 3790 1615}%
\special{pa 3819 1620}%
\special{pa 3848 1627}%
\special{pa 3877 1635}%
\special{pa 3906 1645}%
\special{pa 3935 1657}%
\special{pa 3964 1671}%
\special{pa 3993 1685}%
\special{pa 4022 1701}%
\special{pa 4052 1718}%
\special{pa 4081 1735}%
\special{pa 4110 1754}%
\special{pa 4139 1773}%
\special{pa 4168 1793}%
\special{pa 4197 1812}%
\special{pa 4226 1832}%
\special{pa 4237 1840}%
\special{sp}%
% SPLINE 0 0 3 0
% 4 4237 2230 4750 2040 5282 2278 5282 2278
% 
\special{pn 20}%
\special{pa 4237 1830}%
\special{pa 4267 1813}%
\special{pa 4297 1796}%
\special{pa 4328 1780}%
\special{pa 4358 1763}%
\special{pa 4388 1748}%
\special{pa 4418 1733}%
\special{pa 4448 1718}%
\special{pa 4478 1705}%
\special{pa 4508 1692}%
\special{pa 4538 1680}%
\special{pa 4568 1670}%
\special{pa 4598 1661}%
\special{pa 4628 1653}%
\special{pa 4658 1647}%
\special{pa 4688 1643}%
\special{pa 4718 1641}%
\special{pa 4747 1640}%
\special{pa 4777 1641}%
\special{pa 4806 1645}%
\special{pa 4836 1650}%
\special{pa 4865 1657}%
\special{pa 4895 1665}%
\special{pa 4924 1675}%
\special{pa 4953 1687}%
\special{pa 4982 1700}%
\special{pa 5012 1713}%
\special{pa 5041 1728}%
\special{pa 5070 1744}%
\special{pa 5099 1761}%
\special{pa 5128 1778}%
\special{pa 5157 1796}%
\special{pa 5186 1815}%
\special{pa 5215 1834}%
\special{pa 5244 1853}%
\special{pa 5273 1872}%
\special{pa 5282 1878}%
\special{sp}%
% STR 2 0 3 0
% 3 3734 2547 3734 2642 2 0
% cubic surface 
\put(37.3400,-22.4200){\makebox(0,0)[lb]{cubic surface }}%
% STR 2 0 3 0
% 3 2575 1578 2575 1673 2 0
% $\mathrm{RH}_{z,\kappa}$
\put(25.7500,-12.7300){\makebox(0,0)[lb]{$\mathrm{RH}_{z,\kappa}$}}%
% STR 2 0 3 0
% 3 2470 742 2470 837 2 0
% {\small resolution of}
\put(24.7000,-4.3700){\makebox(0,0)[lb]{{\small resolution of}}}%
% STR 2 0 3 0
% 3 2546 980 2546 1075 2 0
% {\small singularity}
\put(25.4600,-6.7500){\makebox(0,0)[lb]{{\small singularity}}}%
% CIRCLE 0 0 0 0
% 4 4247 672 4247 672 4266 729 4256 729
% 
\special{pn 20}%
\special{ar 4247 272 0 0  1.4141944 6.2831853}%
\special{ar 4247 272 0 0  0.0000000 1.2490458}%
% CIRCLE 0 0 0 0
% 4 4247 672 4275 739 4275 739 4275 739
% 
\special{pn 20}%
\special{sh 0.600}%
\special{ar 4247 272 73 73  0.0000000 6.2831853}%
% STR 2 0 3 0
% 3 3848 1778 3848 1873 2 0
% $\mathcal{S}(\theta)$
\put(38.4800,-14.7300){\makebox(0,0)[lb]{$\mathcal{S}(\theta)$}}%
% VECTOR 0 0 3 0
% 2 2404 1749 3192 1749
% 
\special{pn 20}%
\special{pa 2404 1349}%
\special{pa 3192 1349}%
\special{fp}%
\special{sh 1}%
\special{pa 3192 1349}%
\special{pa 3125 1329}%
\special{pa 3139 1349}%
\special{pa 3125 1369}%
\special{pa 3192 1349}%
\special{fp}%
% STR 2 0 3 0
% 3 1283 2110 1283 2205 2 0
% $\mathcal{M}_z(\kappa)$
\put(12.8300,-18.0500){\makebox(0,0)[lb]{$\mathcal{M}_z(\kappa)$}}%
% LINE 0 0 3 0
% 2 931 828 1482 1303
% 
\special{pn 20}%
\special{pa 931 428}%
\special{pa 1482 903}%
\special{fp}%
% LINE 0 0 3 0
% 2 713 961 1131 1540
% 
\special{pn 20}%
\special{pa 713 561}%
\special{pa 1131 1140}%
\special{fp}%
% LINE 0 0 3 0
% 2 494 1122 713 1711
% 
\special{pn 20}%
\special{pa 494 722}%
\special{pa 713 1311}%
\special{fp}%
% STR 2 0 3 0
% 3 4722 714 4722 809 2 0
% $D_4$
\put(47.2200,-4.0900){\makebox(0,0)[lb]{$D_4$}}%
% VECTOR 2 0 3 0
% 2 4608 695 4389 695
% 
\special{pn 8}%
\special{pa 4608 295}%
\special{pa 4389 295}%
\special{fp}%
\special{sh 1}%
\special{pa 4389 295}%
\special{pa 4456 315}%
\special{pa 4442 295}%
\special{pa 4456 275}%
\special{pa 4389 295}%
\special{fp}%
% VECTOR 2 0 3 0
% 2 2888 1160 2765 1426
% 
\special{pn 8}%
\special{pa 2888 760}%
\special{pa 2765 1026}%
\special{fp}%
\special{sh 1}%
\special{pa 2765 1026}%
\special{pa 2811 974}%
\special{pa 2787 978}%
\special{pa 2775 957}%
\special{pa 2765 1026}%
\special{fp}%
% STR 2 0 3 0
% 3 310 830 310 930 2 0
% $\mathcal{E}_z(\kappa)$
\put(3.1000,-5.3000){\makebox(0,0)[lb]{$\mathcal{E}_z(\kappa)$}}%
\end{picture}%
\end{center}
\caption{Resolution of singularities by Riemann-Hilbert 
correspondence} 
\label{fig:sresol}
\end{figure}

 
If $\k \in \K-\Wall$ then the surface $\Sol(\th)$ is smooth 
and $\RH_{z,\k}$ is a biholomorphism, while if $\k \in \Wall$, 
it is not a biholomorphism but only gives a resolution of 
singularities (proper and surjective, but not injective). 
For example, see Figure \ref{fig:sresol} for the case 
$\k = (0,0,0,0,1)$ where a singularity of type $D_4$ occurs. 
In the latter case, however, if we take a standard 
{\sl algebraic} minimal resolution of Kleinian singularities 
as constructed by Brieskorn \cite{Brieskorn} and others, 
%%%%%%%%%%%%%%%%%%%%%%%% eqn:brieskorn %%%%%%%%%%%%%%%%%%%%%
\begin{equation} \label{eqn:brieskorn}
\varphi : \wt{\Sol}(\th) \to \Sol(\th) 
\end{equation}
%%%%%%%%%%%%%%%%%%%%%%%%%%%%%%%%%%%%%%%%%%%%%%%%%%%%%%%%%%%%
then we can lift the Riemann-Hilbert correspondence 
(\ref{eqn:RHkth}) to have a commutative diagram 
%%%%%%%%%%%%%%%%%%%%%%%%%%%%% cd:lift %%%%%%%%%%%%%%%%%%%%%%
\begin{equation} \label{cd:lift}
\begin{CD}
\M_z(\k) @> \wt{\RH}_{z,\k} >> \wt{\Sol}(\th) \\
@|     @VV \varphi V  \\
\M_z(\k) @> \RH_{z,\k} >> \Sol(\th). 
\end{CD}
\end{equation}
%%%%%%%%%%%%%%%%%%%%%%%%%%%%%%%%%%%%%%%%%%%%%%%%%%%%%%%%%%%%
The lifted Riemann-Hilbert correspondence $\wt{\RH}_{z,\k}$ 
is a biholomorphism and hence gives a strict conjugacy 
between the nonlinear monodromy (\ref{eqn:nm}) of 
$\PVI(\k)$ and a certain automorphism 
%%%%%%%%%%%%%%%%%%%%%%%%%%%%%% eqn:nm2 %%%%%%%%%%%%%%%%%%%%%
\begin{equation} \label{eqn:nm2} 
\tilde{g} : \wt{\Sol}(\th) \to \wt{\Sol}(\th). 
\end{equation}
%%%%%%%%%%%%%%%%%%%%%%%%%%%%%%%%%%%%%%%%%%%%%%%%%%%%%%%%%%%%
This latter map will be described explicitly in Section 
\ref{sec:cubic} (see Theorem \ref{thm:conjugacy}). 
\par 
The singularity structure of the affine cubic surface 
$\Sol(\th)$ is closely related to the Riccati 
solutions to $\PVI(\k)$ \cite{IIS1}, where a Riccati 
solution is a particular solution that arises from the Riccati 
equation associated to a Gauss hypergeometric equation. 
Let $\E_z(\k) \subset \M_z(\k)$ be the exceptional set of 
the resolution of singularities by the Riemann-Hilbert 
correspondence (\ref{eqn:RHkth}). 
%%%%%%%%%%%%%%%%%%%%%%% fig:sresol %%%%%%%%%%%%%%%%%%%%%%%%%
%%%%%%%%%%%%%%%%%%%%%%%%%%%%%%%%%%%%%%%%%%%%%%%%%%%%%%%%%%%%
Similarly, let $\E(\th) \subset \wt{\Sol}(\th)$ be the 
exceptional set of the algebraic resolution of 
singularities (\ref{eqn:brieskorn}). 
%%%%%%%%%%%%%%%%%%%%%%%%% thm:riccati %%%%%%%%%%%%%%%%%%%%%% 
\begin{theorem}[\cite{IIS1,STe}] \label{thm:riccati} 
Equation $\PVI(\k)$ admits Riccati solutions if and only if 
$\k \in \Wall$. 
All Riccati solution germs at time $z \in Z$ are parametrized 
by the exceptional set $\E_z(\k) \subset \M_z(\k)$, which 
precisely corresponds to the exceptional set 
$\E(\th) \subset \wt{\Sol}(\th)$ through the lifted 
Riemann-Hilbert correspondence $(\ref{cd:lift})$. 
\end{theorem}
%%%%%%%%%%%%%%%%%%%%%%%%%%%%%%%%%%%%%%%%%%%%%%%%%%%%%%%%%%%%
Fot this reason we may refer to $\E_z(\k)$ and 
$\M_z^{\circ}(\k) := \M_z(\k) - \E_z(\k)$ as the 
{\sl Riccati locus} and {\sl non-Riccati locus} of 
$\M_z(\k)$ respectively. 
They are invariant under the action of the nonlinear 
monodromy (\ref{eqn:nm}). 
Corresponding to them, let $\mathrm{Sing}(\th)$ and  
$\Sol^{\circ}(\th) := \Sol(\th)-\mathrm{Sing}(\th)$ be the 
singular locus and the smooth locus of the cubic surface 
$\Sol(\th)$ respectively. 
%%%%%%%%%%%%%%%%%%%%%%%%% rem:riccati %%%%%%%%%%%%%%%%%%%%%%
\begin{remark} \label{rem:riccati} 
Two remarks are in order at this stage. 
\begin{enumerate}
\item By Theorem \ref{thm:SolRHP} the Riemann-Hilbert 
correspondence $(\ref{eqn:RHkth})$ restricts to a 
biholomorphism 
%%%%%%%%%%%%%%%%%%%%%%%%% eqn:RHkth2 %%%%%%%%%%%%%%%%%%%%%%
\begin{equation} \label{eqn:RHkth2} 
\RH_{z,\k}^{\circ} : 
\M_z^{\circ}(\k) \to \Sol^{\circ}(\th) 
\end{equation} 
%%%%%%%%%%%%%%%%%%%%%%%%%%%%%%%%%%%%%%%%%%%%%%%%%%%%%%%%%%%
between the non-Riccati locus of $\M_z(\k)$ and the smooth 
locus of $\Sol(\th)$, while it collapses the Riccati locus 
$\E_z(\k)$ to the singular locus $\mathrm{Sing}(\th)$. 
In order to resolve this degeneracy and obtain an 
isomorphism, we had to take the lifted Riemann-Hibert 
correspondence $(\ref{cd:lift})$, which induces an 
isomorphism between the exceptional sets $\E_z(\k)$ and 
$\E(\th)$. 
\item For the Riccati solutions the main problem of this 
article is trivial; if a Riccati solution is a finite branch 
solution around a fixed singular point, then it is an 
algebraic branch solution, because the Riccati solution is 
(essentially) the logarithmic derivative of a Gauss 
hypergeometric function. 
Thus we may restrict our attention to the 
non-Riccati locus. 
\end{enumerate} 
\end{remark}
%%%%%%%%%%%%%%%%%%%%%%%%%%%%%%%%%%%%%%%%%%%%%%%%%%%%%%%%%%%%%%
%%%%%%%%%%%%%%%%%%%%%%%%% sec:cubic %%%%%%%%%%%%%%%%%%%%%%%%%%
\section{Dynamics on Cubic Surface} \label{sec:cubic} 
%%%%%%%%%%%%%%%%%%%%%%%%%%%%%%%%%%%%%%%%%%%%%%%%%%%%%%%%%%%%%%
We shall describe the strict conjugacy (\ref{eqn:nm2}) 
of the nonlinear monodromy (\ref{eqn:nm}). 
For a cyclic permutation $(i,j,k)$ of $(1,2,3)$ we define 
an isomorphism $g_i : \Sol(\th) \to \Sol(\th')$, 
$(x,\th) \mapsto (x',\th')$ by 
%%%%%%%%%%%%%%%%%%%%%%%%%% eqn:gi %%%%%%%%%%%%%%%%%%%%%%%%%%%% 
\begin{equation} \label{eqn:gi}
g_i : (x_i',x_j',x_k', \th_i', \th_j', \th_k') = 
(\th_j - x_j- x_k x_i, x_i, x_k, \th_j, \th_i, \th_k). 
\end{equation}
%%%%%%%%%%%%%%%%%%%%%%%%%%%%%%%%%%%%%%%%%%%%%%%%%%%%%%%%%%%%%
Through the resolution of singularities 
(\ref{eqn:brieskorn}), the map $g_i$ is uniquely lifted 
to an isomorphism 
%%%%% 
\[
\tilde{g}_i : \wt{\Sol}(\th) \to \wt{\Sol}(\th'), 
\qquad (i = 1,2,3). 
\]
%%%%%$
We remark that the square $g_i^2$ is an automorphism of 
$\Sol(\th)$ with $\tilde{g}_i^2$ being its lift to 
$\wt{\Sol}(\th)$. 
%%%%%%%%%%%%%%%%%%%%%%% thm:conjugacy %%%%%%%%%%%%%%%%%%%%%%%%
\begin{theorem}[\cite{IIS1}] \label{thm:conjugacy} 
For each $i \in \{1,2,3\}$ the nonlinear monodromy 
$\ga_{i*} : \M_z(\k) \carl$ along the $i$-th basic 
loop $\ga_i$ is strictly conjugated to the automorphism 
$\tilde{g}_i^2 : \wt{\Sol}(\th) \carl$ via 
the lifted Riemann-Hilbert correspondence $(\ref{cd:lift})$. 
More generally, if $\ga \in \pi_1(Z,z)$ is represented by 
$\ga = \ga_{i_1}^{\ve_1} \ga_{i_2}^{\ve_2} \cdots 
\ga_{i_n}^{\ve_n}$ 
with $(i_1,\dots,i_n) \in \{1,2,3\}^n$ and 
$(\ve_1,\dots,\ve_n) \in \{\pm1\}^n$, then the map 
$(\ref{eqn:nm2})$ is given by 
%%%%
\[
\tilde{g} = \tilde{g}_{i_1}^{2\ve_1} 
\tilde{g}_{i_2}^{2\ve_2} \cdots \tilde{g}_{i_n}^{2\ve_n}. 
\]
\end{theorem} 
%%%%%%%%%%%%%%%%%%%%%%%%%%%%%%%%%%%%%%%%%%%%%%%%%%%%%%%%%%%%%%%
\par 
Let $\wt{\mathrm{Fix}}_j(\th)$ be the set of all 
fixed points of the transformation 
$\tilde{g}_j^2 : \wt{\Sol}(\th) \carl$. 
Moreover, for any integer $n > 1$, let 
$\wt{\mathrm{Per}}_j(\th;n)$ be the set of all periodic 
points of {\sl prime} period $n$ of the transformation 
$\tilde{g}_j^2 : \wt{\Sol}(\th) \carl$. 
Theorem \ref{thm:conjugacy} then implies that 
all single-valued solution germs and all $n$-branch solution 
germs to $\PVI(\k)$ around the fixed singular point $z_j$ are 
parametrized by the sets $\wt{\mathrm{Fix}}_j(\th)$ and 
$\wt{\mathrm{Per}}_j(\th;n)$ respectively. 
By Remark \ref{rem:riccati}, considering 
$\wt{\mathrm{Fix}}_j(\th)$ and $\wt{\mathrm{Per}}_j(\th;n)$ 
upstairs is the same thing as considering $\mathrm{Fix}_j(\th)$ 
and $\mathrm{Per}_j(\th;n)$ downstairs, except for the 
exceptional locus upstairs and the singular locus 
downstairs. 
Here $\mathrm{Fix}_j(\th)$ and $\mathrm{Per}_j(\th;n)$ denote 
the set of all fixed points and the set of all periodic 
points of prime period $n$ of the transformation 
$g_j^2 : \Sol(\th) \carl$ downstairs. 
In order to make the situation more transparent, we 
begin by investigating simultaneous fixed points 
of $g_1^2$, $g_2^2$, $g_3^2$ downstairs. 
%%%%%%%%%%%%%%%%%%%%%%%% thm:SingFix %%%%%%%%%%%%%%%%%%%%%%%%%%
\begin{theorem} \label{thm:SingFix} 
If $\mathrm{Fix}(\th)$ is the set of all simultaneous fixed 
points of $g_1^2$, $g_2^2$, $g_3^2 : \Sol(\th) \carl$, then 
%%%%%%%%%%%%%%%%%%%%%%%%%%% eqn:SingFix %%%%%%%%%%%%%%%%%%%%%%
\begin{equation} \label{eqn:SingFix} 
\mathrm{Fix}(\th) = \mathrm{Sing}(\th). 
\end{equation}
%%%%%%%%%%%%%%%%%%%%%%%%%%%%%%%%%%%%%%%%%%%%%%%%%%%%%%%%%%%%%%
\end{theorem}
%%%%%%%%%%%%%%%%%%%%%%%%%%%%%%%%%%%%%%%%%%%%%%%%%%%%%%%%%%%%%%
{\it Proof}. A point $x \in \Sol(\th)$ is a singular point of 
the surface $\Sol(\th)$ if and only if its gradient vector field 
$y(x,\th) = (y_1(x,\th), y_2(x,\th), y_3(x,\th))$ vanishes 
at the point $x$, where 
%%%%%%%%%%%%%%%%%%%%%%%%% eqn:yi %%%%%%%%%%%%%%%%%%%%%%%%%%%%
\begin{equation} \label{eqn:yi}
y_i(x,\th) := \dfrac{\partial f}{\partial x_i}(x,\th) = 
2 x_i + x_j x_k - \th_i
\end{equation}
%%%%%%%%%%%%%%%%%%%%%%%%%%%%%%%%%%%%%%%%%%%%%%%%%%%%%%%%%%%%%%
On the other hand, an inspection of formula (\ref{eqn:gi}) 
readily shows that $x \in \Sol(\th)$ is a simultaneous 
fixed point of $g_1^2$, $g_2^2$, $g_3^2$ if and only if 
$x$ is a common root of equations 
%%%%%%%%%%%%%%%%%%%%%%%%%%%% eqn:SingFix2 %%%%%%%%%%%%%%%%%%%%
\begin{equation} \label{eqn:SingFix2} 
f(x,\th) = y_1(x,\th) = y_2(x,\th) = y_3(x,\th) = 0. 
\end{equation}
%%%%%%%%%%%%%%%%%%%%%%%%%%%%%%%%%%%%%%%%%%%%%%%%%%%%%%%%%%%%%%
Then the equality (\ref{eqn:SingFix}) immediately follows 
from these observations. \hfill $\Box$  \par\medskip
%%%%%%%%%%%%%%%%%%%%%%%%% end of proof %%%%%%%%%%%%%%%%%%%%%%%
As is announced in \cite{IIS1}, this theorem yields a 
characterization of the rational solutions. 
%%%%%%%%%%%%%%%%%%%%%%%%% cor:ratinal %%%%%%%%%%%%%%%%%%%%%%%%
\begin{corollary} \label{cor:rational}
Any single-valued global solution to $\PVI(\k)$ is a 
rational Riccati solution. 
\end{corollary}
%%%%%%%%%%%%%%%%%%%%%%%%%%%%%%%%%%%%%%%%%%%%%%%%%%%%%%%%%%%%%%
{\it Proof}. 
If a single-valued solution $Q \in \M_z(\k)$ belongs 
to the non-Riccati locus $Q \in \M_z^{\circ}(\k)$, 
then the Riemann-Hilbert correspondence (\ref{eqn:RHkth2}) 
sends $Q$ to a smooth point $x \in \Sol^{\circ}(\th)$. 
Since the single-valued solution $Q$ is a simultaneous 
fixed point of the nonlinear monodromies $\ga_{1*}$, 
$\ga_{2*}$, $\ga_{3*}$, the corresponding point $x$ must 
lie in $\mathrm{Fix}(\th)$. 
Then Theorem \ref{thm:SingFix} implies that 
$x \in \mathrm{Sing}(\th)$, which contradicts the fact 
that $x \in \Sol^{\circ}(\th)$. 
Hence any single-valued solution is a Riccati solution. 
Since any Riccati solution is (essentially) the logarithmic 
derivative of a Gauss hypergeometric function, any 
single-valued Riccati solution must be a rational 
solution. \hfill $\Box$ \par\medskip
%%%%%%%%%%%%%%%%%%%%%%%%% end of proof %%%%%%%%%%%%%%%%%%%%%%%
All the rational solutions to Painlev\'e VI are classified in 
\cite{Mazzocco}. 
We come back to our discussion downstairs and give a 
simple characterization of the sets $\mathrm{Fix}_j(\th)$ 
and $\mathrm{Per}_j(\th;n)$. 
%%%%%%%%%%%%%%%%%%%%%%%%% lem:per %%%%%%%%%%%%%%%%%%%%%%%%%%%%
\begin{lemma} \label{lem:per} 
Let $x = (x_1,x_2,x_3) \in \Sol(\th)$ be any point and let 
$n$ be any integer $> 1$. 
%%%%%%
\begin{enumerate} 
\item $x \in \mathrm{Fix}_{j}(\th)$ if and only if $x$ is a 
root of equations 
%%%%%%%%%%%%%%%%%%%%%%%%% eqn:per1 %%%%%%%%%%%%%%%%%%%%%%%%%%%
\begin{equation} \label{eqn:per1}
f(x,\th) = y_j(x,\th) = y_k(x,\th) = 0. 
\end{equation}
%%%%%%%%%%%%%%%%%%%%%%%%%%%%%%%%%%%%%%%%%%%%%%%%%%%%%%%%%%%%%%
\item $x \in \mathrm{Per}_{j}(\th;n)$ if and only if 
there exists an integer $0 < m < n$ coprime to $n$ such that 
%%%%%%%%%%%%%%%%%%%%%%%%% eqn:pern %%%%%%%%%%%%%%%%%%%%%%%%%%%
\begin{equation} \label{eqn:pern}
f(x,\th) = 0, \qquad x_i = 2 \cos(\pi m/n). 
\end{equation}
%%%%%%%%%%%%%%%%%%%%%%%%%%%%%%%%%%%%%%%%%%%%%%%%%%%%%%%%%%%%%%
\end{enumerate} 
\end{lemma}
%%%%%%%%%%%%%%%%%%%%%%%%%%%%%%%%%%%%%%%%%%%%%%%%%%%%%%%%%%%%%%
{\it Proof}. 
We put $(x',\th') = g_j(x,\th)$ and $y' = y(x',\th')$. 
Then formula (\ref{eqn:gi}) yields 
%%%%%%%%%%%%%%%%%%%%% eqn:action1 %%%%%%%%%%%%%%%%%%%%%%% 
\begin{equation} \label{eqn:action1}
y_i' = y_i - x_j y_k, \qquad y_j' = - y_k, \qquad 
y_k' = y_j - x_i y_k. 
\end{equation} 
%%%%%%%%%%%%%%%%%%%%%%%%%%%%%%%%%%%%%%%%%%%%%%%%%%%%%%%%% 
For each integer $n \in \Z$ we write 
$(x^{(n)}, \th^{(n)}) = g_i^n(x,\th)$ and  
$y^{(n)} = y(x^{(n)},\th^{(n)})$.   
From formulas (\ref{eqn:gi}) and (\ref{eqn:action1}), 
we can easily obtain three recurrence relations    
%%%%%%%%%%%%%%%%%%%%%% eqn:recurrence1 %%%%%%%%%%%%%%%%%%%%%%% 
\begin{eqnarray} 
y_j^{(n+2)} + x_i \, y_j^{(n+1)} + y_j^{(n)} &=& 0,  
\label{eqn:recurrence1} \\[2mm]
%%%%%%%%%%%%%%%%%%%%%% eqn:recurrence2 %%%%%%%%%%%%%%%%%%%%%%%
x_j^{(n+2)} - x_j^{(n)} &=& y_j^{(n+2)}, 
\label{eqn:recurrence2} \\[2mm]
%%%%%%%%%%%%%%%%%%%%%% eqn:recurrence3 %%%%%%%%%%%%%%%%%%%%%%%
x_k^{(n+1)} &=& x_j^{(n)}. \label{eqn:recurrence3}    
\end{eqnarray}
%%%%%%%%%%%%%%%%%%%%%%%%%
The characteristic equation of the recurrence relation 
(\ref{eqn:recurrence1}) is the quadratic equation 
%%%%%%%%%%%%%%%%%%%%%% eqn:quadratic %%%%%%%%%%%%%%%%%%%%%%%%%%
\begin{equation} \label{eqn:quadratic}
\lambda^2 + x_i \, \lambda + 1 = 0,  
\end{equation} 
%%%%%%%%%%%%%%%%%%%%%%%%%%%%%%%%%%%%%%%%%%%%%%%%%%%%%%%%%%%%%%%
the roots of which are denoted by $\alpha$ and 
$\beta = \alpha^{-1}$.   
Since $\alpha \beta = 1$, we may and shall assume that 
$|\alpha| \ge 1 \ge |\beta| > 0$ in the sequel. 
The discussion is divided into two cases.   
\par
%%%%%%%%%%%%%%%%%%%%%%% case 1 %%%%%%%%%%%%%%%%%%%%%%%%%%%%%%%%
{\bf Case} $\mbox{\boldmath $x_i \in \C-\{\pm 2\}$}$:   
In this case, the roots $\alpha$ and $\beta$ are distinct and 
different from $\pm 1$ and the recurrence relation 
(\ref{eqn:recurrence1}) is settled as   
%%%%%
\[
y_j^{(n)} = \dfrac{\beta^n (\alpha y_j + y_k) - 
            \alpha^n (\beta y_j + y_k)}{\alpha - \beta}.  
\]
%%%%%
Then it follows from (\ref{eqn:recurrence2}) and 
(\ref{eqn:recurrence3}) that the sequences $x_j^{(n)}$ and 
$x_k^{(n)}$ are determined as 
%%%%%%%%%%%%%%%%%%%%%%%% eqn:solution %%%%%%%%%%%%%%%%%%%%%%%%
\begin{equation} \label{eqn:solution}
\left\{
\begin{array}{rcccl}  
x_j^{(2n)} &=& x_k^{(2n+1)} &=& 
p \, \alpha^{2n} + q \, \beta^{2n} + r_1, \\[3mm]
x_j^{(2n+1)} &=& x_k^{(2n)} &=& 
p \, \alpha^{2n+1} + q \, \beta^{2n+1} + r_2, \\
\end{array}
\right. 
\end{equation}  
%%%%%%%%%%%%%%%%%%%%%%%%%%%%%%%%%%%%%%%%%%%%%%%%%%%%%%%%%%%%%
where the constants $p$, $q$, $r_1$ and $r_2$ are given by     
%%%%%% 
\[
\begin{array}{rclrcl}
p &=& 
-\dfrac{\alpha^2(\beta y_j+y_k)}{(\alpha-\beta)(\alpha^2-1)}, 
\qquad\, &
q &=& 
\dfrac{\beta^2(\alpha y_j+y_k)}{(\alpha-\beta)(\beta^2-1)}, 
\\[5mm]
r_1 &=& x_j - p - q, \qquad\, & r_2 &=& x_j' 
- \alpha p - \beta q.  
\end{array} 
\] 
%%%%%%
Notice that $p = q = 0$ if and only if $x$ satisfies 
equations (\ref{eqn:per1}). 
Indeed, the condition $p = q = 0$ is equivalent to 
$\alpha y_j + y_k = \beta y_j + y_k = 0$, which is equivalent 
to the condition $y_j = y_k = 0$, because the roots $\alpha$ 
and $\beta$ are distinct. 
\par 
Now we assume that $x$ is a root of equations 
(\ref{eqn:per1}). 
Then (\ref{eqn:solution}) implies that the sequence 
$x^{(n)}$ is periodic of period two, that is, $x$ is 
a fixed point of $g_j^2$. 
Next we assume that $x$ is not a root of equations 
(\ref{eqn:per1}). 
If $x$ is a periodic point of $g_j^2$ of prime period 
$n \ge 1$, then (\ref{eqn:solution}) yields 
%%%%%%%%
\[
\begin{array}{rclcl}
x_j^{(2n)} - x_j &=& 
(\alpha^{2n}-1)(p - q \beta^{2n}) &=& 0, \\[2mm]
x_k^{(2n)} - x_k &=& 
(\alpha^{2n}-1)(p \alpha - q \beta^{2n+1}) &=& 0. 
\end{array} 
\]
%%%%%%%%
Here it cannot happen that 
$p - q \beta^{2n} = p \alpha - q \beta^{2n+1} = 0$. 
Indeed, otherwise, we have $p = q \beta^{2n}$ and 
$q(1-\beta^2) = 0$. 
Since at least one of $p$ and $q$ is nonzero, we have 
$\beta \in \{\pm 1\}$ and hence $x_i \in \{\pm2\}$, 
which contradicts the assumption that 
$x_i \not\in \{\pm2\}$. 
Therefore, $\alpha^{2n} = 1$, that is, $\alpha$ is 
a primitive $2n$-th root of unity. 
Note that $n \ge 2$ since $\alpha \not\in \{\pm1\}$. 
Thus there is an integer $0 < m < n$ comprime to $n$ 
such that $\alpha = \exp(\pi i m/n)$ and so 
$x_i = \alpha + \alpha^{-1} = 2 \cos(\pi m/n)$, which 
leads to condition (\ref{eqn:pern}). 
Conversely, if condition (\ref{eqn:pern}) is 
satisfied, then it is easy to see that $x$ is a periodic 
point of $g_j^2$ of prime period $n$. 
\par
%%%%%%%%%%%%%%%%%%%%% case 2 %%%%%%%%%%%%%%%%%%%%%%%%%% 
{\bf Case} $\mbox{\boldmath $x_i \in \{\pm 2\}$}$:  
In this case we have $x_i = -2 \ve$ for some sign 
$\ve \in \{\pm1\}$ and hence equation (\ref{eqn:quadratic}) 
has a double root $\alpha = \beta = \ve$. 
Then the recurrence equation (\ref{eqn:recurrence1}) is 
settled as 
$y_j^{(n)} = \ve^n \{y_j-n(y_j+\ve y_k) \}$. 
If the sequence $x^{(n)}$ is periodic, then so is the 
sequence $y_j^{(n)}$. 
This is the case if and only if $y_j+\ve y_k = 0$.
Conversely, if this condition is satisfied, then we have 
$y_j^{(n)} = \ve^n y_j$. 
Substituting this equation into 
(\ref{eqn:recurrence2}) yields 
%%%%%%%%%%%%%%%%%%%%%% eqn:solution2 %%%%%%%%%%%%%%%%%%%%%%
\begin{equation} \label{eqn:solution2} 
\left\{
\begin{array}{ccccl}
x_j^{(2n)} &=& x_k^{(2n+1)} &=& x_j + n y_j, \\[3mm]
x_j^{(2n+1)} &=& x_k^{(2n+2)} &=& x_j'+\ve n y_j.  
\end{array}
\right.
\end{equation} 
%%%%%%%%%%%%%%%%%%%%%%%%%%%%%%%%%%%%%%%%%%%%%%%%%%%%%%%%%% 
Hence the sequence $x^{(2n)}$ is periodic if and only 
if $y_j = y_k = 0$, namely, if and only if 
$x$ is a root of (\ref{eqn:per1}). 
In this case $x$ is a fixed point of $g_j^2$. 
\hfill $\Box$ \par\medskip
%%%%%%%%%%%%%%%%%%%%% end of proof %%%%%%%%%%%%%%%%%%%%%%%%
In order to give the relation between the fixed points 
upstairs and those downstairs, we put
%%%%%%%
\[
\mathrm{Fix}_j^{\circ}(\th) := \mathrm{Fix}_j(\th)- 
\mathrm{Sing}(\th), \quad  
\wt{\mathrm{Fix}}_j^{\circ}(\th) := 
\wt{\mathrm{Fix}}_j(\th)- \mathcal{E}(\th), \quad  
\wt{\mathrm{Fix}}_j^{e}(\th) := 
\wt{\mathrm{Fix}}_j(\th) \cap \mathcal{E}(\th). 
\]
%%%%%%%
For the periodic points of prime period $n > 1$, we 
define $\mathrm{Per}_j^{\circ}(\th;n)$, 
$\wt{\mathrm{Per}}_j^{\circ}(\th;n)$ and 
$\wt{\mathrm{Per}}_j^{e}(\th;n)$ in a similar manner.  
Then there exist direct sum decompositions
%%%%%%%%
\[
\wt{\mathrm{Fix}}_j(\th) = 
\wt{\mathrm{Fix}}_j^{\circ}(\th) \amalg 
\wt{\mathrm{Fix}}_j^{e}(\th), 
\qquad 
\wt{\mathrm{Per}}_j(\th;n) = 
\wt{\mathrm{Per}}_j^{\circ}(\th;n) \amalg 
\wt{\mathrm{Per}}_j^{e}(\th;n), 
\]
%%%%%%%%
where the exceptional components 
$\wt{\mathrm{Fix}}_j^{e}(\th)$ and 
$\wt{\mathrm{Per}}_j^{e}(\th;n)$ 
parametrize the single-valued Riccati solutions and 
the $n$-branched Riccati solutions around the fixed 
singular point $z_j$ respectively. 
%%%%%%%%%%%%%%%%%%%%%%% lem:updown %%%%%%%%%%%%%%%%%%%%%%%%
\begin{lemma} \label{lem:updown} 
The minimal resolution $(\ref{eqn:brieskorn})$ induces an 
isomorphism 
%%%%%%%%%%%%%%%%%%%%%%% eqn:updownFix %%%%%%%%%%%%%%%%%%%%
\begin{equation} \label{eqn:updownFix}
\varphi : \wt{\mathrm{Fix}}_j^{\circ}(\th) \to 
\mathrm{Fix}_j^{\circ}(\th). 
\end{equation}
%%%%%%%%%%%%%%%%%%%%%%%%%%%%%%%%%%%%%%%%%%%%%%%%%%%%%%%%%%
For any $n > 1$ we have $\mathrm{Per}(\th;n) \cap 
\mathrm{Sing}(\th) = \emptyset$, that is, 
$\mathrm{Per}_j^{\circ}(\th;n) = \mathrm{Per}_j(\th;n)$, 
and the minimal resolution $(\ref{eqn:brieskorn})$ 
induces an isomorphism 
%%%%%%%%%%%%%%%%%%%%%%% eqn:updownPer %%%%%%%%%%%%%%%%%%%%
\begin{equation} \label{eqn:updownPer}
\varphi : \wt{\mathrm{Per}}_j^{\circ}(\th;n) \to 
\mathrm{Per}_j(\th;n) \qquad (n > 1). 
\end{equation}
%%%%%%%%%%%%%%%%%%%%%%%%%%%%%%%%%%%%%%%%%%%%%%%%%%%%%%%%%%
\end{lemma}
%%%%%%%%%%%%%%%%%%%%%%%%%% proof %%%%%%%%%%%%%%%%%%%%%%%%%%
{\it Proof}. 
The isomorphism (\ref{eqn:updownFix}) is trivial from the 
definition. 
The assertion $\mathrm{Per}(\th;n) \cap \mathrm{Sing}(\th) 
= \emptyset$ follows from (\ref{eqn:SingFix}). 
Then the isomorphism (\ref{eqn:updownPer}) is again trivial 
from the definition. \hfill $\Box$ \par\medskip 
%%%%%%%%%%%%%%%%%%%%%%% end of proof %%%%%%%%%%%%%%%%%%%%%%
The fixed point set and the periodic point set, upstairs 
or downstairs, will be investigated more closely in 
\S\ref{sec:fixed}. 
For this purpose it is convenient to consider the 
symmetric group $S_4$ of degree $4$ acting on $\K$ by 
permuting the entries $\k_1$, $\k_2$, $\k_3$, $\k_4$ of 
$\k \in \K$ and fixing $\k_0$. 
Through the Riemann-Hilbert correspondence in the parameter 
level, $\rh : \K \to \Th$, the action $S_4 \car \K$ 
induces an action of $S_3 \ltimes \Kl$ on $\Th$, 
where $\Kl$ is Klein's $4$-group realized as the group 
of even triple signs, 
$\Kl = \{ \ve = (\ve_1,\ve_2,\ve_3) \in \{\pm1\}^3 \,:\, 
\ve_1\ve_2\ve_3 = 1\}$, acting on $\Th$ by the sign 
changes $(\th_1,\th_2,\th_3,\th_4) \mapsto 
(\ve_1\th_1,\ve_2\th_2,\ve_3\th_3,\th_4)$, while 
$S_3$ acts on $\Th$ by permuting the entries 
$\th_1$, $\th_2$, $\th_3$ of $\th \in \Th$ and fixing 
$\th_4$. 
This construction defines an isomorphism of groups 
%%%%%%%%%%%%%%%%%%%%%%%% eqn:equiv %%%%%%%%%%%%%%%%%%%%
\begin{equation} \label{eqn:equiv} 
S_4 \cong S_3 \ltimes \Kl, \quad \si \mapsto (\tau, \ve), 
\end{equation} 
%%%%%%%%%%%%%%%%%%%%%%%%%%%%%%%%%%%%%%%%%%%%%%%%%%%%%%%
with respect to which the map $\rh : \K \to \Th$ becomes 
$S_4$-equivariant. 
Viewed as a subgroup of $S_4$, Klein's $4$-group is the 
permutation group 
$\Kl = \{1, (14)(23), (24)(31), (34)(12) \}$. 
\par 
Let $\si \in S_4$ act on $x = (x_1,x_2,x_3)$ in 
the same manner as it does on $(\th_1,\th_2,\th_3)$. 
Then the polynomial $f(x,\th)$ is $\si$-invariant and 
hence $\si$ induces an isomorphism of algebraic surfaces, 
$\si : \Sol(\th) \to \Sol(\si(\th))$. 
As for the action $g_j^2 : \Sol(\th) \carl$, 
we have the commutative diagram 
%%%%%%%%%%%%%%%%%%%%%%%%% cd:S4 %%%%%%%%%%%%%%%%%%%%%%%%
\begin{equation} \label{cd:S4}
\begin{CD}
\Sol(\th) @> g_j^2 >> \Sol(\th) \\
 @V \si VV          @VV \si V \\
\Sol(\si(\th)) @>> g_{\tau(j)}^2 > \Sol(\si(\th)), 
\end{CD}
\end{equation}
%%%%%%%%%%%%%%%%%%%%%%%%%%%%%%%%%%%%%%%%%%%%%%%%%%%%%%%%
for any element $\si \in S_4$ with $\tau \in S_3$ 
determined by (\ref{eqn:equiv}). 
It induces isomorphisms 
%%%%%%
\[
\si : \mathrm{Fix}_j(\th) \to \mathrm{Fix}_{\tau(j)}(\si(\th)), 
\qquad 
\si : \mathrm{Per}_j(\th;n) \to \mathrm{Per}_{\tau(j)}(\si(\th);n), 
\]
%%%%%%
which, via the minimal resolution (\ref{eqn:brieskorn}), 
lift up to isomorphisms 
%%%%%%
\[
\wt{\si} : \wt{\mathrm{Fix}}_j(\th) \to 
\wt{\mathrm{Fix}}_{\tau(j)}(\si(\th)), 
\qquad 
\wt{\si} : \wt{\mathrm{Per}}_j(\th;n) \to 
\wt{\mathrm{Per}}_{\tau(j)}(\si(\th);n). 
\]
%%%%%%
\par 
The action of the symmetric group $S_4$ on $\K$ mentioned 
above is just induced from its action on the index set 
$\{0,1,2,3,4\}$ fixing the element $0$, namely, from the 
realization of $S_4$ as the automorphism group of the 
Dynkin diagram $D_4^{(1)}$. 
By taking the semi-direct product by the symmetric group 
$S_4$ or by Klein's $4$-group $\Kl$, we can enlarge the affine 
Weyl group $W(D_4^{(1)})$ to the affine Weyl group of 
type $F_4^{(1)}$ or to the extended affine Weyl group of 
type $D_4^{(1)}$, 
%%%%%%%
\[
W(F_4^{(1)}) = S_4 \ltimes W(D_4^{(1)}) \supset 
\wt{W}(D_4^{(1)}) = \Kl \ltimes W(D_4^{(1)}). 
\]
%%%%%%%
%%%%%%%%%%%%%%%%%%%%%%%%%%% def:Kl %%%%%%%%%%%%%%%%%%%%%%%%%%%
\begin{definition} \label{def:Kl}
Replacing the group $W(D_4^{(1)})$ with $W(F_4^{(1)})$ 
in Definition \ref{def:stratification}, we can define a 
coarser stratification of $\K$ than the 
$W(D_4^{(1)})$-stratification, called the 
$W(F_4^{(1)})$-{\sl stratification}. 
Moreover, replacing $W(D_4^{(1)})$ with $\wt{W}(D_4^{1})$, 
we can also think of a stratification of $\K$ intermediate 
between these two stratifications, 
called the $\wt{W}(D_4^{(1)})$-{\sl stratification}. 
\end{definition}
%%%%%%%%%%%%%%%%%%%%%%%%%%%%%%%%%%%%%%%%%%%%%%%%%%%%%%%%%%%%%%%%
\par 
The following is the classification of the $W(F_4^{(1)})$-strata 
and $\wt{W}(D_4^{(1)})$-strata. 
%%%%%%%%%%%%%%%%%%%%%%% fig:sklequiv %%%%%%%%%%%%%%%%%%%%%%%
\begin{figure}[t] 
\begin{center}
%WinTpicVersion2.15
\unitlength 0.1in
\begin{picture}(38.10,12.60)(2.70,-15.20)
% CIRCLE 0 0 3 0
% 4 448 852 448 900 448 900 448 900
% 
\special{pn 20}%
\special{ar 448 452 48 48  0.0000000 6.2831853}%
% CIRCLE 0 0 3 0
% 4 1248 844 1248 892 1248 892 1248 892
% 
\special{pn 20}%
\special{ar 1248 444 48 48  0.0000000 6.2831853}%
% CIRCLE 0 0 0 0
% 4 456 1652 456 1700 456 1700 456 1700
% 
\special{pn 20}%
\special{sh 0.600}%
\special{ar 456 1252 48 48  0.0000000 6.2831853}%
% CIRCLE 0 0 0 0
% 4 1256 1652 1256 1700 1256 1700 1256 1700
% 
\special{pn 20}%
\special{sh 0.600}%
\special{ar 1256 1252 48 48  0.0000000 6.2831853}%
% CIRCLE 0 0 0 0
% 4 848 1252 848 1300 848 1300 848 1300
% 
\special{pn 20}%
\special{sh 0.600}%
\special{ar 848 852 48 48  0.0000000 6.2831853}%
% LINE 0 2 3 0
% 4 480 900 816 1228 824 1228 824 1228
% 
\special{pn 20}%
\special{pa 480 500}%
\special{pa 816 828}%
\special{dt 0.054}%
\special{pa 816 828}%
\special{pa 816 828}%
\special{dt 0.054}%
\special{pa 824 828}%
\special{pa 824 828}%
\special{dt 0.054}%
% LINE 0 0 3 0
% 4 880 1292 1216 1620 1216 1620 1216 1620
% 
\special{pn 20}%
\special{pa 880 892}%
\special{pa 1216 1220}%
\special{fp}%
\special{pa 1216 1220}%
\special{pa 1216 1220}%
\special{fp}%
% LINE 0 2 3 0
% 2 1216 884 880 1212
% 
\special{pn 20}%
\special{pa 1216 484}%
\special{pa 880 812}%
\special{dt 0.054}%
\special{pa 880 812}%
\special{pa 880 812}%
\special{dt 0.054}%
% LINE 0 0 3 0
% 2 808 1292 504 1604
% 
\special{pn 20}%
\special{pa 808 892}%
\special{pa 504 1204}%
\special{fp}%
% STR 2 0 3 0
% 3 270 750 270 830 2 0
% $\scriptstyle{i}$
\put(2.7000,-4.3000){\makebox(0,0)[lb]{$\scriptstyle{i}$}}%
% STR 2 0 3 0
% 3 1350 760 1350 840 2 0
% $\scriptstyle{4}$
\put(13.5000,-4.4000){\makebox(0,0)[lb]{$\scriptstyle{4}$}}%
% STR 2 0 3 0
% 3 290 1730 290 1810 2 0
% $\scriptstyle{j}$
\put(2.9000,-14.1000){\makebox(0,0)[lb]{$\scriptstyle{j}$}}%
% STR 2 0 3 0
% 3 1350 1730 1350 1810 2 0
% $\scriptstyle{k}$
\put(13.5000,-14.1000){\makebox(0,0)[lb]{$\scriptstyle{k}$}}%
% STR 2 0 3 0
% 3 768 1060 768 1140 2 0
% $\scriptstyle{0}$
\put(7.6800,-7.4000){\makebox(0,0)[lb]{$\scriptstyle{0}$}}%
% CIRCLE 0 0 3 0
% 4 3166 852 3166 900 3166 900 3166 900
% 
\special{pn 20}%
\special{ar 3166 452 48 48  0.0000000 6.2831853}%
% CIRCLE 0 0 3 0
% 4 3966 844 3966 892 3966 892 3966 892
% 
\special{pn 20}%
\special{ar 3966 444 48 48  0.0000000 6.2831853}%
% CIRCLE 0 0 0 0
% 4 3174 1652 3174 1700 3174 1700 3174 1700
% 
\special{pn 20}%
\special{sh 0.600}%
\special{ar 3174 1252 48 48  0.0000000 6.2831853}%
% CIRCLE 0 0 0 0
% 4 3974 1652 3974 1700 3974 1700 3974 1700
% 
\special{pn 20}%
\special{sh 0.600}%
\special{ar 3974 1252 48 48  0.0000000 6.2831853}%
% CIRCLE 0 0 3 0
% 4 3566 1252 3566 1300 3566 1300 3566 1300
% 
\special{pn 20}%
\special{ar 3566 852 48 48  0.0000000 6.2831853}%
% LINE 0 2 3 0
% 4 3198 900 3534 1228 3542 1228 3542 1228
% 
\special{pn 20}%
\special{pa 3198 500}%
\special{pa 3534 828}%
\special{dt 0.054}%
\special{pa 3534 828}%
\special{pa 3534 828}%
\special{dt 0.054}%
\special{pa 3542 828}%
\special{pa 3542 828}%
\special{dt 0.054}%
% LINE 0 2 3 0
% 4 3598 1292 3934 1620 3934 1620 3934 1620
% 
\special{pn 20}%
\special{pa 3598 892}%
\special{pa 3934 1220}%
\special{dt 0.054}%
\special{pa 3934 1220}%
\special{pa 3934 1220}%
\special{dt 0.054}%
\special{pa 3934 1220}%
\special{pa 3934 1220}%
\special{dt 0.054}%
% LINE 0 2 3 0
% 2 3934 884 3598 1212
% 
\special{pn 20}%
\special{pa 3934 484}%
\special{pa 3598 812}%
\special{dt 0.054}%
\special{pa 3598 812}%
\special{pa 3598 812}%
\special{dt 0.054}%
% LINE 0 2 3 0
% 2 3526 1292 3222 1604
% 
\special{pn 20}%
\special{pa 3526 892}%
\special{pa 3222 1204}%
\special{dt 0.054}%
\special{pa 3222 1204}%
\special{pa 3222 1204}%
\special{dt 0.054}%
% STR 2 0 3 0
% 3 2960 750 2960 830 2 0
% $\scriptstyle{i}$
\put(29.6000,-4.3000){\makebox(0,0)[lb]{$\scriptstyle{i}$}}%
% STR 2 0 3 0
% 3 4080 760 4080 840 2 0
% $\scriptstyle{4}$
\put(40.8000,-4.4000){\makebox(0,0)[lb]{$\scriptstyle{4}$}}%
% STR 2 0 3 0
% 3 3020 1710 3020 1790 2 0
% $\scriptstyle{j}$
\put(30.2000,-13.9000){\makebox(0,0)[lb]{$\scriptstyle{j}$}}%
% STR 2 0 3 0
% 3 4070 1720 4070 1800 2 0
% $\scriptstyle{k}$
\put(40.7000,-14.0000){\makebox(0,0)[lb]{$\scriptstyle{k}$}}%
% STR 2 0 3 0
% 3 3494 1060 3494 1140 2 0
% $\scriptstyle{0}$
\put(34.9400,-7.4000){\makebox(0,0)[lb]{$\scriptstyle{0}$}}%
% STR 2 0 3 0
% 3 580 2010 580 2090 2 0
% $(A_3)_i$
\put(5.8000,-16.9000){\makebox(0,0)[lb]{$(A_3)_i$}}%
% STR 2 0 3 0
% 3 3400 2000 3400 2080 2 0
% $(A_1^{\oplus 2})_i$
\put(34.0000,-16.8000){\makebox(0,0)[lb]{$(A_1^{\oplus 2})_i$}}%
\end{picture}%
\end{center}
\caption{$\wt{W}(D_4^{(1)})$-strata $(A_3)_i$ and 
$(A_1^{\oplus 2})_i$} 
\label{fig:sklequiv} 
\end{figure}
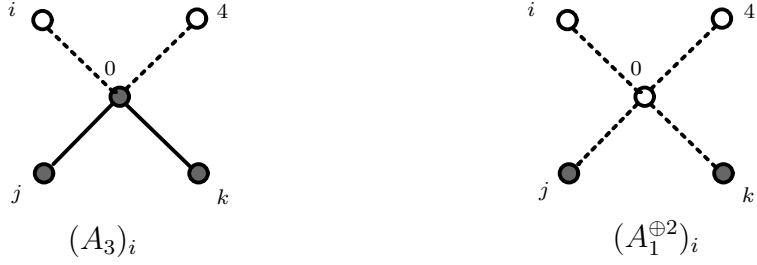
%%%%%%%%%%%%%%%%%%%%%%% lem:equiv %%%%%%%%%%%%%%%%%%%%%%%%%%
\begin{lemma} \label{lem:equiv} 
For each abstract Dynkin type $*$ in Table $\ref{tab:type}$, 
there is a unique $W(F_4^{(1)})$-stratum of type $*$. 
As for the $\wt{W}(D_4^{(1)})$-strata, we have the 
following classification $($see also 
Figure $\ref{fig:sklequiv})$. 
\begin{enumerate}
\item For $* \in \{D_4, A_1^{\oplus 4}, A_1^{\oplus 3}, 
A_2, A_1, \emptyset \}$, there is a unique 
$\wt{W}(D_4^{(1)})$-stratum of abstract type $*$ and 
this unique stratum is denoted by the same symbol $*$.  
\item For $* \in \{A_3, A_1^{\oplus 2}\}$, there are 
exactly three $\wt{W}(D_4^{(1)})$-strata of abstract 
type $*;$ 
\begin{enumerate}
\item for $* = A_3$, the stratum $(A_3)_i$ represented 
by $I = \{0,j,k\}$ with $\{i,j,k\} = \{1,2,3\};$
\item for $* = A_1^{\oplus 2}$, the stratum 
$(A_1^{\oplus 2})_i$ represented by 
$I = \{j,k\}$ with $\{i,j,k\} = \{1,2,3\}$. 
\end{enumerate} 
\end{enumerate} 
\end{lemma} 
%%%%%%%%%%%%%%%%%%%%%%%%%%%%%%%%%%%%%%%%%%%%%%%%%%%%%%%%%%%%%%%%
\par
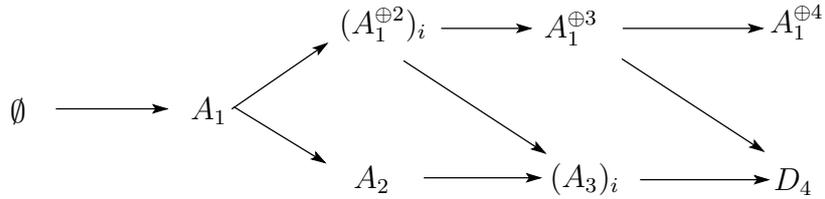
\begin{figure}[t] 
\begin{center}
%WinTpicVersion2.15
\unitlength 0.1in
\begin{picture}(39.55,9.05)(6.70,-13.65)
% STR 2 0 3 0
% 3 1700 1300 1700 1400 5 0
% $A_1$
\put(17.0000,-10.0000){\makebox(0,0){$A_1$}}%
% STR 2 0 3 0
% 3 670 1370 670 1470 2 0
% $\emptyset$
\put(6.7000,-10.7000){\makebox(0,0)[lb]{$\emptyset$}}%
% STR 2 0 3 0
% 3 2410 900 2410 1000 2 0
% 
\put(24.1000,-6.0000){\makebox(0,0)[lb]{}}%
% STR 2 0 3 0
% 3 2450 1730 2450 1830 2 0
% $A_2$
\put(24.5000,-14.3000){\makebox(0,0)[lb]{$A_2$}}%
% STR 2 0 3 0
% 3 2375 940 2375 1040 2 0
% $(A_1^{\oplus 2})_i$
\put(23.7500,-6.4000){\makebox(0,0)[lb]{$(A_1^{\oplus 2})_i$}}%
% STR 2 0 3 0
% 3 3465 1740 3465 1840 2 0
% $(A_3)_i$
\put(34.6500,-14.4000){\makebox(0,0)[lb]{$(A_3)_i$}}%
% STR 2 0 3 0
% 3 3435 960 3435 1060 2 0
% $A_1^{\oplus 3}$
\put(34.3500,-6.6000){\makebox(0,0)[lb]{$A_1^{\oplus 3}$}}%
% STR 2 0 3 0
% 3 4625 1740 4625 1840 2 0
% $D_4$
\put(46.2500,-14.4000){\makebox(0,0)[lb]{$D_4$}}%
% STR 2 0 3 0
% 3 4605 930 4605 1030 2 0
% $A_1^{\oplus 4}$
\put(46.0500,-6.3000){\makebox(0,0)[lb]{$A_1^{\oplus 4}$}}%
% VECTOR 2 0 3 0
% 2 1820 1390 2310 1050
% 
\special{pn 8}%
\special{pa 1820 990}%
\special{pa 2310 650}%
\special{fp}%
\special{sh 1}%
\special{pa 2310 650}%
\special{pa 2244 672}%
\special{pa 2266 680}%
\special{pa 2267 704}%
\special{pa 2310 650}%
\special{fp}%
% VECTOR 2 0 3 0
% 2 1830 1380 2300 1670
% 
\special{pn 8}%
\special{pa 1830 980}%
\special{pa 2300 1270}%
\special{fp}%
\special{sh 1}%
\special{pa 2300 1270}%
\special{pa 2254 1218}%
\special{pa 2255 1242}%
\special{pa 2233 1252}%
\special{pa 2300 1270}%
\special{fp}%
% VECTOR 2 0 3 0
% 2 2815 1750 3415 1750
% 
\special{pn 8}%
\special{pa 2815 1350}%
\special{pa 3415 1350}%
\special{fp}%
\special{sh 1}%
\special{pa 3415 1350}%
\special{pa 3348 1330}%
\special{pa 3362 1350}%
\special{pa 3348 1370}%
\special{pa 3415 1350}%
\special{fp}%
% VECTOR 2 0 3 0
% 2 3935 1760 4595 1760
% 
\special{pn 8}%
\special{pa 3935 1360}%
\special{pa 4595 1360}%
\special{fp}%
\special{sh 1}%
\special{pa 4595 1360}%
\special{pa 4528 1340}%
\special{pa 4542 1360}%
\special{pa 4528 1380}%
\special{pa 4595 1360}%
\special{fp}%
% VECTOR 2 0 3 0
% 2 2905 970 3365 970
% 
\special{pn 8}%
\special{pa 2905 570}%
\special{pa 3365 570}%
\special{fp}%
\special{sh 1}%
\special{pa 3365 570}%
\special{pa 3298 550}%
\special{pa 3312 570}%
\special{pa 3298 590}%
\special{pa 3365 570}%
\special{fp}%
% VECTOR 2 0 3 0
% 2 3845 970 4565 970
% 
\special{pn 8}%
\special{pa 3845 570}%
\special{pa 4565 570}%
\special{fp}%
\special{sh 1}%
\special{pa 4565 570}%
\special{pa 4498 550}%
\special{pa 4512 570}%
\special{pa 4498 590}%
\special{pa 4565 570}%
\special{fp}%
% VECTOR 2 0 3 0
% 2 2710 1140 3440 1620
% 
\special{pn 8}%
\special{pa 2710 740}%
\special{pa 3440 1220}%
\special{fp}%
\special{sh 1}%
\special{pa 3440 1220}%
\special{pa 3395 1167}%
\special{pa 3395 1191}%
\special{pa 3373 1200}%
\special{pa 3440 1220}%
\special{fp}%
% VECTOR 2 0 3 0
% 2 3840 1130 4570 1610
% 
\special{pn 8}%
\special{pa 3840 730}%
\special{pa 4570 1210}%
\special{fp}%
\special{sh 1}%
\special{pa 4570 1210}%
\special{pa 4525 1157}%
\special{pa 4525 1181}%
\special{pa 4503 1190}%
\special{pa 4570 1210}%
\special{fp}%
% VECTOR 2 0 3 0
% 4 910 1390 1480 1390 1480 1390 1480 1390
% 
\special{pn 8}%
\special{pa 910 990}%
\special{pa 1480 990}%
\special{fp}%
\special{sh 1}%
\special{pa 1480 990}%
\special{pa 1413 970}%
\special{pa 1427 990}%
\special{pa 1413 1010}%
\special{pa 1480 990}%
\special{fp}%
\special{pa 1480 990}%
\special{pa 1480 990}%
\special{fp}%
\end{picture}%
\end{center}
\caption{Adjacency relations among $\wt{W}(D_4^{(1)})$-strata 
$(i = 1,2,3)$} 
\label{fig:adjac}
\end{figure}

 
If something about the transformation $g_j^2$ is discussed for 
a {\sl fixed} index $j$, the relevant stratification is the 
$\wt{W}(D_4^{(1)})$-stratification. 
Namely we may discuss the issue on each 
$\wt{W}(D_4^{(1)})$-stratum, choosing any representative of 
each $\wt{W}(D_4^{(1)})$-orbit, since in the commutative 
diagram (\ref{cd:S4}) we have $\tau(j) = j$ and hence 
$g_{\tau(j)}^2 = g_j^2$ for every $\si \in \Kl$ (see 
also Remark \ref{rem:backlund}). 
For two $\wt{W}(D_4^{(1)})$-strata, say $*$ and $**$, we 
write $* \to **$ if the stratum $**$ lies on the boundary 
of the stratum $*$. 
All the possible adjacency relations $* \to **$ are 
depicted in Figure \ref{fig:adjac}. 
Note that there are no adjacency relations between 
$(A_1^{\oplus 2})_i$ and $(A_3)_j$ for any distinct 
$i$, $j \in \{1,2,3\}$. 
%%%%%%%%%%%%%%%%%%%%%%% fig:adjac %%%%%%%%%%%%%%%%%%%%%%%%%%%%%
%%%%%%%%%%%%%%%%%%%%%%%%%%%%%%%%%%%%%%%%%%%%%%%%%%%%%%%%%%%%%%%%
%%%%%%%%%%%%%%%%%%%%%%% sec:backlund %%%%%%%%%%%%%%%%%%%%%%%%%%%
\section{B\"acklund Transformations} \label{sec:backlund} 
%%%%%%%%%%%%%%%%%%%%%%%%%%%%%%%%%%%%%%%%%%%%%%%%%%%%%%%%%%%%%%%%
In this section we briefly discuss B\"acklund transformations, 
especially the characterization of them in terms of 
Riemann-Hilbert correspondence \cite{IIS0,IIS1}. 
This topic is included here in order to confirm that 
our problem may be treated modulo B\"acklund transformations. 
\par 
For each $\si \in S_4$ we define the isomorphism of affine 
cubic surfaces 
%%%%%
\[
\si : \Sol(\th) \to \Sol(\si(\th)), \qquad 
(x_1,x_2,x_3) \mapsto 
(\ve_{\tau(1)} x_{\tau(1)}, \, \ve_{\tau(2)} x_{\tau(2)}, \, 
\ve_{\tau(3)} x_{\tau(3)}), 
\]
%%%%%
where $\si \in S_4$ is identified with 
$(\tau, \ve) \in S_3 \ltimes \Kl$ via the isomorphism 
(\ref{eqn:equiv}). 
Consider the natural homomorphism 
$W(F_4^{(1)}) = S_4 \ltimes W(D_4^{(1)}) \to S_4$, 
$w \mapsto \si$. 
Since the Riemann-Hilbert correspondence (\ref{eqn:RHkth}) 
is an analytic minimal resolution of singularities, 
for each $w \in W(F_4^{(1)})$, there exists an analytic 
isomorphism $w : \M_z(\k) \to \M_{z}(w(\k))$ such that 
the diagram 
%%%%%%%%%%%%%%%%%%%%%%%% cd:backlund %%%%%%%%%%%%%%%%%%%%%%
\begin{equation} \label{cd:backlund}
\begin{CD}
\M_z(\k) @> w >> \M_z(w(\k)) \\
@V \RH_{z,\k} VV  @VV \RH_{z,w(\k)} V \\
\Sol(\th) @>> \si > \Sol(\si(\th))
\end{CD}
\end{equation}
%%%%%%%%%%%%%%%%%%%%%%%%%%%%%%%%%%%%%%%%%%%%%%%%%%%%%%%%%%%%
is commutative, for any fixed $\k \in \K$ with 
$\th = \rh(\k) \in \Th$. 
\par 
The commutative diagram (\ref{cd:backlund}) characterizes 
the B\"acklund transformations of Painlev\'e VI. 
Namely the map $w : \M_z(\k) \to \M_z(w(\k))$ turns out to 
be algebraic and there are suitable affine coordinates on 
$\M_z(\k)$ and $\M_z(w(\k))$ in terms of which the map $w$ 
can be represented by the usual formula for B\"acklund 
transformations known as birational canonical transforamtions 
\cite{Okamoto2} (see \cite{IIS0,IIS1} for the precise 
statement). 
In other words the Riemann-Hilbert correspondence is 
equivariant under the B\"acklund transformations and so 
is our main problem. 
%%%%%%%%%%%%%%%%%%%%%%% rem:backlund %%%%%%%%%%%%%%%%%%%%%%%%
\begin{remark} \label{rem:backlund} 
The $S_4$-factor of $W(F_4^{(1)}) = S_4 \ltimes W(D_4^{(1)})$ 
or more strictly the $S_3$-factor of $S_4 = S_3 \ltimes \Kl$ 
permutes the three fixed singular points $0$, $1$ and 
$\infty$, while they are fixed by $\wt{W}(D_4^{(1)}) 
= \Kl \ltimes W(D_4^{(1)})$. 
Hence we may consider our problem only around the origin 
$z = 0$ and, upon restricting our attention to $z = 0$, 
we may discuss it modulo the B\"acklund action of 
$\wt{W}(D_4^{(1)})$. 
\end{remark} 
%%%%%%%%%%%%%%%%%%%%%%% sec:fixed %%%%%%%%%%%%%%%%%%%%%%%%%%%
\section{Fixed Points and Periodic Points} \label{sec:fixed} 
%%%%%%%%%%%%%%%%%%%%%%%%%%%%%%%%%%%%%%%%%%%%%%%%%%%%%%%%%%%%%
We shall more closely investigate the fixed point 
set $\mathrm{Fix}_j(\th)$, or rather its subset 
$\mathrm{Fix}_j^{\circ}(\th) \cong 
\wt{\mathrm{Fix}}_j^{\circ}(\th)$ of smooth fixed points, 
by solving the system of equations (\ref{eqn:per1}). 
In view of (\ref{eqn:yi}) the last two equations in 
(\ref{eqn:per1}) are expressed as a linear system for 
the unknowns $(x_j,x_k)$, 
%%%%%%%%%%%%%%%%%%%%%% eqn:system %%%%%%%%%%%%%%%%%%%%%%%%%%%
\begin{equation} \label{eqn:system} 
\left\{ 
\begin{array}{rcl} 
2 x_j + x_i x_k &=& \th_j, \\[2mm] 
x_i x_j + 2 x_k &=& \th_k, 
\end{array} 
\right. 
\end{equation} 
%%%%%%%%%%%%%%%%%%%%%%%%%%%%%%%%%%%%%%%%%%%%%%%%%%%%%%%%%%%%%%%% 
If its determinant $4 - x_i^2$ is nonzero, then system 
(\ref{eqn:system}) is uniquely settled as 
%%%%%%%%%%%%%%%%%%%%%% eqn:x %%%%%%%%%%%%%%%%%%%%%%%%%%%%%%%%%%% 
\begin{equation} \label{eqn:x}
x_j = \dfrac{2 \th_j - x_i \th_k}{4 - x_i^2}, \qquad 
x_k = \dfrac{2 \th_k - x_i \th_j}{4 - x_i^2}.  
\end{equation} 
%%%%%%%%%%%%%%%%%%%%%%%%%%%%%%%%%%%%%%%%%%%%%%%%%%%%%%%%%%%%%%%% 
Substituting (\ref{eqn:x}) into equation $f(x,\th) = 0$ yields 
a quartic equation for the unknown $x_i$, 
%%%%%%%%%%%%%%%%%%%%%%%% eqn:quartic %%%%%%%%%%%%%%%%%%%%%%%%%%%%
\begin{equation} \label{eqn:quartic}
x_i^4 - \th_i x_i^3 + (\th_4-4) x_i^2 + 
(4 \th_i - \th_j \th_k) x_i + \th_j^2 + \th_k^2 - 4 \th_4 = 0. 
\end{equation} 
%%%%%%%%%%%%%%%%%%%%%%%%%%%%%%%%%%%%%%%%%%%%%%%%%%%%%%%%%%%%%%%%%%
Conversely, if $x_i$ is a root of equation (\ref{eqn:quartic}) 
with nonzero $x_i^2- 4$, then subsituting 
this into formula (\ref{eqn:x}) yields a root of system 
(\ref{eqn:per1}). 
The four roots of quartic equation (\ref{eqn:quartic}) are 
given by 
%%%%%%
\[
F(b_i,b_4;b_j,b_k), \quad F(b_i,b_4^{-1};b_j,b_k), \quad 
F(b_j,b_k;b_i,b_4), \quad  F(b_j,b_k^{-1};b_i,b_4), 
\]
%%%%%
counted with multiplicities, where $F(b_i,b_4;b_j,b_k)$ 
is defined by 
%%%%%%%%%%%%%%%%%%%%%%%%%% eqn:F %%%%%%%%%%%%%%%%%%%%%%%%%%%%%%%%
\begin{equation} \label{eqn:F} 
F(b_i,b_4;b_j,b_k) = b_i b_4 + b_i^{-1} b_4^{-1}. 
\end{equation}
%%%%%%%%%%%%%%%%%%%%%%%%%%%%%%%%%%%%%%%%%%%%%%%%%%%%%%%%%%%%%%%%%
\par 
We pick up the root $x_i = F(b_i,b_4;b_j,b_k)$. 
Note that $F(b_i,b_4;b_j,b_k)^2 - 4$ is nonzero precisely when 
$b_i^2b_4^2 \neq 1$. 
If this is the case, then substituting $x_i = F(b_i,b_4;b_j,b_k)$ 
into formula (\ref{eqn:x}) yields $x_j = G(b_i,b_4;b_j,b_k)$ and 
$x_k = G(b_i,b_4;b_k,b_j)$, where $G(b_i,b_4;b_j,b_k)$ is 
defined by 
%%%%%%%%%%%%%%%%%%%%%%%%%% eqn:G %%%%%%%%%%%%%%%%%%%%%%%%%%%%%%%%%
\begin{equation} \label{eqn:G} 
G(b_i,b_4;b_j,b_k) = 
\dfrac{(b_i+b_4)(b_j+b_k)(b_j b_k+1)}{2(b_ib_4+1)b_jb_k} 
+ \dfrac{(b_i-b_4)(b_j-b_k)(b_j b_k-1)}{2(b_ib_4-1)b_jb_k}.  
\end{equation} 
%%%%%%%%%%%%%%%%%%%%%%%%%%%%%%%%%%%%%%%%%%%%%%%%%%%%%%%%%%%%%%%%%%
Therefore, if $P(b_i,b_4;b_j,b_k)$ denotes the point defined by 
%%%%%%%%%%%%%%%%%%%%%%%%%%%
\[
x_i = F(b_i,b_4;b_j,b_k), \qquad 
x_j = G(b_i,b_4;b_j,b_k), \qquad 
x_k = G(b_i,b_4;b_k,b_j), 
\]
%%%%%%%%%%%%%%%%%%%%%%%%%%%
then $x = P(b_i,b_4;b_j,b_k)$ gives a root of system 
(\ref{eqn:per1}) with nonzero $x_i^2-4$ provided that 
$b_i^2b_4^2 \neq 1$. 
If $x$ is at this root, then $y_i(x,\th)$ admits the following 
nice factorization 
%%%%%%%%%%%%%%%%%%%%%%%%%%%% eqn:Y %%%%%%%%%%%%%%%%%%%%%%%%%%%%%%%
\begin{equation} \label{eqn:Y}
\begin{array}{rcrl}
y_i(x,\th) &=& 
\dfrac{(b_i-b_i^{-1})(b_4-b_4^{-1})}{(b_i^2b_4^2-1)^2} & 
\displaystyle \prod_{(\ve_j, \ve_k) \in \{\pm1\}^2}
(b_i b_j^{\ve_j} b_k^{\ve_k} b_4-1) \\[8mm]
&=& (b_ib_4-b_i^{-1}b_4^{-1})^{-2} & 
\{F(b_i,b_4;b_j,b_k)-F(b_i,b_4^{-1};b_j,b_k)\}
\\[1mm]
&&& \{F(b_i,b_4;b_j,b_k)-F(b_j,b_k \,\,\,;b_i,b_4)\} 
\\[1mm]
&&& \{F(b_i,b_4;b_j,b_k)-F(b_j,b_k^{-1};b_i,b_4)\}. 
\end{array}
\end{equation} 
%%%%%%%%%%%%%%%%%%%%%%%%%%%%%%%%%%%%%%%%%%%%%%%%%%%%%%%%%%%%%%%%%
Hence $P(b_i,b_4;b_j,b_k)$ is a smooth point of $\Sol(\th)$ if 
and only if $F(b_i,b_4;b_j,b_k)$ is a simple root of equation 
(\ref{eqn:quartic}). 
In terms of $\k \in\K$, the existence and smoothness conditions 
for $P(b_i,b_4;b_j,b_k)$ are given by 
$\k_i + \k_4 \not\in \Z$ and $\k_i \not\in \Z$, 
$\k_4 \not\in \Z$, $\k_i+\k_4 \pm \k_j \pm \k_k \not\in 2\Z+1$, 
respectively. 
%%%%%%%%%%%%%%%%%%%%%%% lem:fixedpt %%%%%%%%%%%%%%%%%%%%%%%%%%%%
\begin{lemma} \label{lem:fixedpt}
The smooth fixed points $x \in \mathrm{Fix}_j^{\circ}(\th)$ 
with nonzero $x_i^2 - 4$ are precisely those points in 
Table $\ref{tab:fixed}$ which satisfy the existence and 
smoothness conditions mentioned there. 
\end{lemma} 
%%%%%%%%%%%%%%%%%%%%%%%%%%%%%%%%%%%%%%%%%%%%%%%%%%%%%%%%%%%%%%%%
%%%%%%%%%%%%%%%%%%%%%%%%%%% tab:fixed %%%%%%%%%%%%%%%%%%%%%%%%%%
\begin{table}[t] 
\begin{center} 
\begin{tabular}{|c||c|c|c|}
\hline
\vspace{-3mm} & & & \\
label & fixed point & existence & smoothness condition \\[2mm]
\hline
\hline 
\vspace{-3mm} & & & \\
1 & $P(b_i,b_4;b_j,b_k)$ &$\k_i+\k_4 \not\in \Z$ & 
$\k_i \not\in \Z$, \, $\k_4 \not\in \Z$, \, 
$\k_i+\k_4 \pm \k_j \pm \k_k \not\in 2 \Z+1$ \\[2mm]
\hline 
\vspace{-3mm} & & & \\
2 & $P(b_i,b_4^{-1};b_j,b_k)$ & $\k_i-\k_4 \not\in \Z$ & 
$\k_i \not\in \Z$, \, $\k_4 \not\in \Z$, \, 
$\k_i-\k_4 \pm \k_j \pm \k_k \not\in 2 \Z+1$ \\[2mm]
\hline 
\vspace{-3mm} & & & \\
3 & $P(b_j,b_k;b_i,b_4)$ & $\k_j+\k_k \not\in \Z$ & 
$\k_j \not\in \Z$, \, $\k_k \not\in \Z$, \, 
$\k_j+\k_k \pm \k_i \pm \k_4 \not\in 2 \Z+1$ \\[2mm]
\hline 
\vspace{-3mm} & & & \\
4 & $P(b_j,b_k^{-1};b_i,b_4)$ & $\k_j-\k_k \not\in \Z$ & 
$\k_j \not\in \Z$, \, $\k_k \not\in \Z$, \, 
$\k_j-\k_k \pm \k_i \pm \k_4 \not\in 2 \Z+1$ \\[2mm]
\hline 
\end{tabular}
\end{center} 
\caption{Smooth fixed points $x \in 
\mathrm{Fix}_j^{\circ}(\th)$ with nonzero $x_i^2-4$}
\label{tab:fixed} 
\end{table}
%%%%%%%%%%%%%%%%%%%%%%%%%%%%%%%%%%%%%%%%%%%%%%%%%%%%%%%%%%%%%
\par 
The fixed points in Table \ref{tab:fixed} is closely 
related to the configuration of lines on the affine 
cubic surface $\Sol(\th)$ or on its compactification  
$\ol{\Sol}(\th)$ by the standard embedding 
%%%%%%%
\[
\Sol(\th) \hookrightarrow \ol{\Sol}(\th) \subset \P^3, 
\qquad x = (x_1,x_2,x_3) \mapsto [1:x_1:x_2:x_3], 
\]
%%%%%%%
where the projective cubic surface $\ol{\Sol}(\th)$ 
is defined by the homogeneous equation
%%%%%%%
\[
F(X,\th) := X_1X_2X_3 + X_0(X_1^2 + X_2^2 + X_3^2) - 
X_0^2(\th_1 X_1 + \th_2 X_2 + \th_3 X_3) + \th_4 X_0^3 = 0. 
\]
%%%%%%%
It is obtained from the affine surface $\Sol(\th)$ by 
adding three lines at infinity 
%%%%%%%%
\[
L_i = \{\, X \in \P^3 \,:\, X_0 = X_i = 0 \,\} 
\qquad (i = 1,2,3), 
\]
%%%%%%%
whose union $L = L_1 \cup L_2 \cup L_3$ is called the 
tritangent lines at infinity. 
\par
\begin{figure}[t]
\begin{center}
%WinTpicVersion2.15
\unitlength 0.1in
\begin{picture}(49.30,35.20)(8.60,-35.60)
% STR 2 0 3 0
% 3 860 3460 860 3560 2 0
% $L_1$
\put(8.6000,-31.6000){\makebox(0,0)[lb]{$L_1$}}%
% STR 2 0 3 0
% 3 5590 3990 5590 4090 2 0
% $L_2$
\put(55.9000,-36.9000){\makebox(0,0)[lb]{$L_2$}}%
% STR 2 0 3 0
% 3 3700 510 3700 610 2 0
% $L_3$
\put(37.0000,-2.1000){\makebox(0,0)[lb]{$L_3$}}%
% LINE 0 0 3 0
% 2 1170 3480 5790 3480
% 
\special{pn 20}%
\special{pa 1170 3080}%
\special{pa 5790 3080}%
\special{fp}%
% LINE 0 0 3 0
% 4 3790 670 1380 3910 3180 670 5570 3870
% 
\special{pn 20}%
\special{pa 3790 270}%
\special{pa 1380 3510}%
\special{fp}%
\special{pa 3180 270}%
\special{pa 5570 3470}%
\special{fp}%
% LINE 1 0 3 0
% 16 2090 3320 2470 3920 2480 3320 2090 3910 2890 3320 3280 3910 3280 3320 2890 3910 3680 3310 4080 3920 4090 3320 3680 3920 4480 3330 4870 3920 4880 3330 4490 3920
% 
\special{pn 13}%
\special{pa 2090 2920}%
\special{pa 2470 3520}%
\special{fp}%
\special{pa 2480 2920}%
\special{pa 2090 3510}%
\special{fp}%
\special{pa 2890 2920}%
\special{pa 3280 3510}%
\special{fp}%
\special{pa 3280 2920}%
\special{pa 2890 3510}%
\special{fp}%
\special{pa 3680 2910}%
\special{pa 4080 3520}%
\special{fp}%
\special{pa 4090 2920}%
\special{pa 3680 3520}%
\special{fp}%
\special{pa 4480 2930}%
\special{pa 4870 3520}%
\special{fp}%
\special{pa 4880 2930}%
\special{pa 4490 3520}%
\special{fp}%
% LINE 1 0 3 0
% 16 3500 1240 4187 1168 3726 1550 3959 873 3963 1875 4648 1818 4189 2184 4422 1508 4412 2508 5118 2446 4657 2826 4887 2128 4891 3131 5575 3072 5122 3447 5355 2771
% 
\special{pn 13}%
\special{pa 3500 840}%
\special{pa 4187 768}%
\special{fp}%
\special{pa 3726 1150}%
\special{pa 3959 473}%
\special{fp}%
\special{pa 3963 1475}%
\special{pa 4648 1418}%
\special{fp}%
\special{pa 4189 1784}%
\special{pa 4422 1108}%
\special{fp}%
\special{pa 4412 2108}%
\special{pa 5118 2046}%
\special{fp}%
\special{pa 4657 2426}%
\special{pa 4887 1728}%
\special{fp}%
\special{pa 4891 2731}%
\special{pa 5575 2672}%
\special{fp}%
\special{pa 5122 3047}%
\special{pa 5355 2371}%
\special{fp}%
% LINE 1 0 3 0
% 16 2970 920 3242 1555 2741 1227 3457 1249 2501 1550 2758 2187 2272 1857 2987 1879 2029 2166 2297 2822 1798 2494 2533 2508 1576 2807 1835 3444 1342 3122 2057 3145
% 
\special{pn 13}%
\special{pa 2970 520}%
\special{pa 3242 1155}%
\special{fp}%
\special{pa 2741 827}%
\special{pa 3457 849}%
\special{fp}%
\special{pa 2501 1150}%
\special{pa 2758 1787}%
\special{fp}%
\special{pa 2272 1457}%
\special{pa 2987 1479}%
\special{fp}%
\special{pa 2029 1766}%
\special{pa 2297 2422}%
\special{fp}%
\special{pa 1798 2094}%
\special{pa 2533 2108}%
\special{fp}%
\special{pa 1576 2407}%
\special{pa 1835 3044}%
\special{fp}%
\special{pa 1342 2722}%
\special{pa 2057 2745}%
\special{fp}%
% STR 2 0 3 0
% 3 5620 3050 5620 3150 2 0
% $L_{21}^{+}$
\put(56.2000,-27.5000){\makebox(0,0)[lb]{$L_{21}^{+}$}}%
% STR 2 0 3 0
% 3 5400 2690 5400 2790 2 0
% $L_{21}^{-}$
\put(54.0000,-23.9000){\makebox(0,0)[lb]{$L_{21}^{-}$}}%
% STR 2 0 3 0
% 3 5160 2410 5160 2510 2 0
% $L_{22}^{+}$
\put(51.6000,-21.1000){\makebox(0,0)[lb]{$L_{22}^{+}$}}%
% STR 2 0 3 0
% 3 4920 2020 4920 2120 2 0
% $L_{22}^{-}$
\put(49.2000,-17.2000){\makebox(0,0)[lb]{$L_{22}^{-}$}}%
% STR 2 0 3 0
% 3 4700 1780 4700 1880 2 0
% $L_{23}^{+}$
\put(47.0000,-14.8000){\makebox(0,0)[lb]{$L_{23}^{+}$}}%
% STR 2 0 3 0
% 3 4440 1440 4440 1540 2 0
% $L_{23}^{-}$
\put(44.4000,-11.4000){\makebox(0,0)[lb]{$L_{23}^{-}$}}%
% STR 2 0 3 0
% 3 4220 1120 4220 1220 2 0
% $L_{24}^{+}$
\put(42.2000,-8.2000){\makebox(0,0)[lb]{$L_{24}^{+}$}}%
% STR 2 0 3 0
% 3 3990 820 3990 920 2 0
% $L_{24}^{-}$
\put(39.9000,-5.2000){\makebox(0,0)[lb]{$L_{24}^{-}$}}%
% STR 2 0 3 0
% 3 2820 800 2820 900 2 0
% $L_{31}^{+}$
\put(28.2000,-5.0000){\makebox(0,0)[lb]{$L_{31}^{+}$}}%
% STR 2 0 3 0
% 3 2500 1200 2500 1300 2 0
% $L_{31}^{-}$
\put(25.0000,-9.0000){\makebox(0,0)[lb]{$L_{31}^{-}$}}%
% STR 2 0 3 0
% 3 2380 1430 2380 1530 2 0
% $L_{32}^{+}$
\put(23.8000,-11.3000){\makebox(0,0)[lb]{$L_{32}^{+}$}}%
% STR 2 0 3 0
% 3 2030 1830 2030 1930 2 0
% $L_{32}^{-}$
\put(20.3000,-15.3000){\makebox(0,0)[lb]{$L_{32}^{-}$}}%
% STR 2 0 3 0
% 3 1810 2040 1810 2140 2 0
% $F_{33}^{+}$
\put(18.1000,-17.4000){\makebox(0,0)[lb]{$F_{33}^{+}$}}%
% STR 2 0 3 0
% 3 1540 2450 1540 2550 2 0
% $L_{33}^{-}$
\put(15.4000,-21.5000){\makebox(0,0)[lb]{$L_{33}^{-}$}}%
% STR 2 0 3 0
% 3 1460 2680 1460 2780 2 0
% $L_{34}^{+}$
\put(14.6000,-23.8000){\makebox(0,0)[lb]{$L_{34}^{+}$}}%
% STR 2 0 3 0
% 3 1110 3100 1110 3200 2 0
% $L_{34}^{-}$
\put(11.1000,-28.0000){\makebox(0,0)[lb]{$L_{34}^{-}$}}%
% STR 2 0 3 0
% 3 1990 4010 1990 4110 2 0
% $L_{11}^+$
\put(19.9000,-37.1000){\makebox(0,0)[lb]{$L_{11}^+$}}%
% STR 2 0 3 0
% 3 2390 4030 2390 4130 2 0
% $L_{11}^-$
\put(23.9000,-37.3000){\makebox(0,0)[lb]{$L_{11}^-$}}%
% STR 2 0 3 0
% 3 2820 4030 2820 4130 2 0
% $L_{12}^+$
\put(28.2000,-37.3000){\makebox(0,0)[lb]{$L_{12}^+$}}%
% STR 2 0 3 0
% 3 3210 4020 3210 4120 2 0
% $L_{12}^-$
\put(32.1000,-37.2000){\makebox(0,0)[lb]{$L_{12}^-$}}%
% STR 2 0 3 0
% 3 3620 4010 3620 4110 2 0
% $L_{13}^+$
\put(36.2000,-37.1000){\makebox(0,0)[lb]{$L_{13}^+$}}%
% STR 2 0 3 0
% 3 4030 4010 4030 4110 2 0
% $L_{13}^-$
\put(40.3000,-37.1000){\makebox(0,0)[lb]{$L_{13}^-$}}%
% STR 2 0 3 0
% 3 4440 4010 4440 4110 2 0
% $L_{14}^+$
\put(44.4000,-37.1000){\makebox(0,0)[lb]{$L_{14}^+$}}%
% STR 2 0 3 0
% 3 4840 4020 4840 4120 2 0
% $L_{14}^-$
\put(48.4000,-37.2000){\makebox(0,0)[lb]{$L_{14}^-$}}%
\end{picture}%
\end{center}
\caption{The $27$ lines on a smooth cubic surface viewed from 
the tritangent lines at infinity} 
\label{fig:lines} 
\end{figure}
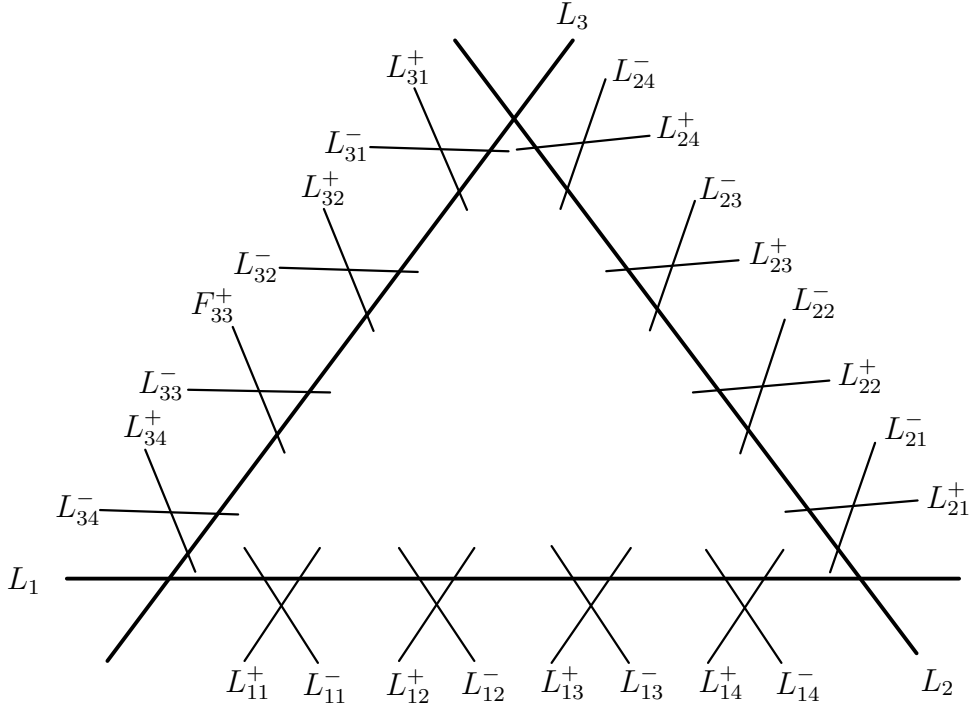

 
%%%%%%%%%%%%%%%%%%%%%%% tab:line %%%%%%%%%%%%%%%%%%%%%%%%%%%%%
\begin{table}[b] 
\begin{center} 
\begin{tabular}{|c|l|l|}
\hline
\vspace{-4mm} & & \\
$1$ & $L_{i1}^{+} = L_i(b_i,b_4;b_j,b_k)$ & 
$L_{i1}^{-} = L_i(b_i^{-1},b_4^{-1};b_j,b_k)$ \\[1mm]
\hline
\hline
\vspace{-4mm} & & \\
$2$ & $L_{i2}^{+} = L_i(b_i,b_4^{-1};b_j,b_k)$ & 
$L_{i2}^{-} = L_i(b_i^{-1},b_4;b_j,b_k)$ \\[1mm]
\hline
\hline
\vspace{-4mm} & & \\
$3$ & $L_{i3}^{+} = L_i(b_j,b_k;b_i,b_4)$ & 
$L_{i3}^{-} = L_i(b_j^{-1},b_k^{-1};b_i,b_4)$ \\[1mm]
\hline
\hline
\vspace{-4mm} & & \\
$4$ & $L_{i4}^{+} =L_i(b_j,b_k^{-1};b_i,b_4)$ & 
$L_{i4}^{-} =L_i(b_j^{-1},b_k;b_i,b_4)$ \\[1mm]
\hline
\end{tabular}
\end{center} 
\caption{Eight lines intersecting the line $L_i$ at infinity, 
divided into four pairs} 
\label{tab:line}
\end{table}
%%%%%%%%%%%%%%%%%%%%%%%%%%%%%%%%%%%%%%%%%%%%%%%%%%%%%%%%%%%%%%%
It is well known that a smooth projective cubic surface 
has exactly $27$ lines on it. 
We describe them in the current situation \cite{IU}. 
Let $L_i(b_i,b_4;b_j,b_k)$ be the line in $\P^3$ defined by 
%%%%%%%%%%%%%%%%%%%%%%%%% eqn:Li %%%%%%%%%%%%%%%%%%%%%%%%%%%%%
\begin{equation} \label{eqn:Li} 
X_i = (b_ib_4+ b_i^{-1}b_4^{-1}) X_0, \qquad 
X_j + (b_ib_4) X_k = 
\{b_i(b_k+b_k^{-1})\} + b_4 (b_j+b_j^{-1}) X_0. 
\end{equation}
%%%%%%%%%%%%%%%%%%%%%%%%%%%%%%%%%%%%%%%%%%%%%%%%%%%%%%%%%%%%%%
For each $i \in \{1,2,3\}$ the eight lines in 
Table \ref{tab:line} are the only lines on $\ol{\Sol}(\th)$ 
that intersect the $i$-th line $L_i$ at infinity, but 
they do not intersect the remaining two lines $L_j$ and $L_k$ 
at infinity. 
These lines are divided into four pairs as in Table 
\ref{tab:line}. 
The surface $\ol{\Sol}(\th)$ is always smooth at infinity 
\cite{IU} and hence, if $\k \in \K-\Wall$, then 
$\ol{\Sol}(\th)$ is smooth everywhere. 
In this case, the two lines in the same pair intersect, 
while two lines from different pairs do not. 
The intersection point of the $i$-th pair is 
exactly the $i$-th fixed point in Table \ref{tab:fixed}. 
See Figure \ref{fig:lines} for a total image of these 
situations. 
{\sl Caution}: for a pair of distinct indices $i$ and $j$, 
the intersection relations between $L_{i\mu}^{\pm}$ and 
$L_{j\nu}^{\pm}$ are not depicted in the Figure \ref{fig:lines}. 
We also remark that in some degenerate cases the lines 
$L_{i\mu}^{+}$ and $L_{i\mu}^{-}$ may meet in a point on 
the line $L_i$ at infinity. 
%%%%%%%%%%%%%%%%%%%%%%%%%% fig:lines %%%%%%%%%%%%%%%%%%%%%%%%%%
%%%%%%%%%%%%%%%%%%%%%%%%%%%%%%%%%%%%%%%%%%%%%%%%%%%%%%%%%%%%%%%
\par 
Next we consider the case where the determinant $4-x_i^2$ 
of system (\ref{eqn:system}) vanishes. 
In other words we ask when the fixed point set 
$\mathrm{Fix}_j(\th)$ contains points $x$ such that 
$x_i \in \{\pm2\}$. 
%%%%%%%%%%%%%%%%%%%%%% lem:two1 %%%%%%%%%%%%%%%%%%%%%%%%%%%%%%%%%
\begin{lemma} \label{lem:two1}
$\mathrm{Fix}_j(\th)$ contains a point $x$ such that 
$x_i = 2\delta$ with $\delta \in \{\pm1\}$ if and only 
if either 
%%%%%%
\begin{enumerate} 
\item $b_ib_4 = b_ib_4^{-1} = \delta;$ or 
\item $b_jb_k = b_jb_k^{-1} = \delta;$ or 
\item $b_ib_4^{\ve_4} = b_jb_k^{\ve_k} = 
\delta$ for some double sign 
$(\ve_k,\ve_4) \in \{\pm1\}^2$. 
\end{enumerate} 
%%%%%%
If this is the case, then $\th_k = \delta \th_j$ and 
all such poins $x$ are exactly those points on the line 
%%%%%%%%%%%%%%%%%%%%%%%% eqn:line0 %%%%%%%%%%%%%%%%%%%%%%%%%%%%% 
\begin{equation} \label{eqn:line0} 
\ell_j^{\delta} := 
\{\, x_i = 2 \delta, \, x_j + \delta x_k = \th_j/2 \,\}. 
\end{equation} 
%%%%%%%%%%%%%%%%%%%%%%%%%%%%%%%%%%%%%%%%%%%%%%%%%%%%%%%%%%%%%%%% 
In particular $\ell_j^{\delta} \subset \mathrm{Fix}_j(\th)$ 
precisely when $x_i = 2\delta$ is a multiple root of the 
quartic equation $(\ref{eqn:quartic})$. 
%%%%%%%%%%%%%%%%%%%%%%%%%%%%%%%%%%%%%%%%%%%%%%%%%%%%%%%%%%%%%%%% 
\end{lemma} 
%%%%%%%%%%%%%%%%%%%%%%%%%%%%%%%%%%%%%%%%%%%%%%%%%%%%%%%%%%%%%%%% 
{\it Proof}. 
If $x_i = 2 \delta$ with $\delta \in \{\pm1\}$ then 
system (\ref{eqn:system}) is linearly dependent, 
so that $\th_k - \delta \th_j = 0$. 
However, since 
$\th_k - \delta \th_j = 
(b_ib_j)^{-1}(b_ib_4-\delta)(b_ib_4^{-1}-\delta)
(b_jb_k-\delta)(b_jb_k^{-1}-\delta)$, 
we have either $b_ib_4^{\ve_4} = \delta$ for some sign 
$\ve_4 \in \{\pm1\}$ or $b_jb_k^{\ve_k} = \delta$ for 
some sign $\ve_k \in \{\pm1\}$. 
Taking the equation $x_j + \delta x_k = \th_j/2$ into 
account, we observe that $f(x,\th)$ factors as 
%%%%%
\[
f(x,\th) = \left\{ 
\begin{array}{ll} 
-(2b_ib_j)^{-2} (b_ib_4^{-\ve_4}-\delta)^2
(b_jb_k-\delta)^2(b_jb_k^{-1}-\delta)^2 
\qquad & (\mbox{if} \quad b_ib_4^{\ve_4} = \delta), \\[2mm]
-(2b_ib_j)^{-2} (b_jb_k^{-\ve_k}-\delta)^2
(b_ib_4-\delta)^2(b_ib_4^{-1}-\delta)^2 
\qquad & (\mbox{if} \quad b_jb_k^{\ve_k} = \delta). 
\end{array}
\right.
\] 
%%%%% 
If $b_ib_4^{\ve_4} = \delta$ then equation $f(x,\th) = 0$ 
yields either $b_ib_4^{-\ve_4} = \delta$ or 
$b_j b_k^{\ve_k} = \delta$ for some sign $\ve_k \in \{\pm1\}$; 
the former case falls into case (1) while the latter 
falls into case (3). 
In a similar manner the other case $b_jb_k^{\ve_k} = \delta$ 
falls into case (2) or case (3). 
\par 
Next, if $\mathrm{Fix}_j(\th)$ contains the line 
$\ell_j^{\delta}$, then what we have just proved implies 
that 
%%%%
\[
\begin{array}{rclcll}
F(b_i,b_4;b_j,b_k) &=& F(b_i,b_4^{-1};b_j,b_k) &=& 
2 \delta \quad & \mbox{if condition (1) is satisfied;} \\[2mm] 
F(b_j,b_k;b_i,b_4) &=& F(b_j,b_k^{-1};b_i,b_4) &=& 
2 \delta \quad & \mbox{if condition (2) is satisfied;} \\[2mm]
F(b_i,b_4^{\ve_4};b_j,b_k) &=& F(b_j,b_k^{\ve_k};b_i,b_4) 
&=& 2 \delta \quad & \mbox{if condition (3) is satisfied.}
\end{array} 
\]
%%%%
Hence $x_i = 2\delta$ is a multiple root of the quartic 
equation (\ref{eqn:quartic}). 
Conversely, if $x_i = 2\delta$ is a multiple root of 
(\ref{eqn:quartic}), then we can trace the argument 
backwards to conclude that the system (\ref{eqn:per1}) 
admits the line solution $\ell_j^{\delta}$, that is, 
$\mathrm{Fix}_j(\th)$ contains $\ell_j^{\delta}$. 
\hfill $\Box$ \par\medskip 
%%%%%%%%%%%%%%%%%%%%%%%%%% end of proof %%%%%%%%%%%%%%%%%%%%%
Summarizing the arguments so far yields a classification of 
the irreducible components of the algebraic 
set $\mathrm{Fix}_j(\th)$ in terms of certain roots of 
quartic equation (\ref{eqn:quartic}). 
%%%%%%%%%%%%%%%%%%%%%%%%%% tab:classify %%%%%%%%%%%%%%%%%%%%%
\begin{table}[t]
\begin{center}
\begin{tabular}{|c|c||c|c|}
\hline
\vspace{-4mm} & & & \\
$x_i \in \{\pm 2\}$ & multiplicity & component & remark \\[1mm]
\hline
\hline 
\vspace{-4mm} & & & \\
no & simple & smooth point & intersection point of 
$L_{i\mu}^{\pm}$ \\[1mm] 
\hline 
\vspace{-4mm} & & & \\
no & multiple & singular point & Riccati locus \\[1mm] 
\hline 
\vspace{-4mm} & & & \\
yes & multiple & line $\ell_j^{+}$ or $\ell_j^{-}$ & 
line contains singular points \\[1mm] 
\hline 
\hline 
\vspace{-4mm} & & & \\
yes & simple & empty & $L_{i\mu}^{\pm}$ intersects at 
infinity \\[1mm] 
\hline 
\end{tabular}
\end{center}
\caption{The roots of quartic equation (\ref{eqn:quartic}) 
and the components of $\mathrm{Fix}_j(\th)$}
\label{tab:classify}
\end{table}
%%%%%%%%%%%%%%%%%%%%%%%%%% thm:quartic %%%%%%%%%%%%%%%%%%%%%%
\begin{theorem} \label{thm:quartic} 
Any irreducible component of $\mathrm{Fix}_j(\th)$ is just 
a single point or a single affine line; the former is called 
a point component and the latter is called a line component 
respectively. 
The irreducible components of $\mathrm{Fix}_j(\th)$ are 
in one-to-one correspondence with those roots of quartic 
equation $(\ref{eqn:quartic})$ which are not a simple 
root $x = (x_1,x_2,x_3)$ such that $x_i \in \{\pm2\}$. 
\begin{enumerate} 
\item A simple root with $x_i \not\in \{\pm2\}$ 
corresponds to a point component that is a smooth point 
of the surface $\Sol(\th)$ and is given in 
Table $\ref{tab:fixed}$. 
\item A multiple root with $x_i \not\in \{\pm2\}$ 
corresponds to a point component that is a singular 
point of the surface $\Sol(\th)$ and is associated 
with Riccati solutions. 
\item A multiple root with $x_i \in \{\pm2\}$ 
corresponds to a line component; either $\ell_j^{+}$ or 
$\ell_j^{-}$. 
\item A simple root with $x_i \in \{\pm2\}$ corresponds to 
no component of $\mathrm{Fix}_j(\th)$. 
\end{enumerate} 
\end{theorem}
%%%%%%%%%%%%%%%%%%%%%%%%%%%%%%%%%%%%%%%%%%%%%%%%%%%%%%%%%%%%%
\par 
A summary of Theorem \ref{thm:quartic} is given in 
Table \ref{tab:classify} and the following remark may be 
helpful. 
%%%%%%%%%%%%%%%%%%%%%%%%%% rem:quartic %%%%%%%%%%%%%%%%%%%%%%
\begin{remark} \label{rem:quartic} 
The assertions (3) and (4) of Theorem \ref{thm:quartic} may 
be well understood through the degeneration of line 
configration on the projective surface $\ol{\Sol}(\th)$ as 
the parameter $\th = \rh (\k)$ tends to a special position. 
For a generic value of $\th$ the lines $L_{i\mu}^{\pm}$ 
intersect in a single (smooth) point on the affine part 
$\Sol(\th)$ of $\ol{\Sol}(\th)$. 
If the parameter $\th$ tends to a special position so that 
a corresponding root $x_i$ of quartic equation 
(\ref{eqn:quartic}) approaches $\{\pm2\}$, then the two 
line $L_{i\mu}^{\pm}$ are getting ``parallel" and 
eventually either coincide completely or meet in a point 
at infinity. 
The former case falls into assertion (3) and the latter 
case falls into assertion (4) respectively. 
\end{remark} 
%%%%%%%%%%%%%%%%%%%%%%%%%%%%%%%%%%%%%%%%%%%%%%%%%%%%%%%%%%%%%
\par 
Let us investigate more closely the case where 
$\mathrm{Fix}_j(\th)$ contains line components. 
%%%%%%%%%%%%%%%%%%%%%%%%%% lem:two2 %%%%%%%%%%%%%%%%%%%%%%%%%
\begin{lemma} \label{lem:two2} 
Let $\th = \rh(\k)$ with $\k \in \K$ and 
$(i,j,k)$ be any cyclic permutation of $(1,2,3)$. 
\begin{enumerate}
\item $\mathrm{Fix}_j(\th)$ contains either 
$\ell_j^+$ or $\ell_j^-$ but not both of them 
if and only if $\k$ lies in a $\wt{W}(D_4^{(1)})$-stratum 
appearing in the following adjacency diagram 
$($see also Figure $\ref{fig:adjac}):$  
%%%%%%%%%%%%%%%%%%%%%%%%% cd:boundary %%%%%%%%%%%%%%%%%%%%%
\begin{equation} \label{cd:boundary}
\begin{CD}
(A_1^{\oplus 2})_i @>>> A_1^{\oplus 3} \\
 @VVV     @VVV  \\
(A_3)_i @>>> D_4
\end{CD}
\end{equation} 
%%%%%%%%%%%%%%%%%%%%%%%%%%%%%%%%%%%%%%%%%%%%%%%%%%%%%%%%%%%%
\item $\mathrm{Fix}_j(\th)$ contains both $\ell_j^+$ and 
$\ell_j^-$ if and only if $\th = (0,0,0,-4)$, that is, 
precisely when $\k$ is in the $\wt{W}(D_4^{(1)})$-stratum 
of type $A_1^{\oplus 4}$. 
In this case one has 
$\mathrm{Fix}_j(\th) = \ell_j^+ \amalg \ell_j^-$. 
\end{enumerate} 
\end{lemma}
%%%%%%%%%%%%%%%%%%%%%%%%%%%%%%%%%%%%%%%%%%%%%%%%%%%%%%%%%%%%%%%%
{\it Proof}. 
Lemma \ref{lem:two1} implies that $\mathrm{Fix}_j(\th)$ 
contains at least one of $\ell_j^{+}$ and $\ell_j^{-}$ if 
and only if either (a) $b_ib_4 = b_ib_4^{-1} \in \{\pm1\}$; 
or (b) $b_jb_k = b_jb_k^{-1} \in \{\pm1\}$; or (c) 
$b_ib_4^{\ve_4} = b_jb_k^{\ve_k} \in \{\pm1\}$ for some 
double sign $(\ve_k,\ve_4) \in \{\pm1\}^2$. 
This property is invariant under 
the action of $\wt{W}(D_4^{(1)}) = \Kl \ltimes W(D_4^{(1)})$ 
on $\K$. 
Using this action we can reduce conditions (a) and (c) to 
condition (b). 
First, observe that the permutation $(i,j)(k4) \in \Kl$ 
induces the map $(b_0,b_i,b_j,b_k,b_4) \mapsto 
(b_0,b_j,b_i,-b_4,-b_k)$, which reduces condition (a) to 
condition (b). 
Next, formula (\ref{eqn:cartan}) implies that the reflection 
$w_i$ induces the multiplicative transformation 
$w_i : B \to B$, $b \mapsto b'$, where 
%%%%%
\[
b_j' = \left\{
\begin{array}{rl}
-b_j b_i^{c_{ij}} \qquad & (i=4, j=0), \\[2mm]
 b_j b_i^{c_{ij}} \qquad & (\mbox{otherwise}). 
\end{array}
\right. 
\]
%%%%%%
Applying $w_4$ or $w_k$ if necessary, we may assume from the 
beginning that $\ve_4 = 1$ and $\ve_k = -1$ in condition (c). 
Then using $w_0$ there yields $b_0^2 b_i b_4 = 
b_j b_k^{-1} \in \{\pm1\}$. 
But since $b_0^2b_ib_4b_jb_k = 1$, 
we have $b_j b_k = b_jb_k^{-1} \in \{\pm1\}$, 
that is, condition (b). 
Note that condition (b) means $\k_j$, $\k_k \in \Z$. 
On the other hand, the extended affine Weyl group 
$\wt{W}(D_4^{(1)})$ contains shifts 
%%%%%%%%%%%%%%%%%%%%%%%% eqn:shift %%%%%%%%%%%%%%%%%%%%%%%%
\begin{equation} \label{eqn:shift}
\begin{array}{rcl}
(\k_0,\k_i,\k_j,\k_k,\k_4) &\mapsto& 
(\k_0,\k_i-1,\k_j+1,\k_k,\k_4), \\[1mm]
(\k_0,\k_i,\k_j,\k_k,\k_4) &\mapsto& 
(\k_0,\k_i,\k_j,\k_k+1,\k_4-1). 
\end{array} 
\end{equation}
%%%%%%%%%%%%%%%%%%%%%%%%%%%%%%%%%%%%%%%%%%%%%%%%%%%%%%%%%%%
Repeated applications of these operations and 
their inverses can shift $\k_j$ and $\k_k$ independently 
by arbitrary integers. 
Thus the condition $\k_j$, $\k_k \in \Z$ can further be 
reduced to $\k_j = \k_k = 0$. 
Thus we have shown that if $\mathrm{Fix}_j(\th)$ 
contains at least one of $\ell_j^{+}$ and $\ell_j^{-}$, 
then $\k$ must lie in the $\wt{W}(D_4^{(1)})$-stratum of 
type $(A_1^{\oplus 2})_i$ or on its boundary strata 
of types $(A_3)_i$, $A_1^{\oplus 3}$, $D_4$, $A_1^{\oplus 4}$. 
Moreover, it is easy to see that the converse is also true. 
\par 
For a sign $\delta \in \{\pm1\}$ the conditions 
(1), (2), (3) in Lemma \ref{lem:two1} are denoted by 
$(1^{\delta})$, $(2^{\delta})$, $(3^{\delta})$, 
respectively. 
Now we assume that $\mathrm{Fix}_j(\th)$ contains both 
$\ell_j^{+}$ and $\ell_j^{-}$. 
Then there exists a pair of conditions, 
one from $\{(1^{+}),(2^{+}),(3^{+})\}$ and the other 
from $\{(1^{-}),(2^{-}),(3^{-})\}$, that are valid at 
the same time. 
Such a pair can be consistent only if it is either 
$(1^{-})+(2^{-})$; or $(2^{+})+(1^{-})$; or $(3^{+})+(3^{-})$ 
where if the sign for $(3^{+})$ is $(\ve_k,\ve_4)$ then 
the sign for $(3^{-})$ must be its antipode $(-\ve_k,-\ve_4)$. 
The first and second pairs lead to  
$b_1^2 = b_2^2 = b_3^2 = b_4^2 = 1$ and to 
$b_1b_2b_3b_4 = -1$, while the third pair 
yields $b_1^2 = b_2^2 = b_3^2 = b_4^2 = -1$. 
These are nothing but the conditions (a) and (b) in 
Example \ref{ex:sing}.\,(1). 
Therefore $\k$ must lie in the stratum of type 
$A_1^{\oplus 4}$. 
Combining this with the discussion in the last paragraph 
establishes the assertion (1), as well as a large part 
of the assertion (2). 
The only thing yet to be proved is the assertion that 
if $\mathrm{Fix}_j(\th)$ contains both $\ell_j^{+}$ and 
$\ell_j^{-}$, then $\mathrm{Fix}_j(\th) = \ell_j^{+} 
\amalg \ell_j^{-}$. 
For this, the last part of Lemma \ref{lem:two1} implies 
that  both $x_i = 2$ and $x_i = -2$ are multiple roots of 
the quartic equation (\ref{eqn:quartic}), so that there 
are no other roots of the equation (\ref{eqn:quartic}). 
Thus $\mathrm{Fix}_j(\th)$ has no elements other than 
those in $\ell_j^{+} \amalg \ell_j^{-}$. 
\hfill $\Box$ 
%%%%%%%%%%%%%%%%%%%%%%%% end of proof %%%%%%%%%%%%%%%%%%%%%%%
\par\medskip 
Now we turn our attention to periodic points and investigate 
the set $\wt{\mathrm{Per}}_j^{\circ}(\th;n)$ of periodic 
points of prime period $n > 1$ on the non-Riccati locus. 
%%%%%%%%%%%%%%%%%%%%%%%%%%% lem:pern %%%%%%%%%%%%%%%%%%%%%%%
\begin{lemma} \label{lem:pern} 
For any integer $n > 1$ the set 
$\wt{\mathrm{Per}}_j^{\circ}(\th;n)$ is biholomorphic to 
the disjoint union of $\varphi(n)$ copies of $\C^{\times}$, 
where $\varphi(n)$ denotes the number of integers $0 < m < n$ 
coprime to $n$. 
\end{lemma}
%%%%%%%%%%%%%%%%%%%%%%%%%%%%%%%%%%%%%%%%%%%%%%%%%%%%%%%%%%%%
{\it Proof}. 
By Lemma \ref{lem:updown} we can identify 
$\wt{\mathrm{Per}}_j^{\circ}(\th;n)$ with 
$\mathrm{Per}_j(\th;n)$ 
and hence may work downstairs. 
For any integer $0< m < n$ coprime to $n$, we consider 
the projective curve $\ol{C}_m$ in $\P^3$ defined by 
%%%%%%%%%%%%%%% eqn:conic1 %%%%% eqn:conic2 %%%%%%%%%%%%%%%
\begin{eqnarray} 
\{4 \cos^2(\pi m/n)-2 \th_i \cos(\pi m/n)+\th_4\} X_0^2 
- X_0(\th_j X_j + \th_k X_k) & & \nonumber \\[2mm]
+ X_j^2 +X_k^2 + 2 \cos(\pi m/n) X_j X_k &=& 0, 
\label{eqn:conic1} \\[3mm] 
X_i- 2 \cos(\pi m/n) X_0 &=& 0, \label{eqn:conic2} 
\end{eqnarray} 
%%%%%%%%%%%%%%%%%%%%%%%%%%%%%%%%%%%%%%%%%%%%%%%%%%%%%%%%%%%%
where (\ref{eqn:conic1}) is obtained from $F(X,\th) = 0$ 
by substituting (\ref{eqn:conic2}) and factoring $X_0$ 
out of it. 
It follows from $-2 < 2 \cos(\pi m/n) < 2$ that 
$\ol{C}_m$ is an irreducible smooth conic curve. 
By equations (\ref{eqn:pern}) of Lemma \ref{lem:per} the 
closure $\ol{\mathrm{Per}}_j(\th;n)$ of 
$\mathrm{Per}_j(\th;n)$ in $\ol{\Sol}(\th)$ is the union 
of these $\varphi(n)$ curves $\ol{C}_m$. 
The curve $\ol{C}_m$ intersects the lines $L = L_i 
\cup L_j \cup L_k$ at infinity in the two points 
%%%%
\[
P_m^{\pm} : \quad [X_0:X_i:X_j:X_k] = 
[0:0:-1:\exp(\pm\pi\sqrt{-1}m/n)] \in L_i. 
\] 
%%%%
If $C_m := \ol{C}_m - \{P_m^{+}, P_m^{-}\}$, then 
$C_m$ is biholomorphic to $\C^{\times}$, since 
$\ol{C}_m \cong \P^1$. 
So $\mathrm{Per}_j(\th;n)$ is the disjoint union 
of the $\varphi(n)$ curves $C_m$ with $0 < m < n$, 
$(m,n)=1$, and hence biholomorphic to the disjoint union 
of $\varphi(n)$ copies of $\C^{\times}$. \hfill $\Box$ 
%%%%%%%%%%%%%%%%%%%%%%%% sec:case %%%%%%%%%%%%%%%%%%%%%%%%%%%
\section{Case-by-Case Study} \label{sec:case}
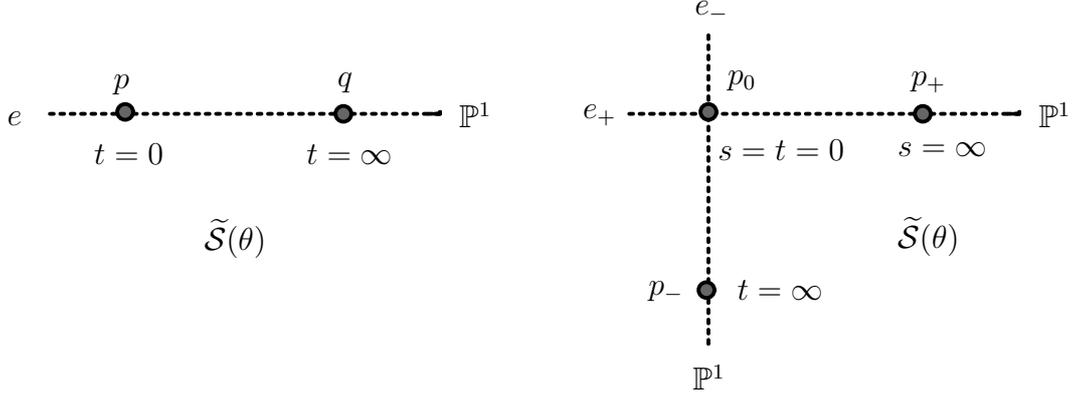
\begin{figure}[t] 
\begin{center}
%WinTpicVersion2.15
\unitlength 0.1in
\begin{picture}(53.40,19.20)(1.80,-24.80)
% LINE 0 2 3 0
% 6 3400 1610 5410 1610 5360 1610 5420 1610 5420 1610 5420 1600
% 
\special{pn 20}%
\special{pa 3400 1210}%
\special{pa 5410 1210}%
\special{dt 0.054}%
\special{pa 5410 1210}%
\special{pa 5409 1210}%
\special{dt 0.054}%
\special{pa 5360 1210}%
\special{pa 5420 1210}%
\special{dt 0.054}%
\special{pa 5420 1210}%
\special{pa 5419 1210}%
\special{dt 0.054}%
\special{pa 5420 1210}%
\special{pa 5420 1200}%
\special{dt 0.054}%
\special{pa 5420 1200}%
\special{pa 5420 1201}%
\special{dt 0.054}%
% LINE 0 2 3 0
% 2 3810 1200 3810 2810
% 
\special{pn 20}%
\special{pa 3810 800}%
\special{pa 3810 2410}%
\special{dt 0.054}%
\special{pa 3810 2410}%
\special{pa 3810 2409}%
\special{dt 0.054}%
% CIRCLE 0 0 0 0
% 4 3800 2530 3810 2570 3810 2570 3810 2570
% 
\special{pn 20}%
\special{sh 0.600}%
\special{ar 3800 2130 41 41  0.0000000 6.2831853}%
% CIRCLE 0 0 0 0
% 4 3810 1600 3820 1640 3820 1640 3820 1640
% 
\special{pn 20}%
\special{sh 0.600}%
\special{ar 3810 1200 41 41  0.0000000 6.2831853}%
% CIRCLE 0 0 0 0
% 4 4920 1610 4930 1650 4930 1650 4930 1650
% 
\special{pn 20}%
\special{sh 0.600}%
\special{ar 4920 1210 41 41  0.0000000 6.2831853}%
% STR 2 0 3 0
% 3 3160 1570 3160 1670 2 0
% $e_{+}$
\put(31.6000,-12.7000){\makebox(0,0)[lb]{$e_{+}$}}%
% STR 2 0 3 0
% 3 3740 1030 3740 1130 2 0
% $e_{-}$
\put(37.4000,-7.3000){\makebox(0,0)[lb]{$e_{-}$}}%
% STR 2 0 3 0
% 3 5520 1580 5520 1680 2 0
% $\mathbb{P}^1$
\put(55.2000,-12.8000){\makebox(0,0)[lb]{$\mathbb{P}^1$}}%
% STR 2 0 3 0
% 3 3730 2950 3730 3050 2 0
% $\mathbb{P}^1$
\put(37.3000,-26.5000){\makebox(0,0)[lb]{$\mathbb{P}^1$}}%
% STR 2 0 3 0
% 3 3910 1390 3910 1490 2 0
% $p_0$
\put(39.1000,-10.9000){\makebox(0,0)[lb]{$p_0$}}%
% STR 2 0 3 0
% 3 4860 1400 4860 1500 2 0
% $p_{+}$
\put(48.6000,-11.0000){\makebox(0,0)[lb]{$p_{+}$}}%
% STR 2 0 3 0
% 3 3500 2490 3500 2590 2 0
% $p_{-}$
\put(35.0000,-21.9000){\makebox(0,0)[lb]{$p_{-}$}}%
% STR 2 0 3 0
% 3 4790 2250 4790 2350 2 0
% $\wt{\mathcal{S}}(\th)$
\put(47.9000,-19.5000){\makebox(0,0)[lb]{$\wt{\mathcal{S}}(\th)$}}%
% LINE 0 2 3 0
% 6 400 1610 2410 1610 2360 1610 2420 1610 2420 1610 2420 1600
% 
\special{pn 20}%
\special{pa 400 1210}%
\special{pa 2410 1210}%
\special{dt 0.054}%
\special{pa 2410 1210}%
\special{pa 2409 1210}%
\special{dt 0.054}%
\special{pa 2360 1210}%
\special{pa 2420 1210}%
\special{dt 0.054}%
\special{pa 2420 1210}%
\special{pa 2419 1210}%
\special{dt 0.054}%
\special{pa 2420 1210}%
\special{pa 2420 1200}%
\special{dt 0.054}%
\special{pa 2420 1200}%
\special{pa 2420 1201}%
\special{dt 0.054}%
% CIRCLE 0 0 0 0
% 4 1920 1610 1930 1650 1930 1650 1930 1650
% 
\special{pn 20}%
\special{sh 0.600}%
\special{ar 1920 1210 41 41  0.0000000 6.2831853}%
% STR 2 0 3 0
% 3 180 1570 180 1670 2 0
% $e$
\put(1.8000,-12.7000){\makebox(0,0)[lb]{$e$}}%
% STR 2 0 3 0
% 3 2520 1580 2520 1680 2 0
% $\mathbb{P}^1$
\put(25.2000,-12.8000){\makebox(0,0)[lb]{$\mathbb{P}^1$}}%
% STR 2 0 3 0
% 3 730 1410 730 1510 2 0
% $p$
\put(7.3000,-11.1000){\makebox(0,0)[lb]{$p$}}%
% STR 2 0 3 0
% 3 1890 1400 1890 1500 2 0
% $q$
\put(18.9000,-11.0000){\makebox(0,0)[lb]{$q$}}%
% CIRCLE 0 0 0 0
% 4 790 1600 800 1640 800 1640 800 1640
% 
\special{pn 20}%
\special{sh 0.600}%
\special{ar 790 1200 41 41  0.0000000 6.2831853}%
% STR 2 0 3 0
% 3 630 1770 630 1870 2 0
% $t=0$
\put(6.3000,-14.7000){\makebox(0,0)[lb]{$t=0$}}%
% STR 2 0 3 0
% 3 1730 1770 1730 1870 2 0
% $t=\infty$
\put(17.3000,-14.7000){\makebox(0,0)[lb]{$t=\infty$}}%
% STR 2 0 3 0
% 3 3860 1750 3860 1850 2 0
% $s=t=0$
\put(38.6000,-14.5000){\makebox(0,0)[lb]{$s=t=0$}}%
% STR 2 0 3 0
% 3 4790 1730 4790 1830 2 0
% $s=\infty$
\put(47.9000,-14.3000){\makebox(0,0)[lb]{$s=\infty$}}%
% STR 2 0 3 0
% 3 3960 2480 3960 2580 2 0
% $t=\infty$
\put(39.6000,-21.8000){\makebox(0,0)[lb]{$t=\infty$}}%
% STR 2 0 3 0
% 3 1200 2260 1200 2360 2 0
% $\wt{\mathcal{S}}(\th)$
\put(12.0000,-19.6000){\makebox(0,0)[lb]{$\wt{\mathcal{S}}(\th)$}}%
\end{picture}%
\end{center}
\caption{Surface of types $A_1$ (left) and $A_2$ (right)} 
\label{fig:A1A2} 
\end{figure}

%%%%%%%%%%%%%%%%%%%%%%%%%%%%%%%%%%%%%%%%%%%%%%%%%%%%%%%%%%%%%
We make case-by-case studies of $\wt{\mathrm{Fix}}_j(\th)$ 
and $\wt{\mathrm{Per}}_j^{e}(\th;n)$ according to the 
adjacency diagram in Figure \ref{fig:adjac}. 
Now we need to introduce some notation. 
Recall that we have the resolution of singularities 
(\ref{eqn:brieskorn}) which restricts to an isomorphism 
$\varphi : \wt{\Sol}^{\circ}(\th) \to \Sol^{\circ}(\th)$ 
and that the smooth fixed points $\mathrm{Fix}_j^{\circ}(\th)$ 
in $\Sol^{\circ}(\th)$ are listed in Table \ref{tab:fixed}. 
For each $P \in \mathrm{Fix}_j^{\circ}(\th)$ let $\wt{P} \in 
\wt{\mathrm{Fix}}_j^{\circ}(\th)$ denote its lift through the 
isomorphism $\varphi$. 
For example, $\wt{P}(b_i,b_4;b_j,b_k)$ denotes the lift of 
$P(b_i,b_4;b_j,b_k)$. 
If $\{\cdots\}$ is a set of expressions 
$\wt{P}(b_i,b_4^{\pm 1};b_j,b_k)$, 
$\wt{P}(b_j,b_k^{\pm1};b_i,b_4)$, then we denote by 
$\{\!\{\cdots \}\!\}$ its subset obtained by discarding 
those expressions which do not satisfy either the existence 
condition or the smoothness condition of Table \ref{tab:fixed}. 
An example is given in (\ref{eqn:FixEmpty}) below. 
%%%%%%%%%%%%%%%%%%%%%%%%% ex:empty %%%%%%%%%%%%%%%%%%%%%%%%%%
\begin{example}[$\mbox{\boldmath $\emptyset$}$] 
\label{ex:empty}
Consider the $\wt{W}(D_4^{(1)})$-stratum of type 
$\emptyset$, namely, the big open $\K-\Wall$. 
%%%%%%%%%%%%%%%%%%%%%%%% eqn:FixEmpty %%%%%%%%%%%%%%%%%%%%%%
\begin{equation} \label{eqn:FixEmpty} 
\wt{\mathrm{Fix}}_j(\th) = 
\{\!\{\, \wt{P}(b_i,b_4;b_j,b_k), \, 
\wt{P}(b_i,b_4^{-1};b_j,b_k), \, 
\wt{P}(b_j,b_k;b_i,b_4), \, 
\wt{P}(b_j,b_k^{-1};b_i,b_4) \, \}\!\}. 
\end{equation}
%%%%%%%%%%%%%%%%%%%%%%%%%%%%%%%%%%%%%%%%%%%%%%%%%%%%%%%%%%%%%
Here we have only to care the existence condition, as we are 
in the big open where the smoothness condition is 
fulfilled by hypothesis. 
If a finer stratification of $\K$ attached to 
the $W(F_4^{(1)})$-action on $\K$ is introduced, then 
a more precise description of (\ref{eqn:FixEmpty}) is feasible, 
detecting how many and which elements are there in 
(\ref{eqn:FixEmpty}), but the details are omitted. 
We only remark that $\wt{\mathrm{Fix}}_j(\th)$ consists of 
four distinct points in the most generic case where none of 
$\k_i \pm \k_4$ and $\k_j \pm \k_k$ are integers. 
As for the periodic points, since there is 
no Riccati locus, we have 
%%%%%%%
\[
\wt{\mathrm{Per}}_j(\th;n) = 
\wt{\mathrm{Per}}_j^{\circ}(\th;n), \qquad  
\wt{\mathrm{Per}}_j^{e}(\th;n) = \emptyset, 
\qquad (n > 1). 
\]
%%%%%%%
\end{example} 
%%%%%%%%%%%%%%%%%%%%%%%%% ex:A1 %%%%%%%%%%%%%%%%%%%%%%%%%%%%%
\begin{example}[$\mbox{\boldmath $A_1$}$] \label{ex:A1}
Consider the $\wt{W}(D_4^{(1)})$-stratum of type $A_1$. 
We may assume that $\k_0 = 0$ so that $b_0 = 1$ and 
$b_ib_jb_kb_4=1$. 
Note that none of $b_i^2$, $b_j^2$, $b_k^2$, $b_4^2$ 
equals $1$. 
We claim that $b_j^2 b_k^2 \neq 1$. 
Otherwise, we would have $\k_j + \k_k \in \Z$. 
Applying a shift as in (\ref{eqn:shift}) to $\k$ repeatedly, 
we may assume that $\k_j + \k_k = 0$ while keeping the 
condition $\k_0 = 0$. 
Then the transformation $w_0 w_j$ sends $\k$ to $\k'$ with  
$\k'_j = 0$ and $\k'_k = \k_j + \k_k = 0$, so that one has 
$\k \in \ol{\K}_{\{j,k\}}$, namely, $\k$ lies in the 
closure of the stratum of type $(A_1^{\oplus 2})_i$. 
This contradicts the assumption that we are in the 
stratum of type $A_1$. 
In this case, the surface $\Sol(\th)$ has a unique singular 
point of type $A_1$ at 
%%%%%%%%
\[
(x_i,x_j,x_k) = (b_ib_4+b_i^{-1}b_4^{-1}, 
b_jb_4+b_j^{-1}b_4^{-1}, b_kb_4+b_k^{-1}b_4^{-1}). 
\]
%%%%%%%%
Blow up $\Sol(\th)$ at this point to obtain a 
minimal resolution (\ref{eqn:brieskorn}). 
Write the blowing-up as 
%%%%%%%%
\[
(x_i,x_j,x_k) = (u_i u_j + b_ib_4+b_i^{-1}b_4^{-1}, 
u_j + b_jb_4+b_j^{-1}b_4^{-1}, 
u_k u_j + b_kb_4+b_k^{-1}b_4^{-1})
\]
%%%%%%%%
in terms of coordinates $(u_i,u_j,u_k)$. 
The exceptional set $e$ is the irreducible quadratic curve 
%%%%%%%%%
\[
u_j = b_ib_jb_k + (1 + b_i^2 b_j^2) b_k u_i + 
(1 + b_k^2  b_j^2) b_i u_k + b_ib_jb_k (u_i^2 + u_k^2) 
+ (1 + b_i^2 b_k^2) b_j u_i u_k = 0, 
\]
%%%%%%%%%
which can be paramatrized as $u_j = 0$ and 
%%%%%%%%%
\[
\begin{array}{rcl}
u_i &=& \dfrac{b_i^2 b_4^2 (b_j^2b_k^2 - 1)^3 t}{b_j 
\{b_k(b_i^2 - 1)(b_j^2 - 1) + b_i(b_j^2b_k^2 - 1)t\}
\{(b_k^2 - 1)(b_4^2-1) + b_i b_k b_4^2(b_j^2 b_k^2 - 1) t\}}, 
\\[6mm]
u_k &=& \dfrac{\{b_j^2 b_k (b_i^2 - 1)
(b_k^2 - 1) - b_i(b_j^2 b_k^2 - 1)t\}
\{b_k(b_j^2 - 1)(1 - b_4^2) + b_ib_4^2(b_j^2b_k^2 - 1)t\}}{b_j 
\{b_k(b_i^2 - 1)(b_j^2 - 1) + b_i(b_j^2b_k^2 - 1)t\}
\{(b_k^2 - 1)(b_4^2-1) + b_ib_k b_4^2(b_j^2 b_k^2 - 1) t\}}. 
\end{array} 
\]
%%%%%%%%%%
In terms of this parametrization, the lifted transformation 
$\tilde{g}_j^2$ acts on the exceptional curve 
$e \simeq \P^1$ by the multiplication 
$t \mapsto b_j^2b_k^2 t$. 
Since $b_j^2 b_k^2 \neq 1$, the set 
$\wt{\mathrm{Fix}}_j^{e}(\th)$ consists of the two points, 
say $p$ and $q$, corresponding to $t = 0$ and $t = \infty$ 
(see Figure \ref{fig:A1A2}, left). 
On the other hand, the possible candidates for the smooth fixed 
points $\wt{\mathrm{Fix}}_j^{\circ}(\th)$ are only the points of 
labels $2$ and $4$ in Table \ref{tab:fixed}, since those of 
labels $1$ and $3$ do not satisfy the smoothness condition. 
Thus, 
%%%%%%%%%%%%%%%%%%%%%%%%% eqn:FixA1 %%%%%%%%%%%%%%%%%%%%%%%%%%%
\begin{equation} \label{eqn:FixA1} 
\wt{\mathrm{Fix}}_j(\th) = 
\{\!\{ \wt{P}(b_i,b_4^{-1};b_j,b_k), \, 
\wt{P}(b_j,b_k^{-1};b_i,b_4) \}\!\} \amalg \{p, q\}. 
\end{equation}
%%%%%%%%%%%%%%%%%%%%%%%%%%%%%%%%%%%%%%%%%%%%%%%%%%%%%%%%%%%%%%%
As for the Riccati periodic points 
$\wt{\mathrm{Per}}_j^{e}(\th;n)$, the discussion above 
implies that for any $n > 1$, 
%%%%%%
\[
\wt{\mathrm{Per}}_j^{e}(\th;n) = \left\{
\begin{array}{ll}
e \qquad & (\mbox{if $b_jb_k$ is a primitive $2n$-th root of unity}), 
\\[2mm]
\emptyset \qquad & (\mbox{otherwise}). 
\end{array}
\right.
\]
%%%%%%
\end{example}
%%%%%%%%%%%%%%%%%%%%%%%%%%%%%%%%%%%%%%%%%%%%%%%%%%%%%%%%%%%%%
%%%%%%%%%%%%%%%%%%%%%%%%% fig:A1A2 %%%%%%%%%%%%%%%%%%%%%%%%%%%
%%%%%%%%%%%%%%%%%%%%%%%%%%%%%%%%%%%%%%%%%%%%%%%%%%%%%%%%%%%%%
%%%%%%%%%%%%%%%%%%%%%%%%% ex:A2 %%%%%%%%%%%%%%%%%%%%%%%%%%%%%
\begin{example}[$\mbox{\boldmath $A_2$}$] \label{ex:A2}
Consider the $\wt{W}(D_4^{(1)})$-stratum of type $A_2$. 
We may assume that $\k_0 = \k_i = 0$ so that $b_0 = b_i = 1$. 
Then the surface $\Sol(\th)$ has a unique singular 
point of type $A_2$ at 
%%%%%%%%
\[
(x_i,x_j,x_k) = (b_4+b_4^{-1}, 
b_jb_4+b_j^{-1}b_4^{-1}, b_kb_4+b_k^{-1}b_4^{-1}). 
\]
%%%%%%%%
Blow up $\Sol(\th)$ at this point to obtain a 
minimal resolution (\ref{eqn:brieskorn}). 
Write the blowing-up as 
%%%%%%%%
\[
(x_i,x_j,x_k) = (u_i u_j + b_4+b_4^{-1}, 
u_j + b_jb_4+b_j^{-1}b_4^{-1}, 
u_k u_j + b_kb_4+b_k^{-1}b_4^{-1})
\]
%%%%%%%%
in terms of coordinates $(u_i,u_j,u_k)$. 
The exceptional set $e$ is the union of two lines 
%%%%%%%%%
\[ 
e^{+}\,\, : \quad u_j = b_k u_i + u_k + b_j b_k = 0, 
\qquad e^{-}\,\, : \quad 
u_j = b_k^{-1} u_i + u_k + b_j^{-1} b_k^{-1}= 0, 
\]
%%%%%%%%%
intersecting in a point. 
These lines are parametrized as 
%%%%%%%%%
\[
\begin{array}{rrcl}
e^{+} \,\, : \quad &(u_i, u_j, u_k) &=& \left(
\dfrac{b_j^2b_k^2-1}{b_j(1-b_k^2)+(b_j^2b_k^2-1)s}, \, 0, \, 
\dfrac{b_k(1-b_j^2)+b_jb_k(1-b_j^2b_k^2)s}{b_j(1-b_k^2)+
(b_j^2b_k^2-1)s}\right), \\[6mm]
e^{-} \,\, : \quad &(u_i, u_j, u_k) &=& \left(
\dfrac{b_j^2b_k^2-1}{b_j(1-b_k^2)+(b_j^2b_k^2-1)t}, \, 0, \, 
\dfrac{b_k(1-b_j^2)+b_j^{-1}b_k^{-1}(1-b_j^2b_k^2)t}{
b_j(1-b_k^2)+(b_j^2b_k^2-1)t}\right), 
\end{array}
\]
%%%%%%%%%
with the intersection point corresponding to $s = t = 0$. 
In terms of these parametrizations, the lifted 
transformation $\tilde{g}_j^2$ acts on $e^{+}$ and 
$e^{-}$ by the multiplications 
$s \mapsto b_j^{-2}b_k^{-2} s$ and $t \mapsto b_j^2b_k^2 t$, 
which are rewritten as $s \mapsto b_4^2 s$ and 
$t \mapsto b_4^{-2} t$, since $b_jb_kb_4 = 1$. 
Note that $b_4^2 \neq 1$, for otherwise $\k$ would be in 
the closure of the $\wt{W}(D_4^{(1)})$-stratum of type 
$(A_3)_i$. 
So $\tilde{g}_j^2$ has exactly two fixed points $p_0$ and 
$p_{+}$ on $e^{+}$ corresponding to $s = 0$ and $s = \infty$. 
Similarly $\tilde{g}_j^2$ has exactly two fixed points 
$p_0$ and $p_{-}$ on $e^{-}$ corresponding to $t = 0$ and 
$t = \infty$, where $p_0$ is the intersection point of 
$e^{+}$ and $e^{-}$ (see Figure \ref{fig:A1A2}, right). 
Thus we have 
$\wt{\mathrm{Fix}}_j^{e}(\th) = \{p_0,p_{+},p_{-}\}$. 
Next we consider the smooth fixed point of $\tilde{g}_j^2$ on 
$\wt{\mathcal{S}}(\th)$. 
Since we are assuming that $\k_i = \k_j+\k_k+\k_4 = 1$, 
the points of labels $1$, $2$, $3$ in Table \ref{tab:fixed} 
do not satisfy the smoothness condition and that of 
label $4$ is the only smooth fixed point. 
Thus $\wt{\mathrm{Fix}}_j^{\circ}(\th) = 
\{\wt{P}(b_j,b_k^{-1};b_i,b_4)\}$ and hence 
%%%%%%%%%%%%%%%%%%%%%%%%%%% eqn:FixA2 %%%%%%%%%%%%%%%%%%%%%%%
\begin{equation} \label{eqn:FixA2} 
\wt{\mathrm{Fix}}_j(\th) = 
\{\wt{P}(b_j,b_k^{-1};b_i,b_4), \, p_0, \, p_{+}, \, p_{-}\}. 
\end{equation} 
%%%%%%%%%%%%%%%%%%%%%%%%%%%%%%%%%%%%%%%%%%%%%%%%%%%%%%%%%%%%%
\end{example} 
%%%%%%%%%%%%%%%%%%%%%%%%%%%%%%%%%%%%%%%%%%%%%%%%%%%%%%%%%%%%%
\par 
In the remaining cases presented below, 
$\mathrm{Fix}_j(\th)$ contains at least one line component.  
%%%%%%%%%%%%%%%%%%%%%%%%% fig:A12 %%%%%%%%%%%%%%%%%%%%%%%%%%%
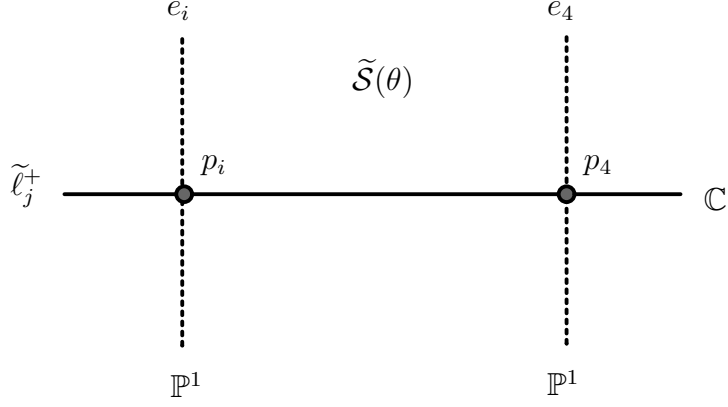
\begin{figure}[t] 
\begin{center}
%WinTpicVersion2.15
\unitlength 0.1in
\begin{picture}(35.90,19.70)(3.20,-23.00)
% LINE 0 0 3 0
% 2 600 1810 3790 1810
% 
\special{pn 20}%
\special{pa 600 1410}%
\special{pa 3790 1410}%
\special{fp}%
% LINE 0 2 3 0
% 2 1210 1000 1210 2600
% 
\special{pn 20}%
\special{pa 1210 600}%
\special{pa 1210 2200}%
\special{dt 0.054}%
\special{pa 1210 2200}%
\special{pa 1210 2199}%
\special{dt 0.054}%
% LINE 0 2 3 0
% 2 3200 990 3200 2590
% 
\special{pn 20}%
\special{pa 3200 590}%
\special{pa 3200 2190}%
\special{dt 0.054}%
\special{pa 3200 2190}%
\special{pa 3200 2189}%
\special{dt 0.054}%
% CIRCLE 0 0 0 0
% 4 1220 1810 1260 1800 1260 1800 1260 1800
% 
\special{pn 20}%
\special{sh 0.600}%
\special{ar 1220 1410 41 41  0.0000000 6.2831853}%
% CIRCLE 0 0 0 0
% 4 3200 1810 3240 1800 3240 1800 3240 1800
% 
\special{pn 20}%
\special{sh 0.600}%
\special{ar 3200 1410 41 41  0.0000000 6.2831853}%
% STR 2 0 3 0
% 3 320 1780 320 1880 2 0
% $\widetilde{\ell}_j^{+}$
\put(3.2000,-14.8000){\makebox(0,0)[lb]{$\widetilde{\ell}_j^{+}$}}%
% STR 2 0 3 0
% 3 3910 1790 3910 1890 2 0
% $\mathbb{C}$
\put(39.1000,-14.9000){\makebox(0,0)[lb]{$\mathbb{C}$}}%
% STR 2 0 3 0
% 3 3290 1610 3290 1710 2 0
% $p_4$
\put(32.9000,-13.1000){\makebox(0,0)[lb]{$p_4$}}%
% STR 2 0 3 0
% 3 1310 1610 1310 1710 2 0
% $p_i$
\put(13.1000,-13.1000){\makebox(0,0)[lb]{$p_i$}}%
% STR 2 0 3 0
% 3 1130 800 1130 900 2 0
% $e_i$
\put(11.3000,-5.0000){\makebox(0,0)[lb]{$e_i$}}%
% STR 2 0 3 0
% 3 3100 800 3100 900 2 0
% $e_4$
\put(31.0000,-5.0000){\makebox(0,0)[lb]{$e_4$}}%
% STR 2 0 3 0
% 3 1150 2770 1150 2870 2 0
% $\mathbb{P}^1$
\put(11.5000,-24.7000){\makebox(0,0)[lb]{$\mathbb{P}^1$}}%
% STR 2 0 3 0
% 3 3100 2760 3100 2860 2 0
% $\mathbb{P}^1$
\put(31.0000,-24.6000){\makebox(0,0)[lb]{$\mathbb{P}^1$}}%
% STR 2 0 3 0
% 3 2090 1210 2090 1310 2 0
% $\widetilde{\mathcal{S}}(\th)$
\put(20.9000,-9.1000){\makebox(0,0)[lb]{$\widetilde{\mathcal{S}}(\th)$}}%
\end{picture}%
\end{center}
\caption{Surface of type $(A_1^{\oplus 2})_i$} 
\label{fig:A12} 
\end{figure}
%%%%%%%%%%%%%%%%%%%%%%%%% ex:A12 %%%%%%%%%%%%%%%%%%%%%%%%%%%%
\begin{example}[$\mbox{\boldmath $A_1^{\oplus 2}$}$]  
\label{ex:A12} 
First we consider $\wt{\mathrm{Fix}}_j(\th)$ and 
$\wt{\mathrm{Per}}_j^{e}(\th)$ on the 
$\wt{W}(D_4^{(1)})$-stratum of type $(A_1^{\oplus 2})_i$. 
We may assume that $\k_j = \k_k = 0$ so that 
$b_j = b_k = 1$. 
Since our stratum is not of type $(A_3)_i$ nor 
of type $D_4$, we have $(b_ib_4-1)(b_ib_4^{-1}-1) \neq 0$ 
or equivalently $b_i+b_i^{-1} \neq b_4+b_4^{-1}$. 
In this case $\mathrm{Fix}_j(\th)$ contains the line 
$\ell_j^{+}$ but does not the line $\ell_j^{-}$ and 
the surface $\Sol(\th)$ has two singular points of type 
$A_1$ at 
$(x_i,x_j,x_k) = (2, b_i+b_i^{-1}, b_4+b_4^{-1})$ and 
$(x_i,x_j,x_k) = (2, b_4+b_4^{-1}, b_i+b_i^{-1})$. 
We denote the former singularity by $q_i$ and the latter 
by $q_4$ respectively; both singularities lie on the line 
$\ell_j^{+}$. 
Blow up $\Sol(\th)$ at these points to obtain a 
minimal resolution as in (\ref{eqn:brieskorn}). 
Let $\wt{\ell}_j^{+}$ be the strict transform of 
$\ell_j^{+}$, and let $e_i$ and $e_4$ be the exceptional 
curves over $q_i$ and $q_4$ respectively. 
Moreover let $p_i$  be the intersection point 
of $\wt{\ell}_j^{+}$ and $e_i$. 
Similarly let $p_4$ be the intersection point of 
$\wt{\ell}_j^{+}$ and $e_4$ (see Figure \ref{fig:A12}). 
Then the blowing-up at the point $q_i$ is represented as 
%%%%%%
\[
(x_i,x_j,x_k) = 
(u_iu_j+2, u_j+b_i+b_i^{-1}, u_ku_j+b_4+b_4^{-1})
\]
%%%%%
in terms of coordinates $(u_i,u_j,u_k)$ around $(0,0,0)$. 
The strict transform $\wt{\ell}_j^{+}$ and the exceptional 
curve $e_i$ are given by $u_i = u_k +1 = 0$ and 
%%%%%
\[ 
u_j = (b_ib_4)(u_i^2+u_k^2) + (b_i^2+1)b_4 (u_iu_k) 
+b_i(b_4^2+1) u_i + 2 (b_ib_4) u_k + (b_ib_4) = 0. 
\]
%%%%%
The exceptional curve $e_i$ admits a parametrization 
%%%%%%%%%%%%%%%%%%%%% eqn:paramA12 %%%%%%%%%%%%%%%%%%%%%%%
\begin{equation} \label{eqn:paramA12}
u_i = \dfrac{(b_ib_4-1)(b_ib_4^{-1}-1)}{(t+b_i)(b_it+1)}, 
\qquad 
u_j = 0, 
\qquad 
u_k = -\dfrac{b_i(t+b_4)(b_4t+1)}{b_4(t+b_i)(b_it+1)}, 
\end{equation}
%%%%%%%%%%%%%%%%%%%%%%%%%%%%%%%%%%%%%%%%%%%%%%%%%%%%%%%%%%
where the intersection point $p_i$ has coordinates 
$(u_i,u_j,u_k) = (0,0,-1)$, which corresponds to 
$t = \infty$. 
The lifted transformation $\tilde{g}_j^2$ acts on $e_i$ as 
a M\"{o}bius transformation fixing $p_i$. 
Some computations show that in terms of the variable $t$ 
this transformation is just the shift 
%%%%%
\[
t \mapsto t + (b_i+b_i^{-1})-(b_4+b_4^{-1}). 
\]
%%%%
and hence a parabolic transformation. 
Thus $\tilde{g}_j^2$ has no periodic points on $e_i$ 
other than the fixed point $p_i$. 
By symmetry, $\tilde{g}_j^2$ also acts on $e_4$ as a 
parabolic M\"obius transformation fixing $p_4$ only. 
Summarizing the arguments, we conclude that 
on the $\wt{W}(D_4^{(1)})$-stratum of type 
$(A_1^{\oplus 2})_i$, 
%%%%%%%%%%%%%%%%%%%%%%% eqn:FixA12i %%%%%%%%%%%%%%%%%%%%%%%%
\begin{equation} \label{eqn:FixA12i} 
\wt{\mathrm{Fix}}_j(\th) = \wt{\ell}_j^{+} \amalg 
\{\, \wt{P}(b_i,b_4;b_j,b_k), \, 
\wt{P}(b_i,b_4^{-1};b_j,b_k) \,\}, \qquad 
\wt{\mathrm{Per}}_j^{e}(\th;n) = \emptyset \quad (n > 1).  
\end{equation} 
%%%%%%%%%%%%%%%%%%%%%%%%%%%%%%%%%%%%%%%%%%%%%%%%%%%%%%%%%%%%
\par 
Next we consider $\wt{\mathrm{Fix}}_i(\th)$ on the 
$\wt{W}(D_4^{(1)})$-stratum of type $(A_1^{\oplus 2})_i$. 
Some calculations show that there are parametrizations 
of $e_i$ and $e_4$ such that $\tilde{g}_i^2$ acts on $e_i$ 
and $e_4$ as the multiplications $t \mapsto b_4^2 t$ and 
$t \mapsto b_i^2 t$ respectively. 
(Modify (\ref{eqn:paramA12}) to get such parametrization.) 
Since $b_4^2 \neq 1$ and $b_i^2 \neq 1$, the transformation 
$\tilde{g}_i^2$ has exactly two fixed points, 
say $p_{ii}$ and $q_{ii}$, on $e_i$, and exactly two fixed 
points, say $p_{i4}$ and $q_{i4}$, on $e_4$. 
There are no smooth fixed points 
$\wt{\mathrm{Fix}}_i^{\circ}(\th)$, because the smoothness 
condition of Table \ref{tab:fixed} with $(i,j,k)$ replaced 
by $(k,i,j)$ is not satisfied for any labels there. 
Thus we have $\wt{\mathrm{Fix}}_i(\th) = 
\wt{\mathrm{Fix}}_i^{e}(\th) = 
\{p_{ii}, \, q_{ii}, \, p_{i4}, \, q_{i4}\}$ and 
$\wt{\mathrm{Fix}}_i^{\circ}(\th) = \emptyset$. 
By symmetry there is a similar characterization of 
$\wt{\mathrm{Fix}}_k(\th)$. 
By permuting the indices $(i,j,k)$, we have 
%%%%%%%%%%%%%%%%%%%%%% eqn:FixA12jk %%%%%%%%%%%%%%%%%%%%%%%%%%
\begin{equation} \label{eqn:FixA12jk}
\wt{\mathrm{Fix}}_j(\th) = \wt{\mathrm{Fix}}_j^{e}(\th) = 
\{\mbox{four points}\}, 
\qquad \wt{\mathrm{Fix}}_j^{\circ}(\th) = \emptyset, 
\end{equation}
%%%%%%%%%%%%%%%%%%%%%%%%%%%%%%%%%%%%%%%%%%%%%%%%%%%%%%%%%%%%%%
on the $\wt{W}(D_4^{(1)})$-strata of types 
$(A_1^{\oplus 2})_j$ and $(A_1^{\oplus 2})_k$. 
A slightly further consideration yields 
%%%%%
\[
\wt{\mathrm{Per}}_j^{e}(\th;n) =
\left\{\begin{array}{ll}
e_j \amalg e_4 \quad & (\mbox{if $b_j$ and $b_4$ are 
primitive $2n$-th roots of unity}), \\[2mm] 
e_j & (\mbox{if $b_4$ is a primitive $2n$-th root of 
unity, but $b_j$ is not}), \\[2mm] 
e_4 & (\mbox{if $b_j$ is a primitive $2n$-th root of 
unity, but $b_4$ is not}), \\[2mm] 
\emptyset \qquad & (\mbox{otherwise}). 
\end{array}\right.
\]
%%%%%%%%%%%%%%%
on the stratum $(A_1^{\oplus 2})_j$ and a similar 
characterization of it on the stratum $(A_1^{\oplus 2})_k$. 
\end{example}
%%%%%%%%%%%%%%%%%%%%%%%%%%%%%%%%%%%%%%%%%%%%%%%%%%%%%%%%%%%%%
%%%%%%%%%%%%%%%%%%%%%%%%%%% fig:A3 %%%%%%%%%%%%%%%%%%%%%%%%%%
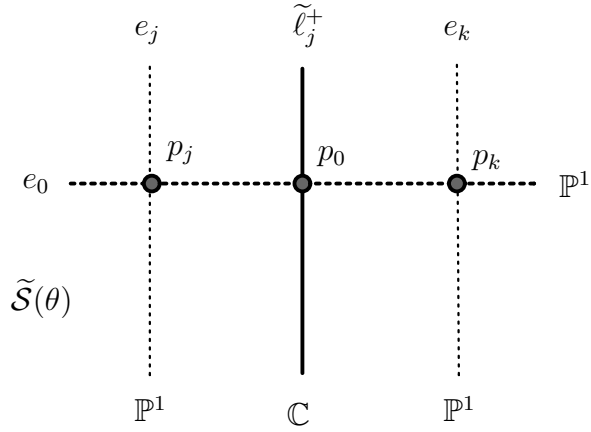
\begin{figure}[t]
\begin{center}
%WinTpicVersion2.15
\unitlength 0.1in
\begin{picture}(28.40,20.10)(2.90,-25.00)
% LINE 0 2 3 0
% 2 600 1810 3000 1810
% 
\special{pn 20}%
\special{pa 600 1410}%
\special{pa 3000 1410}%
\special{dt 0.054}%
\special{pa 3000 1410}%
\special{pa 2999 1410}%
\special{dt 0.054}%
% LINE 0 0 3 0
% 2 1800 1210 1800 2800
% 
\special{pn 20}%
\special{pa 1800 810}%
\special{pa 1800 2400}%
\special{fp}%
% LINE 1 2 3 0
% 2 1010 1200 1010 2810
% 
\special{pn 13}%
\special{pa 1010 800}%
\special{pa 1010 2410}%
\special{dt 0.045}%
\special{pa 1010 2410}%
\special{pa 1010 2409}%
\special{dt 0.045}%
% LINE 1 2 3 0
% 2 2600 1190 2610 2810
% 
\special{pn 13}%
\special{pa 2600 790}%
\special{pa 2610 2410}%
\special{dt 0.045}%
\special{pa 2610 2410}%
\special{pa 2610 2409}%
\special{dt 0.045}%
% CIRCLE 0 0 0 0
% 4 1800 1810 1810 1850 1810 1850 1810 1850
% 
\special{pn 20}%
\special{sh 0.600}%
\special{ar 1800 1410 41 41  0.0000000 6.2831853}%
% CIRCLE 0 0 0 0
% 4 1020 1810 1030 1850 1030 1850 1030 1850
% 
\special{pn 20}%
\special{sh 0.600}%
\special{ar 1020 1410 41 41  0.0000000 6.2831853}%
% CIRCLE 0 0 0 0
% 4 2600 1810 2610 1850 2610 1850 2610 1850
% 
\special{pn 20}%
\special{sh 0.600}%
\special{ar 2600 1410 41 41  0.0000000 6.2831853}%
% STR 2 0 3 0
% 3 350 1760 350 1860 2 0
% $e_0$
\put(3.5000,-14.6000){\makebox(0,0)[lb]{$e_0$}}%
% STR 2 0 3 0
% 3 1750 1000 1750 1100 2 0
% $\widetilde{\ell}_j^{+}$
\put(17.5000,-7.0000){\makebox(0,0)[lb]{$\widetilde{\ell}_j^{+}$}}%
% STR 2 0 3 0
% 3 930 980 930 1080 2 0
% $e_j$
\put(9.3000,-6.8000){\makebox(0,0)[lb]{$e_j$}}%
% STR 2 0 3 0
% 3 2530 960 2530 1060 2 0
% $e_k$
\put(25.3000,-6.6000){\makebox(0,0)[lb]{$e_k$}}%
% STR 2 0 3 0
% 3 3130 1780 3130 1880 2 0
% $\mathbb{P}^1$
\put(31.3000,-14.8000){\makebox(0,0)[lb]{$\mathbb{P}^1$}}%
% STR 2 0 3 0
% 3 2530 2950 2530 3050 2 0
% $\mathbb{P}^1$
\put(25.3000,-26.5000){\makebox(0,0)[lb]{$\mathbb{P}^1$}}%
% STR 2 0 3 0
% 3 1720 2970 1720 3070 2 0
% $\mathbb{C}$
\put(17.2000,-26.7000){\makebox(0,0)[lb]{$\mathbb{C}$}}%
% STR 2 0 3 0
% 3 930 2950 930 3050 2 0
% $\mathbb{P}^1$
\put(9.3000,-26.5000){\makebox(0,0)[lb]{$\mathbb{P}^1$}}%
% STR 2 0 3 0
% 3 1100 1610 1100 1710 2 0
% $p_j$
\put(11.0000,-13.1000){\makebox(0,0)[lb]{$p_j$}}%
% STR 2 0 3 0
% 3 1880 1620 1880 1720 2 0
% $p_0$
\put(18.8000,-13.2000){\makebox(0,0)[lb]{$p_0$}}%
% STR 2 0 3 0
% 3 2680 1640 2680 1740 2 0
% $p_k$
\put(26.8000,-13.4000){\makebox(0,0)[lb]{$p_k$}}%
% STR 2 0 3 0
% 3 290 2410 290 2510 2 0
% $\wt{\mathcal{S}}(\th)$
\put(2.9000,-21.1000){\makebox(0,0)[lb]{$\wt{\mathcal{S}}(\th)$}}%
\end{picture}%
\end{center}
\caption{Surface of type $(A_3)_i$} 
\label{fig:A3} 
\end{figure}
%%%%%%%%%%%%%%%%%%%%%%%%%%%%%%%%%%%%%%%%%%%%%%%%%%%%%%%%%%%%%
%%%%%%%%%%%%%%%%%%%%%%%%%%% ex:A3 %%%%%%%%%%%%%%%%%%%%%%%%%%%
\begin{example}[$\mbox{\boldmath $A_3$}$]  
\label{ex:A3} 
First we consider $\wt{\mathrm{Fix}}_j(\th)$ and 
$\wt{\mathrm{Per}}_j^{e}(\th)$ on the 
$\wt{W}(D_4^{(1)})$-stratum of type $(A_3)_i$. 
We may assume that $\k_0 = \k_j = \k_k = 0$ and 
$\k_i + \k_4 = 1$ so that $b_j = b_k = 1$ and $b_ib_4 = 1$. 
But we have $b_ib_4^{-1} \not\in \{\pm1\}$, since our 
stratum is not of type $D_4$. 
In this case $\mathrm{Fix}_j(\th)$ contains the line 
$\ell_j^{+}$ but does not the line $\ell_j^{-}$. 
The surface $\Sol(\th)$ has only one singular point of type 
$A_3$ at $(x_i,x_j,x_k) = (2, b_4+b_4^{-1}, b_4+b_4^{-1})$, 
which lies on the line $\ell_j^{+}$. 
Blow up the singular point. 
This blowing-up is expressed as $(x_i,x_j,x_k) = 
(u_iu_j+2, u_j+b_4+b_4^{-1}, u_ku_j+b_4+b_4^{-1})$ 
in terms of coordinates $(u_i,u_j,u_k)$ around $(0,0,0)$. 
The strict transform of the surface $\Sol(\th)$ is given by 
%%%%%%%
\[
u_j = b_4 u_iu_ju_k + b_4 u_i^2 + + b_4 u_k^2 + 
(b_4^2+1) u_i u_k + (b_4^2+1) u_i + 2 b_4 u_k + b_4 = 0, 
\]
%%%%%%%
which has yet one singular point, say $q$. 
The exceptional curve consists of two line components 
$u_j = u_i + b_4 u_k + b_4 = 0$ and 
$u_j = b_4 u_i + u_k + 1 = 0$, whose 
intersection point $(u_i,u_j,u_k) = (0,0,-1)$ is exactly the 
singular point $q$. 
The strict transform of $\ell_j^{+}$ is now given by 
$u_i = u_k + 1 = 0$, which also passes through $q$. 
Blow up again the singular point $q$. 
Let $e_0$ be the exceptional curve and let $e_j$, 
$e_k$, $\tilde{\ell}_j^{+}$ be the strict transforms 
of the lines $u_j = u_i + b_4 u_k + b_4 = 0$, 
$u_j = b_4 u_i + u_k + 1 = 0$, $u_i = u_k + 1 = 0$ 
respectively. 
If we express this blowing-up as 
$(u_i,u_j,u_k)= (v_i, v_i v_j, v_i v_k - 1)$, then 
the exceptional curve $e_0$ is given by 
$v_i = b_4 - b_4 v_j + (b_4^2+1) v_k + b_4 v_k^2 = 0$; 
$e_j$ is given by $v_j = 1 + b_4 v_k = 0$; 
and $e_k$ is given by $v_j = b_4 + v_k = 0$.
The intersection point of $e_0$ and $e_j$ is 
$(v_i,v_j,v_k) = (0,0,-b_i)$ and that of $e_0$ and 
$e_k$ is $(v_i,v_j,v_k) = (0,0,-b_4)$. 
If $e_j$ is parametrized as 
$(v_i,v_j,v_k) = ((t+b_i)^{-1},0,-b_i)$, then 
the transformation $\tilde{g}_j^2$ acts on $e_j$ 
as the shift $t \mapsto t + b_4-b_i$. 
Similarly, if $e_k$ is parametrized as 
$(v_i,v_j,v_k) = ((t+b_4)^{-1},0,-b_4)$, then 
$\tilde{g}_j^2$ acts on $e_k$ as the shift 
$t \mapsto t + b_4-b_i$. 
Hence $\tilde{g}_j^2$ acts on $e_j$ and $e_k$ 
as parabolic M\"obius transformations fixing 
only $p_j$ and $q_j$. 
Then $\tilde{g}_j^2$ acts on $e_0$ as the identity, because 
it also fixes the intersection point $p_0$ of $e_0$ and 
$\tilde{\ell}_j^{+}$. 
Summarizing the arguments we see that on the 
stratum of type $(A_3)_i$, 
%%%%%%%%%%%%%%%%%%%%%%% eqn:FixA3i %%%%%%%%%%%%%%%%%%%%%%%%%%
\begin{equation} \label{eqn:FixA3i}
\wt{\mathrm{Fix}}_j(\th) = 
\wt{\ell}_j^{+} \underset{p_0}{\cup} e_0 \amalg 
\{\, \wt{P}(b_i,b_4^{-1};b_j,b_k) \,\}, \qquad 
\wt{\mathrm{Per}}_j^{e}(\th;n) = \emptyset \qquad (n > 1), 
\end{equation}
%%%%%%%%%%%%%%%%%%%%%%%%%%%%%%%%%%%%%%%%%%%%%%%%%%%%%%%%%%%%
where $\wt{\ell}_j^{+} \underset{p_0}{\cup} e_0$ indicates 
that the curves $\wt{\ell}_j^{+}$ and $e_0$ meet in the 
point $p_0$. 
\par 
Next we consider $\wt{\mathrm{Fix}}_i(\th)$ and 
$\wt{\mathrm{Per}}_i^{e}(\th;n)$ on the 
$\wt{W}(D_4^{(1)})$-stratum of type $(A_3)_i$. 
If we take a parametrization of $e_0$ such that 
$t = 0$ and $t = \infty$ correspond to the points 
$p_j$ and $p_k$ respectively, then a simple check shows 
that the transformation $\tilde{g}_i^2$ on $e_j$ 
is expressed as $t \mapsto b_4^{-2} t$. 
There is a parametrization of $e_j$ such that $t = 0$ 
corresponds to $p_j$ and $\tilde{g}_i^2$ is given by 
$t \mapsto b_4^2 t$. 
Since $b_4^2 \neq 1$, the transformation $\tilde{g}_i^2$ has 
exactly two fixed points on $e_j$, one of which is just $p_j$ 
and the other is denoted by $p_{ij}$. 
Similarly, there is a parametrization of $e_k$ such that 
$t = 0$ corresponds to $p_k$ and $\tilde{g}_i^2$ is given 
by $t \mapsto b_4^{-2} t$, and hence $\tilde{g}_i^2$ has 
exactly two fixed points on $e_k$, one of which is just 
$p_k$ and the other is denoted by $p_{ik}$. 
There are no smooth fixed points 
$\wt{\mathrm{Fix}}_i^{\circ}(\th)$, because the smoothness 
condition of Table \ref{tab:fixed} with $(i,j,k)$ replaced 
by $(k,i,j)$ is not satisfied for any labels there. 
So we have $\wt{\mathrm{Fix}}_i(\th) = 
\wt{\mathrm{Fix}}_i^{e}(\th) = 
\{p_j, \, p_{ij}, \, p_k, \, p_{ik}\}$ and 
$\wt{\mathrm{Fix}}_i^{\circ}(\th) = \emptyset$ on the 
$\wt{W}(D_4^{(1)})$-stratum of type $(A_3)_i$. 
By symmetry there is a similar characterization of 
$\wt{\mathrm{Fix}}_k(\th)$ on the same stratum. 
By permuting the indices $(i,j,k)$, we have 
%%%%%%%%%%%%%%%%%%%%%%% eqn:FixA3jk %%%%%%%%%%%%%%%%%%%%%%%%%%
\begin{equation} \label{eqn:FixA3jk}
\wt{\mathrm{Fix}}_j(\th) = 
\wt{\mathrm{Fix}}_j^{e}(\th) = 
\{\mbox{four points}\}, \qquad 
\wt{\mathrm{Fix}}_j^{\circ}(\th) = \emptyset, 
\end{equation} 
%%%%%%%%%%%%%%%%%%%%%%%%%%%%%%%%%%%%%%%%%%%%%%%%%%%%%%%%%%%%%%
on the $\wt{W}(D_4^{(1)})$-strata of types $(A_3)_j$ and 
$(A_3)_k$. 
A slightly further consideration yields 
%%%%%%%
\[
\wt{\mathrm{Per}}_j^{e}(\th;n) = 
\left\{\begin{array}{ll}
e_k \underset{p_k}{\cup} e_0 \underset{p_i}{\cup} e_i 
\quad & (\mbox{if $b_4$ is a primitive $2n$-th root 
of unity}), \\[2mm]
\emptyset \qquad & (\mbox{otherwise}), 
\end{array}\right. 
\]
%%%%%
on the stratum $(A_3)_j$ and a similar 
characterization of it on the stratum $(A_3)_k$. 
\end{example}
%%%%%%%%%%%%%%%%%%%%%%%%%%%%%%%%%%%%%%%%%%%%%%%%%%%%%%%%%%%%%
%%%%%%%%%%%%%%%%%%%%%%%%%%% fig:A13 %%%%%%%%%%%%%%%%%%%%%%%%%
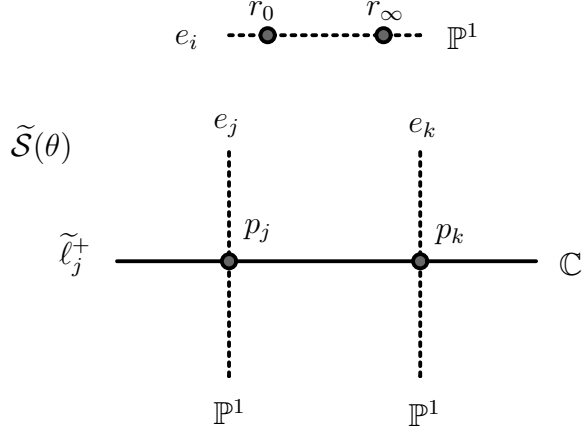
\begin{figure}[t]
\begin{center}
%WinTpicVersion2.15
\unitlength 0.1in
\begin{picture}(28.40,21.20)(0.70,-22.80)
% LINE 0 0 3 0
% 2 620 1990 2790 1990
% 
\special{pn 20}%
\special{pa 620 1590}%
\special{pa 2790 1590}%
\special{fp}%
% LINE 0 2 3 0
% 2 1200 1420 1200 2590
% 
\special{pn 20}%
\special{pa 1200 1020}%
\special{pa 1200 2190}%
\special{dt 0.054}%
\special{pa 1200 2190}%
\special{pa 1200 2189}%
\special{dt 0.054}%
% LINE 0 2 3 0
% 2 2190 1420 2190 2590
% 
\special{pn 20}%
\special{pa 2190 1020}%
\special{pa 2190 2190}%
\special{dt 0.054}%
\special{pa 2190 2190}%
\special{pa 2190 2189}%
\special{dt 0.054}%
% LINE 0 2 3 0
% 2 1200 810 2190 810
% 
\special{pn 20}%
\special{pa 1200 410}%
\special{pa 2190 410}%
\special{dt 0.054}%
\special{pa 2190 410}%
\special{pa 2189 410}%
\special{dt 0.054}%
% STR 2 0 3 0
% 3 320 1970 320 2070 2 0
% $\widetilde{\ell}_j^{+}$
\put(3.2000,-16.7000){\makebox(0,0)[lb]{$\widetilde{\ell}_j^{+}$}}%
% STR 2 0 3 0
% 3 1120 1240 1120 1340 2 0
% $e_j$
\put(11.2000,-9.4000){\makebox(0,0)[lb]{$e_j$}}%
% STR 2 0 3 0
% 3 2130 1240 2130 1340 2 0
% $e_k$
\put(21.3000,-9.4000){\makebox(0,0)[lb]{$e_k$}}%
% STR 2 0 3 0
% 3 920 780 920 880 2 0
% $e_i$
\put(9.2000,-4.8000){\makebox(0,0)[lb]{$e_i$}}%
% STR 2 0 3 0
% 3 1120 2740 1120 2840 2 0
% $\mathbb{P}^1$
\put(11.2000,-24.4000){\makebox(0,0)[lb]{$\mathbb{P}^1$}}%
% STR 2 0 3 0
% 3 2130 2750 2130 2850 2 0
% $\mathbb{P}^1$
\put(21.3000,-24.5000){\makebox(0,0)[lb]{$\mathbb{P}^1$}}%
% STR 2 0 3 0
% 3 2330 780 2330 880 2 0
% $\mathbb{P}^1$
\put(23.3000,-4.8000){\makebox(0,0)[lb]{$\mathbb{P}^1$}}%
% STR 2 0 3 0
% 3 2910 1960 2910 2060 2 0
% $\mathbb{C}$
\put(29.1000,-16.6000){\makebox(0,0)[lb]{$\mathbb{C}$}}%
% CIRCLE 0 0 0 0
% 4 1200 1990 1236 1980 1236 1980 1236 1980
% 
\special{pn 20}%
\special{sh 0.600}%
\special{ar 1200 1590 37 37  0.0000000 6.2831853}%
% CIRCLE 0 0 0 0
% 4 2190 1990 2226 1980 2226 1980 2226 1980
% 
\special{pn 20}%
\special{sh 0.600}%
\special{ar 2190 1590 37 37  0.0000000 6.2831853}%
% CIRCLE 0 0 0 0
% 4 1400 810 1436 800 1436 800 1436 800
% 
\special{pn 20}%
\special{sh 0.600}%
\special{ar 1400 410 37 37  0.0000000 6.2831853}%
% CIRCLE 0 0 0 0
% 4 2000 810 2036 800 2036 800 2036 800
% 
\special{pn 20}%
\special{sh 0.600}%
\special{ar 2000 410 37 37  0.0000000 6.2831853}%
% STR 2 0 3 0
% 3 1280 1790 1280 1890 2 0
% $p_j$
\put(12.8000,-14.9000){\makebox(0,0)[lb]{$p_j$}}%
% STR 2 0 3 0
% 3 2270 1800 2270 1900 2 0
% $p_k$
\put(22.7000,-15.0000){\makebox(0,0)[lb]{$p_k$}}%
% STR 2 0 3 0
% 3 1300 630 1300 730 2 0
% $r_0$
\put(13.0000,-3.3000){\makebox(0,0)[lb]{$r_0$}}%
% STR 2 0 3 0
% 3 1910 630 1910 730 2 0
% $r_{\infty}$
\put(19.1000,-3.3000){\makebox(0,0)[lb]{$r_{\infty}$}}%
% STR 2 0 3 0
% 3 70 1380 70 1480 2 0
% $\widetilde{\mathcal{S}}(\theta)$
\put(0.7000,-10.8000){\makebox(0,0)[lb]{$\widetilde{\mathcal{S}}(\theta)$}}%
\end{picture}%
\end{center}
\caption{Surface of type $A_1^{\oplus 3}$}
\label{fig:A13} 
\end{figure}
%%%%%%%%%%%%%%%%%%%%%%%%%%% ex:A13 %%%%%%%%%%%%%%%%%%%%%%%%%%
\begin{example}[$\mbox{\boldmath $A_1^{\oplus 3}$}$] 
\label{ex:A13} 
Consider the $\wt{W}(D_4^{(1)})$-stratum of type 
$A_1^{\oplus 3}$. 
We may assume that $\k_i = \k_j = \k_k = 0$ so that 
$b_i = b_j = b_k = 1$. 
But we have $b_4 \not\in \{\pm1\}$ since our stratum is not 
of type $D_4$ nor of type $A_1^{\oplus 4}$. 
In this case the surface $\Sol(\th)$ has three singular 
points of type $A_1$ at $(x_i,x_j,x_k) = (b_4+b_4^{-1},2,2)$, 
$(2,b_4+b_4^{-1},2)$, $(2,2,b_4+b_4^{-1})$, which are called 
$q_i$, $q_j$, $q_k$ respectively. 
Note that the two points $q_j$ and $q_k$ lie on the line 
$\ell_j^{+}$ but $q_i$ does not lie on the union 
$\ell_j^{+} \amalg \ell_j^{-}$. 
The minimal resolution (\ref{eqn:brieskorn}) is obtained 
by blowing up these three points (see Figure \ref{fig:A13}). 
First, consider the blowing-up at $q_k$ and represent it 
by $(x_i,x_j,x_k) = (u_iu_j+2,u_j+2,u_ku_j+b_4+b_4^{-1})$. 
Then the strict transform $\wt{\ell}_j^{+}$ of the line 
$\ell_j^{+}$ is given by $u_i = u_k + 1 = 0$, while the 
exceptional curve $e_k$ is given by 
$b_4(u_i + u_k + 1)^2 + (b_4 - 1)^2u_i = 0$. 
The curves $\wt{\ell}_j^{+}$ and $e_k$ intersect in the point 
$(u_i,u_j,u_k) = (0,0,-1)$; this point is called 
$p_k$. 
If we parametrize the curve $e_k$ as 
%%%%%%%%
\[
u_i = - \dfrac{b_4}{(b_4-1)^2 t^2}, \quad 
u_j = 0, \quad 
u_k = - \dfrac{\{(b_4-1)t+1\}\{(b_4-1)t-b_4\}}{(b_4-1)^2 t^2} 
\qquad (t \in \P^1),  
\]
%%%%%%%%
where $t = \infty$ corresponds to the point 
$p_k$, then the lifted transformation $\wt{g}_j^2$ induces 
the shift $t \mapsto t+1$ and hence acts on $e_k$ as a 
parabolic M\"{o}bius transformation fixing $p_k$ only. 
In a similar manner $\wt{g}_j^2$ acts on the exceptional 
curve $e_j$ over $q_j$ as a parabolic M\"{o}bius 
transformation fixing only the intersection point $p_j$ of 
$\wt{\ell}_j^{+}$ and $e_j$. 
Next we consider the blowing-up at $q_i$ and represent it by 
$(x_i,x_j,x_k)=(u_iu_j+b_4+b_4^{-1},u_j+2,u_ku_j+2)$. 
Then the exceptional curve $e_i$ is given by 
$b_4(u_i+u_k+1)^2+(b_4-1)^2 u_k =0$, which can be parametrized 
as 
%%%%%%%%%
\[
u_i = - \dfrac{(b_4+1)^2 t}{(b_4 t+1)^2}, \quad 
u_j = 0, \quad 
u_k = - \dfrac{b_4 (t-1)^2}{(b_4 t+1)^2} \qquad 
(t \in \P^1). 
\]
%%%%%%%%%
In terms of this parametrization, the transformation 
$\wt{g}_j^2$ restricts to the map $t \mapsto b_4^2 t$ on 
the exceptional curve $e_i$. 
Let $r_0$ and $r_{\infty}$ be the points on $e_i$ 
corresponding to $t = 0$ and $t = \infty$ respectively. 
Since $b_4^2 \neq 1$, the map 
$\wt{g}_j^2$ acts on $e_i$ as a M\"obius transformation 
with exactly two fixed points $r_0$ and $r_{\infty}$. 
Hence the set $\wt{\mathrm{Fix}}_j(\th)$ 
contains the line component 
$\wt{\ell}_j^{+}$ and the Riccati component 
$\{r_0, r_{\infty}\}$, but has no smooth-point 
component, since $F(b_i,b_4;b_j,b_k) = 
F(b_i,b_4^{-1};b_j,b_k) = b_4 + b_4^{-1} \not\in \{\pm2\}$ 
is a double root of the quartic equation 
(\ref{eqn:quartic}) (see Theorem \ref{thm:quartic}). 
Thus we have 
%%%%%%%%%%%%%%%%%%%%%%%%% eqn:FixA13 %%%%%%%%%%%%%%%%%%%%%%
\begin{equation} \label{eqn:FixA13} 
\wt{\mathrm{Fix}}_j(\th) = \wt{\ell}_j^{+} \amalg 
\{r_0, r_{\infty}\}. 
\end{equation}
%%%%%%%%%%%%%%%%%%%%%%%%%%%%%%%%%%%%%%%%%%%%%%%%%%%%%%%%%%%
As for the Riccati periodic points, since $\wt{g}_j^2$ acts 
on $e_i$ as $t \mapsto b_4^2 t$, we have for any $n > 0$, 
%%%%%%%%%
\[
\wt{\mathrm{Per}}_j^{e}(\th;n) = 
\left\{\begin{array}{ll}
e_i \quad & (\mbox{if $b_4$ is a primitive 
$2n$-th root of unity}), \\[2mm]
\emptyset \qquad & (\mbox{otherwise}). 
\end{array}\right.
\]
\end{example}
%%%%%%%%%%%%%%%%%%%%%%%%%%%%%%%%%%%%%%%%%%%%%%%%%%%%%%%%%%%%%
%%%%%%%%%%%%%%%%%%%%%%%%%%% fig:D4 %%%%%%%%%%%%%%%%%%%%%%%%%%
\begin{figure}[t]
\begin{center}
%WinTpicVersion2.15
\unitlength 0.1in
\begin{picture}(43.90,27.60)(3.20,-30.90)
% LINE 0 2 3 0
% 2 810 1400 4410 1400
% 
\special{pn 20}%
\special{pa 810 1000}%
\special{pa 4410 1000}%
\special{dt 0.054}%
\special{pa 4410 1000}%
\special{pa 4409 1000}%
\special{dt 0.054}%
% LINE 0 2 3 0
% 2 2600 1000 2600 3400
% 
\special{pn 20}%
\special{pa 2600 600}%
\special{pa 2600 3000}%
\special{dt 0.054}%
\special{pa 2600 3000}%
\special{pa 2600 2999}%
\special{dt 0.054}%
% LINE 1 2 3 0
% 2 1200 1000 1200 3400
% 
\special{pn 13}%
\special{pa 1200 600}%
\special{pa 1200 3000}%
\special{dt 0.045}%
\special{pa 1200 3000}%
\special{pa 1200 2999}%
\special{dt 0.045}%
% LINE 1 2 3 0
% 2 4020 990 4020 3390
% 
\special{pn 13}%
\special{pa 4020 590}%
\special{pa 4020 2990}%
\special{dt 0.045}%
\special{pa 4020 2990}%
\special{pa 4020 2989}%
\special{dt 0.045}%
% LINE 0 0 3 0
% 2 2000 2200 3210 2200
% 
\special{pn 20}%
\special{pa 2000 1800}%
\special{pa 3210 1800}%
\special{fp}%
% LINE 1 0 3 0
% 2 600 3000 1810 3000
% 
\special{pn 13}%
\special{pa 600 2600}%
\special{pa 1810 2600}%
\special{fp}%
% LINE 1 0 3 0
% 2 3420 3000 4630 3000
% 
\special{pn 13}%
\special{pa 3420 2600}%
\special{pa 4630 2600}%
\special{fp}%
% STR 2 0 3 0
% 3 1700 2160 1700 2260 2 0
% $\tilde{\ell}_j^{+}$
\put(17.0000,-18.6000){\makebox(0,0)[lb]{$\tilde{\ell}_j^{+}$}}%
% STR 2 0 3 0
% 3 3120 2960 3120 3060 2 0
% $\tilde{\ell}_k^{+}$
\put(31.2000,-26.6000){\makebox(0,0)[lb]{$\tilde{\ell}_k^{+}$}}%
% STR 2 0 3 0
% 3 320 2970 320 3070 2 0
% $\tilde{\ell}_i^{+}$
\put(3.2000,-26.7000){\makebox(0,0)[lb]{$\tilde{\ell}_i^{+}$}}%
% STR 2 0 3 0
% 3 2510 810 2510 910 2 0
% $e_j$
\put(25.1000,-5.1000){\makebox(0,0)[lb]{$e_j$}}%
% STR 2 0 3 0
% 3 1080 800 1080 900 2 0
% $e_i$
\put(10.8000,-5.0000){\makebox(0,0)[lb]{$e_i$}}%
% STR 2 0 3 0
% 3 3930 810 3930 910 2 0
% $e_k$
\put(39.3000,-5.1000){\makebox(0,0)[lb]{$e_k$}}%
% CIRCLE 0 0 0 0
% 4 2600 2200 2623 2166 2623 2166 2623 2166
% 
\special{pn 20}%
\special{sh 0.600}%
\special{ar 2600 1800 41 41  0.0000000 6.2831853}%
% CIRCLE 0 0 0 0
% 4 2600 1400 2623 1366 2623 1366 2623 1366
% 
\special{pn 20}%
\special{sh 0.600}%
\special{ar 2600 1000 41 41  0.0000000 6.2831853}%
% CIRCLE 0 0 0 0
% 4 1200 1400 1223 1366 1223 1366 1223 1366
% 
\special{pn 20}%
\special{sh 0.600}%
\special{ar 1200 1000 41 41  0.0000000 6.2831853}%
% CIRCLE 0 0 0 0
% 4 4020 1400 4043 1366 4043 1366 4043 1366
% 
\special{pn 20}%
\special{sh 0.600}%
\special{ar 4020 1000 41 41  0.0000000 6.2831853}%
% CIRCLE 0 0 0 0
% 4 4020 2990 4043 2956 4043 2956 4043 2956
% 
\special{pn 20}%
\special{sh 0.600}%
\special{ar 4020 2590 41 41  0.0000000 6.2831853}%
% CIRCLE 0 0 0 0
% 4 1200 3000 1223 2966 1223 2966 1223 2966
% 
\special{pn 20}%
\special{sh 0.600}%
\special{ar 1200 2600 41 41  0.0000000 6.2831853}%
% STR 2 0 3 0
% 3 2690 1980 2690 2080 2 0
% $p_j$
\put(26.9000,-16.8000){\makebox(0,0)[lb]{$p_j$}}%
% STR 2 0 3 0
% 3 4110 2780 4110 2880 2 0
% $p_k$
\put(41.1000,-24.8000){\makebox(0,0)[lb]{$p_k$}}%
% STR 2 0 3 0
% 3 1300 2780 1300 2880 2 0
% $p_i$
\put(13.0000,-24.8000){\makebox(0,0)[lb]{$p_i$}}%
% STR 2 0 3 0
% 3 1310 1210 1310 1310 2 0
% $q_i$
\put(13.1000,-9.1000){\makebox(0,0)[lb]{$q_i$}}%
% STR 2 0 3 0
% 3 2710 1220 2710 1320 2 0
% $q_j$
\put(27.1000,-9.2000){\makebox(0,0)[lb]{$q_j$}}%
% STR 2 0 3 0
% 3 4120 1230 4120 1330 2 0
% $q_k$
\put(41.2000,-9.3000){\makebox(0,0)[lb]{$q_k$}}%
% STR 2 0 3 0
% 3 390 2100 390 2200 2 0
% $\wt{\Sol}(\th)$
\put(3.9000,-18.0000){\makebox(0,0)[lb]{$\wt{\Sol}(\th)$}}%
% STR 2 0 3 0
% 3 540 1370 540 1470 2 0
% $e_0$
\put(5.4000,-10.7000){\makebox(0,0)[lb]{$e_0$}}%
% STR 2 0 3 0
% 3 2500 3560 2500 3660 2 0
% $\mathbb{P}^1$
\put(25.0000,-32.6000){\makebox(0,0)[lb]{$\mathbb{P}^1$}}%
% STR 2 0 3 0
% 3 3300 2170 3300 2270 2 0
% $\mathbb{C}$
\put(33.0000,-18.7000){\makebox(0,0)[lb]{$\mathbb{C}$}}%
% STR 2 0 3 0
% 3 1870 2970 1870 3070 2 0
% $\mathbb{C}$
\put(18.7000,-26.7000){\makebox(0,0)[lb]{$\mathbb{C}$}}%
% STR 2 0 3 0
% 3 4710 2950 4710 3050 2 0
% $\mathbb{C}$
\put(47.1000,-26.5000){\makebox(0,0)[lb]{$\mathbb{C}$}}%
% STR 2 0 3 0
% 3 1110 3560 1110 3660 2 0
% $\mathbb{P}^1$
\put(11.1000,-32.6000){\makebox(0,0)[lb]{$\mathbb{P}^1$}}%
% STR 2 0 3 0
% 3 3940 3530 3940 3630 2 0
% $\mathbb{P}^1$
\put(39.4000,-32.3000){\makebox(0,0)[lb]{$\mathbb{P}^1$}}%
% STR 2 0 3 0
% 3 4530 1350 4530 1450 2 0
% $\mathbb{P}^1$
\put(45.3000,-10.5000){\makebox(0,0)[lb]{$\mathbb{P}^1$}}%
\end{picture}%
\end{center}
\caption{Surface of type $D_4$}
\label{fig:D4} 
\end{figure}
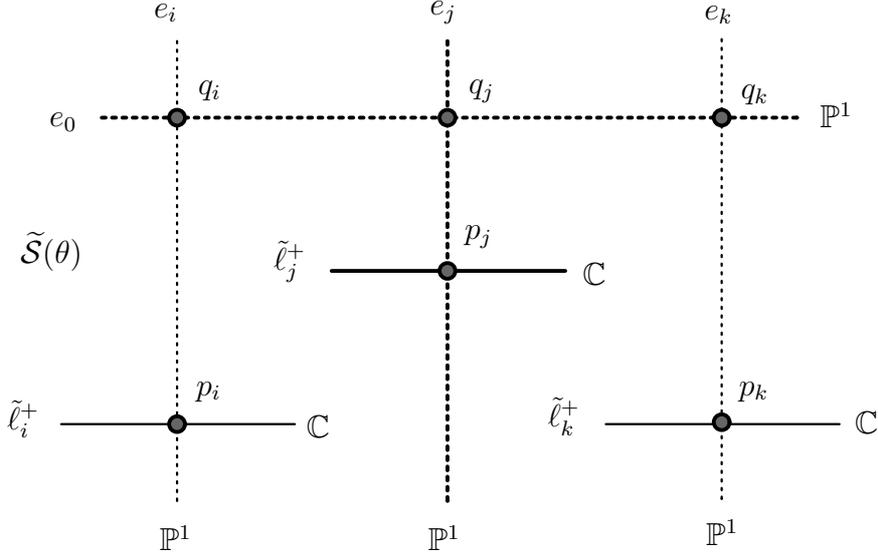
%%%%%%%%%%%%%%%%%%%%%%%%%%%%%%%%%%%%%%%%%%%%%%%%%%%%%%%%%%%%
%%%%%%%%%%%%%%%%%%%%%%%%%%% ex:D4 %%%%%%%%%%%%%%%%%%%%%%%%%%%
\begin{example}[$\mbox{\boldmath $D_4$}$]  \label{ex:D4} 
Consider the $\wt{W}(D_4^{(1)})$-stratum of type $D_4$, say, 
the $W(D_4^{(1)})$-stratum with value $\th = (8,8,8,28)$. 
In this case the surface $\Sol(\th)$ has only one 
singular point of type $D_4$ at $(x_i,x_j,x_k) = (2,2,2)$. 
The minimal resolution (\ref{eqn:brieskorn}) is obtained 
by successive blowing-ups: 
Blow up the singular point. 
If we express the blowing-up as 
$(x_i,x_j,x_k) = (u_iu_j+2, u_j+2, u_ku_j+2)$ in terms of 
coordinates $(u_i,u_j,u_k)$, then the strict transform of 
the surface $\Sol(\th)$ is represented as 
$u_iu_ju_k + (u_i+u_k+1)^2 = 0$. 
The exceptional curve $e$ is given by 
$u_j = u_i+u_k+1 = 0$. 
The strict transforms of $\ell_i^{+}$ and $\ell_j^{+}$ are 
given by $u_i+1= u_k = 0$ and $u_i = u_k+1 = 0$, while 
the strict transform of $\ell_k^{+}$ is at infinity 
and not expressible in terms of the coordinates 
$(u_i,u_j,u_k)$. 
The blow-up surface has three singularities, all of which 
are of type $A_1$ and located at the points in which the 
exceptional curve $e$ intersects the strict transforms of 
$\ell_i^{+}$, $\ell_j^{+}$, $\ell_k^{+}$. 
The lifts of the transformations $g_i^2$, $g_j^2$, $g_k^2$ 
fix the curve $e$ pointwise, since they fix the three 
singular points on it. 
Again blow up these points. 
Then we obtain a minimal resolution (\ref{eqn:brieskorn}) 
of the surface $\Sol(\th)$ as depicted in Figure \ref{fig:D4}, 
where $e_i$, $e_j$, $e_k$ are the exceptional curves over 
the singular points and $e_0$, $\wt{\ell}_i^+$, 
$\wt{\ell}_j^+$, $\wt{\ell}_k^+$ are the strict transforms 
of $e$, $\ell_i^+$, $\ell_j^+$, $\ell_k^+$, respectively. 
Being the strict transform of $e$, the exceptional curve 
$e_0$ is fixed pointwise by the lifts $\tilde{g}_i^2$, 
$\tilde{g}_j^2$, $\tilde{g}_k^2$ of $g_i^2$, $g_j^2$, 
$g_k^2$, and hence carries rational solutions. 
Moreover the lift $\tilde{g}_j^2$ fixes $e_j$ pointwise. 
This can be seen without computation. 
Since $\tilde{g}_j^2$ is area-preserving and fixes 
$\wt{\ell}_j^+$ pointwise, it has derivative $1$ at $p_j$ 
along the curve $e_j$. 
So the M\"obius transformation on $e_j$ induced by 
$\tilde{g}_j^2$ is either identity or a map of parabolic type. 
But the latter is impossible because it has at least 
two fixed points at $p_j$ and $q_j$ (see Figure \ref{fig:D4}). 
Hence $\wt{g}_j^2$ acts on $e_j$ as the identity. 
Next we shall observe that $\wt{g}_j^2$ acts on $e_i$ 
as a parabolic M\"obius transformation fixing $q_i$ only. 
If we express the blowing-up at $(u_i,u_j,u_k)=(-1,0,0)$ 
as $(u_i,u_j,u_k) = (v_iv_k-1,v_jv_k,v_k)$, then the 
exceptional curve $e_i$ is given by 
$v_k = v_j - (v_i+1)^2 = 0$. 
Parametrize $e_i$ as $(v_i,v_j,v_k) = (-(t+1)/t, t^{-2}, 0)$, 
where $q_i$ corresponds to $t = \infty$. 
Then $\wt{g}_j^2$ acts on $e_i$ by the shift $t \mapsto t+1$. 
Similarly $\wt{g}_j^2$ acts on $e_k$ as a parabolic 
transformation fixing $q_k$ only. 
By symmetry, $\wt{g}_i^2$ and $\wt{g}_k^2$ act on $e_j$ as 
parabolic transformations fixing $q_j$ only. 
Notice that the exceptional curve $e_0$ carries rational 
Riccati solutions, while $e_j-\{q_j\}$ carries Riccati 
solutions of infinte period. 
Thus we have 
%%%%%%%%%%%%%%%%%%%%%%%%% eqn:FixD4 %%%%%%%%%%%%%%%%%%%%%%%%
\begin{equation} \label{eqn:FixD4}
\wt{\mathrm{Fix}}_j(\th) = \wt{\ell}_j^{+} 
\underset{p_j}{\cup} e_j \underset{q_j}{\cup} e_0, 
\qquad \wt{\mathrm{Per}}_j^{e}(\th;n) = \emptyset 
\quad (n > 1). 
\end{equation}
%%%%%%%%%%%%%%%%%%%%%%%%%%%%%%%%%%%%%%%%%%%%%%%%%%%%%%%%%%%%%
\end{example}
%%%%%%%%%%%%%%%%%%%%%%%%%%%%%%%%%%%%%%%%%%%%%%%%%%%%%%%%%%%%%
%%%%%%%%%%%%%%%%%%%%%%%%% fig:A14 %%%%%%%%%%%%%%%%%%%%%%%%%%%
\begin{figure}[t]
\begin{center}
%WinTpicVersion2.15
\unitlength 0.1in
\begin{picture}(29.00,24.00)(1.80,-27.80)
% LINE 0 0 3 0
% 4 790 1420 810 1420 820 1420 2980 1420
% 
\special{pn 20}%
\special{pa 790 1020}%
\special{pa 810 1020}%
\special{fp}%
\special{pa 820 1020}%
\special{pa 2980 1020}%
\special{fp}%
% LINE 0 0 3 0
% 4 810 2740 830 2740 840 2740 3000 2740
% 
\special{pn 20}%
\special{pa 810 2340}%
\special{pa 830 2340}%
\special{fp}%
\special{pa 840 2340}%
\special{pa 3000 2340}%
\special{fp}%
% LINE 0 2 3 0
% 2 1200 1020 1200 1810
% 
\special{pn 20}%
\special{pa 1200 620}%
\special{pa 1200 1410}%
\special{dt 0.054}%
\special{pa 1200 1410}%
\special{pa 1200 1409}%
\special{dt 0.054}%
% LINE 0 2 3 0
% 2 2600 1030 2600 1820
% 
\special{pn 20}%
\special{pa 2600 630}%
\special{pa 2600 1420}%
\special{dt 0.054}%
\special{pa 2600 1420}%
\special{pa 2600 1419}%
\special{dt 0.054}%
% LINE 0 2 3 0
% 2 2610 2360 2610 3150
% 
\special{pn 20}%
\special{pa 2610 1960}%
\special{pa 2610 2750}%
\special{dt 0.054}%
\special{pa 2610 2750}%
\special{pa 2610 2749}%
\special{dt 0.054}%
% LINE 0 2 3 0
% 2 1210 2360 1210 3150
% 
\special{pn 20}%
\special{pa 1210 1960}%
\special{pa 1210 2750}%
\special{dt 0.054}%
\special{pa 1210 2750}%
\special{pa 1210 2749}%
\special{dt 0.054}%
% STR 2 0 3 0
% 3 390 1390 390 1490 2 0
% $\tilde{\ell}_j^+$
\put(3.9000,-10.9000){\makebox(0,0)[lb]{$\tilde{\ell}_j^+$}}%
% STR 2 0 3 0
% 3 410 2730 410 2830 2 0
% $\tilde{\ell}_j^-$
\put(4.1000,-24.3000){\makebox(0,0)[lb]{$\tilde{\ell}_j^-$}}%
% CIRCLE 0 0 0 0
% 4 1200 1420 1240 1430 1240 1410 1240 1410
% 
\special{pn 20}%
\special{sh 0.600}%
\special{ar 1200 1020 41 41  0.0000000 6.2831853}%
% CIRCLE 0 0 0 0
% 4 2600 1420 2640 1430 2640 1410 2640 1410
% 
\special{pn 20}%
\special{sh 0.600}%
\special{ar 2600 1020 41 41  0.0000000 6.2831853}%
% CIRCLE 0 0 0 0
% 4 2610 2740 2650 2750 2650 2730 2650 2730
% 
\special{pn 20}%
\special{sh 0.600}%
\special{ar 2610 2340 41 41  0.0000000 6.2831853}%
% CIRCLE 0 0 0 0
% 4 1210 2740 1250 2750 1250 2730 1250 2730
% 
\special{pn 20}%
\special{sh 0.600}%
\special{ar 1210 2340 41 41  0.0000000 6.2831853}%
% STR 2 0 3 0
% 3 2480 850 2480 950 2 0
% $e^{++-}$
\put(24.8000,-5.5000){\makebox(0,0)[lb]{$e^{++-}$}}%
% STR 2 0 3 0
% 3 1040 860 1040 960 2 0
% $e^{+-+}$
\put(10.4000,-5.6000){\makebox(0,0)[lb]{$e^{+-+}$}}%
% STR 2 0 3 0
% 3 1040 2210 1040 2310 2 0
% $e^{---}$
\put(10.4000,-19.1000){\makebox(0,0)[lb]{$e^{---}$}}%
% STR 2 0 3 0
% 3 2500 2200 2500 2300 2 0
% $e^{-++}$
\put(25.0000,-19.0000){\makebox(0,0)[lb]{$e^{-++}$}}%
% STR 2 0 3 0
% 3 2670 1250 2670 1350 2 0
% $p^{++-}$
\put(26.7000,-9.5000){\makebox(0,0)[lb]{$p^{++-}$}}%
% STR 2 0 3 0
% 3 2690 2560 2690 2660 2 0
% $p^{-++}$
\put(26.9000,-22.6000){\makebox(0,0)[lb]{$p^{-++}$}}%
% STR 2 0 3 0
% 3 1290 1250 1290 1350 2 0
% $p^{+-+}$
\put(12.9000,-9.5000){\makebox(0,0)[lb]{$p^{+-+}$}}%
% STR 2 0 3 0
% 3 1300 2560 1300 2660 2 0
% $p^{---}$
\put(13.0000,-22.6000){\makebox(0,0)[lb]{$p^{---}$}}%
% STR 2 0 3 0
% 3 180 2060 180 2160 2 0
% $\wt{\Sol}(\th)$
\put(1.8000,-17.6000){\makebox(0,0)[lb]{$\wt{\Sol}(\th)$}}%
% STR 2 0 3 0
% 3 2510 1920 2510 2020 2 0
% $\mathbb{P}^1$
\put(25.1000,-16.2000){\makebox(0,0)[lb]{$\mathbb{P}^1$}}%
% STR 2 0 3 0
% 3 1110 1910 1110 2010 2 0
% $\mathbb{P}^1$
\put(11.1000,-16.1000){\makebox(0,0)[lb]{$\mathbb{P}^1$}}%
% STR 2 0 3 0
% 3 1120 3250 1120 3350 2 0
% $\mathbb{P}^1$
\put(11.2000,-29.5000){\makebox(0,0)[lb]{$\mathbb{P}^1$}}%
% STR 2 0 3 0
% 3 2530 3250 2530 3350 2 0
% $\mathbb{P}^1$
\put(25.3000,-29.5000){\makebox(0,0)[lb]{$\mathbb{P}^1$}}%
% STR 2 0 3 0
% 3 3070 1380 3070 1480 2 0
% $\mathbb{C}$
\put(30.7000,-10.8000){\makebox(0,0)[lb]{$\mathbb{C}$}}%
% STR 2 0 3 0
% 3 3080 2710 3080 2810 2 0
% $\mathbb{C}$
\put(30.8000,-24.1000){\makebox(0,0)[lb]{$\mathbb{C}$}}%
\end{picture}%
\end{center}
\caption{Surface of type $A_1^{\oplus 4}$} 
\label{fig:A14}
\end{figure}
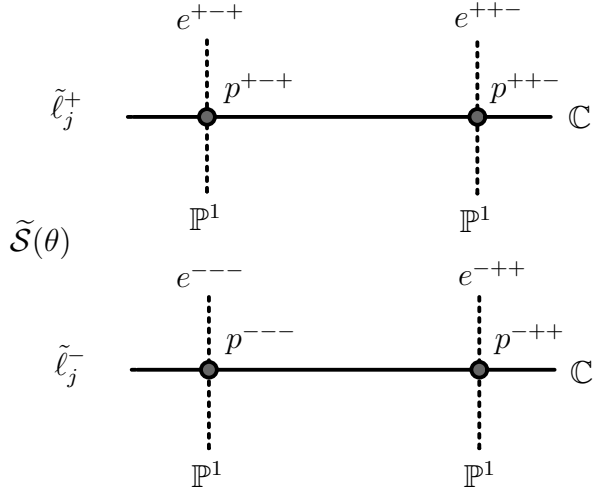
%%%%%%%%%%%%%%%%%%%%%%%%%%%%%%%%%%%%%%%%%%%%%%%%%%%%%%%%%%%%%
%%%%%%%%%%%%%%%%%%%%%%%% ex:A14 %%%%%%%%%%%%%%%%%%%%%%%%%%%%%
\begin{example}[$\mbox{\boldmath $A_1^{\oplus 4}$}$] 
\label{ex:A14} 
Consider the $\wt{W}(D_4^{(1)})$-stratum of type 
$A_1^{\oplus 4}$, where $\th = (0,0,0,-4)$ and 
$\mathrm{Fix}_j(\th) = \ell_j^{+} \amalg \ell_j^{-}$. 
In this case the surface $\Sol(\th)$ has four 
singularities of type $A_1$ at 
$(x_i,x_j,x_k) = (2\ve_i,2\ve_j,2\ve_k) \in \{\pm2\}^3$ with 
$\ve_i\ve_j\ve_k = -1$. 
Blow up at these points to obtain a minimal resolution as in 
(\ref{eqn:brieskorn}). 
Let $e^{\ve_i\ve_j\ve_k}$ be the exceptional line 
over $(x_i,x_j,x_k) = (2\ve_i,2\ve_j,2\ve_k)$ and 
$\tilde{\ell}_j^{\ve_i}$ be the strict transform 
of $\ell_j^{\ve_i}$. 
Moreover let $p^{\ve_i\ve_j\ve_k}$ denote the intersection 
point of the lines $e^{\ve_i\ve_j\ve_k}$ and 
$\wt{\ell}_j^{\ve_i}$ (see Figure \ref{fig:A14}). 
Then the lifted transformation $\tilde{g}_j^2 : 
\wt{\Sol}(\th) \carl$ acts on the exceptional line 
$e^{\ve_i\ve_j\ve_k} \cong \P^1$ as a M\"{o}bius 
transformation. 
It is a parabolic transformation with the only 
fixed point $p^{\ve_i\ve_j\ve_k}$. 
Let us check this for $(\ve_i,\ve_j,\ve_k) = (-1,-1,-1)$. 
The blowing-up of $\C^3$ at $(x_i,x_j,x_k)=(-2,-2,-2)$ is 
described by $x_i = u_i u_j-2$, $x_j = u_j-2$, 
$x_k = u_j u_k-2$, 
in terms of coordinates $(u_i,u_j,u_k)$ around $(0,0,0)$. 
Then the exceptional line $e^{---}$ is represented by the 
equations $u_j = 0$ and $(u_i-u_k)^2-2(u_i+u_k)+1 = 0$ 
and hence it is parametrized as 
%%%%%
\[
u_i = \left(\dfrac{2}{t+1}\right)^2, \qquad u_j = 0, 
\qquad u_k = \left(\dfrac{t-1}{t+1}\right)^2, 
\]
%%%%% 
where the fixed point $p^{---}$ corresponds to $t = \infty$. 
Then we can check that $\tilde{g}_j^2$ acts on the line 
$e^{---}$ as the translation $t \mapsto t+4$, as desired. 
Thus the only fixed points of $\tilde{g}_j^2$ on 
the exceptional set $\E(\th)$ are the four points 
$p^{\ve_i\ve_j\ve_k}$ with $\ve_i\ve_j\ve_k = -1$ and there 
are no periodic points, so that 
%%%%%%%%%%%%%%%%%%%%%%%%%% eqn:FixA14 %%%%%%%%%%%%%%%%%%%%%%%%
\begin{equation} \label{eqn:FixA14}
\wt{\mathrm{Fix}}_j(\th) = 
\wt{\ell}_j^+ \amalg \wt{\ell}_j^-, 
\qquad 
\wt{\mathrm{Per}}_j^{e}(\th;n) = \emptyset \quad (n > 1). 
\end{equation}
%%%%%%%%%%%%%%%%%%%%%%%%%%%%%%%%%%%%%%%%%%%%%%%%%%%%%%%%%%%%%%
\end{example} 
%%%%%%%%%%%%%%%%%%%%%%%%%%%%%%%%%%%%%%%%%%%%%%%%%%%%%%%%%%%%%%
%%%%%%%%%%%%%%%%%%%%%%%% sec:power %%%%%%%%%%%%%%%%%%%%%%%%%%%
\section{Power Geometry} \label{sec:power} 
%%%%%%%%%%%%%%%%%%%%%%%%%%%%%%%%%%%%%%%%%%%%%%%%%%%%%%%%%%%%%%
We apply the method of power geometry \cite{Bruno1,Bruno2,BG} 
to construct as many algebraic branch solutions to $\PVI(\k)$ 
as possible around each fixed singular point. 
Basically we can follow the arguments of \cite{BG}. 
However, while the attention of \cite{BG} is restricted to 
generic parameters, we require a thorough treatment of all 
parameters, where much ampler varieties of patterns 
are present. 
Moreover, the way in \cite{BG} of representing the parameters 
of Painlev\'e VI is not convenient for our purpose. 
So we have to redevelop the necessary arguments on power 
geometry from scratch. 
\par 
In view of Remark \ref{rem:backlund}, it is sufficient to 
work around the origin $z = 0$. 
In order to apply the method in \cite{Bruno1,Bruno2,BG}, 
we reduce the system (\ref{eqn:PVI}) into a single 
second-order equation. 
If $(q,p) = (q(z), p(z))$ is a solution to 
system (\ref{eqn:PVI}) such that $q \not\equiv 0$, $1$, $z$, 
$\infty$, then we solve the first equation of system 
(\ref{eqn:PVI}) with respect to $p = p(z)$ to obtain 
%%%%%%%%%%%%%%%%%%%%%%%%% eqn:p(z) %%%%%%%%%%%%%%%%%%%%%%%%%%%
\begin{equation} \label{eqn:p(z)}
p = \dfrac{z(z-1)q'+\k_1q_1q_z+(\k_2-1)q_0q_1+
\k_3 q_0 q_z}{2q_0q_1q_z}. 
\end{equation}
%%%%%%%%%%%%%%%%%%%%%%%%%%%%%%%%%%%%%%%%%%%%%%%%%%%%%%%%%%%%%%
Substituting this into the second equation yields the 
single second-order equation 
%%%%%%%%%%%%%%%%%%%%%%%% eqn:PVI2 %%%%%%%%%%%%%%%%%%%%%%%%%%%%
\begin{equation} \label{eqn:PVI2} 
\begin{array}{rcl}
\dfrac{d^2 q}{dz^2} 
&=& \dfrac{1}{2} 
\left(\dfrac{1}{q_0} + \dfrac{1}{q_1} + \dfrac{1}{q_z} 
\right) \left(\dfrac{dq}{dz}\right)^2 
- \left(\dfrac{1}{z}+\dfrac{1}{z-1}+\dfrac{1}{q_z} \right) 
\left(\dfrac{dq}{dz}\right) \\[6mm] 
&&+ \dfrac{q_0q_1q_z}{2z^2(z-1)^2} 
\left\{\k_4^2 - \k_1^2 \dfrac{z}{q_0^2} + 
\k_3^2 \dfrac{z-1}{q_1^2} 
+(1-\k_2^2) \dfrac{z(z-1)}{q_z^2}\right\}. 
\end{array}
\end{equation}
%%%%%%%%%%%%%%%%%%%%%%%%%%%%%%%%%%%%%%%%%%%%%%%%%%%%%%%%%%%%%%%
Multiply equation (\ref{eqn:PVI2}) by $2z^2(z-1)^2 q_0 q_1 q_z$ 
and move its right-hand side to the left to obtain 
%%%%%%%%%%%%%%%%%%%%%%%%% eqn:PVI3 %%%%%%%%%%%%%%%%%%%%%%%%%%%%
\begin{equation} \label{eqn:PVI3}
P(z,q) = 0, 
\end{equation}
%%%%%%%%%%%%%%%%%%%%%%%%%%%%%%%%%%%%%%%%%%%%%%%%%%%%%%%%%%%%%%%
where $P(z,q)$ is a polynomial of $(z,q,q',q'')$, that 
is, a differential sum of $(z,q)$, whose explicit formula is 
omitted here but can be found in \cite{BG}. 
Therefore system (\ref{eqn:PVI}) is equivalent to equation 
(\ref{eqn:PVI3}) together with (\ref{eqn:p(z)}) except 
for the possible solutions such that $q \equiv 0$, $1$, $z$, 
$\infty$. 
A simple check shows that the Newton polygon of equation 
(\ref{eqn:PVI3}) is given as in Figure \ref{fig:snewtonPVI}, 
where there are four patterns according as the parameters 
$\k_1$ and $\k_4$ are zero or not. 
%%%%%%%%%%%%%%%%%%%%%%% fig:snewtonPVI %%%%%%%%%%%%%%%%%%%%%%%
\begin{figure}[p]
\begin{center}
\begin{tabular}{cc}
%WinTpicVersion2.15
\unitlength 0.1in
\begin{picture}(23.20,32.40)(5.00,-36.40)
% CIRCLE 0 0 0 0
% 4 1000 1200 1030 1230 1020 1200 1040 1200
% 
\special{pn 20}%
\special{sh 0.600}%
\special{ar 1000 800 42 42  0.0000000 6.2831853}%
% CIRCLE 0 0 0 0
% 4 1010 1610 1040 1640 1030 1610 1050 1610
% 
\special{pn 20}%
\special{sh 0.600}%
\special{ar 1010 1210 42 42  0.0000000 6.2831853}%
% CIRCLE 0 0 0 0
% 4 1010 2000 1040 2030 1030 2000 1050 2000
% 
\special{pn 20}%
\special{sh 0.600}%
\special{ar 1010 1600 42 42  0.0000000 6.2831853}%
% CIRCLE 0 0 0 0
% 4 1010 2400 1040 2430 1030 2400 1050 2400
% 
\special{pn 20}%
\special{sh 0.600}%
\special{ar 1010 2000 42 42  0.0000000 6.2831853}%
% CIRCLE 0 0 0 0
% 4 1400 1610 1430 1640 1420 1610 1440 1610
% 
\special{pn 20}%
\special{sh 0.600}%
\special{ar 1400 1210 42 42  0.0000000 6.2831853}%
% CIRCLE 0 0 3 0
% 4 1390 2000 1420 2030 1410 2000 1430 2000
% 
\special{pn 20}%
\special{ar 1390 1600 42 42  0.0000000 6.2831853}%
% CIRCLE 0 0 3 0
% 4 1390 2400 1420 2430 1410 2400 1430 2400
% 
\special{pn 20}%
\special{ar 1390 2000 42 42  0.0000000 6.2831853}%
% CIRCLE 0 0 0 0
% 4 1390 2790 1420 2820 1410 2790 1430 2790
% 
\special{pn 20}%
\special{sh 0.600}%
\special{ar 1390 2390 42 42  0.0000000 6.2831853}%
% CIRCLE 0 0 0 0
% 4 1810 3200 1840 3230 1830 3200 1850 3200
% 
\special{pn 20}%
\special{sh 0.600}%
\special{ar 1810 2800 42 42  0.0000000 6.2831853}%
% CIRCLE 0 0 3 0
% 4 1820 2790 1850 2820 1840 2790 1860 2790
% 
\special{pn 20}%
\special{ar 1820 2390 42 42  0.0000000 6.2831853}%
% CIRCLE 0 0 3 0
% 4 1810 2400 1840 2430 1830 2400 1850 2400
% 
\special{pn 20}%
\special{ar 1810 2000 42 42  0.0000000 6.2831853}%
% CIRCLE 0 0 0 0
% 4 1800 2000 1830 2030 1820 2000 1840 2000
% 
\special{pn 20}%
\special{sh 0.600}%
\special{ar 1800 1600 42 42  0.0000000 6.2831853}%
% CIRCLE 0 0 0 0
% 4 2200 3600 2230 3630 2220 3600 2240 3600
% 
\special{pn 20}%
\special{sh 0.600}%
\special{ar 2200 3200 42 42  0.0000000 6.2831853}%
% CIRCLE 0 0 0 0
% 4 2200 3200 2230 3230 2220 3200 2240 3200
% 
\special{pn 20}%
\special{sh 0.600}%
\special{ar 2200 2800 42 42  0.0000000 6.2831853}%
% CIRCLE 0 0 0 0
% 4 2200 2800 2230 2830 2220 2800 2240 2800
% 
\special{pn 20}%
\special{sh 0.600}%
\special{ar 2200 2400 42 42  0.0000000 6.2831853}%
% CIRCLE 0 0 0 0
% 4 2200 2400 2230 2430 2220 2400 2240 2400
% 
\special{pn 20}%
\special{sh 0.600}%
\special{ar 2200 2000 42 42  0.0000000 6.2831853}%
% LINE 1 0 3 0
% 2 1000 1210 2200 2400
% 
\special{pn 13}%
\special{pa 1000 810}%
\special{pa 2200 2000}%
\special{fp}%
% LINE 1 0 3 0
% 2 2200 2400 2210 3610
% 
\special{pn 13}%
\special{pa 2200 2000}%
\special{pa 2210 3210}%
\special{fp}%
% LINE 1 0 3 0
% 2 1010 1210 1010 2400
% 
\special{pn 13}%
\special{pa 1010 810}%
\special{pa 1010 2000}%
\special{fp}%
% LINE 1 0 3 0
% 2 1010 2400 2210 3620
% 
\special{pn 13}%
\special{pa 1010 2000}%
\special{pa 2210 3220}%
\special{fp}%
% LINE 2 0 3 0
% 2 610 3600 2820 3600
% 
\special{pn 8}%
\special{pa 610 3200}%
\special{pa 2820 3200}%
\special{fp}%
% LINE 2 0 3 0
% 2 1000 800 1000 3860
% 
\special{pn 8}%
\special{pa 1000 400}%
\special{pa 1000 3460}%
\special{fp}%
% STR 2 0 3 0
% 3 770 3720 770 3820 2 0
% $0$
\put(7.7000,-34.2000){\makebox(0,0)[lb]{$0$}}%
% STR 2 0 3 0
% 3 2020 3780 2020 3880 2 0
% $(3,0)$
\put(20.2000,-34.8000){\makebox(0,0)[lb]{$(3,0)$}}%
% STR 2 0 3 0
% 3 510 2360 510 2460 2 0
% $(0,3)$
\put(5.1000,-20.6000){\makebox(0,0)[lb]{$(0,3)$}}%
% STR 2 0 3 0
% 3 500 1160 500 1260 2 0
% $(0,6)$
\put(5.0000,-8.6000){\makebox(0,0)[lb]{$(0,6)$}}%
% STR 2 0 3 0
% 3 2280 2210 2280 2310 2 0
% $(3,3)$
\put(22.8000,-19.1000){\makebox(0,0)[lb]{$(3,3)$}}%
% STR 2 0 3 0
% 3 690 1810 690 1910 2 0
% $\varGamma_0$
\put(6.9000,-15.1000){\makebox(0,0)[lb]{$\varGamma_0$}}%
% STR 2 0 3 0
% 3 1380 3110 1380 3210 2 0
% $\varGamma_1$
\put(13.8000,-28.1000){\makebox(0,0)[lb]{$\varGamma_1$}}%
% STR 2 0 3 0
% 3 1000 4110 1000 4210 2 0
% Case $\k_1 \neq 0$, $\k_4 \neq 0$.
\put(10.0000,-38.1000){\makebox(0,0)[lb]{Case $\k_1 \neq 0$, $\k_4 \neq 0$.}}%
\end{picture}%
\qquad\qquad  & 
%WinTpicVersion2.15
\unitlength 0.1in
\begin{picture}(23.20,32.40)(5.00,-36.40)
% CIRCLE 0 0 0 0
% 4 1000 1200 1030 1230 1020 1200 1040 1200
% 
\special{pn 20}%
\special{sh 0.600}%
\special{ar 1000 800 42 42  0.0000000 6.2831853}%
% CIRCLE 0 0 0 0
% 4 1010 1610 1040 1640 1030 1610 1050 1610
% 
\special{pn 20}%
\special{sh 0.600}%
\special{ar 1010 1210 42 42  0.0000000 6.2831853}%
% CIRCLE 0 0 0 0
% 4 1010 2000 1040 2030 1030 2000 1050 2000
% 
\special{pn 20}%
\special{sh 0.600}%
\special{ar 1010 1600 42 42  0.0000000 6.2831853}%
% CIRCLE 0 0 0 0
% 4 1010 2400 1040 2430 1030 2400 1050 2400
% 
\special{pn 20}%
\special{sh 0.600}%
\special{ar 1010 2000 42 42  0.0000000 6.2831853}%
% CIRCLE 0 0 0 0
% 4 1400 1610 1430 1640 1420 1610 1440 1610
% 
\special{pn 20}%
\special{sh 0.600}%
\special{ar 1400 1210 42 42  0.0000000 6.2831853}%
% CIRCLE 0 0 3 0
% 4 1390 2000 1420 2030 1410 2000 1430 2000
% 
\special{pn 20}%
\special{ar 1390 1600 42 42  0.0000000 6.2831853}%
% CIRCLE 0 0 3 0
% 4 1390 2400 1420 2430 1410 2400 1430 2400
% 
\special{pn 20}%
\special{ar 1390 2000 42 42  0.0000000 6.2831853}%
% CIRCLE 0 0 0 0
% 4 1390 2790 1420 2820 1410 2790 1430 2790
% 
\special{pn 20}%
\special{sh 0.600}%
\special{ar 1390 2390 42 42  0.0000000 6.2831853}%
% CIRCLE 0 0 0 0
% 4 1820 2790 1850 2820 1840 2790 1860 2790
% 
\special{pn 20}%
\special{sh 0.600}%
\special{ar 1820 2390 42 42  0.0000000 6.2831853}%
% CIRCLE 0 0 3 0
% 4 1810 2400 1840 2430 1830 2400 1850 2400
% 
\special{pn 20}%
\special{ar 1810 2000 42 42  0.0000000 6.2831853}%
% CIRCLE 0 0 0 0
% 4 1800 2000 1830 2030 1820 2000 1840 2000
% 
\special{pn 20}%
\special{sh 0.600}%
\special{ar 1800 1600 42 42  0.0000000 6.2831853}%
% CIRCLE 0 0 0 0
% 4 2200 2800 2230 2830 2220 2800 2240 2800
% 
\special{pn 20}%
\special{sh 0.600}%
\special{ar 2200 2400 42 42  0.0000000 6.2831853}%
% CIRCLE 0 0 0 0
% 4 2200 2400 2230 2430 2220 2400 2240 2400
% 
\special{pn 20}%
\special{sh 0.600}%
\special{ar 2200 2000 42 42  0.0000000 6.2831853}%
% LINE 1 0 3 0
% 2 1000 1210 2200 2400
% 
\special{pn 13}%
\special{pa 1000 810}%
\special{pa 2200 2000}%
\special{fp}%
% LINE 1 0 3 0
% 2 2200 2400 2200 2790
% 
\special{pn 13}%
\special{pa 2200 2000}%
\special{pa 2200 2390}%
\special{fp}%
% LINE 1 0 3 0
% 2 1010 1210 1010 2400
% 
\special{pn 13}%
\special{pa 1010 810}%
\special{pa 1010 2000}%
\special{fp}%
% LINE 1 0 3 0
% 2 1010 2400 1390 2790
% 
\special{pn 13}%
\special{pa 1010 2000}%
\special{pa 1390 2390}%
\special{fp}%
% LINE 2 0 3 0
% 2 610 3600 2820 3600
% 
\special{pn 8}%
\special{pa 610 3200}%
\special{pa 2820 3200}%
\special{fp}%
% LINE 2 0 3 0
% 2 1000 800 1000 3860
% 
\special{pn 8}%
\special{pa 1000 400}%
\special{pa 1000 3460}%
\special{fp}%
% STR 2 0 3 0
% 3 770 3720 770 3820 2 0
% $0$
\put(7.7000,-34.2000){\makebox(0,0)[lb]{$0$}}%
% STR 2 0 3 0
% 3 510 2360 510 2460 2 0
% $(0,3)$
\put(5.1000,-20.6000){\makebox(0,0)[lb]{$(0,3)$}}%
% STR 2 0 3 0
% 3 500 1160 500 1260 2 0
% $(0,6)$
\put(5.0000,-8.6000){\makebox(0,0)[lb]{$(0,6)$}}%
% STR 2 0 3 0
% 3 2280 2210 2280 2310 2 0
% $(3,3)$
\put(22.8000,-19.1000){\makebox(0,0)[lb]{$(3,3)$}}%
% STR 2 0 3 0
% 3 690 1810 690 1910 2 0
% $\varGamma_0$
\put(6.9000,-15.1000){\makebox(0,0)[lb]{$\varGamma_0$}}%
% STR 2 0 3 0
% 3 1070 2680 1070 2780 2 0
% $\varGamma_1$
\put(10.7000,-23.8000){\makebox(0,0)[lb]{$\varGamma_1$}}%
% LINE 0 0 3 0
% 2 1390 2800 2200 2790
% 
\special{pn 20}%
\special{pa 1390 2400}%
\special{pa 2200 2390}%
\special{fp}%
% STR 2 0 3 0
% 3 2280 2940 2280 3040 2 0
% $(3,2)$
\put(22.8000,-26.4000){\makebox(0,0)[lb]{$(3,2)$}}%
% STR 2 0 3 0
% 3 1220 2970 1220 3070 2 0
% $(1,2)$
\put(12.2000,-26.7000){\makebox(0,0)[lb]{$(1,2)$}}%
% STR 2 0 3 0
% 3 1000 4110 1000 4210 2 0
% Case $\k_1 = 0$, $\k_4 \neq 0$.
\put(10.0000,-38.1000){\makebox(0,0)[lb]{Case $\k_1 = 0$, $\k_4 \neq 0$.}}%
\end{picture}%
\\[1.5cm]
%WinTpicVersion2.15
\unitlength 0.1in
\begin{picture}(22.70,32.50)(5.50,-36.50)
% CIRCLE 0 0 0 0
% 4 1010 2000 1040 2030 1030 2000 1050 2000
% 
\special{pn 20}%
\special{sh 0.600}%
\special{ar 1010 1600 42 42  0.0000000 6.2831853}%
% CIRCLE 0 0 0 0
% 4 1010 2400 1040 2430 1030 2400 1050 2400
% 
\special{pn 20}%
\special{sh 0.600}%
\special{ar 1010 2000 42 42  0.0000000 6.2831853}%
% CIRCLE 0 0 0 0
% 4 1390 2000 1420 2030 1410 2000 1430 2000
% 
\special{pn 20}%
\special{sh 0.600}%
\special{ar 1390 1600 42 42  0.0000000 6.2831853}%
% CIRCLE 0 0 3 0
% 4 1390 2400 1420 2430 1410 2400 1430 2400
% 
\special{pn 20}%
\special{ar 1390 2000 42 42  0.0000000 6.2831853}%
% CIRCLE 0 0 0 0
% 4 1390 2790 1420 2820 1410 2790 1430 2790
% 
\special{pn 20}%
\special{sh 0.600}%
\special{ar 1390 2390 42 42  0.0000000 6.2831853}%
% CIRCLE 0 0 0 0
% 4 1810 3200 1840 3230 1830 3200 1850 3200
% 
\special{pn 20}%
\special{sh 0.600}%
\special{ar 1810 2800 42 42  0.0000000 6.2831853}%
% CIRCLE 0 0 3 0
% 4 1820 2790 1850 2820 1840 2790 1860 2790
% 
\special{pn 20}%
\special{ar 1820 2390 42 42  0.0000000 6.2831853}%
% CIRCLE 0 0 3 0
% 4 1810 2400 1840 2430 1830 2400 1850 2400
% 
\special{pn 20}%
\special{ar 1810 2000 42 42  0.0000000 6.2831853}%
% CIRCLE 0 0 0 0
% 4 1800 2000 1830 2030 1820 2000 1840 2000
% 
\special{pn 20}%
\special{sh 0.600}%
\special{ar 1800 1600 42 42  0.0000000 6.2831853}%
% CIRCLE 0 0 0 0
% 4 2200 3600 2230 3630 2220 3600 2240 3600
% 
\special{pn 20}%
\special{sh 0.600}%
\special{ar 2200 3200 42 42  0.0000000 6.2831853}%
% CIRCLE 0 0 0 0
% 4 2200 3200 2230 3230 2220 3200 2240 3200
% 
\special{pn 20}%
\special{sh 0.600}%
\special{ar 2200 2800 42 42  0.0000000 6.2831853}%
% CIRCLE 0 0 0 0
% 4 2200 2800 2230 2830 2220 2800 2240 2800
% 
\special{pn 20}%
\special{sh 0.600}%
\special{ar 2200 2400 42 42  0.0000000 6.2831853}%
% CIRCLE 0 0 0 0
% 4 2200 2400 2230 2430 2220 2400 2240 2400
% 
\special{pn 20}%
\special{sh 0.600}%
\special{ar 2200 2000 42 42  0.0000000 6.2831853}%
% LINE 1 0 3 0
% 2 2200 2400 2210 3610
% 
\special{pn 13}%
\special{pa 2200 2000}%
\special{pa 2210 3210}%
\special{fp}%
% LINE 1 0 3 0
% 2 1010 2400 2210 3620
% 
\special{pn 13}%
\special{pa 1010 2000}%
\special{pa 2210 3220}%
\special{fp}%
% LINE 2 0 3 0
% 2 610 3600 2820 3600
% 
\special{pn 8}%
\special{pa 610 3200}%
\special{pa 2820 3200}%
\special{fp}%
% LINE 2 0 3 0
% 2 1000 800 1000 3860
% 
\special{pn 8}%
\special{pa 1000 400}%
\special{pa 1000 3460}%
\special{fp}%
% STR 2 0 3 0
% 3 770 3720 770 3820 2 0
% $0$
\put(7.7000,-34.2000){\makebox(0,0)[lb]{$0$}}%
% STR 2 0 3 0
% 3 2020 3780 2020 3880 2 0
% $(3,0)$
\put(20.2000,-34.8000){\makebox(0,0)[lb]{$(3,0)$}}%
% STR 2 0 3 0
% 3 550 2510 550 2610 2 0
% $(0,3)$
\put(5.5000,-22.1000){\makebox(0,0)[lb]{$(0,3)$}}%
% STR 2 0 3 0
% 3 2280 2210 2280 2310 2 0
% $(3,3)$
\put(22.8000,-19.1000){\makebox(0,0)[lb]{$(3,3)$}}%
% STR 2 0 3 0
% 3 670 2200 670 2300 2 0
% $\varGamma_0$
\put(6.7000,-19.0000){\makebox(0,0)[lb]{$\varGamma_0$}}%
% STR 2 0 3 0
% 3 1380 3110 1380 3210 2 0
% $\varGamma_1$
\put(13.8000,-28.1000){\makebox(0,0)[lb]{$\varGamma_1$}}%
% LINE 0 0 3 0
% 2 1000 2000 1000 2400
% 
\special{pn 20}%
\special{pa 1000 1600}%
\special{pa 1000 2000}%
\special{fp}%
% LINE 0 0 3 0
% 2 1000 2000 1810 1990
% 
\special{pn 20}%
\special{pa 1000 1600}%
\special{pa 1810 1590}%
\special{fp}%
% LINE 0 0 3 0
% 2 1810 2000 2210 2400
% 
\special{pn 20}%
\special{pa 1810 1600}%
\special{pa 2210 2000}%
\special{fp}%
% STR 2 0 3 0
% 3 560 1800 560 1900 2 0
% $(0,4)$
\put(5.6000,-15.0000){\makebox(0,0)[lb]{$(0,4)$}}%
% STR 2 0 3 0
% 3 1780 1780 1780 1880 2 0
% $(2,4)$
\put(17.8000,-14.8000){\makebox(0,0)[lb]{$(2,4)$}}%
% STR 2 0 3 0
% 3 1000 4120 1000 4220 2 0
% Case $\k_1 \neq 0$, $\k_4 = 0$.
\put(10.0000,-38.2000){\makebox(0,0)[lb]{Case $\k_1 \neq 0$, $\k_4 = 0$.}}%
\end{picture}%
\qquad\qquad  & 
%WinTpicVersion2.15
\unitlength 0.1in
\begin{picture}(22.60,32.50)(5.60,-36.50)
% CIRCLE 0 0 0 0
% 4 1010 2000 1040 2030 1030 2000 1050 2000
% 
\special{pn 20}%
\special{sh 0.600}%
\special{ar 1010 1600 42 42  0.0000000 6.2831853}%
% CIRCLE 0 0 0 0
% 4 1010 2400 1040 2430 1030 2400 1050 2400
% 
\special{pn 20}%
\special{sh 0.600}%
\special{ar 1010 2000 42 42  0.0000000 6.2831853}%
% CIRCLE 0 0 0 0
% 4 1390 2000 1420 2030 1410 2000 1430 2000
% 
\special{pn 20}%
\special{sh 0.600}%
\special{ar 1390 1600 42 42  0.0000000 6.2831853}%
% CIRCLE 0 0 3 0
% 4 1390 2400 1420 2430 1410 2400 1430 2400
% 
\special{pn 20}%
\special{ar 1390 2000 42 42  0.0000000 6.2831853}%
% CIRCLE 0 0 0 0
% 4 1390 2790 1420 2820 1410 2790 1430 2790
% 
\special{pn 20}%
\special{sh 0.600}%
\special{ar 1390 2390 42 42  0.0000000 6.2831853}%
% CIRCLE 0 0 0 0
% 4 1820 2790 1850 2820 1840 2790 1860 2790
% 
\special{pn 20}%
\special{sh 0.600}%
\special{ar 1820 2390 42 42  0.0000000 6.2831853}%
% CIRCLE 0 0 3 0
% 4 1810 2400 1840 2430 1830 2400 1850 2400
% 
\special{pn 20}%
\special{ar 1810 2000 42 42  0.0000000 6.2831853}%
% CIRCLE 0 0 0 0
% 4 1800 2000 1830 2030 1820 2000 1840 2000
% 
\special{pn 20}%
\special{sh 0.600}%
\special{ar 1800 1600 42 42  0.0000000 6.2831853}%
% CIRCLE 0 0 0 0
% 4 2200 2800 2230 2830 2220 2800 2240 2800
% 
\special{pn 20}%
\special{sh 0.600}%
\special{ar 2200 2400 42 42  0.0000000 6.2831853}%
% CIRCLE 0 0 0 0
% 4 2200 2400 2230 2430 2220 2400 2240 2400
% 
\special{pn 20}%
\special{sh 0.600}%
\special{ar 2200 2000 42 42  0.0000000 6.2831853}%
% LINE 0 0 3 0
% 2 2200 2400 2200 2800
% 
\special{pn 20}%
\special{pa 2200 2000}%
\special{pa 2200 2400}%
\special{fp}%
% LINE 0 0 3 0
% 2 1010 2400 1410 2800
% 
\special{pn 20}%
\special{pa 1010 2000}%
\special{pa 1410 2400}%
\special{fp}%
% LINE 2 0 3 0
% 2 610 3600 2820 3600
% 
\special{pn 8}%
\special{pa 610 3200}%
\special{pa 2820 3200}%
\special{fp}%
% LINE 2 0 3 0
% 2 1000 800 1000 3860
% 
\special{pn 8}%
\special{pa 1000 400}%
\special{pa 1000 3460}%
\special{fp}%
% STR 2 0 3 0
% 3 770 3720 770 3820 2 0
% $0$
\put(7.7000,-34.2000){\makebox(0,0)[lb]{$0$}}%
% STR 2 0 3 0
% 3 560 2530 560 2630 2 0
% $(0,3)$
\put(5.6000,-22.3000){\makebox(0,0)[lb]{$(0,3)$}}%
% STR 2 0 3 0
% 3 2280 2210 2280 2310 2 0
% $(3,3)$
\put(22.8000,-19.1000){\makebox(0,0)[lb]{$(3,3)$}}%
% STR 2 0 3 0
% 3 730 2200 730 2300 2 0
% $\varGamma_0$
\put(7.3000,-19.0000){\makebox(0,0)[lb]{$\varGamma_0$}}%
% STR 2 0 3 0
% 3 1070 2680 1070 2780 2 0
% $\varGamma_1$
\put(10.7000,-23.8000){\makebox(0,0)[lb]{$\varGamma_1$}}%
% LINE 0 0 3 0
% 2 1000 2000 1000 2400
% 
\special{pn 20}%
\special{pa 1000 1600}%
\special{pa 1000 2000}%
\special{fp}%
% LINE 0 0 3 0
% 2 1000 2000 1810 1990
% 
\special{pn 20}%
\special{pa 1000 1600}%
\special{pa 1810 1590}%
\special{fp}%
% LINE 0 0 3 0
% 2 1810 2000 2210 2400
% 
\special{pn 20}%
\special{pa 1810 1600}%
\special{pa 2210 2000}%
\special{fp}%
% STR 2 0 3 0
% 3 560 1800 560 1900 2 0
% $(0,4)$
\put(5.6000,-15.0000){\makebox(0,0)[lb]{$(0,4)$}}%
% STR 2 0 3 0
% 3 1780 1780 1780 1880 2 0
% $(2,4)$
\put(17.8000,-14.8000){\makebox(0,0)[lb]{$(2,4)$}}%
% LINE 0 0 3 0
% 2 1390 2800 2200 2800
% 
\special{pn 20}%
\special{pa 1390 2400}%
\special{pa 2200 2400}%
\special{fp}%
% STR 2 0 3 0
% 3 2280 2980 2280 3080 2 0
% $(3,2)$
\put(22.8000,-26.8000){\makebox(0,0)[lb]{$(3,2)$}}%
% STR 2 0 3 0
% 3 1200 2980 1200 3080 2 0
% $(1,2)$
\put(12.0000,-26.8000){\makebox(0,0)[lb]{$(1,2)$}}%
% STR 2 0 3 0
% 3 1010 4120 1010 4220 2 0
% Case $\k_1 = 0$, $\k_4 = 0$.
\put(10.1000,-38.2000){\makebox(0,0)[lb]{Case $\k_1 = 0$, $\k_4 = 0$.}}%
\end{picture}%
\end{tabular}
\end{center}
\vspace{1.5cm}
\caption{Newton polygon for Painlev\'e VI} 
\label{fig:snewtonPVI}
\end{figure}
%%%%%%%%%%%%%%%%%%%%%%%%%%%%%%%%%%%%%%%%%%%%%%%%%%%%%%%%%%%%%%%
\par 
First we search for a {\sl holomorphic} solution germ 
$q = q(z)$ to equation (\ref{eqn:PVI3}) around $z = 0$. 
We have only to construct formal power series solutions of 
the form
%%%%%%%%%%%%%%%%%%%%%%%%%% eqn:asymptotic %%%%%%%%%%%%%%%%%%%%%
\begin{equation} \label{eqn:asymptotic}
q = c z^r + (\mbox{higher order terms}), \qquad 
(r,c) \in \Z \times \C^{\times}, 
\end{equation}
%%%%%%%%%%%%%%%%%%%%%%%%%%%%%%%%%%%%%%%%%%%%%%%%%%%%%%%%%%%%%%%
since any formal power series solution to equation 
(\ref{eqn:PVI3}) is convergent \cite{Gerard,GS}. 
Then it follows from (\ref{eqn:p(z)}) that the associated 
formal Laurent series for $p = p(z)$ is also convergent. 
\par 
In order to construct formal solutions (\ref{eqn:asymptotic}), 
we consider the truncations along the edges $\varGamma_1$ and 
$\varGamma_0$ of the Newton polygons in Figure 
\ref{fig:snewtonPVI}. 
We see that $\varGamma_1$ and $\varGamma_0$ have outer normal 
vectors $(p_1,p_2) = (-1,-1)$ and $(p_1,p_2) = (-1,0)$, whose 
slopes are $p_2/p_1 = 1$ and $p_2/p_1 = 0$ respectively. 
Thus the edges $\varGamma_1$ and $\varGamma_0$ correspond to 
the exponents $r = 1$ and $r = 0$ respectively. 
The truncation of $P = P(z,q)$ along the edge 
$\varGamma_1$ is given by 
%%%%%%%
\[
\begin{array}{rcl}
P_1 &=& - 2 z q^2 q' + 2 z^2 q (q')^2 - 2 z^2 q^2 q'' 
+ (\k_1^2 - \k_2^2 + 1) z q^2 
- z^3 (q')^2 + 2 z^3 q q'' \\[3mm]
&& - 2 \k_1^2 z^2 q + \k_1^2 z^3. 
\end{array}
\]
%%%%%%%%
Substituting $q = cz$ into equation $P_1 = 0$ yields 
$c = \k_1/(\k_1 + \ve \k_2)$ with any sign 
$\ve \in \{\pm1\}$. 
Similarly, the truncation of $P = P(z,q)$ along the edge 
$\varGamma_0$ is given by 
%%%%%
\[
\begin{array}{rcl}
P_0 &=& - 2 z q^2 q' + 2 z^2 q (q')^2 - 2 z^2 q^2 q'' 
+ (\k_3^2 - \k_4^2) q^4 + 2 z q^3 q' 
- 3 z^2 q^2 (q')^2  \\[3mm]
&& + 2 z^2 q^3 q'' + 2 \k_4^2 q^5 - \k_4^2 q^6. 
\end{array} 
\]
%%%%
Substituting $q = c$ into equation $P_0 = 0$ yields 
$c = (\k_4+\ve\k_3)/\k_4$ with any sign $\ve \in \{\pm1\}$.  
%%%%%%%%%%%%%%%%%%%%% fig:snewtonk1k2k3k4 %%%%%%%%%%%%%%%%%%%%%%%%
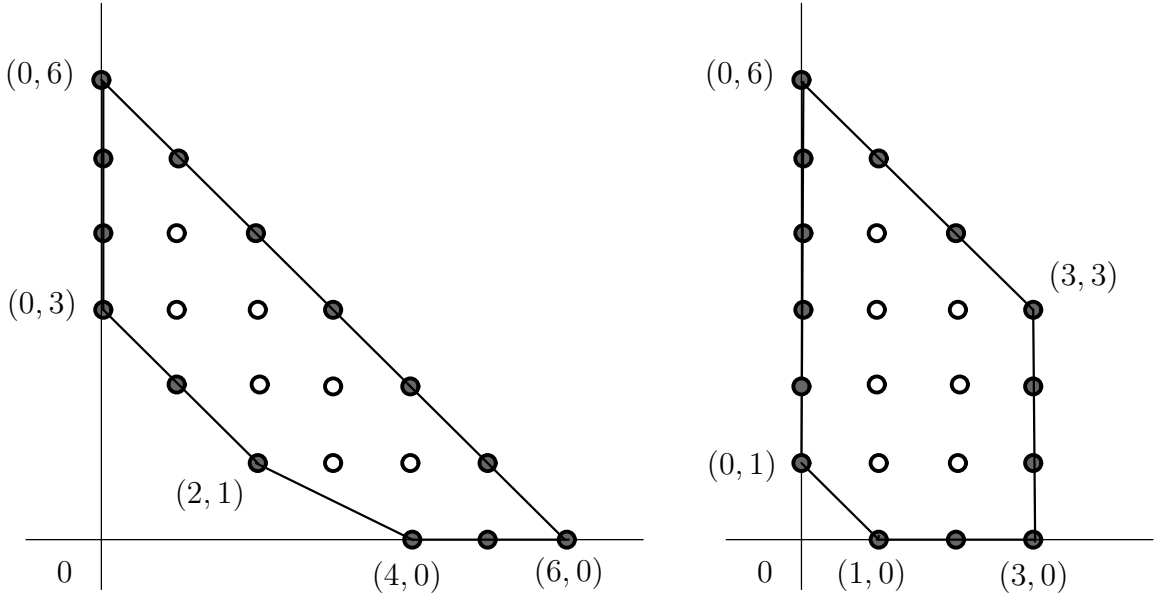
\begin{figure}[t] 
\begin{center}
%WinTpicVersion2.15
\unitlength 0.1in
\begin{picture}(33.00,30.60)(5.00,-34.60)
% CIRCLE 0 0 0 0
% 4 1000 1200 1030 1230 1020 1200 1040 1200
% 
\special{pn 20}%
\special{sh 0.600}%
\special{ar 1000 800 42 42  0.0000000 6.2831853}%
% CIRCLE 0 0 0 0
% 4 1010 1610 1040 1640 1030 1610 1050 1610
% 
\special{pn 20}%
\special{sh 0.600}%
\special{ar 1010 1210 42 42  0.0000000 6.2831853}%
% CIRCLE 0 0 0 0
% 4 1010 2000 1040 2030 1030 2000 1050 2000
% 
\special{pn 20}%
\special{sh 0.600}%
\special{ar 1010 1600 42 42  0.0000000 6.2831853}%
% CIRCLE 0 0 0 0
% 4 1010 2400 1040 2430 1030 2400 1050 2400
% 
\special{pn 20}%
\special{sh 0.600}%
\special{ar 1010 2000 42 42  0.0000000 6.2831853}%
% CIRCLE 0 0 0 0
% 4 1400 1610 1430 1640 1420 1610 1440 1610
% 
\special{pn 20}%
\special{sh 0.600}%
\special{ar 1400 1210 42 42  0.0000000 6.2831853}%
% CIRCLE 0 0 3 0
% 4 1390 2000 1420 2030 1410 2000 1430 2000
% 
\special{pn 20}%
\special{ar 1390 1600 42 42  0.0000000 6.2831853}%
% CIRCLE 0 0 3 0
% 4 1390 2400 1420 2430 1410 2400 1430 2400
% 
\special{pn 20}%
\special{ar 1390 2000 42 42  0.0000000 6.2831853}%
% CIRCLE 0 0 0 0
% 4 1390 2790 1420 2820 1410 2790 1430 2790
% 
\special{pn 20}%
\special{sh 0.600}%
\special{ar 1390 2390 42 42  0.0000000 6.2831853}%
% CIRCLE 0 0 0 0
% 4 1810 3200 1840 3230 1830 3200 1850 3200
% 
\special{pn 20}%
\special{sh 0.600}%
\special{ar 1810 2800 42 42  0.0000000 6.2831853}%
% CIRCLE 0 0 3 0
% 4 1820 2790 1850 2820 1840 2790 1860 2790
% 
\special{pn 20}%
\special{ar 1820 2390 42 42  0.0000000 6.2831853}%
% CIRCLE 0 0 3 0
% 4 1810 2400 1840 2430 1830 2400 1850 2400
% 
\special{pn 20}%
\special{ar 1810 2000 42 42  0.0000000 6.2831853}%
% CIRCLE 0 0 0 0
% 4 1800 2000 1830 2030 1820 2000 1840 2000
% 
\special{pn 20}%
\special{sh 0.600}%
\special{ar 1800 1600 42 42  0.0000000 6.2831853}%
% CIRCLE 0 0 0 0
% 4 2610 3600 2640 3630 2630 3600 2650 3600
% 
\special{pn 20}%
\special{sh 0.600}%
\special{ar 2610 3200 42 42  0.0000000 6.2831853}%
% CIRCLE 0 0 3 0
% 4 2200 3200 2230 3230 2220 3200 2240 3200
% 
\special{pn 20}%
\special{ar 2200 2800 42 42  0.0000000 6.2831853}%
% CIRCLE 0 0 3 0
% 4 2200 2800 2230 2830 2220 2800 2240 2800
% 
\special{pn 20}%
\special{ar 2200 2400 42 42  0.0000000 6.2831853}%
% CIRCLE 0 0 0 0
% 4 2200 2400 2230 2430 2220 2400 2240 2400
% 
\special{pn 20}%
\special{sh 0.600}%
\special{ar 2200 2000 42 42  0.0000000 6.2831853}%
% LINE 1 0 3 0
% 2 1010 1210 1010 2400
% 
\special{pn 13}%
\special{pa 1010 810}%
\special{pa 1010 2000}%
\special{fp}%
% LINE 1 0 3 0
% 2 1010 2400 1810 3200
% 
\special{pn 13}%
\special{pa 1010 2000}%
\special{pa 1810 2800}%
\special{fp}%
% LINE 2 0 3 0
% 4 610 3600 3800 3600 3800 3600 3800 3600
% 
\special{pn 8}%
\special{pa 610 3200}%
\special{pa 3800 3200}%
\special{fp}%
\special{pa 3800 3200}%
\special{pa 3800 3200}%
\special{fp}%
% LINE 2 0 3 0
% 2 1000 800 1000 3860
% 
\special{pn 8}%
\special{pa 1000 400}%
\special{pa 1000 3460}%
\special{fp}%
% STR 2 0 3 0
% 3 770 3720 770 3820 2 0
% $0$
\put(7.7000,-34.2000){\makebox(0,0)[lb]{$0$}}%
% STR 2 0 3 0
% 3 2400 3780 2400 3880 2 0
% $(4,0)$
\put(24.0000,-34.8000){\makebox(0,0)[lb]{$(4,0)$}}%
% STR 2 0 3 0
% 3 510 2360 510 2460 2 0
% $(0,3)$
\put(5.1000,-20.6000){\makebox(0,0)[lb]{$(0,3)$}}%
% STR 2 0 3 0
% 3 500 1160 500 1260 2 0
% $(0,6)$
\put(5.0000,-8.6000){\makebox(0,0)[lb]{$(0,6)$}}%
% STR 2 0 3 0
% 3 1380 3350 1380 3450 2 0
% $(2,1)$
\put(13.8000,-30.5000){\makebox(0,0)[lb]{$(2,1)$}}%
% CIRCLE 0 0 0 0
% 4 3410 3600 3440 3630 3430 3600 3450 3600
% 
\special{pn 20}%
\special{sh 0.600}%
\special{ar 3410 3200 42 42  0.0000000 6.2831853}%
% CIRCLE 0 0 0 0
% 4 2600 2800 2630 2830 2620 2800 2640 2800
% 
\special{pn 20}%
\special{sh 0.600}%
\special{ar 2600 2400 42 42  0.0000000 6.2831853}%
% CIRCLE 0 0 0 0
% 4 3000 3200 3030 3230 3020 3200 3040 3200
% 
\special{pn 20}%
\special{sh 0.600}%
\special{ar 3000 2800 42 42  0.0000000 6.2831853}%
% CIRCLE 0 0 0 0
% 4 3000 3600 3030 3630 3020 3600 3040 3600
% 
\special{pn 20}%
\special{sh 0.600}%
\special{ar 3000 3200 42 42  0.0000000 6.2831853}%
% CIRCLE 0 0 3 0
% 4 2600 3200 2630 3230 2620 3200 2640 3200
% 
\special{pn 20}%
\special{ar 2600 2800 42 42  0.0000000 6.2831853}%
% LINE 1 0 3 0
% 2 1800 3200 2610 3600
% 
\special{pn 13}%
\special{pa 1800 2800}%
\special{pa 2610 3200}%
\special{fp}%
% LINE 1 0 3 0
% 2 2610 3600 3420 3600
% 
\special{pn 13}%
\special{pa 2610 3200}%
\special{pa 3420 3200}%
\special{fp}%
% LINE 1 0 3 0
% 2 1000 1200 3400 3600
% 
\special{pn 13}%
\special{pa 1000 800}%
\special{pa 3400 3200}%
\special{fp}%
% STR 2 0 3 0
% 3 3240 3760 3240 3860 2 0
% $(6,0)$
\put(32.4000,-34.6000){\makebox(0,0)[lb]{$(6,0)$}}%
\end{picture}%
\qquad 
%WinTpicVersion2.15
\unitlength 0.1in
\begin{picture}(23.20,30.60)(5.00,-34.60)
% CIRCLE 0 0 0 0
% 4 1000 1200 1030 1230 1020 1200 1040 1200
% 
\special{pn 20}%
\special{sh 0.600}%
\special{ar 1000 800 42 42  0.0000000 6.2831853}%
% CIRCLE 0 0 0 0
% 4 1010 1610 1040 1640 1030 1610 1050 1610
% 
\special{pn 20}%
\special{sh 0.600}%
\special{ar 1010 1210 42 42  0.0000000 6.2831853}%
% CIRCLE 0 0 0 0
% 4 1010 2000 1040 2030 1030 2000 1050 2000
% 
\special{pn 20}%
\special{sh 0.600}%
\special{ar 1010 1600 42 42  0.0000000 6.2831853}%
% CIRCLE 0 0 0 0
% 4 1010 2400 1040 2430 1030 2400 1050 2400
% 
\special{pn 20}%
\special{sh 0.600}%
\special{ar 1010 2000 42 42  0.0000000 6.2831853}%
% CIRCLE 0 0 0 0
% 4 1400 1610 1430 1640 1420 1610 1440 1610
% 
\special{pn 20}%
\special{sh 0.600}%
\special{ar 1400 1210 42 42  0.0000000 6.2831853}%
% CIRCLE 0 0 3 0
% 4 1390 2000 1420 2030 1410 2000 1430 2000
% 
\special{pn 20}%
\special{ar 1390 1600 42 42  0.0000000 6.2831853}%
% CIRCLE 0 0 3 0
% 4 1390 2400 1420 2430 1410 2400 1430 2400
% 
\special{pn 20}%
\special{ar 1390 2000 42 42  0.0000000 6.2831853}%
% CIRCLE 0 0 3 0
% 4 1390 2790 1420 2820 1410 2790 1430 2790
% 
\special{pn 20}%
\special{ar 1390 2390 42 42  0.0000000 6.2831853}%
% CIRCLE 0 0 3 0
% 4 1810 3200 1840 3230 1830 3200 1850 3200
% 
\special{pn 20}%
\special{ar 1810 2800 42 42  0.0000000 6.2831853}%
% CIRCLE 0 0 3 0
% 4 1820 2790 1850 2820 1840 2790 1860 2790
% 
\special{pn 20}%
\special{ar 1820 2390 42 42  0.0000000 6.2831853}%
% CIRCLE 0 0 3 0
% 4 1810 2400 1840 2430 1830 2400 1850 2400
% 
\special{pn 20}%
\special{ar 1810 2000 42 42  0.0000000 6.2831853}%
% CIRCLE 0 0 0 0
% 4 1800 2000 1830 2030 1820 2000 1840 2000
% 
\special{pn 20}%
\special{sh 0.600}%
\special{ar 1800 1600 42 42  0.0000000 6.2831853}%
% CIRCLE 0 0 0 0
% 4 2200 3600 2230 3630 2220 3600 2240 3600
% 
\special{pn 20}%
\special{sh 0.600}%
\special{ar 2200 3200 42 42  0.0000000 6.2831853}%
% CIRCLE 0 0 0 0
% 4 2200 3200 2230 3230 2220 3200 2240 3200
% 
\special{pn 20}%
\special{sh 0.600}%
\special{ar 2200 2800 42 42  0.0000000 6.2831853}%
% CIRCLE 0 0 0 0
% 4 2200 2800 2230 2830 2220 2800 2240 2800
% 
\special{pn 20}%
\special{sh 0.600}%
\special{ar 2200 2400 42 42  0.0000000 6.2831853}%
% CIRCLE 0 0 0 0
% 4 2200 2400 2230 2430 2220 2400 2240 2400
% 
\special{pn 20}%
\special{sh 0.600}%
\special{ar 2200 2000 42 42  0.0000000 6.2831853}%
% LINE 1 0 3 0
% 2 1000 1210 2200 2400
% 
\special{pn 13}%
\special{pa 1000 810}%
\special{pa 2200 2000}%
\special{fp}%
% LINE 1 0 3 0
% 2 2200 2400 2210 3610
% 
\special{pn 13}%
\special{pa 2200 2000}%
\special{pa 2210 3210}%
\special{fp}%
% LINE 1 0 3 0
% 2 1010 1210 1000 3200
% 
\special{pn 13}%
\special{pa 1010 810}%
\special{pa 1000 2800}%
\special{fp}%
% LINE 2 0 3 0
% 2 610 3600 2820 3600
% 
\special{pn 8}%
\special{pa 610 3200}%
\special{pa 2820 3200}%
\special{fp}%
% LINE 2 0 3 0
% 2 1000 800 1000 3860
% 
\special{pn 8}%
\special{pa 1000 400}%
\special{pa 1000 3460}%
\special{fp}%
% STR 2 0 3 0
% 3 770 3720 770 3820 2 0
% $0$
\put(7.7000,-34.2000){\makebox(0,0)[lb]{$0$}}%
% STR 2 0 3 0
% 3 2020 3780 2020 3880 2 0
% $(3,0)$
\put(20.2000,-34.8000){\makebox(0,0)[lb]{$(3,0)$}}%
% STR 2 0 3 0
% 3 510 3200 510 3300 2 0
% $(0,1)$
\put(5.1000,-29.0000){\makebox(0,0)[lb]{$(0,1)$}}%
% STR 2 0 3 0
% 3 500 1160 500 1260 2 0
% $(0,6)$
\put(5.0000,-8.6000){\makebox(0,0)[lb]{$(0,6)$}}%
% STR 2 0 3 0
% 3 2280 2210 2280 2310 2 0
% $(3,3)$
\put(22.8000,-19.1000){\makebox(0,0)[lb]{$(3,3)$}}%
% CIRCLE 0 0 0 0
% 4 1000 2800 1030 2830 1020 2800 1040 2800
% 
\special{pn 20}%
\special{sh 0.600}%
\special{ar 1000 2400 42 42  0.0000000 6.2831853}%
% CIRCLE 0 0 0 0
% 4 1000 3200 1030 3230 1020 3200 1040 3200
% 
\special{pn 20}%
\special{sh 0.600}%
\special{ar 1000 2800 42 42  0.0000000 6.2831853}%
% CIRCLE 0 0 3 0
% 4 1400 3200 1430 3230 1420 3200 1440 3200
% 
\special{pn 20}%
\special{ar 1400 2800 42 42  0.0000000 6.2831853}%
% CIRCLE 0 0 0 0
% 4 1400 3600 1430 3630 1420 3600 1440 3600
% 
\special{pn 20}%
\special{sh 0.600}%
\special{ar 1400 3200 42 42  0.0000000 6.2831853}%
% CIRCLE 0 0 0 0
% 4 1800 3600 1830 3630 1820 3600 1840 3600
% 
\special{pn 20}%
\special{sh 0.600}%
\special{ar 1800 3200 42 42  0.0000000 6.2831853}%
% LINE 1 0 3 0
% 4 1000 3200 1400 3600 1400 3590 1400 3580
% 
\special{pn 13}%
\special{pa 1000 2800}%
\special{pa 1400 3200}%
\special{fp}%
\special{pa 1400 3190}%
\special{pa 1400 3180}%
\special{fp}%
% LINE 1 0 3 0
% 2 1390 3600 2210 3600
% 
\special{pn 13}%
\special{pa 1390 3200}%
\special{pa 2210 3200}%
\special{fp}%
% STR 2 0 3 0
% 3 1180 3770 1180 3870 2 0
% $(1,0)$
\put(11.8000,-34.7000){\makebox(0,0)[lb]{$(1,0)$}}%
\end{picture}%
\end{center}
\caption{Newton polygons for Lemma \ref{lem:solk1k2} 
(left) and Lemma \ref{lem:solk3k4} (right)} 
\label{fig:snewtonk1k2k3k4}
\end{figure}
%%%%%%%%%%%%%%%%%%%%%%% lem:solk1k2 %%%%%%%%%%%%%%%%%%%%%%
\begin{lemma} \label{lem:solk1k2} 
If $\k_1+\k_2 \not\in \Z$, then there exists a 
holomorphic solution around the origin $z = 0$, 
%%%%%%%%%%%%%%%%%%%%%%% eqn:solk1k21 %%%%%%%%%%%%%%%%%%%%%%
\begin{equation} \label{eqn:solk1k21} 
q = \dfrac{\k_1 z}{\k_1+\k_2} + \k_1 \k_2 
\sum_{k=2}^{\infty} a_{k,+}(\k) \, z^k, \qquad 
p = \k_0(\k_0 + \k_4) 
\sum_{k=0}^{\infty} b_{k,+}(\k) \, z^k, 
\end{equation}
%%%%%%%%%%%%%%%%%%%%%%%%%%%%%%%%%%%%%%%%%%%%%%%%%%%%%%%%%%
depending holomorphically on $\k \in \K$ with 
$\k_1 + \k_2 \not\in \Z$. 
Similarly, if $\k_1-\k_2 \not\in \Z$, then there exists 
a meromorphic solution around the origin $z = 0$, 
%%%%%%%%%%%%%%%%%%%%%%% eqn:solk1k22 %%%%%%%%%%%%%%%%%%%%%%
\begin{equation} \label{eqn:solk1k22} 
q = \dfrac{\k_1 z}{\k_1-\k_2} + \k_1 \k_2 
\sum_{k=2}^{\infty} a_{k,-}(\k) \, z^k, \qquad 
p = \dfrac{\k_1-\k_2}{z} + 
\sum_{k=0}^{\infty} b_{k,-}(\k) \, z^k, 
\end{equation}
%%%%%%%%%%%%%%%%%%%%%%%%%%%%%%%%%%%%%%%%%%%%%%%%%%%%%%%%%%
depending holomorphically on $\k \in \K$ with 
$\k_1 - \k_2 \not\in \Z$. 
\end{lemma} 
%%%%%%%%%%%%%%%%%%%%%%% proof %%%%%%%%%%%%%%%%%%%%%%%%%%%%
{\it Proof}. 
Substituting $q = \k_1 z (\k_1+\ve \k_2)^{-1}+ \k_1\k_2 Q$ 
with $\ve \in \{\pm1\}$ into equation (\ref{eqn:PVI3}) 
yields 
%%%%
\[
(\k_1 + \ve\k_2)^6 P\left(z, \k_1 z(\k_1+\ve\k_2)^{-1} 
+ \k_1\k_2 Q \right) = \k_1^2 \k_2^2 \, p(z, Q), 
\]
%%%%
where $p(z;Q)$ is a differential sum of $(z,Q)$ with 
coefficients in $\C[\k]$ whose Newton polygon is given 
as in Figure \ref{fig:snewtonk1k2k3k4} (left). 
The vertex $(2,1)$ carries the linear differential 
expression 
%%%%%
\[
\mathcal{L}_{\ve}Q = 2 \ve (\k_1 +\ve \k_2)^4 x^2 
\{x^2 Q''-x Q'-(\k_1+\ve\k_2+1)(\k_1+\ve\k_2-1)Q\}, 
\]
%%%%%
while the vertex $(4,0)$ carries the monomial 
$(\k_1+\ve \k_2)^2\{(\k_1+\ve\k_2)^2+\k_3^2-\k_4^2-1\} x^4$. 
The corresponding characteristic polynomial is given by 
%%%%%%
\[
v_{\ve}(k) = 
2 \ve (\k_1+\ve \k_2)^4 
(k-1-\k_1-\ve \k_2)(k-1+\k_1+\ve\k_2). 
\]
%%%%%%
Hence $1+|\k_1+\ve\k_2|$ is the unique critical value 
of the problem. 
If it is not an integer, then the coefficients 
$a_{k,\ve}(\k)$ of the expansions (\ref{eqn:solk1k21}) 
and (\ref{eqn:solk1k22}) are determined uniquely and 
recursively. 
By substituting the resulting power series $q = q(z)$ 
into equation (\ref{eqn:p(z)}), the Laurent 
series for $p = p(z)$ is uniquely determined as in 
(\ref{eqn:solk1k21}) and (\ref{eqn:solk1k22}). 
As is mentioned earlier, the formal solutions 
(\ref{eqn:solk1k21}) and (\ref{eqn:solk1k22}) so 
obtained are convergent. 
\hfill $\Box$ 
%%%%%%%%%%%%%%%%%%%%%% end of proof %%%%%%%%%%%%%%%%%%%%%%%%
%%%%%%%%%%%%%%%%%%%%%%% lem:solk3k4 %%%%%%%%%%%%%%%%%%%%%%%
\begin{lemma} \label{lem:solk3k4} 
Assume that $\k_4$ is nonzero. 
If $\k_4+\k_3 \not\in\Z$, then there exists a holomorphic 
solution germ around the origin $z = 0$, 
%%%%%%%%%%%%%%%%%%%%%%% eqn:solk3k41 %%%%%%%%%%%%%%%%%%%%%%%
\begin{equation} \label{eqn:solk3k41}
q = \dfrac{\k_4+\k_3}{\k_4} + \dfrac{\k_3}{\k_4}
\sum_{k=1}^{\infty} a_{k,+}(\k) \, z^k, 
\qquad 
p = -\k_4\k_0 \sum_{k=0}^{\infty} b_{k,+}(\k) \, z^k, 
\end{equation}
%%%%%%%%%%%%%%%%%%%%%%%%%%%%%%%%%%%%%%%%%%%%%%%%%%%%%%%%%%%%
depending holomorphically on $\k \in \K$ with 
$\k_4 \neq 0$ and $\k_4+\k_3 \not\in \Z$. 
Similarly, if $\k_4-\k_3 \not\in\Z$ then there exists 
a holomorphic solution germ around the origin $z = 0$, 
%%%%%%%%%%%%%%%%%%%%%%% eqn:solk3k42 %%%%%%%%%%%%%%%%%%%%%%%
\begin{equation} \label{eqn:solk3k42}
q = \dfrac{\k_4-\k_3}{\k_4} + \dfrac{\k_3}{\k_4}
\sum_{k=1}^{\infty} a_{k,-}(\k) \, z^k, 
\qquad 
p = -\k_4(\k_0+\k_4 ) \sum_{k=0}^{\infty} b_{k,-}(\k) \, z^k, 
\end{equation}
%%%%%%%%%%%%%%%%%%%%%%%%%%%%%%%%%%%%%%%%%%%%%%%%%%%%%%%%%%%%
depending holomorphically on $\k \in \K$ with 
$\k_4 \neq 0$ and $\k_4-\k_3 \not\in \Z$. 
\end{lemma} 
%%%%%%%%%%%%%%%%%%%%%%% proof %%%%%%%%%%%%%%%%%%%%%%%%%%%%%%
{\it Proof}. 
Substituting $q = \k_4^{-1}(\k_4+\ve \k_3)+\k_4^{-1}\k_3 Q$ 
into equation (\ref{eqn:PVI3}) yields 
%%%%
\[
\k_4^4 P\left(z, \k_4^{-1}(\k_4+\ve\k_3) + 
\k_4^{-1}\k_3 Q \right) = \k_3^2 \, p(z;Q), 
\]
%%%%
where $p(z;Q)$ is a differential sum of $(z,Q)$ with 
coefficients in $\C[\k]$. 
The vertex $(0,1)$ carries the linear differential 
expression $\mathcal{L}_{\ve}Q = 2(\k_4+\ve\k_3)^2 
\{x^2 Q''+x Q'-(\k_4+\ve\k_3)^2 Q \}$, while the vertex 
$(1,0)$ carries the monomial 
$(\k_4+\ve\k_3)^2 \{1+\k_1^2-\k_2^2-(\k_4+\ve\k_3)^2\}x$. 
The corresponding characteristic polynomial is given by 
$v_{\ve}(k) = 2 (\k_4+\ve\k_3)^2 
\{k-(\k_4+\ve\k_3)\}\{k+(\k_4+\ve\k_3)\}$. 
Hence $|\k_4+\ve\k_3|$ is the unique critical value of 
the problem. 
If it is not an integer, then the coefficients 
$a_{k,\ve}(\k)$ of expansions (\ref{eqn:solk3k41}) and 
(\ref{eqn:solk3k42}) are determined uniquely and recursively. 
Then substituting the resulting series for $q = q(z)$ into 
equation (\ref{eqn:p(z)}) yields the Laurent series for 
$p = p(z)$ as in (\ref{eqn:solk3k41}) and 
(\ref{eqn:solk3k42}). 
The formal solutions so obtained are convergent. 
\hfill $\Box$ \par\medskip
%%%%%%%%%%%%%%%%%%%%%%%% end of proof %%%%%%%%%%%%%%%%%%%%%%%
The solutions in Lemmas \ref{lem:solk1k2} and 
\ref{lem:solk3k4} are essentially constructed in 
\cite{BG,Kaneko}. 
We construct more particular solutions for the parameters 
on various strata of higher codimensions. 
%%%%%%%%%%%%%%%%%%%%%%% fig:snewtonA12 %%%%%%%%%%%%%%%%%%%%%%
\begin{figure}[t]
\begin{center}
%WinTpicVersion2.15
\unitlength 0.1in
\begin{picture}(33.00,30.60)(5.00,-34.60)
% CIRCLE 0 0 0 0
% 4 1000 1200 1030 1230 1020 1200 1040 1200
% 
\special{pn 20}%
\special{sh 0.600}%
\special{ar 1000 800 42 42  0.0000000 6.2831853}%
% CIRCLE 0 0 0 0
% 4 1010 1610 1040 1640 1030 1610 1050 1610
% 
\special{pn 20}%
\special{sh 0.600}%
\special{ar 1010 1210 42 42  0.0000000 6.2831853}%
% CIRCLE 0 0 0 0
% 4 1010 2000 1040 2030 1030 2000 1050 2000
% 
\special{pn 20}%
\special{sh 0.600}%
\special{ar 1010 1600 42 42  0.0000000 6.2831853}%
% CIRCLE 0 0 0 0
% 4 1010 2400 1040 2430 1030 2400 1050 2400
% 
\special{pn 20}%
\special{sh 0.600}%
\special{ar 1010 2000 42 42  0.0000000 6.2831853}%
% CIRCLE 0 0 3 0
% 4 1400 1610 1430 1640 1420 1610 1440 1610
% 
\special{pn 20}%
\special{ar 1400 1210 42 42  0.0000000 6.2831853}%
% CIRCLE 0 0 3 0
% 4 1390 2000 1420 2030 1410 2000 1430 2000
% 
\special{pn 20}%
\special{ar 1390 1600 42 42  0.0000000 6.2831853}%
% CIRCLE 0 0 3 0
% 4 1390 2400 1420 2430 1410 2400 1430 2400
% 
\special{pn 20}%
\special{ar 1390 2000 42 42  0.0000000 6.2831853}%
% CIRCLE 0 0 0 0
% 4 1390 2790 1420 2820 1410 2790 1430 2790
% 
\special{pn 20}%
\special{sh 0.600}%
\special{ar 1390 2390 42 42  0.0000000 6.2831853}%
% CIRCLE 0 0 0 0
% 4 1810 3200 1840 3230 1830 3200 1850 3200
% 
\special{pn 20}%
\special{sh 0.600}%
\special{ar 1810 2800 42 42  0.0000000 6.2831853}%
% CIRCLE 0 0 3 0
% 4 1820 2790 1850 2820 1840 2790 1860 2790
% 
\special{pn 20}%
\special{ar 1820 2390 42 42  0.0000000 6.2831853}%
% CIRCLE 0 0 3 0
% 4 1810 2400 1840 2430 1830 2400 1850 2400
% 
\special{pn 20}%
\special{ar 1810 2000 42 42  0.0000000 6.2831853}%
% CIRCLE 0 0 3 0
% 4 1800 2000 1830 2030 1820 2000 1840 2000
% 
\special{pn 20}%
\special{ar 1800 1600 42 42  0.0000000 6.2831853}%
% CIRCLE 0 0 0 0
% 4 2610 3600 2640 3630 2630 3600 2650 3600
% 
\special{pn 20}%
\special{sh 0.600}%
\special{ar 2610 3200 42 42  0.0000000 6.2831853}%
% CIRCLE 0 0 3 0
% 4 2200 3200 2230 3230 2220 3200 2240 3200
% 
\special{pn 20}%
\special{ar 2200 2800 42 42  0.0000000 6.2831853}%
% CIRCLE 0 0 3 0
% 4 2200 2800 2230 2830 2220 2800 2240 2800
% 
\special{pn 20}%
\special{ar 2200 2400 42 42  0.0000000 6.2831853}%
% CIRCLE 0 0 3 0
% 4 2200 2400 2230 2430 2220 2400 2240 2400
% 
\special{pn 20}%
\special{ar 2200 2000 42 42  0.0000000 6.2831853}%
% LINE 1 0 3 0
% 2 1010 1210 1010 2400
% 
\special{pn 13}%
\special{pa 1010 810}%
\special{pa 1010 2000}%
\special{fp}%
% LINE 1 0 3 0
% 2 1010 2400 1810 3200
% 
\special{pn 13}%
\special{pa 1010 2000}%
\special{pa 1810 2800}%
\special{fp}%
% LINE 2 0 3 0
% 4 610 3600 3800 3600 3800 3600 3800 3600
% 
\special{pn 8}%
\special{pa 610 3200}%
\special{pa 3800 3200}%
\special{fp}%
\special{pa 3800 3200}%
\special{pa 3800 3200}%
\special{fp}%
% LINE 2 0 3 0
% 2 1000 800 1000 3860
% 
\special{pn 8}%
\special{pa 1000 400}%
\special{pa 1000 3460}%
\special{fp}%
% STR 2 0 3 0
% 3 770 3720 770 3820 2 0
% $0$
\put(7.7000,-34.2000){\makebox(0,0)[lb]{$0$}}%
% STR 2 0 3 0
% 3 2400 3780 2400 3880 2 0
% $(4,0)$
\put(24.0000,-34.8000){\makebox(0,0)[lb]{$(4,0)$}}%
% STR 2 0 3 0
% 3 510 2360 510 2460 2 0
% $(0,3)$
\put(5.1000,-20.6000){\makebox(0,0)[lb]{$(0,3)$}}%
% STR 2 0 3 0
% 3 500 1160 500 1260 2 0
% $(0,6)$
\put(5.0000,-8.6000){\makebox(0,0)[lb]{$(0,6)$}}%
% STR 2 0 3 0
% 3 1380 3350 1380 3450 2 0
% $(2,1)$
\put(13.8000,-30.5000){\makebox(0,0)[lb]{$(2,1)$}}%
% CIRCLE 0 0 0 0
% 4 3410 3600 3440 3630 3430 3600 3450 3600
% 
\special{pn 20}%
\special{sh 0.600}%
\special{ar 3410 3200 42 42  0.0000000 6.2831853}%
% CIRCLE 0 0 3 0
% 4 2600 2800 2630 2830 2620 2800 2640 2800
% 
\special{pn 20}%
\special{ar 2600 2400 42 42  0.0000000 6.2831853}%
% CIRCLE 0 0 3 0
% 4 3000 3200 3030 3230 3020 3200 3040 3200
% 
\special{pn 20}%
\special{ar 3000 2800 42 42  0.0000000 6.2831853}%
% CIRCLE 0 0 0 0
% 4 3000 3600 3030 3630 3020 3600 3040 3600
% 
\special{pn 20}%
\special{sh 0.600}%
\special{ar 3000 3200 42 42  0.0000000 6.2831853}%
% CIRCLE 0 0 3 0
% 4 2600 3200 2630 3230 2620 3200 2640 3200
% 
\special{pn 20}%
\special{ar 2600 2800 42 42  0.0000000 6.2831853}%
% LINE 1 0 3 0
% 2 1800 3200 2610 3600
% 
\special{pn 13}%
\special{pa 1800 2800}%
\special{pa 2610 3200}%
\special{fp}%
% LINE 1 0 3 0
% 4 2600 3600 3780 3610 3780 3610 3780 3610
% 
\special{pn 13}%
\special{pa 2600 3200}%
\special{pa 3780 3210}%
\special{fp}%
\special{pa 3780 3210}%
\special{pa 3780 3210}%
\special{fp}%
% CIRCLE 0 0 3 0
% 4 3010 2800 3040 2830 3030 2800 3050 2800
% 
\special{pn 20}%
\special{ar 3010 2400 42 42  0.0000000 6.2831853}%
% CIRCLE 0 0 3 0
% 4 2600 2400 2630 2430 2620 2400 2640 2400
% 
\special{pn 20}%
\special{ar 2600 2000 42 42  0.0000000 6.2831853}%
% CIRCLE 0 0 3 0
% 4 3010 2400 3040 2430 3030 2400 3050 2400
% 
\special{pn 20}%
\special{ar 3010 2000 42 42  0.0000000 6.2831853}%
% CIRCLE 0 0 3 0
% 4 2200 1990 2230 2020 2220 1990 2240 1990
% 
\special{pn 20}%
\special{ar 2200 1590 42 42  0.0000000 6.2831853}%
% CIRCLE 0 0 3 0
% 4 1800 1600 1830 1630 1820 1600 1840 1600
% 
\special{pn 20}%
\special{ar 1800 1200 42 42  0.0000000 6.2831853}%
% CIRCLE 0 0 0 0
% 4 1400 1200 1430 1230 1420 1200 1440 1200
% 
\special{pn 20}%
\special{sh 0.600}%
\special{ar 1400 800 42 42  0.0000000 6.2831853}%
% CIRCLE 0 0 0 0
% 4 1800 1190 1830 1220 1820 1190 1840 1190
% 
\special{pn 20}%
\special{sh 0.600}%
\special{ar 1800 790 42 42  0.0000000 6.2831853}%
% CIRCLE 0 0 3 0
% 4 2190 1600 2220 1630 2210 1600 2230 1600
% 
\special{pn 20}%
\special{ar 2190 1200 42 42  0.0000000 6.2831853}%
% CIRCLE 0 0 0 0
% 4 2200 1200 2230 1230 2220 1200 2240 1200
% 
\special{pn 20}%
\special{sh 0.600}%
\special{ar 2200 800 42 42  0.0000000 6.2831853}%
% CIRCLE 0 0 3 0
% 4 2610 1990 2640 2020 2630 1990 2650 1990
% 
\special{pn 20}%
\special{ar 2610 1590 42 42  0.0000000 6.2831853}%
% CIRCLE 0 0 3 0
% 4 3000 2000 3030 2030 3020 2000 3040 2000
% 
\special{pn 20}%
\special{ar 3000 1600 42 42  0.0000000 6.2831853}%
% CIRCLE 0 0 3 0
% 4 2600 1600 2630 1630 2620 1600 2640 1600
% 
\special{pn 20}%
\special{ar 2600 1200 42 42  0.0000000 6.2831853}%
% CIRCLE 0 0 3 0
% 4 3000 1600 3030 1630 3020 1600 3040 1600
% 
\special{pn 20}%
\special{ar 3000 1200 42 42  0.0000000 6.2831853}%
% CIRCLE 0 0 0 0
% 4 2610 1200 2640 1230 2630 1200 2650 1200
% 
\special{pn 20}%
\special{sh 0.600}%
\special{ar 2610 800 42 42  0.0000000 6.2831853}%
% CIRCLE 0 0 0 0
% 4 3010 1200 3040 1230 3030 1200 3050 1200
% 
\special{pn 20}%
\special{sh 0.600}%
\special{ar 3010 800 42 42  0.0000000 6.2831853}%
% CIRCLE 0 0 0 0
% 4 3400 1200 3430 1230 3420 1200 3440 1200
% 
\special{pn 20}%
\special{sh 0.600}%
\special{ar 3400 800 42 42  0.0000000 6.2831853}%
% CIRCLE 0 0 3 0
% 4 3400 1600 3430 1630 3420 1600 3440 1600
% 
\special{pn 20}%
\special{ar 3400 1200 42 42  0.0000000 6.2831853}%
% CIRCLE 0 0 3 0
% 4 3400 1990 3430 2020 3420 1990 3440 1990
% 
\special{pn 20}%
\special{ar 3400 1590 42 42  0.0000000 6.2831853}%
% CIRCLE 0 0 3 0
% 4 3400 2400 3430 2430 3420 2400 3440 2400
% 
\special{pn 20}%
\special{ar 3400 2000 42 42  0.0000000 6.2831853}%
% CIRCLE 0 0 3 0
% 4 3400 2800 3430 2830 3420 2800 3440 2800
% 
\special{pn 20}%
\special{ar 3400 2400 42 42  0.0000000 6.2831853}%
% CIRCLE 0 0 3 0
% 4 3400 3200 3430 3230 3420 3200 3440 3200
% 
\special{pn 20}%
\special{ar 3400 2800 42 42  0.0000000 6.2831853}%
% LINE 1 0 3 0
% 2 1000 1200 3790 1190
% 
\special{pn 13}%
\special{pa 1000 800}%
\special{pa 3790 790}%
\special{fp}%
\end{picture}%
\end{center}
\caption{Newton polygon for Lemma \ref{lem:solA12}}  
\label{fig:snewtonA12}
\end{figure}
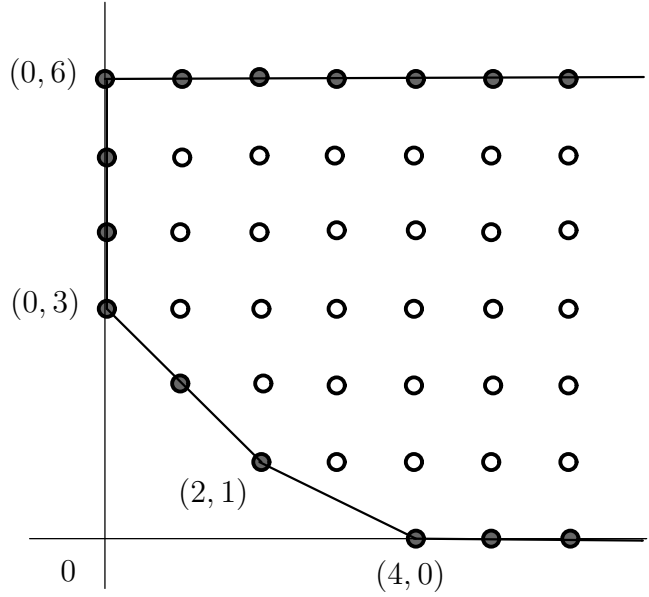
%%%%%%%%%%%%%%%%%%%%%%%%%%%%%%%%%%%%%%%%%%%%%%%%%%%%%%%%%%%
%%%%%%%%%%%%%%%%%%%%%%%% lem:solA12 %%%%%%%%%%%%%%%%%%%%%%%
\begin{lemma}[$\mbox{\boldmath $A_1^{\oplus 2}$ and 
$A_1^{\oplus 3}$}$] \label{lem:solA12} 
If $\k_1 = \k_2 = 0$, then there exists a $1$-parameter 
family of holomorphic solution around the origin $z = 0$, 
%%%%%%%%%%%%%%%%%%%%%%% eqn:solA12 %%%%%%%%%%%%%%%%%%%%%%%
\begin{equation} \label{eqn:solA12} 
\begin{array}{rcl}
q &=& \dfrac{t z}{t+(1-t)(1-z)^{\k_4}} + t(1-t)
\k_0(\k_0+\k_3) \displaystyle 
\sum_{k=2}^{\infty} a_k(t;\k) \, z^k, 
\\[3mm]
p &=& \k_0(\k_0+\k_4) \displaystyle  
\sum_{k=0}^{\infty} b_k(t;\k) \, z^k, 
\end{array}
\end{equation}
%%%%%%%%%%%%%%%%%%%%%%%%%%%%%%%%%%%%%%%%%%%%%%%%%%%%%%%%%%
depending on $t \in \C$, where $a_2(t;\k) = 2$, 
$b_0(t;\k) = 1$ and the remaining coefficients $a_k(t;\k)$, 
$k \ge 3$, and $b_k(t;\k)$, $k \ge 1$, are polynomials of 
$(t,\k_3,\k_4)$ determined uniquely and recursively. 
\end{lemma}
%%%%%%%%%%%%%%%%%%%%%%%%%%%%%%%%%%%%%%%%%%%%%%%%%%%%%%%%%%
{\it Proof}. 
We put $R(z;t) = tz\{t + (1-t)(1-z)^{\k_4}\}^{-1}$. 
Substituting $q = R(z;t) + t(1-t) \k_0 (\k_0 + \k_3) Q$ 
into equation (\ref{eqn:PVI3}) and multiplying the result 
by $\{t+(1-t) (1-z)^{\k_4} \}^4$ yield 
%%%%%%%%%%%%%%%%%%%%%%%% eqn:p(z,Q;t) %%%%%%%%%%%%%%%%%%%%
\begin{equation} \label{eqn:p(z,Q;t)} 
\begin{array}{rrl}
p(z,Q;t) &:=& \{t+(1-t) (1-z)^{\k_4} \}^4 
P(z, R(z;t)+t(1-t) \k_0 (\k_0+\k_3) Q) \\[2mm] 
&=& 2 t^2(1-t)^2 \k_0(\k_0+\k_3) 
\{\mathcal{L}Q + g(z,Q;t) + h(z)\} = 0, 
\end{array}
\end{equation}
%%%%%%%%%%%%%%%%%%%%%%%%%%%%%%%%%%%%%%%%%%%%%%%%%%%%%%%%%%
where $\mathcal{L}Q = z^2 \{z^2 Q''- z Q' + Q\}$ 
and $h(z) = -2z^4 (1-z)^{2\k_4+1}$. 
The Newton polygon of (\ref{eqn:p(z,Q;t)}) is given as 
in Figure \ref{fig:snewtonA12}, where the terms 
$\mathcal{L}Q$ and $h(z)$ correspond to the vertex 
$(2,1)$ and the horizontal infinite edge emanating from 
the vertex $(4,0)$ respectively, and the remaining term 
$g(z,Q;t)$ corresponds to the remaining part of the polygon. 
Since the characteristic equation of $\mathcal{L}Q$ is 
$(k-1)^2 = 0$ having the unique root $k = 1$, the coefficients 
$a_k(t;\k)$, $k \ge 2$, in (\ref{eqn:solA12}) are determined 
uniquely and recursively. 
Here the leading coefficient $a_2(t;\k)$ is found to be
$a_2(t;\k) = 2$ by substituting $Q = a_2(t;\k) z^2$ into 
the truncation $\mathcal{L}Q - 2 z^4 = 0$ of equation 
(\ref{eqn:p(z,Q;t)}) along the edge connecting the vertices 
$(2,1)$ and $(4,0)$. 
Substituting the resulting series $q = q(z)$ into 
(\ref{eqn:p(z)}) we have $p = p(z)$ as in 
(\ref{eqn:solA12}). \hfill $\Box$ 
%%%%%%%%%%%%%%%%%%%%%%% end of proof %%%%%%%%%%%%%%%%%%%%%%%
%%%%%%%%%%%%%%%%%%%%% fig:snewtonA3 %%%%%%%%%%%%%%%%%%%%%%%%
\begin{figure}[t]
\begin{center}
%WinTpicVersion2.15
\unitlength 0.1in
\begin{picture}(33.00,30.60)(5.00,-34.60)
% CIRCLE 0 0 0 0
% 4 1000 1200 1030 1230 1020 1200 1040 1200
% 
\special{pn 20}%
\special{sh 0.600}%
\special{ar 1000 800 42 42  0.0000000 6.2831853}%
% CIRCLE 0 0 0 0
% 4 1010 1610 1040 1640 1030 1610 1050 1610
% 
\special{pn 20}%
\special{sh 0.600}%
\special{ar 1010 1210 42 42  0.0000000 6.2831853}%
% CIRCLE 0 0 0 0
% 4 1010 2000 1040 2030 1030 2000 1050 2000
% 
\special{pn 20}%
\special{sh 0.600}%
\special{ar 1010 1600 42 42  0.0000000 6.2831853}%
% CIRCLE 0 0 0 0
% 4 1010 2400 1040 2430 1030 2400 1050 2400
% 
\special{pn 20}%
\special{sh 0.600}%
\special{ar 1010 2000 42 42  0.0000000 6.2831853}%
% CIRCLE 0 0 3 0
% 4 1400 1610 1430 1640 1420 1610 1440 1610
% 
\special{pn 20}%
\special{ar 1400 1210 42 42  0.0000000 6.2831853}%
% CIRCLE 0 0 3 0
% 4 1390 2000 1420 2030 1410 2000 1430 2000
% 
\special{pn 20}%
\special{ar 1390 1600 42 42  0.0000000 6.2831853}%
% CIRCLE 0 0 3 0
% 4 1390 2400 1420 2430 1410 2400 1430 2400
% 
\special{pn 20}%
\special{ar 1390 2000 42 42  0.0000000 6.2831853}%
% CIRCLE 0 0 3 0
% 4 1400 2800 1430 2830 1420 2800 1440 2800
% 
\special{pn 20}%
\special{ar 1400 2400 42 42  0.0000000 6.2831853}%
% CIRCLE 0 0 3 0
% 4 1810 3200 1840 3230 1830 3200 1850 3200
% 
\special{pn 20}%
\special{ar 1810 2800 42 42  0.0000000 6.2831853}%
% CIRCLE 0 0 3 0
% 4 1800 2800 1830 2830 1820 2800 1840 2800
% 
\special{pn 20}%
\special{ar 1800 2400 42 42  0.0000000 6.2831853}%
% CIRCLE 0 0 3 0
% 4 1810 2400 1840 2430 1830 2400 1850 2400
% 
\special{pn 20}%
\special{ar 1810 2000 42 42  0.0000000 6.2831853}%
% CIRCLE 0 0 3 0
% 4 1800 2000 1830 2030 1820 2000 1840 2000
% 
\special{pn 20}%
\special{ar 1800 1600 42 42  0.0000000 6.2831853}%
% CIRCLE 0 0 0 0
% 4 1000 2800 1030 2830 1020 2800 1040 2800
% 
\special{pn 20}%
\special{sh 0.600}%
\special{ar 1000 2400 42 42  0.0000000 6.2831853}%
% CIRCLE 0 0 3 0
% 4 2200 3200 2230 3230 2220 3200 2240 3200
% 
\special{pn 20}%
\special{ar 2200 2800 42 42  0.0000000 6.2831853}%
% CIRCLE 0 0 3 0
% 4 2200 2810 2230 2840 2220 2810 2240 2810
% 
\special{pn 20}%
\special{ar 2200 2410 42 42  0.0000000 6.2831853}%
% CIRCLE 0 0 3 0
% 4 2200 2400 2230 2430 2220 2400 2240 2400
% 
\special{pn 20}%
\special{ar 2200 2000 42 42  0.0000000 6.2831853}%
% LINE 1 0 3 0
% 2 1010 1210 1010 3200
% 
\special{pn 13}%
\special{pa 1010 810}%
\special{pa 1010 2800}%
\special{fp}%
% LINE 2 0 3 0
% 4 610 3600 3800 3600 3800 3600 3800 3600
% 
\special{pn 8}%
\special{pa 610 3200}%
\special{pa 3800 3200}%
\special{fp}%
\special{pa 3800 3200}%
\special{pa 3800 3200}%
\special{fp}%
% LINE 2 0 3 0
% 2 1000 800 1000 3860
% 
\special{pn 8}%
\special{pa 1000 400}%
\special{pa 1000 3460}%
\special{fp}%
% STR 2 0 3 0
% 3 770 3720 770 3820 2 0
% $0$
\put(7.7000,-34.2000){\makebox(0,0)[lb]{$0$}}%
% STR 2 0 3 0
% 3 500 1160 500 1260 2 0
% $(0,6)$
\put(5.0000,-8.6000){\makebox(0,0)[lb]{$(0,6)$}}%
% STR 2 0 3 0
% 3 500 3220 500 3320 2 0
% $(0,1)$
\put(5.0000,-29.2000){\makebox(0,0)[lb]{$(0,1)$}}%
% CIRCLE 0 0 0 0
% 4 1000 3210 1030 3240 1020 3210 1040 3210
% 
\special{pn 20}%
\special{sh 0.600}%
\special{ar 1000 2810 42 42  0.0000000 6.2831853}%
% CIRCLE 0 0 3 0
% 4 2600 2800 2630 2830 2620 2800 2640 2800
% 
\special{pn 20}%
\special{ar 2600 2400 42 42  0.0000000 6.2831853}%
% CIRCLE 0 0 3 0
% 4 3000 3200 3030 3230 3020 3200 3040 3200
% 
\special{pn 20}%
\special{ar 3000 2800 42 42  0.0000000 6.2831853}%
% CIRCLE 0 0 3 0
% 4 1400 3200 1430 3230 1420 3200 1440 3200
% 
\special{pn 20}%
\special{ar 1400 2800 42 42  0.0000000 6.2831853}%
% CIRCLE 0 0 3 0
% 4 2600 3200 2630 3230 2620 3200 2640 3200
% 
\special{pn 20}%
\special{ar 2600 2800 42 42  0.0000000 6.2831853}%
% CIRCLE 0 0 3 0
% 4 3010 2800 3040 2830 3030 2800 3050 2800
% 
\special{pn 20}%
\special{ar 3010 2400 42 42  0.0000000 6.2831853}%
% CIRCLE 0 0 3 0
% 4 2600 2400 2630 2430 2620 2400 2640 2400
% 
\special{pn 20}%
\special{ar 2600 2000 42 42  0.0000000 6.2831853}%
% CIRCLE 0 0 3 0
% 4 3010 2400 3040 2430 3030 2400 3050 2400
% 
\special{pn 20}%
\special{ar 3010 2000 42 42  0.0000000 6.2831853}%
% CIRCLE 0 0 3 0
% 4 2200 1990 2230 2020 2220 1990 2240 1990
% 
\special{pn 20}%
\special{ar 2200 1590 42 42  0.0000000 6.2831853}%
% CIRCLE 0 0 3 0
% 4 1800 1600 1830 1630 1820 1600 1840 1600
% 
\special{pn 20}%
\special{ar 1800 1200 42 42  0.0000000 6.2831853}%
% CIRCLE 0 0 0 0
% 4 1400 1200 1430 1230 1420 1200 1440 1200
% 
\special{pn 20}%
\special{sh 0.600}%
\special{ar 1400 800 42 42  0.0000000 6.2831853}%
% CIRCLE 0 0 0 0
% 4 1800 1190 1830 1220 1820 1190 1840 1190
% 
\special{pn 20}%
\special{sh 0.600}%
\special{ar 1800 790 42 42  0.0000000 6.2831853}%
% CIRCLE 0 0 3 0
% 4 2190 1600 2220 1630 2210 1600 2230 1600
% 
\special{pn 20}%
\special{ar 2190 1200 42 42  0.0000000 6.2831853}%
% CIRCLE 0 0 0 0
% 4 2200 1200 2230 1230 2220 1200 2240 1200
% 
\special{pn 20}%
\special{sh 0.600}%
\special{ar 2200 800 42 42  0.0000000 6.2831853}%
% CIRCLE 0 0 3 0
% 4 2610 1990 2640 2020 2630 1990 2650 1990
% 
\special{pn 20}%
\special{ar 2610 1590 42 42  0.0000000 6.2831853}%
% CIRCLE 0 0 3 0
% 4 3000 2000 3030 2030 3020 2000 3040 2000
% 
\special{pn 20}%
\special{ar 3000 1600 42 42  0.0000000 6.2831853}%
% CIRCLE 0 0 3 0
% 4 2600 1600 2630 1630 2620 1600 2640 1600
% 
\special{pn 20}%
\special{ar 2600 1200 42 42  0.0000000 6.2831853}%
% CIRCLE 0 0 3 0
% 4 3000 1600 3030 1630 3020 1600 3040 1600
% 
\special{pn 20}%
\special{ar 3000 1200 42 42  0.0000000 6.2831853}%
% CIRCLE 0 0 0 0
% 4 2610 1200 2640 1230 2630 1200 2650 1200
% 
\special{pn 20}%
\special{sh 0.600}%
\special{ar 2610 800 42 42  0.0000000 6.2831853}%
% CIRCLE 0 0 0 0
% 4 3010 1200 3040 1230 3030 1200 3050 1200
% 
\special{pn 20}%
\special{sh 0.600}%
\special{ar 3010 800 42 42  0.0000000 6.2831853}%
% CIRCLE 0 0 0 0
% 4 3400 1200 3430 1230 3420 1200 3440 1200
% 
\special{pn 20}%
\special{sh 0.600}%
\special{ar 3400 800 42 42  0.0000000 6.2831853}%
% CIRCLE 0 0 3 0
% 4 3400 1600 3430 1630 3420 1600 3440 1600
% 
\special{pn 20}%
\special{ar 3400 1200 42 42  0.0000000 6.2831853}%
% CIRCLE 0 0 3 0
% 4 3400 1990 3430 2020 3420 1990 3440 1990
% 
\special{pn 20}%
\special{ar 3400 1590 42 42  0.0000000 6.2831853}%
% CIRCLE 0 0 3 0
% 4 3400 2400 3430 2430 3420 2400 3440 2400
% 
\special{pn 20}%
\special{ar 3400 2000 42 42  0.0000000 6.2831853}%
% CIRCLE 0 0 3 0
% 4 3400 2800 3430 2830 3420 2800 3440 2800
% 
\special{pn 20}%
\special{ar 3400 2400 42 42  0.0000000 6.2831853}%
% CIRCLE 0 0 3 0
% 4 3400 3200 3430 3230 3420 3200 3440 3200
% 
\special{pn 20}%
\special{ar 3400 2800 42 42  0.0000000 6.2831853}%
% LINE 1 0 3 0
% 2 1000 1200 3790 1190
% 
\special{pn 13}%
\special{pa 1000 800}%
\special{pa 3790 790}%
\special{fp}%
% CIRCLE 0 0 0 0
% 4 1810 3600 1840 3630 1830 3600 1850 3600
% 
\special{pn 20}%
\special{sh 0.600}%
\special{ar 1810 3200 42 42  0.0000000 6.2831853}%
% CIRCLE 0 0 0 0
% 4 2200 3600 2230 3630 2220 3600 2240 3600
% 
\special{pn 20}%
\special{sh 0.600}%
\special{ar 2200 3200 42 42  0.0000000 6.2831853}%
% CIRCLE 0 0 0 0
% 4 2600 3600 2630 3630 2620 3600 2640 3600
% 
\special{pn 20}%
\special{sh 0.600}%
\special{ar 2600 3200 42 42  0.0000000 6.2831853}%
% CIRCLE 0 0 0 0
% 4 3000 3600 3030 3630 3020 3600 3040 3600
% 
\special{pn 20}%
\special{sh 0.600}%
\special{ar 3000 3200 42 42  0.0000000 6.2831853}%
% CIRCLE 0 0 0 0
% 4 3400 3600 3430 3630 3420 3600 3440 3600
% 
\special{pn 20}%
\special{sh 0.600}%
\special{ar 3400 3200 42 42  0.0000000 6.2831853}%
% LINE 1 0 3 0
% 2 1010 3210 1800 3600
% 
\special{pn 13}%
\special{pa 1010 2810}%
\special{pa 1800 3200}%
\special{fp}%
% LINE 1 0 3 0
% 2 1800 3600 3800 3600
% 
\special{pn 13}%
\special{pa 1800 3200}%
\special{pa 3800 3200}%
\special{fp}%
% STR 2 0 3 0
% 3 1660 3770 1660 3870 2 0
% $(2,0)$
\put(16.6000,-34.7000){\makebox(0,0)[lb]{$(2,0)$}}%
\end{picture}%
\end{center}
\caption{Newton polygon for Lemma \ref{lem:solA3}}  
\label{fig:snewtonA3}
\end{figure}
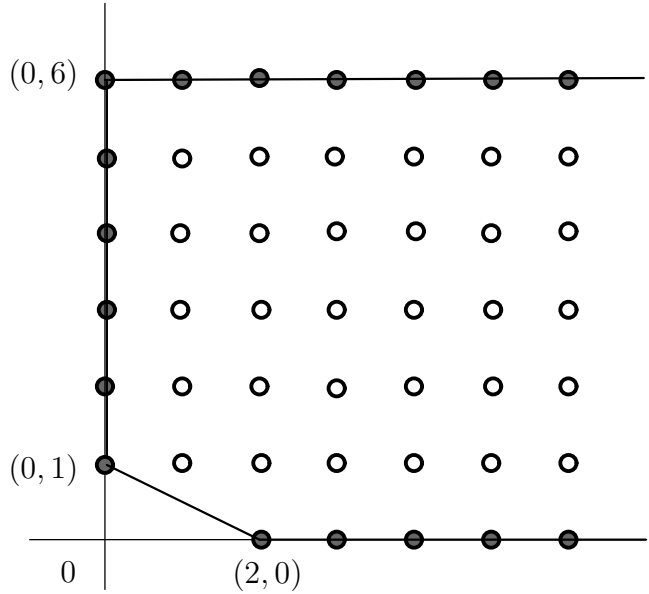
%%%%%%%%%%%%%%%%%%%%%%%%%%%%%%%%%%%%%%%%%%%%%%%%%%%%%%%%%%%
%%%%%%%%%%%%%%%%%%%%%% lem:solA3 %%%%%%%%%%%%%%%%%%%%%%%%%%
\begin{lemma}[$\mbox{\boldmath $A_3$}$] 
\label{lem:solA3} 
If $\k_0 = \k_1 = \k_2 = 0$, $\k_3 + \k_4 = 1$, then there 
exists a $1$-parameter family of holomorphic solutions 
around the origin $z = 0$, 
%%%%%%%%%%%%%%%%%%%%%%% eqn:solA3 %%%%%%%%%%%%%%%%%%%%%%%%%
\begin{equation} \label{eqn:solA3}
q = \dfrac{z}{1-(1-z)^{\k_4}} + t \k_3 z + 
t \k_3 \k_4 \sum_{k=2}^{\infty} a_{k}(t;\k) \, z^k, 
\qquad 
p = - t \k_4^3 \sum_{k=1}^{\infty} b_k(t;\k) \, z^k, 
\end{equation}
%%%%%%%%%%%%%%%%%%%%%%%%%%%%%%%%%%%%%%%%%%%%%%%%%%%%%%%%%%%%
depending on $t \in \C$, where $a_2(t;\k) = t \k_3$, 
$b_1(t;\k) = 1$ and the coeffcients $a_{k}(t;\k)$, $k \ge 2$, 
and $b_k(t;\k)$, $k \ge 1$, are polynomials of $(t,\k_4)$ 
determined uniquely and recursively. 
\end{lemma}
%%%%%%%%%%%%%%%%%%%%%%%%%%%%%%%%%%%%%%%%%%%%%%%%%%%%%%%%%%%%%
{\it Proof}. 
Put $R(z) = z \{1-(1-z)^{\k_4}\}^{-1}$. 
Substituting $q = R(z) + t \k_3 z + t \k_3 \k_4 Q$ into 
(\ref{eqn:PVI3}) yields 
%%%%%%
\[
P(z, R(z) + t \k_3 z + t \k_3 \k_4 Q) = 
t \k_3^2 \k_4^2 R(x)^4 \, p(z,Q;t),
\]
%%%%%%
where $p(z,Q;t)$ is a differential sum of $(z,Q)$ with 
coefficients in $\C[t,\k_4]$ whose Newton polygon is given 
as in Figure \ref{fig:snewtonA3}. 
Especially the vertex $(0,1)$ carries the linear 
differential expression 
$\mathcal{L}Q = 2 (z^2 Q'' + z Q' - Q)$, whose characteristic 
polynomial is $2(k-1)(k+1)$, while the vertex $(2,0)$ carries 
the monomial $-6 t \k_3 z^2$. 
Since the critical values $k = \pm 1$ are smaller than $2$, 
the coefficients $a_k(t;\k)$, $k \ge 2$, in (\ref{eqn:solA3}) 
are determined uniquely and recursively, where the leading 
coefficient $a_2(t;\k)$ is found to be $t \k_3$. 
The rest of the proof is similar to that in 
Lemma \ref{lem:solA12}. \hfill $\Box$ 
%%%%%%%%%%%%%%%%%%%%%%% end of proof %%%%%%%%%%%%%%%%%%%%%%%%%
%%%%%%%%%%%%%%%%%%%%%%% fig:snewtonD4A14 %%%%%%%%%%%%%%%%%%%%%
\begin{figure}[t]
\begin{center}
%WinTpicVersion2.15
\unitlength 0.1in
\begin{picture}(33.00,30.60)(5.00,-34.60)
% CIRCLE 0 0 0 0
% 4 1000 1200 1030 1230 1020 1200 1040 1200
% 
\special{pn 20}%
\special{sh 0.600}%
\special{ar 1000 800 42 42  0.0000000 6.2831853}%
% CIRCLE 0 0 0 0
% 4 1010 1610 1040 1640 1030 1610 1050 1610
% 
\special{pn 20}%
\special{sh 0.600}%
\special{ar 1010 1210 42 42  0.0000000 6.2831853}%
% CIRCLE 0 0 0 0
% 4 1010 2000 1040 2030 1030 2000 1050 2000
% 
\special{pn 20}%
\special{sh 0.600}%
\special{ar 1010 1600 42 42  0.0000000 6.2831853}%
% CIRCLE 0 0 0 0
% 4 1010 2400 1040 2430 1030 2400 1050 2400
% 
\special{pn 20}%
\special{sh 0.600}%
\special{ar 1010 2000 42 42  0.0000000 6.2831853}%
% CIRCLE 0 0 0 0
% 4 1400 1610 1430 1640 1420 1610 1440 1610
% 
\special{pn 20}%
\special{sh 0.600}%
\special{ar 1400 1210 42 42  0.0000000 6.2831853}%
% CIRCLE 0 0 3 0
% 4 1390 2000 1420 2030 1410 2000 1430 2000
% 
\special{pn 20}%
\special{ar 1390 1600 42 42  0.0000000 6.2831853}%
% CIRCLE 0 0 3 0
% 4 1390 2400 1420 2430 1410 2400 1430 2400
% 
\special{pn 20}%
\special{ar 1390 2000 42 42  0.0000000 6.2831853}%
% CIRCLE 0 0 3 0
% 4 1390 2790 1420 2820 1410 2790 1430 2790
% 
\special{pn 20}%
\special{ar 1390 2390 42 42  0.0000000 6.2831853}%
% CIRCLE 0 0 3 0
% 4 1810 3200 1840 3230 1830 3200 1850 3200
% 
\special{pn 20}%
\special{ar 1810 2800 42 42  0.0000000 6.2831853}%
% CIRCLE 0 0 3 0
% 4 1820 2790 1850 2820 1840 2790 1860 2790
% 
\special{pn 20}%
\special{ar 1820 2390 42 42  0.0000000 6.2831853}%
% CIRCLE 0 0 3 0
% 4 1810 2400 1840 2430 1830 2400 1850 2400
% 
\special{pn 20}%
\special{ar 1810 2000 42 42  0.0000000 6.2831853}%
% CIRCLE 0 0 0 0
% 4 1800 2000 1830 2030 1820 2000 1840 2000
% 
\special{pn 20}%
\special{sh 0.600}%
\special{ar 1800 1600 42 42  0.0000000 6.2831853}%
% CIRCLE 0 0 0 0
% 4 2610 3600 2640 3630 2630 3600 2650 3600
% 
\special{pn 20}%
\special{sh 0.600}%
\special{ar 2610 3200 42 42  0.0000000 6.2831853}%
% CIRCLE 0 0 3 0
% 4 2200 3200 2230 3230 2220 3200 2240 3200
% 
\special{pn 20}%
\special{ar 2200 2800 42 42  0.0000000 6.2831853}%
% CIRCLE 0 0 3 0
% 4 2200 2800 2230 2830 2220 2800 2240 2800
% 
\special{pn 20}%
\special{ar 2200 2400 42 42  0.0000000 6.2831853}%
% CIRCLE 0 0 0 0
% 4 2200 2400 2230 2430 2220 2400 2240 2400
% 
\special{pn 20}%
\special{sh 0.600}%
\special{ar 2200 2000 42 42  0.0000000 6.2831853}%
% LINE 1 0 3 0
% 2 1010 1210 990 2800
% 
\special{pn 13}%
\special{pa 1010 810}%
\special{pa 990 2400}%
\special{fp}%
% LINE 2 0 3 0
% 4 610 3600 3800 3600 3800 3600 3800 3600
% 
\special{pn 8}%
\special{pa 610 3200}%
\special{pa 3800 3200}%
\special{fp}%
\special{pa 3800 3200}%
\special{pa 3800 3200}%
\special{fp}%
% LINE 2 0 3 0
% 2 1000 800 1000 3860
% 
\special{pn 8}%
\special{pa 1000 400}%
\special{pa 1000 3460}%
\special{fp}%
% STR 2 0 3 0
% 3 770 3720 770 3820 2 0
% $0$
\put(7.7000,-34.2000){\makebox(0,0)[lb]{$0$}}%
% STR 2 0 3 0
% 3 2000 3760 2000 3860 2 0
% $(3,0)$
\put(20.0000,-34.6000){\makebox(0,0)[lb]{$(3,0)$}}%
% STR 2 0 3 0
% 3 500 2780 500 2880 2 0
% $(0,2)$
\put(5.0000,-24.8000){\makebox(0,0)[lb]{$(0,2)$}}%
% STR 2 0 3 0
% 3 500 1160 500 1260 2 0
% $(0,6)$
\put(5.0000,-8.6000){\makebox(0,0)[lb]{$(0,6)$}}%
% STR 2 0 3 0
% 3 1070 3380 1070 3480 2 0
% $(1,1)$
\put(10.7000,-30.8000){\makebox(0,0)[lb]{$(1,1)$}}%
% CIRCLE 0 0 0 0
% 4 3410 3600 3440 3630 3430 3600 3450 3600
% 
\special{pn 20}%
\special{sh 0.600}%
\special{ar 3410 3200 42 42  0.0000000 6.2831853}%
% CIRCLE 0 0 0 0
% 4 2600 2800 2630 2830 2620 2800 2640 2800
% 
\special{pn 20}%
\special{sh 0.600}%
\special{ar 2600 2400 42 42  0.0000000 6.2831853}%
% CIRCLE 0 0 0 0
% 4 3000 3200 3030 3230 3020 3200 3040 3200
% 
\special{pn 20}%
\special{sh 0.600}%
\special{ar 3000 2800 42 42  0.0000000 6.2831853}%
% CIRCLE 0 0 0 0
% 4 3000 3600 3030 3630 3020 3600 3040 3600
% 
\special{pn 20}%
\special{sh 0.600}%
\special{ar 3000 3200 42 42  0.0000000 6.2831853}%
% CIRCLE 0 0 3 0
% 4 2600 3200 2630 3230 2620 3200 2640 3200
% 
\special{pn 20}%
\special{ar 2600 2800 42 42  0.0000000 6.2831853}%
% LINE 1 0 3 0
% 2 1000 1200 3400 3600
% 
\special{pn 13}%
\special{pa 1000 800}%
\special{pa 3400 3200}%
\special{fp}%
% STR 2 0 3 0
% 3 3240 3760 3240 3860 2 0
% $(6,0)$
\put(32.4000,-34.6000){\makebox(0,0)[lb]{$(6,0)$}}%
% CIRCLE 0 0 0 0
% 4 2200 3600 2230 3630 2220 3600 2240 3600
% 
\special{pn 20}%
\special{sh 0.600}%
\special{ar 2200 3200 42 42  0.0000000 6.2831853}%
% CIRCLE 0 0 0 0
% 4 1400 3200 1430 3230 1420 3200 1440 3200
% 
\special{pn 20}%
\special{sh 0.600}%
\special{ar 1400 2800 42 42  0.0000000 6.2831853}%
% CIRCLE 0 0 0 0
% 4 1000 2800 1030 2830 1020 2800 1040 2800
% 
\special{pn 20}%
\special{sh 0.600}%
\special{ar 1000 2400 42 42  0.0000000 6.2831853}%
% LINE 1 0 3 0
% 2 1000 2800 1390 3190
% 
\special{pn 13}%
\special{pa 1000 2400}%
\special{pa 1390 2790}%
\special{fp}%
% LINE 1 0 3 0
% 2 1390 3200 2200 3600
% 
\special{pn 13}%
\special{pa 1390 2800}%
\special{pa 2200 3200}%
\special{fp}%
% LINE 1 0 3 0
% 2 2210 3600 3400 3600
% 
\special{pn 13}%
\special{pa 2210 3200}%
\special{pa 3400 3200}%
\special{fp}%
\end{picture}%
\qquad 
%WinTpicVersion2.15
\unitlength 0.1in
\begin{picture}(16.00,22.20)(3.00,-25.10)
% CIRCLE 0 0 0 0
% 4 810 1000 840 1030 830 1000 850 1000
% 
\special{pn 20}%
\special{sh 0.600}%
\special{ar 810 600 42 42  0.0000000 6.2831853}%
% CIRCLE 0 0 0 0
% 4 810 1400 840 1430 830 1400 850 1400
% 
\special{pn 20}%
\special{sh 0.600}%
\special{ar 810 1000 42 42  0.0000000 6.2831853}%
% CIRCLE 0 0 0 0
% 4 1190 1000 1220 1030 1210 1000 1230 1000
% 
\special{pn 20}%
\special{sh 0.600}%
\special{ar 1190 600 42 42  0.0000000 6.2831853}%
% CIRCLE 0 0 3 0
% 4 1190 1400 1220 1430 1210 1400 1230 1400
% 
\special{pn 20}%
\special{ar 1190 1000 42 42  0.0000000 6.2831853}%
% CIRCLE 0 0 3 0
% 4 1190 1790 1220 1820 1210 1790 1230 1790
% 
\special{pn 20}%
\special{ar 1190 1390 42 42  0.0000000 6.2831853}%
% CIRCLE 0 0 0 0
% 4 1610 2200 1640 2230 1630 2200 1650 2200
% 
\special{pn 20}%
\special{sh 0.600}%
\special{ar 1610 1800 42 42  0.0000000 6.2831853}%
% CIRCLE 0 0 0 0
% 4 1620 1790 1650 1820 1640 1790 1660 1790
% 
\special{pn 20}%
\special{sh 0.600}%
\special{ar 1620 1390 42 42  0.0000000 6.2831853}%
% CIRCLE 0 0 0 0
% 4 1610 1400 1640 1430 1630 1400 1650 1400
% 
\special{pn 20}%
\special{sh 0.600}%
\special{ar 1610 1000 42 42  0.0000000 6.2831853}%
% LINE 2 0 3 0
% 2 410 2600 1900 2600
% 
\special{pn 8}%
\special{pa 410 2200}%
\special{pa 1900 2200}%
\special{fp}%
% LINE 2 0 3 0
% 2 800 690 790 2910
% 
\special{pn 8}%
\special{pa 800 290}%
\special{pa 790 2510}%
\special{fp}%
% STR 2 0 3 0
% 3 570 2720 570 2820 2 0
% $0$
\put(5.7000,-24.2000){\makebox(0,0)[lb]{$0$}}%
% STR 2 0 3 0
% 3 1650 2230 1650 2330 2 0
% $(2,1)$
\put(16.5000,-19.3000){\makebox(0,0)[lb]{$(2,1)$}}%
% STR 2 0 3 0
% 3 310 2200 310 2300 2 0
% $(0,1)$
\put(3.1000,-19.0000){\makebox(0,0)[lb]{$(0,1)$}}%
% STR 2 0 3 0
% 3 300 920 300 1020 2 0
% $(0,4)$
\put(3.0000,-6.2000){\makebox(0,0)[lb]{$(0,4)$}}%
% STR 2 0 3 0
% 3 1630 1300 1630 1400 2 0
% $(2,3)$
\put(16.3000,-10.0000){\makebox(0,0)[lb]{$(2,3)$}}%
% CIRCLE 0 0 0 0
% 4 800 1800 830 1830 820 1800 840 1800
% 
\special{pn 20}%
\special{sh 0.600}%
\special{ar 800 1400 42 42  0.0000000 6.2831853}%
% CIRCLE 0 0 0 0
% 4 800 2200 830 2230 820 2200 840 2200
% 
\special{pn 20}%
\special{sh 0.600}%
\special{ar 800 1800 42 42  0.0000000 6.2831853}%
% CIRCLE 0 0 3 0
% 4 1200 2200 1230 2230 1220 2200 1240 2200
% 
\special{pn 20}%
\special{ar 1200 1800 42 42  0.0000000 6.2831853}%
% CIRCLE 0 0 0 0
% 4 1200 2600 1230 2630 1220 2600 1240 2600
% 
\special{pn 20}%
\special{sh 0.600}%
\special{ar 1200 2200 42 42  0.0000000 6.2831853}%
% LINE 1 0 3 0
% 4 800 2200 1200 2600 1200 2590 1200 2580
% 
\special{pn 13}%
\special{pa 800 1800}%
\special{pa 1200 2200}%
\special{fp}%
\special{pa 1200 2190}%
\special{pa 1200 2180}%
\special{fp}%
% STR 2 0 3 0
% 3 980 2770 980 2870 2 0
% $(1,0)$
\put(9.8000,-24.7000){\makebox(0,0)[lb]{$(1,0)$}}%
% LINE 1 0 3 0
% 2 800 1000 800 2200
% 
\special{pn 13}%
\special{pa 800 600}%
\special{pa 800 1800}%
\special{fp}%
% LINE 1 0 3 0
% 2 800 1000 1200 1000
% 
\special{pn 13}%
\special{pa 800 600}%
\special{pa 1200 600}%
\special{fp}%
% LINE 1 0 3 0
% 4 1610 1410 1620 2130 1610 2200 1610 2200
% 
\special{pn 13}%
\special{pa 1610 1010}%
\special{pa 1620 1730}%
\special{fp}%
\special{pa 1610 1800}%
\special{pa 1610 1800}%
\special{fp}%
% LINE 1 0 3 0
% 2 1200 1000 1600 1400
% 
\special{pn 13}%
\special{pa 1200 600}%
\special{pa 1600 1000}%
\special{fp}%
% LINE 1 0 3 0
% 2 1610 2200 1200 2600
% 
\special{pn 13}%
\special{pa 1610 1800}%
\special{pa 1200 2200}%
\special{fp}%
\end{picture}%
\end{center}
\caption{Newton polygons for Lemma \ref{lem:solD4} (left) 
and Lemma \ref{lem:solA14} (right)} 
\label{fig:snewtonD4A14} 
\end{figure}
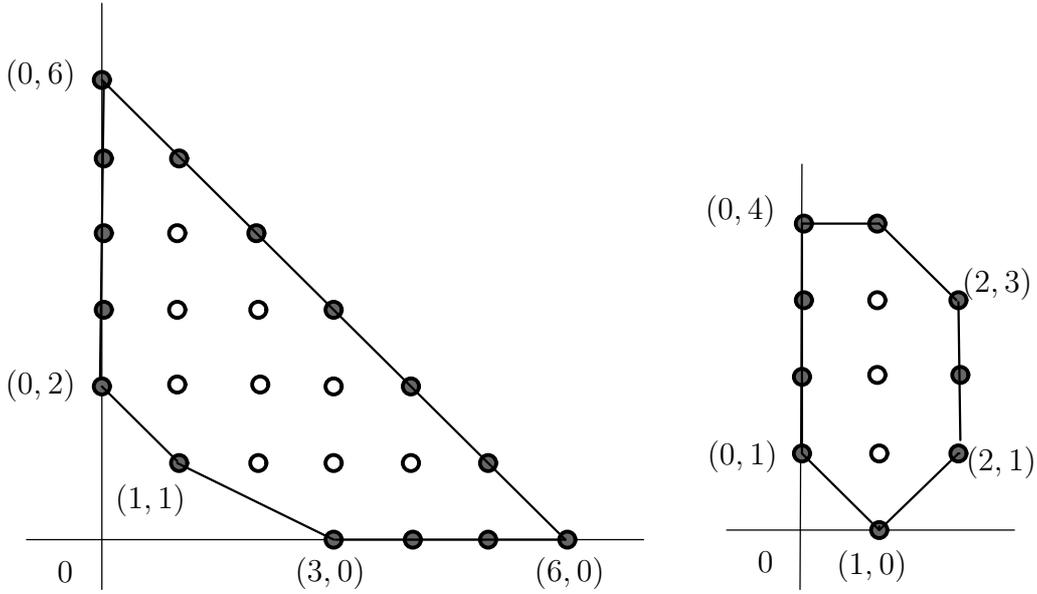
%%%%%%%%%%%%%%%%%%%%%%% lem:solD4 %%%%%%%%%%%%%%%%%%%%%%%%%%%
\begin{lemma}[$\mbox{\boldmath $D_4$}$]  
\label{lem:solD4} 
If $\k_0 = \k_1 = \k_2 = \k_3 = 0$ and $\k_4 = 1$, then there 
exists a $1$-parameter family of holomorphic solution germs 
around the origin $z = 0$, 
%%%%%%%%%%%%%%%%%%%%%%% eqn:solD4 %%%%%%%%%%%%%%%%%%%%%%%%%%%
\begin{equation} \label{eqn:solD4}
q = 1 + t \sum_{k=1}^{\infty} a_k(t) \, z^k, \qquad 
p = \dfrac{z}{(1-z) \log(1-z)} + t \sum_{k=1}^{\infty} 
b_k(t) \, z^k, 
\end{equation}
%%%%%%%%%%%%%%%%%%%%%%%%%%%%%%%%%%%%%%%%%%%%%%%%%%%%%%%%%%%%%
depending on a parameter $t \in \C$, where 
$a_1(t) = b_1(t) = 1$ and the remaining coeffieicents 
$a_k(t)$ and $b_k(t)$, $k \ge 2$, are polynomials of $t$ 
determined uniquely and recursively. 
\end{lemma}
%%%%%%%%%%%%%%%%%%%%%%%%%%%%%%%%%%%%%%%%%%%%%%%%%%%%%%%%%%%%%
%%%%%%%%%%%%%%%%%%%%%%%%%% proof %%%%%%%%%%%%%%%%%%%%%%%%%%%%
{\it Proof}. 
Substituting $q = 1 + tz + t Q$ into equation 
(\ref{eqn:PVI3}) yields $P(z,1+tQ) = t^2 p(z,Q;t)$, where 
$p(z,Q;t)$ is a differential sum of $(z,Q)$ with coefficients 
in $\C[t]$ whose Newton polygon is given as in 
Figure \ref{fig:snewtonD4A14} (left). 
Consider the edge connecting $(1,1)$ and $(3,0)$. 
The vertex $(1,1)$ carries the differential monomial 
$\mathcal{L}Q = 2 z^3 Q''$, whose characteristic 
polynomial is $k(k-1)$, while the vertex $(3,0)$ 
carries the monomial $-4 t z^3$. 
Put $a_1(t) = 1$. 
Since the critical values $k = 0$, $1$ are smaller 
than $2$, the coefficients $a_k(t)$, $k \ge 2$, in 
(\ref{eqn:solD4}) are determined uniquely and 
recursively. 
Substituting the resulting series $q = q(z)$ into 
(\ref{eqn:p(z)}) we have $p = p(z)$ as in 
(\ref{eqn:solD4}). 
Here the term $z/\{(1-z) \log(1-z)\}$ is singled 
out from $p = p(z)$, because putting $t = 0$ yields the 
special solution $q \equiv 1$ and 
$p = z/\{(1-z) \log(1-z)\}$ (see also 
Lemma \ref{lem:d4ric}). \hfill $\Box$ 
%%%%%%%%%%%%%%%%%%%%%%% end of proof %%%%%%%%%%%%%%%%%%%%%%%%
%%%%%%%%%%%%%%%%%%%%%%% lem:solA14 %%%%%%%%%%%%%%%%%%%%%%%%%%
\begin{lemma}[$\mbox{\boldmath $A_1^{\oplus 4}$}$]  
\label{lem:solA14} 
Let $\k_0 = 1/2$ and $\k_1=\k_2=\k_3=\k_4=0$. 
Then there exists a $1$-parameter family of holomorphic 
solutions around the origin $z = 0$ depending on a parameter 
$t \in \C$, 
%%%%%%%%%%%%%%%%%%%%%%% eqn:solA141 %%%%%%%%%%%%%%%%%%%%%%%%%
\begin{equation} \label{eqn:solA141} 
q = t z + t(1-t) \sum_{k=2}^{\infty} a_k(t) \, z^k, \qquad 
p = \sum_{k=0}^{\infty} b_k(t) \, z^k, 
\end{equation}
%%%%%%%%%%%%%%%%%%%%%%%%%%%%%%%%%%%%%%%%%%%%%%%%%%%%%%%%%%%%%%
where the coefficients $a_k(t)$ and $b_k(t)$ are polynomials 
of $t$ beginning with $a_2(t) = 1/2$ and $b_0(t) = 1/4$. 
Moreover there is another $1$-parameter family of solutions 
around $z = 0$, 
%%%%%%%%%%%%%%%%%%%%%%% eqn:solA142 %%%%%%%%%%%%%%%%%%%%%%%%%
\begin{equation} \label{eqn:solA142} 
q =\dfrac{1}{t} + \dfrac{t-1}{t} 
\sum_{k=1}^{\infty} c_k(t) \, z^k, \qquad 
p = t \sum_{k=0}^{\infty} d_k(t) \, z^k, 
\end{equation}
%%%%%%%%%%%%%%%%%%%%%%%%%%%%%%%%%%%%%%%%%%%%%%%%%%%%%%%%%%%%%%
depending on $t \in \C$, where $c_k(t)$ and $d_k(t)$ are 
polynomials of $t$ beginning with $c_1(t) = 1/2$ and 
$d_0(t) = -1/2$. 
For $t = 0$, formula $(\ref{eqn:solA142})$ represents 
the solution such that $q \equiv \infty$ and $p \equiv 0$. 
\end{lemma} 
%%%%%%%%%%%%%%%%%%%%%%% proof %%%%%%%%%%%%%%%%%%%%%%%%%%%%%%%
{\it Proof}. 
We only derive (\ref{eqn:solA142}), as (\ref{eqn:solA141}) 
is derived in a similar manner. 
Substituting $q = t^{-1} + t^{-1}(t-1) Q$ into 
(\ref{eqn:PVI3}) yields 
$P(z, t^{-1}+ t^{-1}(t-1) Q) = 
t^{-4}(t-1)^2 z^{-1} p(z,Q;t)$, 
where $p(z,Q;t)$ is a differential sum of $(z,Q)$ with 
coefficients in $\C[t]$ whose Newton polygon is given as 
in Figure \ref{fig:snewtonD4A14} (right). 
Consider the edge connecting $(0,1)$ and $(1,0)$. 
The vertex $(0,1)$ carries the differential sum 
$\mathcal{L}Q = -2(z Q' + z^2 Q'')$, whose 
characteristic polynomial is $-2 k^2$, while the vertex 
$(1,0)$ carries the monomial $z$. 
Since the critical value $k = 0$ is smaller than $1$, 
the coefficients $c_k(t)$, $k \ge 1$, in (\ref{eqn:solA142}) 
are determined uniquely and recursively, where the leading 
term is $c_1(t) = 1/2$.
Substituting the resulting series $q = q(z)$ into 
(\ref{eqn:p(z)}) we have $p = p(z)$ as in 
(\ref{eqn:solA142}). 
\hfill $\Box$ \par\medskip 
%%%%%%%%%%%%%%%%%%%%%%% end of proof %%%%%%%%%%%%%%%%%%%%%%%%
%%%%%%%%%%%%%%%%%%%%%%% vertex (0,3) %%%%%%%%%%%%%%%%%%%%%%%%
So far we have considered the edges $\varGamma_1$ and 
$\varGamma_0$ of the Newton polygon in 
Figure \ref{fig:snewtonPVI} and constructed meromorphic 
solutions around the origin $z = 0$. 
Now let us consider the vertex $(0,3)$ of the polygon, 
which gives rise to algebraic branch solutions around 
$z = 0$. 
The truncation of the differential sum $P = P(z,q)$ at the 
vertex $(0,3)$ is given by 
%%%%%%%%%%%%%
\[
P_3 = -2z q^2 q' + 2z^2 q (q')^2 - 2z^2 q^2 q''. 
\]
%%%%%%%%%%%%%
Its charactersistic polynomial $\chi(r)$ is defined by 
substituting $q = z^r$ into $P_3(z,q)$ and dividing 
the result by $q^3 = z^{3r}$. 
In the present situation we see that $\chi(r)$ is 
identically zero; $\chi(r) \equiv 0$. 
The normal cone of the vertex $(0,3)$ is 
$U_3 = \{\, (p_1, p_2) \in \bR^2 \,:\, 
p_1 < 0, \, 0 < r = p_2/p_1 < 1 \,\}$. 
Thus the truncated solutions at the vertex $(0,3)$ are 
$q = t z^r$ for an arbitrary $0 < r < 1$ and 
$t \in \C^{\times}$. 
The Fr\'{e}chet derivative with respect to $q$ at the 
truncated solution $q = t z^r$ is given by 
%%%%%%%%
\[
\mathcal{L}_3 Q = 
-2 t^2 z^{2r} \{ z^2 Q''+(1-2r)z Q' + r^2 Q \}. 
\]
%%%%%%%%
The corresponding characteristic equation is given by 
$v_3(k) := -2 t^2 (k - r)^2 = 0$, which has the 
only root $k = r$. 
Let $n$ be any integer greater than $1$. 
In order to search for an algebraic $n$-branch solution 
around $z = 0$, we take $r = m/n$ for any integer 
$0 < m < n$ coprime to $n$, and consider a formal 
Puiseux series solution of the form 
%%%%%%% 
\[
q = t z^{m/n} + \sum_{\nu=m+1}^{\infty} a_{\nu}(t) 
\, z^{\nu/n}. 
\] 
%%%%%%
Since the characteristic equation $v_3(k) = 0$ has 
no roots such that $k > m/n$, the coefficients 
$a_{\nu} = a_{\nu}(t)$ can be determined uniquely and 
recursively for any given initial coefficient 
$t \in \C^{\times}$. 
The convergence of the formal solution and its 
holomorphic dependence on parameters follow 
easily if we rewrite the equation (\ref{eqn:PVI3}) 
in terms of the new independent variable 
$\zeta = z^{1/n}$ and apply the convergence 
arguments in \cite{Gerard,GS}. 
Thus we have established the following lemma. 
%%%%%%%%%%%%%%%%%%%%%%%%% lem:algsol %%%%%%%%%%%%%%%%%%%%%%%%%
\begin{lemma} \label{lem:algsol} 
For any integer $n > 1$, there exist $\varphi(n)$ mutually 
disjoint $1$-parameter families of $n$-branch solution 
germs to $\PVI(\k)$ around the origin $z = 0$, 
%%%%%%%%%%%%%%%%%%%%%%%%% eqn:algsol %%%%%%%%%%%%%%%%%%%%%%%%%
\begin{equation} \label{eqn:algsol} 
\begin{array}{rcl}
q(z) &=& t z^{m/n} + 
\displaystyle \sum_{\nu=m+1}^{\infty} 
a_{\nu}(n,m,t;\k) \, z^{\nu/n}, \\[5mm]
p(z) &=& \dfrac{m+n(\k_1+\k_2-1)}{2 n t} \, 
z^{-m/n} + \displaystyle \sum_{\nu=-m+1}^{\infty} 
b_{\nu}(n,m,t;\k) \, z^{\nu/n},
\end{array} 
\end{equation}
%%%%%%%%%%%%%%%%%%%%%%%%%%%%%%%%%%%%%%%%%%%%%%%%%%%%%%%%%%%%%%
where the discrete parameter $m$ ranges over all integers 
$0 < m < n$ coprime to $n$ and the continuous parameter 
$t$ takes any value of the punctured complex line 
$\C^{\times}$. 
\end{lemma}
%%%%%%%%%%%%%%%%%%%%%%%%%%%%%%%%%%%%%%%%%%%%%%%%%%%%%%%%%%%%%%
%%%%%%%%%%%%%%%%%%%%%%%%% rem:algsol %%%%%%%%%%%%%%%%%%%%%%%%%
\begin{remark} \label{rem:algsol} 
The family (\ref{eqn:algsol}) contains no Riccati solutions, 
even if $\k \in \Wall$. 
This will be shown in the proof of Lemma \ref{lem:proofPer}. 
\end{remark}
%%%%%%%%%%%%%%%%%%%% sec:surjection %%%%%%%%%%%%%%%%%%%%%%%%%%
\section{Injection Implies Surjection} \label{sec:surjection} 
%%%%%%%%%%%%%%%%%%%%%%%%%%%%%%%%%%%%%%%%%%%%%%%%%%%%%%%%%%%%%%
We establish Theorem \ref{thm:main} based on the main idea 
described in \S\ref{sec:intro}. 
Due to the $S_3$-symmetry permuting the three fixed 
singular points, it suffices to work around $z = 0$ 
(see Remark \ref{rem:backlund}). 
\par 
As a preliminary we begin by constructing some Riccati 
solutions to equation (\ref{eqn:PVI}). 
Assume that $\k_0 = 0$ so that $\k_1+\k_2+\k_3+\k_4 = 1$. 
Then the second equation of system (\ref{eqn:PVI}) has the 
null solution $p(z) \equiv 0$. 
Substituting this into the first equation yields the 
Riccati equation 
%%%%%%%
\[
z(z-1) q' + \k_1 q_1 q_z + (\k_2-1) q_0 q_1 
+ \k_3 q_0 q_z = 0. 
\]
%%%%%%
If $\k_4$ is nonzero, then the change of dependent variable 
%%%%%
\[
q = \dfrac{z(1-z)}{\k_4} \dfrac{d}{dz} 
\log \{(1-z)^{-\k_4} f \}
\]
%%%%%
transfers the Riccati equation to the Gauss hypergeometric 
equation 
%%%%%%%%%%%%%%%%%%%%%%%%% eqn:gauss1 %%%%%%%%%%%%%%%%%%%%%%%%%
\begin{equation} \label{eqn:gauss1}
z(1-z) f'' + \{(1-\k_3-\k_4)-(\k_2-\k_4+1)z \} f'+ 
\k_2\k_4 f = 0. 
\end{equation}
%%%%%%%%%%%%%%%%%%%%%%%%%%%%%%%%%%%%%%%%%%%%%%%%%%%%%%%%%%%%%%
\par 
Next assume that $\k_0 = \k_1 = 0$ so that $\k_2+\k_3+\k_4=1$. 
In this case there is another type of Riccati solution to the 
system (\ref{eqn:PVI}). 
The first equation of the system (\ref{eqn:PVI}) has the null 
solution $q(z) \equiv 0$. 
Substituting this into the second equation yields the Riccati 
equation 
%%%%%%%%%
\[
z(z-1) p' + z p^2 + (\k_2-1 + \k_3 z) p = 0. 
\]
%%%%%%%%
Then change of independent variable 
$p = (z-1) \dfrac{d}{dz} \log g$ takes it to the linear 
equation 
%%%%%%%%%%%%%%%%%%%%%%%%%% eqn:gauss2 %%%%%%%%%%%%%%%%%%%%%%%%
\begin{equation} \label{eqn:gauss2} 
z(1-z) g'' + \{(1-\k_2) -(\k_3 + 1) z\} g' = 0. 
\end{equation}
%%%%%%%%%%%%%%%%%%%%%%%%%%%%%%%%%%%%%%%%%%%%%%%%%%%%%%%%%%%%%%
%%%%%%%%%%%%%%%%%%%%%%%%%% lem:a1ric %%%%%%%%%%%%%%%%%%%%%%%%%
\begin{lemma}[$\mbox{\boldmath $A_1$}$] \label{lem:a1ric} 
Assume that $\k_0 = 0$, $\k_1+\k_2+\k_3+\k_4 = 1$, 
$\k_1 + \k_2 \not\in \Z$, $\k_3 + \k_4 \not\in \Z$ and 
$\k_1 \k_4 \neq 0$. 
Then system $(\ref{eqn:PVI})$ has two single-valued 
Riccati solutions around the origin $z = 0$, 
%%%%%%%%%
\[
\begin{array}{crclrcl}
(\mathrm{a}) \qquad & q &=& \dfrac{\k_1 z}{\k_1+\k_2} 
+ O(\k_1 z^2), 
\qquad & p &\equiv& 0, \\[4mm] 
(\mathrm{b}) \qquad & q &=& \dfrac{\k_3 + \k_4}{\k_4} 
+ O(\k_3 z), 
\qquad & p &\equiv& 0. 
\end{array}
\]
%%%%%%%%%
\end{lemma}
%%%%%%%%%%%%%%%%%%%%%%%%%%%%%%%%%%%%%%%%%%%%%%%%%%%%%%%%%%%%%%
{\it Proof}. 
The solutions (a) and (b) are obtained from two linearly 
independent solutions 
%%%%%
\[
{}_2F_1(\k_2,-\k_4,1-\k_3-\k_4;z), \qquad  
z^{\k_3+\k_4} {}_2F_1(\k_2+\k_3+\k_4,\k_3,\k_3+\k_4+1;z), 
\]
%%%%% 
of equation (\ref{eqn:gauss1}) repsectively. 
They also come from solutions (\ref{eqn:solk1k21}) 
and (\ref{eqn:solk3k41}) respectively. 
\hfill $\Box$
%%%%%%%%%%%%%%%%%%%%%%%%%% lem:a2ric %%%%%%%%%%%%%%%%%%%%%%%%%
\begin{lemma}[$\mbox{\boldmath $A_2$}$] \label{lem:a2ric} 
Assume that $\k_0 = \k_1 = 0$, $\k_2 + \k_3 + \k_4 = 1$, 
$\k_2 \not\in \Z$, $\k_3 \neq 1$ and $\k_4 \neq 0$. 
Then system $(\ref{eqn:PVI})$ has three single-valued 
Riccati solutions around the origin $z = 0$, 
%%%%%%%%%%%%%%
\[ 
\begin{array}{crclrcl}
(\mathrm{a}) \qquad & q &\equiv& 0, \qquad & p &\equiv& 0,  
\\[2mm] 
(\mathrm{b}) \qquad & q &=& \dfrac{\k_3+\k_4}{\k_4} + O(z), 
\qquad & p &\equiv& 0, \\[2mm]
(\mathrm{c}) \qquad & q &\equiv& 0, \qquad 
& p &=& -\dfrac{\k_2(\k_3-1)}{\k_2+1} + O(z). 
\end{array}
\]
%%%%%%%%%%%%%%
\end{lemma}
%%%%%%%%%%%%%%%%%%%%%%%%%%%%%%%%%%%%%%%%%%%%%%%%%%%%%%%%%%%%%%
{\it Proof}. The solutions (a) and (b) just come from 
(a) and (b) of Lemma \ref{lem:a1ric} respectively, 
while solution (c) is obtained from the solution 
$z^{\k_2} {}_2F_1(\k_2,\k_2+\k_3,\k_2+1;z)$ of 
equation (\ref{eqn:gauss2}). \hfill $\Box$ 
%%%%%%%%%%%%%%%%%%%%%%%%%% lem:a3ric %%%%%%%%%%%%%%%%%%%%%%%%%
\begin{lemma}[$\mbox{\boldmath $A_3$}$] \label{lem:a3ric} 
Assume that $\k_0 = \k_1 = \k_2 = 0$, $\k_3+\k_4 = 1$ 
and $\k_4 \neq 0$. 
Then system $(\ref{eqn:PVI})$ admits a $1$-parameter 
family of single-valued Riccati solutions around the origin 
$z = 0$, 
%%%%%%%%%%%%%%%%%%%%%%%%%% eqn:a3single %%%%%%%%%%%%%%%%%%%%%%
\begin{equation} \label{eqn:a3single} 
q(z;s) = \dfrac{s_0 z}{s_0 + s_1 (1-z)^{\k_4}}, \quad 
p(z;s) \equiv 0, \qquad (s = [s_0:s_1] \in \P^1). 
\end{equation}
%%%%%%%%%%%%%%%%%%%%%%%%%%%%%%%%%%%%%%%%%%%%%%%%%%%%%%%%%%%%%%
\end{lemma}
%%%%%%%%%%%%%%%%%%%%%%%%%%%%%%%%%%%%%%%%%%%%%%%%%%%%%%%%%%%%%%
{\it Proof}. 
The hypergeometric equation (\ref{eqn:gauss1}) becomes 
$(1-z) f''-\k_3 f' = 0$, whose nontrivial solutions 
are given by $f(z) = s_0 + s_1 (1-z)^{\k_4}$ with 
$(s_0,s_1) \in \C^2-\{(0,0)\}$. 
The corresponding Riccati solutions are the $1$-parameter 
family of single-valued solutions $q = q(z;s)$ as in 
(\ref{eqn:a3single}). \hfill $\Box$ 
%%%%%%%%%%%%%%%%%%%%%%%%% end of proof %%%%%%%%%%%%%%%%%%%%%%%
%%%%%%%%%%%%%%%%%%%%%%%%%% lem:d4ric %%%%%%%%%%%%%%%%%%%%%%%%%
\begin{lemma}[$\mbox{\boldmath $D_4$}$] \label{lem:d4ric} 
Assume that $\k_0 = \k_1 = \k_2 = \k_3 = 0$ and $\k_4 = 1$. 
\begin{enumerate}
\item System $(\ref{eqn:PVI})$ admits a $1$-parameter 
family of rational Riccati solutions 
%%%%%%%%%%%%%%%%%%%%%%%%%% eqn:d4rat %%%%%%%%%%%%%%%%%%%%%%%%%
\begin{equation} \label{eqn:d4rat} 
q(z;s) = \dfrac{s_0 z}{s_0 + s_1 (1-z)}, 
\quad p(z;s) \equiv 0, \qquad (s = [s_0:s_1] \in \P^1). 
\end{equation}
%%%%%%%%%%%%%%%%%%%%%%%%%%%%%%%%%%%%%%%%%%%%%%%%%%%%%%%%%%%%%%
\item System $(\ref{eqn:PVI})$ admits a $1$-parameter 
family of single-valued Riccati solutions around $z = 0$, 
%%%%%%%%%%%%%%%%%%%%%%%%% eqn:d4single %%%%%%%%%%%%%%%%%%%%%%%
\begin{equation} \label{eqn:d4single}
q(z;t) \equiv 1, \quad 
p(z;t) = \dfrac{t_0 z}{(1-z) \{t_0 \log(1-z) + t_1\}}, 
\qquad (t = [t_0:t_1] \in \P^1). 
\end{equation}
%%%%%%%%%%%%%%%%%%%%%%%%%%%%%%%%%%%%%%%%%%%%%%%%%%%%%%%%%%%%%%
\end{enumerate}
\end{lemma}
%%%%%%%%%%%%%%%%%%%%%%%%%%%%%%%%%%%%%%%%%%%%%%%%%%%%%%%%%%%%%%
{\it Proof}. 
In this case (\ref{eqn:a3single}) gives the $1$-parameter 
family of rational solutions (\ref{eqn:d4rat}). 
Moreover the first equation of system (\ref{eqn:PVI}) admits 
a constant solution $q(z) \equiv 1$. 
Substituting this into the second 
equation yelds the Riccati equation 
$z(z-1) p' + (1-z) p^2 + p = 0$. 
Change of dependent variable $p = -z f'/f$ takes this 
into the linear equation $(1-z)f''-f'=0$, whose nontrivial 
solutions are given by $f = t_0 \log(1-z) + t_1$ with 
$(t_0,t_1) \in \C^2-\{(0,0)\}$. 
Thus the Riccati equation has the $1$-parameter family of 
single-valued solutions $p = p(z;t)$ as in 
(\ref{eqn:d4single}). \hfill $\Box$ \par\medskip 
%%%%%%%%%%%%%%%%%%%%%%%%% end of proof %%%%%%%%%%%%%%%%%%%%%%% 
Now we proceed to the proof of Theorem \ref{thm:main}. 
From now on we fix the indices as $(i,j,k) = (3,1,2)$ in 
accordance with the choice of indices in \S\ref{sec:power}. 
First we treat the fixed point case. 
%%%%%%%%%%%%%%%%%%%%%%%% lem:proofFix %%%%%%%%%%%%%%%%%%%%%%%%
\begin{lemma} \label{lem:proofFix}
The set $\wt{\mathrm{Fix}}_j(\th)$ is exhausted by meromorphic 
solutions around $z = 0$. 
\end{lemma}
%%%%%%%%%%%%%%%%%%%%%%%% proof %%%%%%%%%%%%%%%%%%%%%%%%%%%%%%%
{\it Proof}. Case-by-case check based on the 
``injection-implies-surjection" principle described in 
\S\ref{sec:intro}. 
%%%%%%%%%%%%%%%%%%%%%%%%% ex:surjEmpty %%%%%%%%%%%%%%%%%%%%%%%
\begin{example}[$\mbox{\boldmath $\emptyset$}$] 
\label{ex:surjEmpty} 
We combine the results of Example \ref{ex:empty}, Lemmas 
\ref{lem:solk1k2} and \ref{lem:solk3k4}. 
A key observation is that Lemmas \ref{lem:solk1k2} and 
\ref{lem:solk3k4} give us as many meromorphic solutions 
around $z = 0$ as the cardinality of the set 
$\wt{\mathrm{Fix}}_j(\th)$ in (\ref{eqn:FixEmpty}). 
For example, if $\wt{P}(b_i,b_4;b_j,b_k) \in 
\wt{\mathrm{Fix}}_j(\th)$, then the existing and 
smoothness conditions for it (see Table \ref{tab:fixed}) 
makes it possible to apply Lemma \ref{lem:solk3k4} to 
conclude that the meromorphic solution (\ref{eqn:solk3k41}) 
exists corresponding to the fixed point 
$\wt{P}(b_i,b_4;b_j,b_k)$. 
Thus the set $\wt{\mathrm{Fix}}_j(\th)$ is exhausted by 
meromorphic solutions around $z = 0$. 
\end{example}
%%%%%%%%%%%%%%%%%%%%%%%%% ex:surjA1 %%%%%%%%%%%%%%%%%%%%%%%%%%
\begin{example}[$\mbox{\boldmath $A_1$}$] 
\label{ex:surjA1} 
We combine the results of Example \ref{ex:A1}, Lemmas 
\ref{lem:solk1k2}, \ref{lem:solk3k4} and \ref{lem:a1ric}. 
First we notice that the two single-valued Riccati solutions 
in Lemma \ref{lem:a1ric} correspond to the two Riccati fixed 
points $\wt{\mathrm{Fix}}_j^{e}(\th) = \{p,q\}$ in 
(\ref{eqn:FixA1}). 
On the other hand, for the same reason as in 
Example \ref{ex:surjEmpty}, formulas (\ref{eqn:solk1k22}) 
and (\ref{eqn:solk3k42}) in Lemmas \ref{lem:solk1k2} and 
\ref{lem:solk3k4} give us as many meromorphic solutions 
around $z = 0$ as the cardinality of smooth fixed points 
$\wt{\mathrm{Fix}}_j^{\circ}(\th) = 
\{\!\{ \wt{P}(b_i,b_4^{-1};b_j,b_k), \, 
\wt{P}(b_j,b_k^{-1};b_i,b_4)\}\!\}$ in (\ref{eqn:FixA1}). 
Thus the set $\wt{\mathrm{Fix}}_j(\th)$ is exhausted by 
meromorphic solutions around $z = 0$. 
\end{example}
%%%%%%%%%%%%%%%%%%%%%%%%% ex:surjA2 %%%%%%%%%%%%%%%%%%%%%%%%%%
\begin{example}[$\mbox{\boldmath $A_2$}$] \label{ex:surjA2} 
We combine the results of Example \ref{ex:A2}, Lemmas 
\ref{lem:solk3k4} and \ref{lem:a2ric}. 
As (\ref{eqn:FixA2}) shows, the set $\wt{\mathrm{Fix}}_j(\th)$ 
consists of the four points $p_0$, $p_+$, $p_-$ and 
$\wt{P}(b_j,b_k^{-1};b_i,b_4)$. 
On the other hand, we have the three single-valued Riccati 
solutions of Lemma \ref{lem:a2ric} and one 
non-Riccati holomorphic solution (\ref{eqn:solk3k42}). 
Clearly, the three Riccati solutions correspond to 
the points $p_0$, $p_+$ and $p_-$, while the non-Riccati 
solution corresponds to the remaining point 
$\wt{P}(b_j,b_k^{-1};b_i,b_4)$. 
Thus any single-valued solution around $z = 0$ is a 
meromorphic solution. 
\end{example}
%%%%%%%%%%%%%%%%%%%%%%%%% ex:surjA12 %%%%%%%%%%%%%%%%%%%%%%%%%
\begin{example}[$\mbox{\boldmath $A_1^{\oplus 2}$}$] 
\label{ex:surjA12} 
We combine the results of Example \ref{ex:A12}, Lemmas 
\ref{lem:solk3k4} and \ref{lem:solA12}. 
First we consider the $\wt{W}(D_4^{(1)})$-stratum of type 
$(A_1^{\oplus 2})_i$. 
The $\C$-parameter family (\ref{eqn:solA12}) of 
holomorphic solutions injects into the line 
$\wt{\ell}_j^{+} \simeq \C$ in (\ref{eqn:FixA12i}), so 
that we have an injection $\C \hookrightarrow \C$. 
Since this injection is holomorphic, it must be a 
surjection. 
Thus the line $\wt{\ell}_j^{+}$ is exhausted by the 
family (\ref{eqn:solA12}). 
Moreover we have the two holomorphic solutions 
(\ref{eqn:solk3k41}) and (\ref{eqn:solk3k42}), 
which do not lie in the family (\ref{eqn:solA12}). 
Thus they must correspond to the points 
$\wt{P}(b_i,b_4;b_j,b_k)$ and 
$\wt{P}(b_i,b_4^{-1};b_j,b_k)$ in (\ref{eqn:FixA12i}). 
So on the stratum of type $(A_1^{\oplus 2})_i$ the set 
$\wt{\mathrm{Fix}}_j(\th)$ is exhausted by 
meromorphic solutions around $z = 0$. 
Next we consider the $\wt{W}(D_4^{(1)})$-strata of types 
$(A_1^{\oplus 2})_j$ and $(A_1^{\oplus 2})_k$. 
On these strata the equality 
$\wt{\mathrm{Fix}}_j(\th) = \wt{\mathrm{Fix}}_j^{e}(\th)$ 
in (\ref{eqn:FixA12jk}) implies that any single-valued 
solution around $z = 0$ is a Riccati and hence meromorphic 
solution. 
\end{example} 
%%%%%%%%%%%%%%%%%%%%%%%%%%%%%%%%%%%%%%%%%%%%%%%%%%%%%%%%%%%%%%
\begin{example}[$\mbox{\boldmath $A_3$}$] \label{ex:surjA3} 
We combine the results of Example \ref{ex:A3}, Lemmas 
\ref{lem:solk3k4}, \ref{lem:solA3} and \ref{lem:a3ric}. 
First we consider the $\wt{W}(D_4^{(1)})$-stratum of type 
$(A_1^{\oplus 2})_i$. 
The $\P^1$-family (\ref{eqn:a3single}) of single-valued 
Riccati solutions exactly corresponds to the exceptional 
curve $e_0 \simeq \P^1$ in (\ref{eqn:FixA3i}). 
So the $\C$-family (\ref{eqn:solA3}) of holomorphic 
solutions must inject into the line 
$\wt{\ell}_j^{+} \simeq \C$ in (\ref{eqn:FixA3i}). 
Since this injection $\C \hookrightarrow \C$ is 
holomorphic, it must be a surjection. 
Hence $\wt{\ell}_j^{+}$ is exhausted by the family 
(\ref{eqn:solA3}). 
Moreover there is the holomorphic solution 
(\ref{eqn:solk3k42}), which must correspond to the point 
$\wt{P}(b_i,b_4^{-1};b_j,b_k)$ in (\ref{eqn:FixA3i}). 
Thus on the stratum of type $(A_1^{\oplus 2})_i$ the set 
$\wt{\mathrm{Fix}}_j(\th)$ is exhausted by meromorphic 
solutions around $z = 0$. 
Next we consider the $\wt{W}(D_4^{(1)})$-strata of types 
$(A_1^{\oplus 2})_j$ and $(A_1^{\oplus 2})_k$. 
On these strata the equality 
$\wt{\mathrm{Fix}}_j(\th) = \wt{\mathrm{Fix}}_j^{e}(\th)$ 
in (\ref{eqn:FixA3jk}) implies that any single-valued 
solution around $z = 0$ is a Riccati and hence meromorphic 
solution. 
\end{example}
%%%%%%%%%%%%%%%%%%%%%%%%% ex:surjA13 %%%%%%%%%%%%%%%%%%%%%%%%%
\begin{example}[$\mbox{\boldmath $A_1^{\oplus 3}$}$] 
\label{ex:surjA13} 
We combine the results of Example \ref{ex:A13} and 
Lemma \ref{lem:solA12}. 
As (\ref{eqn:FixA13}) shows, $\wt{\mathrm{Fix}}_j(\th)$ has 
only one line component $\wt{\ell}_j^{+} \simeq \C$. 
Hence the $\C$-family (\ref{eqn:solA12}) of holomorphic 
solutions must inject into this line, so that we 
have an inclusion $\C \hookrightarrow \C$. 
Since this injection is holomorphic, it must be a 
surjection. 
Thus $\wt{\ell}_j^{+}$ is exhausted by the family 
(\ref{eqn:solA12}). 
\end{example}
%%%%%%%%%%%%%%%%%%%%%%%%% ex:surjD4 %%%%%%%%%%%%%%%%%%%%%%%%%%
\begin{example}[$\mbox{\boldmath $D_4$}$] \label{ex:surjD4} 
We combine the results of Example \ref{ex:D4}, Lemmas 
\ref{lem:solD4} and \ref{lem:d4ric}. 
The $\P^1$-family (\ref{eqn:d4rat}) of rational Riccati 
solutions corresponds to the exceptional curve $e_0$ in 
(\ref{eqn:FixD4}), while the $\P^1$-family (\ref{eqn:d4single}) 
of single-valued Riccati solutions corresponds to the 
exceptional curve $e_j$ there. 
Hence $\C$-family (\ref{eqn:solD4}) of holomorphic 
solutions, which is different from (\ref{eqn:d4rat}) and 
(\ref{eqn:d4single}), must inject into the line 
$\wt{\ell}_j^{+} \simeq \C$ in (\ref{eqn:FixD4}). 
Since this injection $\C \hookrightarrow \C$ is 
holomorphic, it must be a surjection. 
Thus $\wt{\ell}_j^{+}$ is exhausted by the family 
(\ref{eqn:solD4}). 
\end{example}
%%%%%%%%%%%%%%%%%%%%%%%%% ex:surjA14 %%%%%%%%%%%%%%%%%%%%%%%%%
\begin{example}[$\mbox{\boldmath $A_1^{\oplus 4}$}$] 
\label{ex:surjA14} 
We combine the results of Example \ref{ex:A14} and 
Lemma \ref{lem:solA14}. 
In view of $\wt{\mathrm{Fix}}_j(\th) = \wt{\ell}_j^{+} 
\amalg \wt{\ell}_j^{-}$, the $\C$-family (\ref{eqn:solA141}) 
of holomorphic solutions injects into the line 
$\wt{\ell}_j^{\ve} \simeq \C$ for some sign 
$\ve \in \{\pm 1\}$. 
So we have an injection $\C \hookrightarrow \C$. 
Since this injection is holomorphic, it must be a 
surjection, so that $\wt{\ell}_j^{\ve}$ is exhausted by 
the family (\ref{eqn:solA141}). 
Then the other $\C$-family (\ref{eqn:solA142}) of 
holomorphic solutions injects into the remaining line 
$\wt{\ell}_j^{-\ve} \simeq \C$. 
So we have another injection $\C \hookrightarrow \C$. 
Since this injection is holomorphic, it must be a 
surjection. 
Hence $\wt{\mathrm{Fix}}_j(\th) = 
\wt{\ell}_j^{+} \amalg \wt{\ell}_j^{-}$ is exhausted by 
the families (\ref{eqn:solA141}) and (\ref{eqn:solA142}). 
The proof of Lemma \ref{lem:proofFix} is now complete. 
\hfill $\Box$ 
\end{example}
%%%%%%%%%%%%%%%%%%%%%%%%%%%%%%%%%%%%%%%%%%%%%%%%%%%%%%%%%%%%%%
\par 
Finally we argue the periodic point case 
using the ``injection-implies-surjection" principle. 
%%%%%%%%%%%%%%%%%%%%%%%% lem:proofPer %%%%%%%%%%%%%%%%%%%%%%%%
\begin{lemma} \label{lem:proofPer}
For any $n > 1$ the set $\wt{\mathrm{Per}}_j(\th;n)$ is 
exhausted by algebraic $n$-branch solutions around $z = 0$. 
\end{lemma} 
%%%%%%%%%%%%%%%%%%%%%%%% proof %%%%%%%%%%%%%%%%%%%%%%%%%%%%%%%
{\it Proof}. 
We combine Lemmas \ref{lem:pern} and \ref{lem:algsol}. 
First we consider the generic case where $\k \in \K-\Wall$, 
namely, where $\th = \rh(\k)$ is such that 
$\varDelta(\th) \neq 0$. 
In this case there is no Riccati locus and hence 
$\wt{\mathrm{Per}}_j(\th;n) = 
\wt{\mathrm{Per}}_j^{\circ}(\th;n)$, which is biholomorphic 
to the disjoint union of $\varphi(n)$ copies of 
$\C^{\times}$ by Lemma \ref{lem:pern}. 
On the other hand, by Lemma \ref{lem:algsol}, there are 
$\varphi(n)$ mutually disjoint $\C^{\times}$-parameter 
families of {\sl algebraic} $n$-branch 
solutions around $z = 0$ as in (\ref{eqn:algsol}). 
Number these families from $1$ to $\varphi(n)$. 
The first family injects into a (unique) connected 
component ($\simeq \C^{\times}$) of 
$\wt{\mathrm{Per}}_j(\th;n)$, which we call the first 
component, and we have an injection $\C^{\times} 
\hookrightarrow \C^{\times}$. 
Since this injection is holomorphic, it must be a 
surjection and hence the first component is exhausted 
by the first family. 
Consider the second family of solutions and the 
corresponding second component of 
$\wt{\mathrm{Per}}_j(\th;n)$. 
Notice that the second component is different from 
the first one, because the first component is already 
occupied by the first family and so it cannot contain 
the second family. 
For the same reason as above, the second component is 
exhausted by the second family. 
Since the families and the components have the same 
cardinality $\varphi(n)$, we can repeat this argument 
to conclude that $\wt{\mathrm{Per}}_j(\th;n)$ is 
exhausted by the $\varphi(n)$ families of algebraic 
$n$-branch solutions. 
\par 
Next we consider the case where $\k \in \Wall$, namely, 
where the Riccati part $\wt{\mathrm{Per}}_j^{e}(\th;n)$ 
may appear. 
Since the lemma is trivial for the Riccati part, we 
have only to consider the non-Riccati part 
$\wt{\mathrm{Per}}_j^{\circ}(\th;n)$. 
The argument proceeds just in the same manner as in 
the last paragraph, once we show that the family of 
solutions in (\ref{eqn:algsol}) contains no Riccati 
solutions (see Remark \ref{rem:algsol}). 
To see this, we consider the family $\wt{\Sol} \to \Th$ 
of surfaces $\wt{\Sol}(\th)$ parametrized by 
$\th \in\Th$ and put 
%%%%%%%%
\[
\wt{\mathrm{Per}}_j^{\circ}(n) = \coprod_{\th \in \Th} 
\wt{\mathrm{Per}}_j^{\circ}(\th;n), 
\qquad 
\mathcal{E} = \coprod_{\th \in \Th} \mathcal{E}(\th), 
\]
%%%%%%%%%
where $\mathcal{E}(\th)$ is the exceptional set in 
$\wt{\Sol}(\th)$. 
(Precisely speaking, the parameter space $\Th$ should 
be replaced by a finite covering of it to get a 
simultaneous minimal resolution.) 
Then $\wt{\mathrm{Per}}_j^{\circ}(n)$ and 
$\mathcal{E}$ are closed subsets of $\wt{\Sol}$ 
which are {\sl disjoint} by Lemma \ref{lem:updown}. 
Now we look at the family of solutions in 
(\ref{eqn:algsol}). 
It depends continuously on $\k \in \K$. 
Take any point $\k^* \in \Wall$ and let 
$\K-\Wall \ni \k \to \k^*$. 
For any $\k \in \K-\Wall$, the family at $\k$ is 
contained in $\wt{\mathrm{Per}}_j^{\circ}(\th;n)$ 
with $\th = \rh(\k)$ and hence in 
$\wt{\mathrm{Per}}_j^{\circ}(n)$. 
Taking the limit $\k \to \k^*$, we see that 
the family at $\k^*$ is contained in 
$\wt{\mathrm{Per}}_j^{\circ}(n)$, hence in 
$\wt{\mathrm{Per}}_j^{\circ}(\th^*;n)$ with 
$\th^* = \rh(\k^*)$. 
Since $\wt{\mathrm{Per}}_j^{\circ}(\th^*;n)$ 
is disjoint from $\mathcal{E}(\th^*)$, 
the family at $\k^*$ contains no Riccati solutions. 
Therefore the proof is complete. 
\hfill $\Box$\par\medskip 
%%%%%%%%%%%%%%%%%%%%%%%% end of proof %%%%%%%%%%%%%%%%%%%%%%%%
Now the local statement of Theorem \ref{thm:main} around a 
fixed singular point, say $z = 0$, is an immediate 
consequence of Lemmas \ref{lem:proofFix} and \ref{lem:proofPer}. 
At the same time all the finite branch solutions around 
$z = 0$ have been classified up to B\"acklund transformations. 
The global statement about algebraic solutions follows readily 
from the local statements around $z = 0$, $1$, $\infty$, 
together with the analytic Painlev\'e property on 
$Z = \P^1-\{0,1,\infty\}$. 
The proof of Theorem \ref{thm:main} is complete. 
%%%%%%%%%%%%%%%%%%%%%%%%% references %%%%%%%%%%%%%%%%%%%%%%%%%

%%%%%%%%%%%%%%
\end{document}